\input amstex
\documentstyle{amsppt}
\nopagenumbers
\nologo
%
%
%
\catcode`@=11
\redefine\output@{%
  \def\break{\penalty-\@M}\let\par\endgraf
  \ifnum\pageno=1\global\voffset=90pt\else\global\voffset=25pt\fi
  \ifodd\pageno\global\hoffset=105pt\else\global\hoffset=8pt\fi  
  \shipout\vbox{%
    \ifplain@
      \let\makeheadline\relax \let\makefootline\relax
    \else
      \iffirstpage@ \global\firstpage@false
        \let\rightheadline\frheadline
        \let\leftheadline\flheadline
      \else
        \ifrunheads@ 
        \else \let\makeheadline\relax
        \fi
      \fi
    \fi
    \makeheadline \pagebody \makefootline}%
  \advancepageno \ifnum\outputpenalty>-\@MM\else\dosupereject\fi
}
\font\cpr=cmr7
\newcount\xnumber
\footline={\xnumber=\pageno
\divide\xnumber by 7
\multiply\xnumber by -7
\advance\xnumber by\pageno
\ifnum\xnumber>0\hfil\else\vtop{\vskip 0.5cm
\noindent\cpr CopyRight \copyright\ Sharipov R.A., 1996, 
2004.}\hfil\fi}
\def\setfirstpage{\global\firstpage@true}
\catcode`\@=\active

\fontdimen3\tenrm=3pt
\fontdimen4\tenrm=0.7pt

\def\leaderfill{\leaders\hbox to 0.3em{\hss.\hss}\hfill}
\font\tvbf=cmbx12
\font\tvrm=cmr12
\font\etbf=cmbx8
\font\tencyr=wncyr10
\font\eightcyr=wncyr8

  \let\tg=\tan
  \let\arctg=\arctan
  \let\ctg=\cot

\Monograph
\def\negskp{\hskip -2pt}
\def\compos{\,\raise 1pt\hbox{$\sssize\circ$} \,}
\def\tr{\operatorname{tr}}
\def\const{\operatorname{const}}
\def\rank{\operatorname{rank}}
\def\sign{\operatorname{sign}}
\def\grad{\operatorname{grad}}
\def\divr{\operatorname{div}}
\def\rot{\operatorname{rot}}
\newcount\chapternum
\def\blue#1{#1}
\catcode`#=11\def\diez{#}\catcode`#=6
\def\mycite#1{\cite{\blue{#1}}\immediate\special{ps:
     ShrHPSdict begin /ShrBORDERthickness 0 def}}
\def\mytag#1{%
    \tag#1}
\def\mythetag#1{\thetag{\blue{#1}}\immediate\special{ps:
     ShrHPSdict begin /ShrBORDERthickness 0 def}}
\def\mythetagchapter#1#2{\thetag{\blue{#1}}\immediate\special{ps:
     ShrHPSdict begin /ShrBORDERthickness 0 def}}
\def\myhref#1#2{\blue{#2}\immediate\special{ps:
     ShrHPSdict begin /ShrBORDERthickness 0 def}}
\def\myEarXivlink{\myhref{http://arXiv.org}{http:/\negskp/arXiv.org}}
\def\mysection#1{%
    #1.}

\def\mythesectionchaptertext#1#2#3{\blue{#3}\immediate\special{ps:
     ShrHPSdict begin /ShrBORDERthickness 0 def}}
\def\mydefinition#1{\definition{Definition #1}}
\def\mythedefinition#1{\blue{#1}\immediate\special{ps:
     ShrHPSdict begin /ShrBORDERthickness 0 def}}
\def\mythedefinitionchapter#1#2{\blue{#1}\immediate\special{ps:
     ShrHPSdict begin /ShrBORDERthickness 0 def}}
\def\mytheorem#1{\csname proclaim\endcsname{Theorem #1}}
\def\mythetheorem#1{\blue{#1}\immediate\special{ps:
     ShrHPSdict begin /ShrBORDERthickness 0 def}}
\def\mythetheoremchapter#1#2{\blue{#1}\immediate\special{ps:
     ShrHPSdict begin /ShrBORDERthickness 0 def}}
\def\mylemma#1{\csname proclaim\endcsname{Lemma #1}}
\def\mythelemma#1{\blue{#1}\immediate\special{ps:
     ShrHPSdict begin /ShrBORDERthickness 0 def}}
\def\mythelemmachapter#1#2{\blue{#1}\immediate\special{ps:
     ShrHPSdict begin /ShrBORDERthickness 0 def}}
\def\myref#1#2{\line{\vtop{\hsize 11pt\noindent
#1.}
\vtop{\advance\hsize -10pt\noindent #2}\hss}\vskip 4pt}
\pagewidth{360pt}
\pageheight{606pt}
\loadbold
\TagsOnRight
\document
\vbox to\vsize{

\centerline{\etbf RUSSIAN FEDERAL COMMITTEE}
\centerline{\etbf FOR HIGHER EDUCATION}
\bigskip
\centerline{\etbf BASHKIR STATE UNIVERSITY}
\vskip 3cm
\centerline{SHARIPOV\ R.\,A.}
\vskip 1.5cm
\centerline{\tvbf COURSE \ OF \ DIFFERENTIAL \ GEOMETRY}
\vskip 1.3cm
\centerline{\tvrm The Textbook}
\vfill
\centerline{Ufa 1996}

}
\vbox to\vsize{
MSC 97U20\par
UDC 514.7\par
\medskip
Sharipov R. A. {\bf Course of Differential Geometry}: the textbook / 
Publ\. of Bashkir State University
--- Ufa, 1996. --- pp\.~132. --- ISBN 5-7477-0129-0.
\bigskip
\bigskip
This book is a textbook for the basic course of differential
geometry. It is recommended as an introductory material for 
this subject.\par
     In preparing Russian edition of this book I used the computer
typesetting on the base of the \AmSTeX\ package and I used Cyrillic
fonts of the Lh-family distributed by the CyrTUG association of Cyrillic
\TeX\ users. English edition of this book is also typeset by means of
the \AmSTeX\ package.\par
\medskip
Referees:\ \ \ \
\vtop{\hsize 9.5cm\noindent Mathematics group of Ufa State University for 
Aircraft and Technology ({\tencyr UGATU});
\vskip 0.1cm
\noindent Prof\.~V.~V.~Sokolov, Mathematical Institute of Ural Branch
of Russian Academy of Sciences ({\tencyr IM UrO RAN}).}
\medskip
\noindent {\bf Contacts to author}.
\medskip
\line{\vtop to 150pt{\hsize=300pt\settabs\+\indent Office:\ &\cr
\+ Office:\hss &Mathematics Department, Bashkir State University,\cr
\+\hss &32 Frunze street, 450074 Ufa, Russia\cr
\+ Phone:\hss &7-(3472)-23-67-18\cr
\+ Fax:\hss   &7-(3472)-23-67-74\cr
\medskip
\+ Home:\hss &5 Rabochaya street, 450003 Ufa, Russia\cr
\+ Phone:\hss &7-(917)-75-55-786\cr
\+ E-mails:\hss &\myhref{mailto:R_Sharipov\@ic.bashedu.ru}
   {R\_\hskip 1pt Sharipov\@ic.bashedu.ru}\cr
\+\hss &\myhref{mailto:r-sharipov\@mail.ru}{r-sharipov\@mail.ru}\cr
\+\hss &\myhref{mailto:ra_sharipov\@lycos.com}{ra\_\hskip 1pt
   sharipov\@lycos.com}\cr
\+\hss &\myhref{mailto:ra_sharipov\@hotmail.com}{ra\_\hskip 1pt
   sharipov\@hotmail.com}\cr
\+ URL:\hss &\myhref{http://www.geocities.com/r-sharipov}{http:/\negskp
   /www.geocities.com/r-sharipov}\cr
\vfil}\hss}
\vskip -100pt
\vfil
\line{\vtop{\hsize=300pt\settabs\+\indent\kern 180pt &\cr
\+ISBN 5-7477-0129-0\hss
  &\copyright\ Sharipov R.A., 1996\cr
\+English translation\hss &\copyright\ Sharipov R.A., 2004\cr}\hss}
\vskip 1pt plus 1pt minus 1pt}
\topmatter
\title
CONTENTS.
\endtitle
\endtopmatter
\document
\vskip 30pt
\line{CONTENTS.\ \leaderfill\ 3.}
\setfirstpage
\medskip
\line{PREFACE.\ \leaderfill\ 5.}
\medskip
\line{CHAPTER \uppercase\expandafter{\romannumeral 1}.
CURVES IN THREE-DIMENSIONAL SPACE.\ \leaderfill\ 6.}
\medskip
\line{\S~1. Curves. Methods of defining a curve. Regular and 
singular points\hss}
\line{\qquad of a curve.\ \leaderfill\ \mythesectionchaptertext{1}{1}{6.}}
\line{\S~2. The length integral and the natural parametrization 
of a curve.\ \leaderfill\ \mythesectionchaptertext{2}{1}{10.}}
\line{\S~3. Frenet frame. The dynamics of Frenet frame.
Curvature and torsion\hss}
\line{\qquad of a spacial curve.\
\leaderfill\ \mythesectionchaptertext{3}{1}{12.}}
\line{\S~4. The curvature center and the curvature radius of
a spacial curve.\hss}
\line{\qquad  The evolute and the evolvent of a curve.\ \leaderfill\ 
\mythesectionchaptertext{4}{1}{14.}}
\line{\S~5. Curves as trajectories of material points
in mechanics.\ \leaderfill\ \mythesectionchaptertext{5}{1}{16.}}
\medskip
\line{CHAPTER \uppercase\expandafter{\romannumeral 2}.
ELEMENTS OF VECTORIAL\hss}
\line{\qquad AND TENSORIAL ANALYSIS.\ \leaderfill\ 18.}
\medskip
\line{\S~1. Vectorial and tensorial fields in the space.\ \leaderfill\ 
\mythesectionchaptertext{1}{2}{18.}}
\line{\S~2. Tensor product and contraction.\ \leaderfill\ 
\mythesectionchaptertext{2}{2}{20.}}
\line{\S~3. The algebra of tensor fields.\ \leaderfill\ 
\mythesectionchaptertext{3}{2}{24.}}
\line{\S~4. Symmetrization and alternation.\ \leaderfill\ 
\mythesectionchaptertext{4}{2}{26.}}
\line{\S~5. Differentiation of tensor fields.\ \leaderfill\ 
\mythesectionchaptertext{5}{2}{28.}}
\line{\S~6. The metric tensor and the volume pseudotensor.\ \leaderfill\ \mythesectionchaptertext{6}{2}{31.}}
\line{\S~7. The properties of pseudotensors.\ \leaderfill\ 
\mythesectionchaptertext{7}{2}{34.}}
\line{\S~8. A note on the orientation.\ \leaderfill\ 
\mythesectionchaptertext{8}{2}{35.}}
\line{\S~9. Raising and lowering indices.\ \leaderfill\ 
\mythesectionchaptertext{9}{2}{36.}} 
\line{\S~10. Gradient, divergency and rotor. Some identities\hss} 
\line{\qquad\ \ of the vectorial analysis. \ \leaderfill\ 
\mythesectionchaptertext{10}{2}{38.}}
\line{\S~11. Potential and vorticular vector fields.\ \leaderfill\ 
\mythesectionchaptertext{11}{2}{41.}}
\medskip
\line{CHAPTER \uppercase\expandafter{\romannumeral 3}.
CURVILINEAR COORDINATES.\ \leaderfill\ 45.}
\medskip
\line{\S~1. Some examples of curvilinear coordinate systems.\ 
\leaderfill\ \mythesectionchaptertext{1}{3}{45.}}
\line{\S~2. Moving frame of a curvilinear coordinate system.\ \leaderfill\ \mythesectionchaptertext{2}{3}{48.}}
\line{\S~3. Change of curvilinear coordinates.\ \leaderfill\ 
\mythesectionchaptertext{3}{3}{52.}}
\line{\S~4. Vectorial and tensorial fields
in curvilinear coordinates.\ \leaderfill\
\mythesectionchaptertext{4}{3}{55.}}
\line{\S~5. Differentiation of tensor fields
in curvilinear coordinates.\ \leaderfill\ 
\mythesectionchaptertext{5}{3}{57.}}
\line{\S~6. Transformation of the connection components\hss}
\line{\qquad under a change of a coordinate system.\ \leaderfill\ 
\mythesectionchaptertext{6}{3}{62.}}
\line{\S~7. Concordance of metric and connection. Another 
formula\hss}
\line{\qquad for Christoffel symbols.\ \leaderfill\ 
\mythesectionchaptertext{7}{3}{63.}}
\line{\S~8. Parallel translation. The equation of a straight line\hss}
\line{\qquad in curvilinear coordinates.\ \leaderfill\ 
\mythesectionchaptertext{8}{3}{65.}}
\line{\S~9. Some calculations in polar, cylindrical, and spherical
coordinates.\ \leaderfill\ \mythesectionchaptertext{9}{3}{70.}}
\medskip\newpage
\line{CHAPTER \uppercase\expandafter{\romannumeral 4}.
GEOMETRY OF SURFACES.\ \leaderfill\ 74.}
\medskip
\line{\S~1. Parametric surfaces. Curvilinear coordinates on a surface.\
\leaderfill\ \mythesectionchaptertext{1}{4}{74.}}
\line{\S~2. Change of curvilinear coordinates on a surface.
\ \leaderfill\ \mythesectionchaptertext{2}{4}{78.}}
\line{\S~3. The metric tensor and the area tensor.\ \leaderfill\ 
\mythesectionchaptertext{3}{4}{80.}}
\line{\S~4. Moving frame of a surface. Veingarten's derivational 
formulas.\ \leaderfill\ \mythesectionchaptertext{4}{4}{82.}}
\line{\S~5. Christoffel symbols and the second quadratic form.\ 
\leaderfill\ \mythesectionchaptertext{5}{4}{84.}}
\line{\S~6. Covariant differentiation of inner tensorial fields
of a surface.\ \leaderfill\ \mythesectionchaptertext{6}{4}{88.}}
\line{\S~7. Concordance of metric and connection on a surface.
\ \leaderfill\ \mythesectionchaptertext{7}{4}{94.}}
\line{\S~8. Curvature tensor.\ \leaderfill\ 
\mythesectionchaptertext{8}{4}{97.}}
\line{\S~9. Gauss equation and Peterson-Codazzi equation.\ \leaderfill\
\mythesectionchaptertext{9}{4}{103.}}
\medskip
\line{CHAPTER \uppercase\expandafter{\romannumeral 5}.
CURVES ON SURFACES.\ \leaderfill\ 106.}
\medskip
\line{\S~1. Parametric equations of a curve on a surface.\ 
\leaderfill\ \mythesectionchaptertext{1}{5}{106.}}
\line{\S~2. Geodesic and normal curvatures of a curve.\ \leaderfill\ 
\mythesectionchaptertext{2}{5}{107.}}
\line{\S~3. Extremal property of geodesic lines.\ \leaderfill\ 
\mythesectionchaptertext{3}{5}{110.}}
\line{\S~4. Inner parallel translation on a surface.\ \leaderfill\ 
\mythesectionchaptertext{4}{5}{114.}}
\line{\S~5. Integration on surfaces. Green's formula. \ \leaderfill\ 
\mythesectionchaptertext{5}{5}{120.}}
\line{\S~6. Gauss-Bonnet theorem.\ \leaderfill\ 
\mythesectionchaptertext{6}{5}{124.}}
\medskip
\line{REFERENCES.\ \leaderfill\ 132.}
\medskip
\newpage
\topmatter
\title
PREFACE.
\endtitle
\endtopmatter
\document
\setfirstpage
    This book was planned as the third book in the series of three
textbooks for three basic geometric disciplines of the university 
education. These are 
\roster
\item"--" {\tencyr\char '074}Course of analytical geometry\footnotemark%
{\tencyr\char '076};
\item"--" {\tencyr\char '074}Course of linear algebra and multidimensional
geometry{\tencyr\char '076};
\item"--" {\tencyr\char '074}Course of differential 
geometry{\tencyr\char '076}.
\endroster
\footnotetext{\ Russian versions of the second and the third books were
written in 1096, but the first book is not yet written. I understand it
as my duty to complete the series, but I had not enough time all these
years since 1996.}
\adjustfootnotemark{-1}
This book is devoted to the first acquaintance with the differential
geometry. Therefore it begins with the theory of curves in
three-dimensional Euclidean space $\Bbb E$. Then the vectorial analysis
in $\Bbb E$ is stated both in Cartesian and curvilinear coordinates,
afterward the theory of surfaces in the space $\Bbb E$ is considered.
\par
     The newly fashionable approach starting with the concept of a
differentiable manifold, to my opinion, is not suitable for the
introduction to the subject. In this way too many efforts are spent 
for to assimilate this rather abstract notion and the rather special
methods associated with it, while the the essential content of the 
subject is postponed for a later time. I think it is more important 
to make faster acquaintance with other elements of modern geometry
such as the vectorial and tensorial analysis, covariant differentiation,
and the theory of Riemannian curvature. The restriction of the
dimension to the cases $n=2$ and $n=3$ is not an essential obstacle
for this purpose. The further passage from surfaces to 
higher-dimensional manifolds becomes more natural and simple.\par
     I am grateful to D.~N.~Karbushev, R.~R.~Bakhitov,
S.~Yu.~Ubiyko, D.~I.~Borisov, and Yu.~N.~Polyakov for reading and
correcting the manuscript of the Russian edition of this book.
\bigskip\bigskip
\line{\vbox{\hsize 7.5cm\noindent November, 1996;\newline December,
2004.}\hss R.~A.~Sharipov.}
\newpage
\setfirstpage
\topmatter
\title\chapter{1}
CURVES IN THREE-DIMENSIONAL SPACE.
\endtitle
\endtopmatter
\chapternum=1
\document
\head
\S~\mysection{1} Curves. Methods of defining a curve. Regular and 
singular points of a curve.
\endhead
\rightheadtext{\S~1. Curves. Methods of defining a curve \dots}
\leftheadtext{CHAPTER \uppercase\expandafter{\romannumeral 1}.
CURVES IN THREE-DIMENSIONAL SPACE.}
     Let $\Bbb E$ be a three-dimensional Euclidean point space. 
The strict mathematical definition of such a space can be found 
in \mycite{1}. However, knowing this definition is not so urgent. 
The matter is that $\Bbb E$ can be understood as the regular 
three-dimensional space (that in which we live). The properties
of the space $\Bbb E$ are studied in elementary mathematics and
in analytical geometry on the base intuitively clear visual
forms. The concept of a {\it line} or a {\it curve} is also
related to some visual form. A {\it curve} in the space $\Bbb E$ 
is a spatially extended one-dimensional geometric form. The
one-dimensionality of a curve reveals when we use the
{\it vectorial-parametric} method of defining it:
$$
\hskip -2em
\bold r=\bold r(t)=\Vmatrix x^1(t)\\x^2(t)\\x^3(t)\endVmatrix.
\mytag{1.1}
$$
We have one degree of freedom when choosing a point on the curve
\mythetag{1.1}, our choice is determined by the value of the numeric 
parameter $t$ taken from some interval, e\.\,g\. from the unit 
interval $[0,\,1]$ on the real axis $\Bbb R$. Points of the curve
\mythetag{1.1} are given by their radius-vectors\footnotemark\ 
$\bold r=\bold r(t)$ whose components $x^1(t)$, $x^2(t)$, $x^3(t)$
are functions of the parameter $t$.
\par
\footnotetext{\ Here we assume that some Cartesian coordinate system in
$\Bbb E$ is taken.}
\adjustfootnotemark{-1}
     The continuity of the curve \mythetag{1.1} means that the functions
$x^1(t)$, $x^2(t)$, $x^3(t)$ should be continuous. However, this 
condition is too weak. Among continuous curves there are some instances
which do not agree with our intuitive understanding of a curve. In the
course of mathematical analysis the Peano curve is often considered 
as an example (see \mycite{2}). This is a continuous parametric curve on 
a plane such that it is enclosed within a unit square, has no self
intersections, and passes through each point of this square. In order
to avoid such unusual curves the functions $x^i(t)$ in \mythetag{1.1}
are assumed to be continuously differentiable ($C^1$ class) functions or, 
at least, piecewise continuously differentiable functions.\par
     Now let's consider another method of defining a curve. An arbitrary
point of the space $\Bbb E$ is given by three arbitrary parameters
$x^1$, $x^2$, $x^3$ --- its coordinates. We can restrict the degree of
arbitrariness by considering a set of points whose coordinates 
$x^1,\,x^2,\,x^3$ satisfy an equation of the form
$$
\pagebreak
\hskip -2em
F(x^1,x^2,x^3)=0,
\mytag{1.2}
$$
where $F$ is some continuously differentiable function of three
variables. In a typical situation formula \mythetag{1.2} still admits 
two-parametric arbitrariness: choosing arbitrarily two coordinates 
of a point, we can determine its third coordinate by solving the
equation \mythetag{1.2}. Therefore, \mythetag{1.2} is an equation of
a surface. In the intersection of two surfaces usually a curve 
arises. Hence, a system of two equations of the form \mythetag{1.2}
defines a curve in $\Bbb E$:
$$
\hskip -2em
\cases
F(x^1,x^2,x^3)=0,\\
G(x^1,x^2,x^3)=0.
\endcases
\mytag{1.3}
$$
If a curve lies on a plane, we say that it is a {\it plane} curve.
For a plane curve one of the equations \mythetag{1.3} can be replaced
by the equation of a plane: $A\,x^1+B\,x^2+C\,x^3+D=0$.\par
     Suppose that a curve is given by the equations \mythetag{1.3}. 
Let's choose one of the variables $x^1$, $x^2$, or $x^3$ for a parameter,
e\.\,g\. we can take $x^1=t$ to make certain. Then, writing the system of
the equations \mythetag{1.3} as
$$
\cases
F(t,x^2,x^3)=0,\\
G(t,x^2,x^3)=0,
\endcases
$$
and solving them with respect to $x^2$ and $x^3$, we get two 
functions $x^2(t)$ and $x^3(t)$. Hence, the same curve can be given 
in vectorial-parametric form:
$$
\bold r=\bold r(t)=\Vmatrix t\\x^2(t)\\x^3(t)\endVmatrix.
$$
\par
     Conversely, assume that a curve is initially given in
vectorial-parametric form by means of vector-functions \mythetag{1.1}. 
Then, using the functions $x^1(t)$, $x^2(t)$,  $x^3(t)$, we construct 
the following two systems of equations:
$$
\xalignat 2
&\hskip -2em
\cases
x^1-x^1(t)=0,\\
x^2-x^2(t)=0,
\endcases
&&\cases
x^1-x^1(t)=0,\\
x^3-x^3(t)=0.
\endcases
\mytag{1.4}
\endxalignat
$$
Excluding the parameter $t$ from the first system of equations
\mythetag{1.4}, we obtain some functional relation for two variable 
$x^1$ and $x^2$. We can write it as $F(x^1,x^2)=0$. Similarly, the 
second system reduces to the equation $G(x^1,x^3)=0$. Both these 
equations constitute a system, which is a special instance of
\mythetag{1.3}:
$$
\cases
F(x^1,x^2)=0,\\
G(x^1,x^3)=0.
\endcases
$$
This means that the vectorial-parametric representation of a curve can 
be transformed to the form of a system of equations \mythetag{1.3}.\par
\parshape 28 0pt 360pt 0pt 360pt 0pt 360pt 0pt 360pt 0pt 360pt
0pt 360pt 0pt 360pt 0pt 360pt 0pt 360pt 180pt 180pt 180pt 180pt 
180pt 180pt 180pt 180pt 180pt 180pt 180pt 180pt 180pt 180pt 180pt 180pt 
180pt 180pt 180pt 180pt 180pt 180pt 180pt 180pt 180pt 180pt 180pt 180pt 
180pt 180pt 180pt 180pt 180pt 180pt 180pt 180pt 
0pt 360pt 
     None of the above two methods of defining a curve in $\Bbb E$
is absolutely preferable. In some cases the first method is 
better, in other cases the second one is used. However, for constructing
the theory of curves the vectorial-parametric method is more suitable.
Suppose that we have a parametric curve $\gamma$ of the {\it smoothness
class\/} $C^1$. This is a curve with the coordinate functions $x^1(t)$,
$x^2(t)$, $x^3(t)$ being continuously differentiable. Let's choose two 
different values of the parameter: $t$ and $\tilde t=t+\triangle t$, where
$\triangle t$ is an increment of the parameter. Let $A$ and $B$ be two
points on the curve corresponding to that two values of the parameter
$t$. We draw the straight line passing through these points $A$ and $B$; 
this is a secant for the curve $\gamma$. Directing vectors of this
secant are collinear to the vector $\overrightarrow{AB\ }$. We choose one
of them:
$$
\bold a=\frac{\overrightarrow{AB\ }}{\triangle t}=
\frac{\bold r(t+\triangle t)-\bold r(t)}{\triangle t}.\qquad
\mytag{1.5}
$$
Tending $\triangle t$ to zero, we find that the point $B$ moves toward
the point $A$. Then the secant tends to its limit position and becomes
the tangent line of the curve at the point $A$. Therefore limit value of
the vector \mythetag{1.5} is a tangent vector of the curve $\gamma$ at the
point $A$:
$$
\boldsymbol\tau(t)=\lim_{\triangle t\to\infty}\bold a=
\frac{d\bold r(t)}{dt}=\dot\bold r(t).\quad
\mytag{1.6}
$$
The components of the tangent vector \mythetag{1.6} are evaluated by 
differentiating the components of the radius-vector $\bold r(t)$ with
respect to the variable $t$.
\vadjust{\vskip 0pt\hbox to 0pt{\kern -10pt
\includegraphics{ris01.eps}\hss}\vskip 0pt}\par
     The tangent vector $\dot\bold r(t)$ determines the direction of the 
instantaneous displacement of the point $\bold r(t)$ for the given value
of the parameter $t$. Those points, at which the derivative $\dot\bold r(t)$ vanishes, are
special ones. They are {\tencyr\char '074}stopping points{\tencyr\char '076}.
Upon stopping, the point can begin moving in quite different direction.
For example, let's consider the following two plane curves:
$$
\xalignat 2
&\hskip -2em
\bold r(t)=\Vmatrix t^2\\t^3\endVmatrix,
&&\bold r(t)=\Vmatrix t^4\\t^3\endVmatrix.
\mytag{1.7}
\endxalignat
$$
At $t=0$ both curves \mythetag{1.7} pass through the origin and the tangent
vectors of both curves at the origin are equal to zero. However, the
behavior of these curves near the origin is quite different:
\vadjust{\vskip 5pt\hbox to 0pt{\kern -10pt
\includegraphics{ris02.eps}\hss}\vskip 157pt}
\noindent
the first curve has a beak-like fracture at the origin, while the
second one is smooth. Therefore, vanishing of the derivative
$$
\hskip -2em
\boldsymbol\tau(t)=\dot\bold r(t)=0
\mytag{1.8}
$$
is only the necessary, but not sufficient condition for a parametric curve 
to have a singularity at the point $\bold r(t)$. The opposite condition
$$
\hskip -2em
\boldsymbol\tau(t)=\dot\bold r(t)\neq 0
\mytag{1.9}
$$
guaranties that the point $\bold r(t)$ is free of singularities.
Therefore, those points of a parametric curve, where the condition
\mythetag{1.9} is fulfilled, are called {\it regular points}.\par
     Let's study the problem of separating regular and singular 
points on a curve given by a system of equations \mythetag{1.3}. Let
$A=(a^1,a^2,a^3)$ be a point of such a curve. The functions 
$F(x^1,x^2,x^3)$ and $G(x^1,x^2,x^3)$ in \mythetag{1.3} are assumed
to be continuously differentiable. The matrix
$$
\hskip -2em
J=\left\Vert
\vrule height 5.5ex depth 6.0ex width 0ex
\matrix
\dsize\frac{\partial F}{\partial x^1} &
\dsize\frac{\partial F}{\partial x^2} &
\dsize\frac{\partial F}{\partial x^3} \\
\vspace{1.5ex}
\dsize\frac{\partial G}{\partial x^1} &
\dsize\frac{\partial G}{\partial x^2} &
\dsize\frac{\partial G}{\partial x^3}
\endmatrix
\right\Vert
\mytag{1.10}
$$
composed of partial derivatives of $F$ and $G$ at the point $A$
is called the {\it Jacobi} matrix or the {\it Jacobian}  of the
system of equations \mythetag{1.3}. If the minor
$$
M_1=\det\left\vert
\vrule height 5.5ex depth 6.0ex width 0ex
\matrix
\dsize\frac{\partial F}{\partial x^2} &
\dsize\frac{\partial F}{\partial x^3} \\
\vspace{1.5ex}
\dsize\frac{\partial G}{\partial x^2} &
\dsize\frac{\partial G}{\partial x^3}
\endmatrix
\right\vert\neq 0
$$
in Jacobi matrix is nonzero, the equations \mythetag{1.3} can be
resolved with respect to $x^2$ and $x^3$ in some neighborhood of the
point $A$. Then we have three functions $x^1=t$, $x^2=x^2(t)$, 
$x^3=x^3(t)$ which determine the parametric representation of our 
curve. This fact follows from the theorem {\it on implicit functions} 
(see \mycite{2}). Note that the tangent vector of the curve in this 
parametrization
$$
\boldsymbol\tau=\Vmatrix 1\\ \dot x^2\\ \dot x^3\endVmatrix
\neq 0
$$
is nonzero because of its first component. This means that the
condition $M_1\neq 0$ is sufficient for the point $A$ to be a 
regular point of a curve given by the system of equations
\mythetag{1.3}. Remember that the Jacobi matrix \mythetag{1.10}
has two other minors:
$$
\pagebreak
\xalignat 2
&M_2=\det\left\vert
\vrule height 5.5ex depth 6.0ex width 0ex
\matrix
\dsize\frac{\partial F}{\partial x^3} &
\dsize\frac{\partial F}{\partial x^1} \\
\vspace{1.5ex}
\dsize\frac{\partial G}{\partial x^3} &
\dsize\frac{\partial G}{\partial x^1}
\endmatrix
\right\vert,
&&M_3=\det\left\vert
\vrule height 5.5ex depth 6.0ex width 0ex
\matrix
\dsize\frac{\partial F}{\partial x^1} &
\dsize\frac{\partial F}{\partial x^2} \\
\vspace{1.5ex}
\dsize\frac{\partial G}{\partial x^1} &
\dsize\frac{\partial G}{\partial x^2}
\endmatrix
\right\vert.
\endxalignat
$$
For both of them the similar propositions are fulfilled. Therefore,
we can formulate the following theorem.
\mytheorem{1.1} A curve given by a system of equations
\mythetag{1.3} is regular at all points, where the rank of its Jacobi 
matrix \mythetag{1.10} is equal to $2$.
\endproclaim
     A plane curve lying on the plane $x^3=0$ can be defined by
one equation $F(x^1,x^2)=0$. The second equation here reduces to
$x^3=0$. Therefore, $G(x^1,x^2,x^3)=x^3$. The Jacoby matrix for
the system \mythetag{1.3} in this case is 
$$
\hskip -2em
J=\left\Vert
\vrule height 5.0ex depth 3.0ex width 0ex
\matrix
\dsize\frac{\partial F}{\partial x^1} &
\dsize\frac{\partial F}{\partial x^2} &
0\\
\vspace{1.5ex}
0 & 0 & 1
\endmatrix
\right\Vert.
\mytag{1.11}
$$
If $\rank J=2$, this means that at least one of two partial derivatives
in the matrix \mythetag{1.11} is nonzero. These derivatives form the
{\it gradient} vector for the function $F$:
$$
\grad F=\left(\frac{\partial F}{\partial x^1},\
\frac{\partial F}{\partial x^2}\right).
$$
\mytheorem{1.2} A plane curve given by an equation
$F(x^1,x^2)=0$ is regular at all points where $\grad F\neq 0$.
\endproclaim
     This theorem~\mythetheorem{1.2} is a simple corollary from the
theorem~\mythetheorem{1.1} and the relationship \mythetag{1.11}. Note 
that the theorems~\mythetheorem{1.1} and \mythetheorem{1.2} yield only
sufficient conditions for regularity of curve points. Therefore, some
points where these theorems are not applicable can also be regular 
points of a curve.\par
\head
\S~\mysection{2} The length integral\\and the natural parametrization 
of a curve.
\endhead
\rightheadtext{\S~2. The length integral \dots}
     Let $\bold r=\bold r(t)$ be a parametric curve of smoothness
class $C^1$, where the parameter $t$ runs over the interval $[a,\,b]$.
Let's consider a monotonic increasing continuously differentiable
function $\varphi(\tilde t\kern 1pt)$ on a segment $[\tilde a,\,\tilde 
b]$ such that $\varphi(\tilde a)=a$ and $\varphi(\tilde b)=b$. Then it
takes each value from the segment $[a,\,b]$ exactly once. Substituting
$t=\varphi(\tilde t\kern 1pt)$ into $\bold r(t)$, we define the new
vector-function $\tilde\bold r(\tilde t\kern 1pt)=\bold r(\varphi(
\tilde t\kern 1pt))$, it describes the same curve as the original
vector-function $\bold r(t)$. This procedure is called the {\it
reparametrization} of a curve. We can calculate the tangent vector
in the new parametrization by means of the chain rule:
$$
\hskip -2em
\tilde{\boldsymbol\tau}(\tilde t\kern 1pt)=
\varphi'(\tilde t\kern 1pt)\cdot
\boldsymbol\tau(\varphi(\tilde t\kern 1pt)).
\mytag{2.1}
$$
Here $\varphi'(\tilde t\kern 1pt)$ is the derivative of the function
$\varphi(\tilde t\kern 1pt)$. The formula \mythetag{2.1} is known as 
the {\it transformation rule} for the tangent vector of a curve 
under a change of parametrization.\par
     A monotonic decreasing function $\varphi(\tilde t\kern 1pt)$ can 
also be used for the reparametrization of curves. In this case
$\varphi(\tilde a)=b$ and $\varphi(\tilde b)=a$, i\.\,e\. the beginning
point and the ending point of a curve are exchanged. Such 
reparametrizations are called {\it changing the orientation of a curve}.
\par
      From the formula \mythetag{2.1}, we see that the tangent vector
$\tilde{\boldsymbol\tau}(\tilde t\kern 1pt)$ can vanish at some points
of the curve due to the derivative $\varphi'(\tilde t\kern 1pt)$ even 
when $\boldsymbol\tau(\varphi(\tilde t\kern 1pt))$ is nonzero.
Certainly, such points are not actually the singular points of a curve.
In order to exclude such formal singularities, only those reparametrizations
of a curve are admitted for which the function $\varphi(\tilde t\kern 1pt)$
is a strictly monotonic function, i\.\,e\. $\varphi'(\tilde t\kern 1pt)>0$
or $\varphi'(\tilde t\kern 1pt)<0$.\par
      The formula \mythetag{2.1} means that the tangent vector of a curve
at its regular point depends not only on the geometry of the curve, but
also on its parametrization. However, the effect of parametrization is not
so big, it can yield a numeric factor to the vector $\boldsymbol\tau$ only.
Therefore, the natural question arises: is there some preferable
parametrization on a curve\,? The answer to this question is given by the
length integral.\par
\parshape 26 0pt 360pt 0pt 360pt 0pt 360pt 0pt 360pt 0pt 360pt 0pt 360pt 
180pt 180pt 180pt 180pt 180pt 180pt 180pt 180pt 180pt 180pt 180pt 180pt 
180pt 180pt 180pt 180pt 180pt 180pt 180pt 180pt 180pt 180pt 180pt 180pt 
180pt 180pt 180pt 180pt 180pt 180pt 180pt 180pt 180pt 180pt 180pt 180pt 
180pt 180pt 0pt 360pt 
     Let's consider a segment of a parametric curve of the smoothness class
$C^1$ with the parameter $t$ running over the segment $[a,\,b]$ of real numbers.
Let 
$$
\hskip -2em
a=t_0<t_1<\ldots<t_n=b
\mytag{2.2}
$$
be a series of points breaking this segment into $n$ parts. The points
$\bold r(t_0),\ldots,\,\bold r(t_n)$ on the curve define a polygonal line
with $n$ segments. Denote $\triangle t_k=t_k-t_{k-1}$ and let $\varepsilon$
be the maximum of $\triangle t_k$:
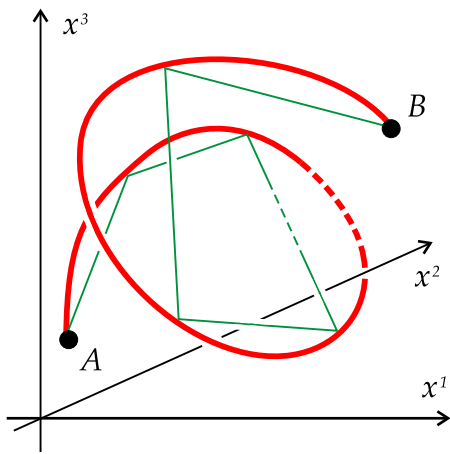
$$
\varepsilon=\max_{k=1,\,\ldots,\,n}\triangle t_k.
$$
The quantity $\varepsilon$ is the fineness of the partition
\mythetag{2.2}. The length of $k$-th segment of the polygonal line $AB$ 
is calculated by the formula $L_k=|\bold r(t_k)-\bold r(t_{k-1})|$. Using
the continuous differentiability of the vector-function $\bold r(t)$, from
the Taylor expansion of $\bold r(t)$ at the point $t_{k-1}$ we
get $L_k=|\boldsymbol\tau(t_{k-1})|\cdot\triangle t_k+o(\varepsilon)$.
\vadjust{\vskip 80pt\hbox to 0pt{\kern -30pt
\includegraphics{ris03.eps}\hss}\vskip -80pt}There\-fore, as the
fineness $\varepsilon$ of the partition \mythetag{2.2} tends to zero, the
length of the polygonal line $AB$ has the limit equal to the integral
of the modulus of tangent vector $\boldsymbol\tau(t)$ along the curve:
$$
\hskip -2em
L=\lim_{\varepsilon\to 0}\,\sum^n_{k=1}L_k=
\int\limits^{\ \,b}_{\!a}|\boldsymbol\tau(t)|\,dt.
\mytag{2.3}
$$
It is natural to take the quantity $L$ in \mythetag{2.3} for the length
of the curve $AB$. Note that if we reparametrize a curve according to the 
formula \mythetag{2.1}, this leads to a change of variable in the integral.
Nevertheless, the value of the integral $L$ remains unchanged. Hence, the
length of a curve is its geometric invariant which does not depend on the
way how it is parameterized.\par
     The length integral \mythetag{2.3} defines the preferable way for
parameterizing a curve in the Euclidean space $\Bbb E$. Let's denote by
$s(t)$ an antiderivative \pagebreak of the function $\psi(t)
=|\boldsymbol\tau(t)|$ being under integration in the formula
\mythetag{2.3}:
$$
\hskip -2em
s(t)=\int\limits^t_{t_0}|\boldsymbol\tau(t)|\,dt.
\mytag{2.4}
$$
\mydefinition{2.1} The quantity $s$ determined by the
integral \mythetag{2.4} is called the {\it natural parameter} of a curve in the Euclidean space $\Bbb E$.
\enddefinition
     Note that once the reference point $\bold r(t_0)$ and some direction
(orientation) on a curve have been chosen, the value of natural parameter
depends on the point of the curve only. Then the change of $s$ for $-s$ 
means the change of orientation of the curve for the opposite one.\par
     Let's differentiate the integral \mythetag{2.4} with respect to its 
upper limit $t$. As a result we obtain the following relationship:
$$
\hskip -2em
\frac{ds}{dt}=|\boldsymbol\tau(t)|.
\mytag{2.5}
$$
Now, using the formula \mythetag{2.5},
we can calculate the tangent vector of a curve in its natural
parametrization, i\.\,e\. when $s$ is used instead of $t$ as a parameter:
$$
\hskip -2em
\frac{d\bold r}{ds}=\frac{d\bold r}{dt}\cdot\frac{dt}{ds}=
\frac{d\bold r}{dt}\bigg/\frac{ds}{dt}=
\frac{\boldsymbol\tau}{|\boldsymbol\tau|}.
\mytag{2.6}
$$
From the formula \mythetag{2.6}, we see that in the tangent vector of
a curve in natural parametrization is a unit vector at all regular
points. In singular points this vector is not defined at all.
\head
\S~\mysection{3} Frenet frame. The dynamics of Frenet frame.
Curvature and torsion of a spacial curve.
\endhead
\rightheadtext{\S~3. Frenet frame. The dynamics of Frenet frame\dots}
     Let's consider a smooth parametric curve $\bold r(s)$ in natural
parametrization. The components of the radius-vector $\bold r(s)$ for such
a curve are smooth functions of $s$ (smoothness class $C^\infty$). They
are differentiable unlimitedly many times with respect to $s$.
The unit vector $\boldsymbol\tau(s)$ is obtained as the derivative
of $\bold r(s)$:
$$
\hskip -2em
\boldsymbol\tau(s)=\frac{d\bold r}{ds}.
\mytag{3.1}
$$
Let's differentiate the vector $\boldsymbol\tau(s)$ with respect to $s$ 
and then apply the following lemma to its derivative $\boldsymbol\tau'(s)$.
\mylemma{3.1} The derivative of a vector of a constant length is
a vector perpen\-dicular to the original one.
\endproclaim
\demo{Proof} In order to prove the lemma we choose some standard 
rectangular Cartesian coordinate system in $\Bbb E$. Then 
$$
|\boldsymbol\tau(s)|^2=(\boldsymbol\tau(s)\,|\,\boldsymbol\tau(s))=
(\tau^1)^2+(\tau^2)^2+(\tau^3)^2=\const.
$$
Let's differentiate this expression with respect to $s$. As a result 
we get the following relationship:
$$
\align
\frac{d}{ds}\left(|\boldsymbol\tau(s)|^2\right)&=
\frac{d}{ds}\left((\tau^1)^2+(\tau^2)^2+(\tau^3)^2\right)=\\
&=2\,\tau^1\,(\tau^1)'+2\,\tau^2\,(\tau^2)'+
2\,\tau^3\,(\tau^3)'=0.
\endalign
$$
One can easily see that this relationship is equivalent to 
$(\boldsymbol\tau(s)\,|\,\boldsymbol\tau'(s))=0$. Hence,
$\boldsymbol\tau(s)\perp\boldsymbol\tau'(s)$. The lemma is
proved.\qed\enddemo
     Due to the above lemma the vector $\boldsymbol\tau'(s)$
is perpendicular to the unit vector $\boldsymbol\tau(s)$. If
the length of $\boldsymbol\tau'(s)$ is nonzero, one can represent
it as 
$$
\hskip -2em
\boldsymbol\tau'(s)=k(s)\cdot\bold n(s),
\mytag{3.2}
$$
where $k(s)=|\boldsymbol\tau'(s)|$ and $|\bold n(s)|=1$. The scalar
quantity $k(s)=|\boldsymbol\tau'(s)|$ in formula \mythetag{3.2} is called
the {\it curvature} of a curve, while the unit vector $\bold n(s)$
is called its {\it primary normal vector} or simply the {\it normal
vector} of a curve at the point $\bold r(s)$.
The unit vectors $\boldsymbol\tau(s)$ and $\bold n(s)$ are orthogonal to 
each other. We can complement them by the third unit vector $\bold b(s)$
so that $\boldsymbol\tau$, $\bold n$, $\bold b$ become a {\it right 
triple\footnotemark}:
\footnotetext{\ A non-coplanar ordered triple of vectors $\bold a_1$,
$\bold a_2$, $\bold a_3$ is called a {\it right triple} if, upon moving
these vectors to a common origin, when looking from the end of the
third vector $\bold a_3$, we see the shortest rotation from $\bold a_1$ 
to $\bold a_2$ as a counterclockwise rotation.}
\adjustfootnotemark{-1}
$$
\bold b(s)=[\boldsymbol\tau(s),\,\bold n(s)].
\mytag{3.3}
$$
The vector $\bold b(s)$ defined by the formula \mythetag{3.3} is called 
the {\it secondary normal} vector or the {\it binormal} vector of
a curve. Vectors $\boldsymbol\tau(s)$, $\bold n(s)$, $\bold b(s)$ 
compose an orthonormal right basis attached to the point $\bold r(s)$.
\par
    Bases, which are attached to some points, are usually called 
{\it frames}. One should distinguish frames from coordinate systems.
Cartesian coordinate systems are also defined by choosing some point
(an origin) and some basis. However, coordinate systems are used for
describing the points of the space through their coordinates. The
purpose of frames is different. They are used for to expand the
vectors which, by their nature, are attached to the same points as
the vectors of the frame.\par
     The isolated frames are rarely considered, frames usually
arise within families of frames: typically at each point of some
set (a curve, a surface, or even the whole space) there arises some
frame attached to this point. The frame $\boldsymbol\tau(s)$, $\bold n(s)$, 
$\bold b(s)$ is an example of such frame. It is called the {\it Frenet
frame} of a curve. This is the moving frame: in typical situation the 
vectors of this frame change when we move the attachment point along the
curve.\par
     Let's consider the derivative $\bold n'(s)$. This vector attached
to the point $\bold r(s)$ can be expanded in the Frenet frame at that
point. Due to the lemma~\mythelemma{3.1} the vector $\bold n'(s)$ is
orthogonal to the vector $\bold n(s)$. Therefore its expansion has the 
form
$$
\hskip -2em
\bold n'(s)=\alpha\cdot\boldsymbol\tau(s)+\varkappa\cdot
\bold b(s).
\mytag{3.4}
$$
The quantity $\alpha$ in formula \mythetag{3.4} \pagebreak can be expressed
through the curvature of the curve. Indeed, as a result of the following
calculations we derive
$$
\hskip -2em
\aligned
\alpha(s)&=(\boldsymbol\tau(s)\,|\,\bold n'(s))=
(\boldsymbol\tau(s)\,|\,\bold n(s))'-\\
&-(\boldsymbol\tau'(s)\,|\,\bold n(s))=
-(k(s)\cdot\bold n(s)\,|\,\bold n(s))=-k(s).
\endaligned
\mytag{3.5}
$$
The quantity $\varkappa=\varkappa(s)$ cannot be expressed 
through the curvature. This is an additional parameter
characterizing a curve in the space $\Bbb E$. It is called
the {\it torsion} of the curve at the point $\bold r=\bold r(s)$.
The above expansion \mythetag{3.4} of the vector $\bold n'(s)$ now 
is written in the following form:
$$
\hskip -2em
\bold n'(s)=-k(s)\cdot\boldsymbol\tau(s)+\varkappa(s)
\cdot\bold b(s).
\mytag{3.6}
$$\par
     Let's consider the derivative of the binormal vector
$\bold b'(s)$. It is perpendicular to $\bold b(s)$. This derivative
can also be expanded in the Frenet frame. Due to 
$\bold b'(s)\perp \bold b(s)$ we have $\bold b'(s)=\beta\cdot\bold n(s)
+\gamma\cdot\boldsymbol
\tau(s)$. The coefficients $\beta$ and $\gamma$ in this expansion
can be found by means of the calculations similar to \mythetag{3.5}:
$$
\align
&\aligned
 \beta(s)&=(\bold n(s)\,|\,\bold b'(s))=
 (\bold n(s)\,|\,\bold b(s))'-(\bold n'(s)\,|\,\bold b(s))=\\
 &=-(-k(s)\cdot\boldsymbol\tau(s)+\varkappa(s)\cdot\bold b(s)\,
 |\,\bold b(s))=-\varkappa(s).
 \endaligned\\
\vspace{1ex}
&\aligned
 \gamma(s)&=(\boldsymbol\tau(s)\,|\,\bold b'(s))=
 (\boldsymbol\tau(s)\,|\,\bold b(s))'-(\boldsymbol\tau'(s)\,|
 \,\bold b(s))=\\
 &=-(k(s)\cdot\bold n(s)\,|\,\bold b(s))=0.
 \endaligned
\endalign
$$
Hence, for the expansion of the vector $\bold b'(s)$ in the Frenet frame
we get
$$
\hskip -2em
\bold b'(s)=-\varkappa(s)\cdot\bold n(s).
\mytag{3.7}
$$
Let's gather the equations \mythetag{3.2}, \mythetag{3.6}, and \mythetag{3.7} 
into a system:
$$
\cases
\boldsymbol\tau'(s)=k(s)\cdot\bold n(s),\\
\bold n'(s)=-k(s)\cdot\boldsymbol\tau(s)+\varkappa(s)
\cdot\bold b(s),\\
\bold b'(s)=-\varkappa(s)\cdot\bold n(s).
\endcases
\mytag{3.8}
$$
The equations \mythetag{3.8} relate the vectors $\boldsymbol\tau(s)$,
$\bold n(s)$, $\bold b(s)$ and their derivatives with respect to $s$. 
These differential equations describe the dynamics of the Frenet frame.
They are called the {\it Frenet equations}. The equations \mythetag{3.8}
should be complemented with the equation \mythetag{3.1} which describes
the dynamics of the point $\bold r(s)$ (the point to which the vectors
of the Frenet frame are attached).
\par
\head
\S~\mysection{4} The curvature center and the curvature radius\\ of
a spacial curve. The evolute and the evolvent of a curve.
\endhead
\rightheadtext{\S~4. The curvature center and the curvature radius \dots}
     In the case of a planar curve the vectors $\boldsymbol\tau(s)$ and
$\bold n(s)$ lie in the same plane as the curve itself. Therefore, 
binormal vector \mythetag{3.3} in this case coincides with the unit normal
vector of the plane. Its derivative $\bold b'(s)$ is equal to zero. Hence,
due to the third Frenet equation \mythetag{3.7} we find that for a planar
curve $\varkappa(s)\equiv 0$. The Frenet equations \mythetag{3.8} then are
reduced to
$$
\hskip -2em
\cases
\boldsymbol\tau'(s)=k(s)\cdot\bold n(s),\\
\bold n'(s)=-k(s)\cdot\boldsymbol\tau(s).
\endcases
\mytag{4.1}
$$
Let's consider the circle of the radius $R$ with the center at the origin
lying in the coordinate plane $x^3=0$. It is convenient to define this
circle as follows:
$$
\hskip -2em
\bold r(s)=
\Vmatrix
R\,\cos(s/R)\\
R\,\sin(s/R)
\endVmatrix,
\mytag{4.2}
$$
here $s$ is the natural parameter. Substituting \mythetag{4.2} into
\mythetag{3.1} and then into \mythetag{3.2}, we find the unit tangent vector
$\boldsymbol\tau(s)$ and the primary normal vector $\bold n(s)$:
$$
\xalignat 2
&\hskip -2em
\boldsymbol\tau(s)=
\Vmatrix
-\sin(s/R)\\
\ \ \cos(s/R)
\endVmatrix,
&&\bold n(s)=
\Vmatrix
-\cos(s/R)\\
-\sin(s/R)
\endVmatrix.
\mytag{4.3}
\endxalignat
$$
Now, substituting \mythetag{4.3} into the formula \mythetag{4.1}, we
calculate the curvature of a circle $k(s)=1/R=\const$. The curvature $k$
of a circle is constant, the inverse curvature $1/k$ coincides with its
radius.\par
     Let's make a step from the point $\bold r(s)$ on a circle to the
distance $1/k$ in the direction of its primary normal vector $\bold n(s)$.
It is easy to see that we come to the center of a circle. Let's make the 
same step for an arbitrary spacial curve. As a result of this step we come
from the initial point $\bold r(s)$ on the curve to the point with the 
following radius-vector:
$$
\hskip -2em
\boldsymbol\rho(s)=\bold r(s)+\frac{\bold n(s)}{k(s)}.
\mytag{4.4}
$$
Certainly, this can be done only for that points of a curve, where $k(s)
\neq 0$. The analogy with a circle induces the following terminology:
the quantity $R(s)=1/k(s)$ is called the {\it curvature radius}, the point 
with the radius-vector \mythetag{4.4} is called the {\it curvature center} 
of a curve at the point $\bold r(s)$.\par
     In the case of an arbitrary curve its curvature center is not a fixed
point. When parameter $s$ is varied, the curvature center of the curve
moves in the space drawing another  curve, which is called the {\it evolute}
of the original curve. The formula \mythetag{4.4} is a vectorial-parametric
equation of the evolute. However, note that the natural parameter $s$ of
the original curve is not a natural parameter for its evolute.\par
     Suppose that some spacial curve $\bold r(t)$ is given. A curve 
$\tilde\bold r(\tilde s)$ whose evolute $\tilde{\boldsymbol\rho}(\tilde s)$
coincides with the curve $\bold r(t)$ is called an {\it evolvent} of the
curve $\bold r(t)$. The problem of constructing the evolute of a given
curve is solved by the formula \mythetag{4.4}. The inverse problem of 
constructing an evolvent for a given curve appears to be more complicated.
It is effectively solved only in the case of a planar curve.\par
     Let $\bold r(s)$ be a vector-function defining some planar curve in 
natural parametri\-zation and let $\tilde\bold r(\tilde s)$ be the evolvent
in its own natural parametrization. Two natural parameters $s$ and 
$\tilde s$ are related to each other by some function $\varphi$ in form
of the relationship $\tilde s=\varphi(s)$. Let $\psi=\varphi^{-1}$ be the
inverse function for $\varphi$, then $s=\psi(\tilde s)$. Using the formula
\mythetag{4.4}, now we obtain
$$
\hskip -2em
\bold r(\psi(\tilde s))=\tilde\bold r(\tilde s)+
\frac{\tilde\bold n(\tilde s)}{\tilde k(\tilde s)}.
\mytag{4.5}
$$
Let's differentiate the relationship \mythetag{4.5} with respect to 
$\tilde s$ and then let's apply the formula \mythetag{3.1} and the Frenet
equations written in form of \mythetag{4.1}:
$$
\psi'(\tilde s)\cdot\boldsymbol\tau(\psi(\tilde s))=
\frac{d}{d\tilde s}\left(\frac{1}{\tilde k(\tilde s)}\right)
\cdot\tilde\bold n(\tilde s).
$$
Here $\boldsymbol\tau(\psi(\tilde s))$ and $\tilde\bold n
(\tilde s)$ both are unit vectors which are collinear due to the 
above relationship. Hence, we have the following two equalities:
$$
\xalignat 2
&\hskip -2em
\tilde\bold n(\tilde s)=\pm\boldsymbol\tau(\psi(\tilde s)),
&&\psi'(\tilde s)=\pm\frac{d}{d\tilde s}\left(\frac{1}{\tilde k(\tilde s)}
\right).
\mytag{4.6}
\endxalignat
$$
The second equality \mythetag{4.6} can be integrated:
$$
\hskip -2em
\frac{1}{\tilde k(\tilde s)}=\pm(\psi(\tilde s)-C).
\mytag{4.7}
$$
Here $C$ is a constant of integration. Let's combine \mythetag{4.7} with 
the first relationship \mythetag{4.6} and substitute it into the formula
\mythetag{4.5}:
$$
\tilde\bold r(\tilde s)=\bold r(\psi(\tilde s))+
(C-\psi(\tilde s))\cdot\boldsymbol\tau(\psi(\tilde s)).
$$
Then we substitute $\tilde s=\varphi(s)$ into the above formula
and denote $\boldsymbol\rho(s)=\tilde\bold r(\varphi(s))$. As a result 
we obtain the following equality:
$$
\hskip -2em
\boldsymbol\rho(s)=\bold r(s)+(C-s)\cdot\boldsymbol\tau(s).
\mytag{4.8}
$$
The formula \mythetag{4.8} is a parametric equation for the evolvent
of a planar curve $\bold r(s)$. The entry of an arbitrary constant
in the equation \mythetag{4.8} means the evolvent is not unique. Each
curve has the family of evolvents. This fact is valid for non-planar
curves either. However, we should emphasize that the formula \mythetag{4.8}
cannot be applied to general spacial curves. 
\head
\S~\mysection{5} Curves as trajectories of material points in mechanics.
\endhead
\rightheadtext{\S~5. Curves as trajectories of material 
points \dots}
     The presentation of classical mechanics traditionally begins
with considering the motion of material points. Saying {\it material 
point}, we understand any material object whose sizes are much smaller
than its displacement in the space. The position of such an object
can be characterized by its radius-vector in some Cartesian coordinate
system, while its motion is described by a vector-function $\bold r(t)$. 
The curve $\bold r(t)$ is called the {\it trajectory} of a material point.
Unlike to purely geometric curves, the trajectories of material points
possess preferable parameter $t$, which is usually distinct from the 
natural parameter $s$. This preferable parameter is the time variable 
$t$.\par
     The tangent vector of a trajectory, when computed in the time parametrization, is called the {\it velocity} of a material point:
$$
\hskip -2em
\bold v(t)=\frac{d\bold r}{dt}=\dot\bold r(t)=
\Vmatrix
v^1(t)\\v^2(t)\\v^3(t)
\endVmatrix.
\mytag{5.1}
$$
The time derivative of the velocity vector is called the {\it acceleration 
vector}:
$$
\hskip -2em
\bold a(t)=\frac{d\bold v}{dt}=\dot\bold v(t)=
\Vmatrix
a^1(t)\\a^2(t)\\a^3(t)
\endVmatrix.
\mytag{5.2}
$$
The motion of a material point in mechanics is described by Newton's 
second law:
$$
m\,\bold a=\bold F(\bold r,\bold v).
\mytag{5.3}
$$
Here $m$ is the mass of a material point. This is a constant characterizing
the amount of matter enclosed in this material object. The vector $\bold F$
is the {\it force vector}. By means of the force vector in mechanics one
describes the action of ambient objects (which are sometimes very far apart)
upon the material point under consideration. The magnitude of this action 
usually depends on the position of a point relative to the ambient objects,
but sometimes it can also depend on the velocity of the point itself.
Newton's second law in form of \mythetag{5.3} shows that the external
action immediately affects the acceleration of a material point, but neither
the velocity nor the coordinates of a point.\par
     Let $s=s(t)$ be the natural parameter on the trajectory of a material
point expressed through the time variable. Then the formula \mythetag{2.5} 
yields
$$
\hskip -2em
\dot s(t)=|\bold v(t)|=v(t).
\mytag{5.4}
$$
Through $v(t)$ in \mythetag{5.4} we denote the modulus of the velocity
vector.\par
     Let's consider a trajectory of a material point in natural 
parametrization: $\bold r=\bold r(s)$. Then for the velocity vector 
\mythetag{5.1} and for the acceleration vector \mythetag{5.2} we get the
following expressions:
$$
\aligned
&\bold v(t)=\dot s(t)\cdot\boldsymbol\tau(s(t)),\\
&\bold a(t)=\ddot s(t)\cdot\boldsymbol\tau(s(t))+
(\dot s(t))^2\cdot\boldsymbol\tau'(s(t)).
\endaligned
$$
Taking into account the formula \mythetag{5.4} and the first Frenet equation,
these expressions can be rewritten as 
$$
\hskip -2em
\aligned
&\bold v(t)=v(t)\cdot\boldsymbol\tau(s(t)),\\
&\bold a(t)=\dot v(t)\cdot\boldsymbol\tau(s(t))+
\left(k(s(t))\,v(t)^2\right)\cdot\bold n(s(t)).
\endaligned
\mytag{5.5}
$$
The second formula \mythetag{5.5} determines the expansion of the
acceleration vector into two components. The first component is 
tangent to the trajectory, it is called the {\it tangential
acceleration}. The second component is perpendicular to the 
trajectory and directed toward the curvature center. It is called
the {\it centripetal acceleration}. It is important to note that
the centripetal acceleration is determined by the modulus of
the velocity and by the geometry of the trajectory (by its
curvature).\par
\newpage
\topmatter
\title\chapter{2}
Elements of vectorial and tensorial analysis.
\endtitle
\endtopmatter
\chapternum=2
\document
\head
\S~\mysection{1} Vectorial and tensorial fields in the space.
\endhead
\setfirstpage
\leftheadtext{Chapter \uppercase\expandafter{\romannumeral 2}.
ELEMENTS OF TENSORIAL ANALYSIS.}
\rightheadtext{\S~1. Vectorial and tensorial fields in the space.}
     Let again $\Bbb E$ be a three-dimensional Euclidean point space.
We say that in $\Bbb E$ a {\it vectorial field} is given if at each
point of the space $\Bbb E$ some vector attached to this point is given.
Let's choose some Cartesian coordinate system in $\Bbb E$; in general,
this system is skew-angular. Then we can define the points of the space
by their coordinates $x^1$, $x^2$, $x^3$, and, simultaneously, we get 
the basis $\bold e_1$, $\bold e_2$, $\bold e_3$ for expanding the vectors
attached to these points. In this case we can present any vector field
$\bold F$ by three numeric functions
$$
\hskip -2em
\bold F=
\Vmatrix
F^1(\bold x)\\
F^2(\bold x)\\
F^3(\bold x)
\endVmatrix,
\mytag{1.1}
$$
where $\bold x=(x^1,x^2,x^3)$ are the components of the radius-vector 
of an arbitrary point of the space $\Bbb E$. Writing $\bold F(\bold x)$
instead of $\bold F(x^1,x^2,x^3)$, we make all formulas more compact.
\par
     The vectorial nature of the field $\bold F$ reveals when we replace
one coordinate system by another. Let \mythetag{1.1} be the coordinates of
a vector field in some coordinate system $O,\,\bold e_1,\,\bold e_2,\,
\bold e_3$ and let $\tilde O,\,\tilde\bold e_1,\,\tilde\bold e_2,\,
\tilde\bold e_3$ be some other coordinate system. The transformation rule
for the components of a vectorial field under a change of a Cartesian
coordinate system is written as follows:
$$
\hskip -2em
\aligned
&F^i(\bold x)=\sum^3_{j=1} S^i_j\,\,
\tilde F^j(\tilde\bold x),\\
\vspace{-0.7ex}
&x^i=\sum^3_{j=1} S^i_j\,\tilde x^j+a^i.
\endaligned
\mytag{1.2}
$$
Here $S^i_j$ are the components of the transition matrix relating the
basis $\bold e_1,\,\bold e_2,\,\bold e_3$ with the new basis 
$\tilde\bold e_1,\,\tilde\bold e_2,\,\tilde\bold e_3$, while 
$a^1,a^2,a^3$ are the components of the vector 
$\overrightarrow{O\tilde O\ }$ in the basis $\bold e_1,\,\bold e_2,\,
\bold e_3$.\par
    The formula \mythetag{1.2} combines the transformation rule 
for the components of a vector under a change of a basis and the
transformation rule for the coordinates of a point under a change 
of a Cartesian coordinate system (see  \mycite{1}). The arguments 
$\bold x$ and $\tilde\bold x$ beside the vector components $F^i$ and
$\tilde F^i$ in \mythetag{1.2} is an important novelty as compared to
\mycite{1}. It is due to the fact that here we deal with vector fields, 
not with separate vectors.\par
      Not only vectors can be associated with the points of the space
$\Bbb E$. In linear algebra along with vectors one considers covectors,
linear operators, bilinear forms and quadratic forms. Associating 
some covector with each point of $\Bbb E$, we get a {\it covector}
field. If we associate some linear operator with each point of the
space, we get an {\it operator field}. An finally, associating
a bilinear (quadratic) form with each point of $\Bbb E$, we obtain
a {\it field of bilinear (quadratic) forms}. Any choice of a Cartesian
coordinate system $O,\,\bold e_1,\,\bold e_2,\,\bold e_3$ assumes the
choice of a basis $\bold e_1,\bold e_2,\bold e_3$, while the basis
defines the numeric representations for all of the above objects:
for a covector this is the list of its components, for linear operators,
bilinear and quadratic forms these are their matrices. Therefore defining
a covector field $\bold F$ is equivalent to defining three functions
$F_1(\bold x)$, $F_2(\bold x)$, $F_3(\bold x)$ that transform according
to the following rule under a change of a coordinate system:
$$
\hskip -2em
\aligned
&F_i(\bold x)=\sum^3_{j=1} T^j_i\,\,
\tilde F_j(\tilde\bold x),\\
\vspace{-0.7ex}
&\dsize x^i=\sum^3_{j=1} S^i_j\,\tilde x^j+a^i.
\endaligned
\mytag{1.3}
$$
In the case of operator field $\bold F$ the transformation formula
for the components of its matrix under a change of a coordinate system
has the following form:
$$
\hskip -2em
\aligned
&F^i_j(\bold x)=\sum^3_{p=1}\sum^3_{q=1} S^i_p\,
T^q_j\,\,\tilde F^p_q(\tilde\bold x),\\
\vspace{-0.7ex}
&\dsize x^i=\sum^3_{p=1} S^i_p\,\tilde x^p+a^i.
\endaligned
\mytag{1.4}
$$
For a field of bilinear (quadratic) forms $\bold F$ the transformation
rule for its components under a change of Cartesian coordinates looks
like
$$
\aligned
&F_{ij}(\bold x)=\sum^3_{p=1}\sum^3_{q=1} T^p_i\,
T^q_j\,\,\tilde F_{p\,q}(\tilde\bold x),\\
\vspace{-0.7ex}
&\dsize x^i=\sum^3_{p=1} S^i_p\,\tilde x^p+a^i.
\endaligned
\mytag{1.5}
$$
Each of the relationships \mythetag{1.2}, \mythetag{1.3}, \mythetag{1.4}, 
and \mythetag{1.5} consists of two formulas. The first formula relates
the components of a field, which are the functions of two different sets
of arguments $\bold x=(x^1,x^2,x^3)$ and $\tilde\bold x=(\tilde x^1,
\tilde x^2,\tilde x^3)$. The second formula establishes the functional
dependence of these two sets of arguments.\par
     The first formulas in \mythetag{1.2}, \mythetag{1.3}, and \mythetag{1.4} are different. However, one can see some regular pattern
in them. The number of summation signs and the number of summation indices in their right hand sides are determined by the number of indices in
the components of a field $\bold F$. The total number of transition 
matrices used in the right hand sides of these formulas is also determined
by the number of indices in the components of $\bold F$. Thus, each
upper index of $\bold F$ implies the usage of the transition matrix $S$,
while each lower index of $\bold F$ means that the inverse matrix $T=
S^{-1}$ is used.\par
     The number of indices of the field $\bold F$ in the above examples
doesn't exceed two. However, the regular pattern detected in the transformation
rules for the components of $\bold F$ can be generalized for the case of an 
arbitrary number of indices:
$$
\hskip -2em
\ F^{i_1\ldots i_r}_{j_1\ldots j_s}=
\sum\Sb p_1\ldots p_r\\ q_1\ldots q_s\endSb
S^{i_1}_{p_1}\ldots\,S^{i_r}_{p_r}\,\,
T^{q_1}_{j_1}\ldots\,T^{q_s}_{j_s}\,\,
\tilde F^{p_1\ldots p_r}_{q_1\ldots q_s}
\mytag{1.6}
$$
The formula \mythetag{1.6} comprises the multiple summation with respect
to $(r+s)$ indices $p_1,\ldots,p_r$ and $q_1,\ldots,q_s$ each of which 
runs from 1 to 3.
\mydefinition{1.1} A tensor of the type $(r,s)$ is a geometric object 
$\bold F$ whose components in each basis are enumerated by $(r+s)$
indices and obey the transformation rule \mythetag{1.6} under a change of
basis.
\enddefinition
     Lower indices in the components of a tensor are called {\it covariant 
indices}, upper indices are called {\it contravariant indices} respectively. 
Generalizing the concept of a vector field, we can attach some tensor
of the type $(r,s)$, to each point of the space. As a result we get the
concept of a {\it tensor field}. This concept is convenient because it
describes in the unified way any vectorial and covectorial fields, operator
fields, and arbitrary fields of bilinear (quadratic) forms. Vectorial fields
are fields of the type $(1,0)$, covectorial fields have the type $(0,1)$,
operator fields are of the type $(1,1)$, and finally, any field of bilinear
(quadratic) forms are of the type $(0,2)$. Tensor fields of some other types
are also meaningful. In Chapter \uppercase\expandafter{\romannumeral 4}
we consider the curvature field with four indices.\par
     Passing from separate tensors to tensor fields, we acquire the 
arguments in formula \mythetag{1.6}. Now this formula should be written 
as the couple of two relationships similar to \mythetag{1.2}, \mythetag{1.3}, \mythetag{1.4}, or \mythetag{1.5}:
$$
\aligned
&F^{i_1\ldots\,i_r}_{j_1\ldots\,j_s}(\bold x)=
 \sum\Sb p_1\ldots\,p_r\\ q_1\ldots\,q_s\endSb
   S^{i_1}_{p_1}\ldots\,S^{i_r}_{p_r}\,\,
   T^{q_1}_{j_1}\ldots\,T^{q_s}_{j_s}\,\,
   \tilde F^{p_1\ldots\,p_r}_{q_1\ldots\,q_s}
   (\tilde\bold x),\\
&x^i=\sum^3_{j=1} S^i_j\,\tilde x^j+a^i.
\endaligned
\mytag{1.7}
$$
The formula \mythetag{1.7} expresses the transformation rule 
for the components of a tensorial field of the type $(r,s)$ 
under a change of Cartesian coordinates.\par
      The most simple type of tensorial fields is the type $(0,0)$.
Such fields are called {\it scalar fields}. Their components have
no indices at all, i.\,e\. they are numeric functions in the space 
$\Bbb E$.\par
\head
\S~\mysection{2} Tensor product and contraction.
\endhead
\rightheadtext{\S~2. Tensor product and contraction.}
     Let's consider two covectorial fields $\bold a$ and $\bold b$.
In some Cartesian coordinate system they are given by their 
components $a_i(\bold x)$ and $b_j(\bold x)$. These are two sets
of functions with three functions in each set. Let's form a new set
of nine functions by multiplying the functions of initial sets:
$$
\hskip -2em
c_{ij}(\bold x)=a_i(\bold x)\,b_j(\bold x).
\mytag{2.1}
$$
Applying the formula \mythetag{1.3} we can express the right hand side
of \mythetag{2.1} through the components of the fields $\bold a$ and 
$\bold b$ in the other coordinate system:
$$
\pagebreak
c_{ij}(\bold x)=\left(\shave{\sum^3_{p=1}}T^p_i\,\tilde a_p
\right)\left(\shave{\sum^3_{q=1}}T^q_j\,\tilde b_q\right)=
\sum^3_{p=1}\sum^3_{q=1} T^p_i\,T^q_j\,(\tilde a_p\,\tilde
b_q).
$$
If we denote by $\tilde c_{pq}(\tilde\bold x)$ the product of
$\tilde a_i(\tilde\bold x)$ and $\tilde b_j(\tilde\bold x)$, then
we find that the quantities $c_{ij}(\bold x)$ and $\tilde c_{pq}(\tilde
\bold x)$ are related by the formula \mythetag{1.5}. This means that
taking two covectorial fields one can compose a field of bilinear
forms by multiplying the components of these two covectorial fields 
in an arbitrary Cartesian coordinate system. This operation is called
the {\it tensor product} of the fields $\bold a$ and $\bold b$.
Its result is denoted as $\bold c=\bold a\otimes\bold b$.\par
     The above trick of multiplying components can be applied to an
arbitrary pair of tensor fields. Suppose we have a tensorial field
$\bold A$ of the type $(r,s)$ and another tensorial field $\bold B$
of the type $(m,n)$. Denote
$$
\hskip -2em
C^{\,i_1\ldots\,i_ri_{r+1}\ldots\,i_{r+m}}_{j_1\ldots\,j_s
j_{s+1}\ldots\,j_{s+n}}(\bold x)=
A^{i_1\ldots\,i_r}_{j_1\ldots\,j_s}(\bold x)\,\,
B^{i_{r+1}\ldots\,i_{r+m}}_{j_{s+1}\ldots\,j_{s+n}}(\bold x).
\mytag{2.2}
$$
\mydefinition{2.1} The tensor field $\bold C$ of the type
$(r+m,s+n)$ whose components are determined by the formula \mythetag{2.2}
is called the {\it tensor product} of the fields $\bold A$ and $\bold B$.
It is denoted  $\bold C=\bold A \otimes\bold B$.
\enddefinition
     This definition should be checked for correctness. We should make
sure that the components of the field $\bold C$ are transformed according
to the rule \mythetag{1.7} when we pass from one Cartesian coordinate 
system to another. The transformation rule \mythetag{1.7}, when applied
to the fields $\bold A$ and $\bold B$, yields
$$
\aligned
&A^{i_1\ldots\,i_r}_{j_1\ldots\,j_s}=
   \sum_{p..q}
   S^{i_1}_{p_1}\ldots\,S^{i_r}_{p_r}\,\,
   T^{q_1}_{j_1}\ldots\,T^{q_s}_{j_s}\,\,
   \tilde A^{p_1\ldots\,p_r}_{q_1\ldots\,q_s},\\
&B^{i_{r+1}\ldots\,i_{r+m}}_{j_{s+1}\ldots\,j_{s+n}}=
   \sum_{p..q}
   S^{i_{r+1}}_{p_{r+1}}\ldots\,S^{i_{r+m}}_{p_{r+m}}\,\,
   T^{q_{s+1}}_{j_{s+1}}\ldots\,T^{q_{s+n}}_{j_{s+n}}\,\,
   \tilde B^{p_{r+1}\ldots\,p_{r+m}}_{q_{s+1}\ldots\,q_{s+n}}.
\endaligned
$$
The summation in right hand sides of this formulas is carried out
with respect to each double index which enters the formula twice
--- once as an upper index and once as a lower index. Multiplying
these two formulas, we get exactly the transformation rule 
\mythetag{1.7} for the components of $\bold C$.\par
\mytheorem{2.1} The operation of tensor product is associative, 
this means that $(\bold A\otimes\bold B)\otimes\bold C=
\bold A\otimes(\bold B\otimes\bold C)$.
\endproclaim
\demo{Proof} Let $\bold A$ be a tensor of the type $(r,s)$,
let $\bold B$ be a tensor of the type $(m,n)$, and let $\bold C$ 
be a tensor of the type $(p,q)$. Then one can write the following
obvious numeric equality for their components:
$$
\align
&\hskip -2em
\left(A^{i_1\ldots\,i_r}_{j_1\ldots\,j_s}\,\,
B^{i_{r+1}\ldots\,i_{r+m}}_{j_{s+1}\ldots\,j_{s+n}}\right)\,
C^{\,i_{r+m+1}\ldots\,i_{r+m+p}}_{j_{s+n+1}\ldots\,j_{s+n+q}}=\\
\vspace{-0.4ex}
&\hskip -2em
\mytag{2.3}\\
\vspace{-1.5ex}
&\hskip -2em
\quad=A^{i_1\ldots\,i_r}_{j_1\ldots\,j_s}\,\left(
B^{i_{r+1}\ldots\,i_{r+m}}_{j_{s+1}\ldots\,j_{s+n}}\,\,
C^{\,i_{r+m+1}\ldots\,i_{r+m+p}}_{j_{s+n+1}\ldots\,j_{s+n+q}}
\right).
\endalign
$$
As we see in \mythetag{2.3}, the associativity of the tensor product follows
from the associativity of the multiplication of numbers.
\qed\enddemo
     The tensor product is not commutative. One can easily construct
an example illustrating this fact. Let's consider two covectorial
fields $\bold a$ and $\bold b$ with the following components in some
coordinate system: $\bold a=(1,0,0)$ and $\bold b=(0,1,0)$. Denote
$\bold c=\bold a\otimes\bold b$ and $\bold d=\bold b\otimes\bold a$. 
Then for $c_{12}$ and $d_{12}$ with the use of the formula \mythetag{2.2} 
we derive: $c_{12}=1$ and $d_{12}=0$. Hence, $\bold c\neq\bold d$ and
$\bold a\otimes\bold b\neq\bold b\otimes\bold a$.\par
    Let's consider an operator field $\bold F$. Its components
$F^i_j(\bold x)$ are the components of the operator $\bold F(\bold x)$ 
in the basis $\bold e_1,\,\bold e_2,\,\bold e_3$. It is known that 
the trace of the matrix $F^i_j(\bold x)$ is a scalar invariant of the
operator $\bold F(\bold x)$ (see \mycite{1}). Therefore, the formula
$$
\hskip -2em
f(\bold x)=\tr\bold F(\bold x)=\sum^3_{i=1}
F^i_i(\bold x)
\mytag{2.4}
$$
determines a scalar field $f(\bold x)$ in the space $\Bbb E$. The sum
similar to \mythetag{2.4} can be written for an arbitrary tensorial field $\bold F$ with at least one upper index and at least one lower index in its 
components:
$$
\hskip -2em
H^{i_1\ldots\,i_{r-1}}_{j_1\ldots\,j_{s-1}}(\bold x)=
\sum^3_{k=1}
F^{i_1\ldots\,i_{m-1}\,k\,i_m\ldots\,i_{r-1}}_{j_1\ldots\,
j_{n-1}\,k\,j_n\ldots\,j_{s-1}}(\bold x).
\mytag{2.5}
$$
In the formula \mythetag{2.5} the summation index $k$ is placed to $m$-th
upper position and to $n$-th lower position. The succeeding indices 
$i_m,\,\ldots\,\,i_{r-1}$ and $j_n,\,\ldots\,\,j_{s-1}$ in writing 
the components of the field $\bold F$ are shifted one 
position to the right as compared to their positions in left hand side
of the equality \mythetag{2.5}:
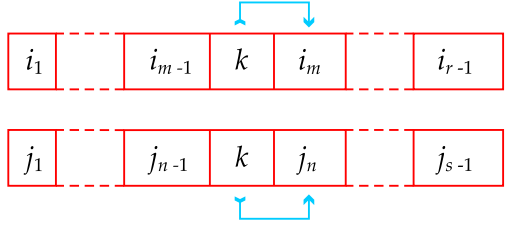
\newline
\vadjust{\vskip 40pt\hbox to 0pt{\kern 80pt
\includegraphics{ris04.eps}\hss}\vskip 55pt}
\mydefinition{2.2} The tensor field $\bold H$ whose components
are calculated according to the formula \mythetag{2.5} from the components
of the tensor field $\bold F$ is called the {\it contraction} of the
field $\bold F$ with respect to $m$-th and $n$-th indices.
\enddefinition
     Like the definition~\mythedefinition{2.1}, this definition should be
tested for correctness. Let's verify that the components of the field
$\bold H$ are transformed according to the formula \mythetag{1.7}. For this
purpose we write the transformation rule \mythetag{1.7} applied to the
components of the field $\bold F$ in right hand side of the formula
\mythetag{2.5}:
$$
\align
&F^{i_1\ldots\,i_{m-1}\,k\,i_m\ldots\,i_{r-1}}_{j_1\ldots\,
j_{n-1}\,k\,j_n\ldots\,j_{s-1}}=
\sum\Sb \alpha\,p_1\ldots p_{r-1}\\
\beta\,q_1\ldots q_{s-1}\endSb
 S^{i_1}_{p_1}\ldots\,S^{i_{m-1}}_{p_{m-1}}\,S^k_{\alpha}\,
 S^{i_m}_{p_m}\ldots\,S^{i_{r-1}}_{p_{r-1}}\times\\
&\qquad\times T^{q_1}_{j_1}\ldots\,T^{q_{n-1}}_{j_{n-1}}\,T^\beta_k\,
   T^{q_n}_{j_n}\ldots\,T^{q_{s-1}}_{j_{s-1}}\,\,
   \tilde F^{p_1\ldots\,p_{m-1}\,\alpha\,p_m\ldots\,p_{r-1}}_{q_1
   \ldots\,q_{n-1}\,\beta\,q_n\ldots\,q_{s-1}}.
\endalign
$$
In order to derive this formula from \mythetag{1.7} we substitute
the index $k$ into the $m$-th and $n$-th positions, then we shift
all succeeding indices one position to the right. In order to have
more similarity of left and right hand sides of this formula we
shift summation indices as well. It is clear that such redesignation
of summation indices does not change the value of the sum.\par
     Now in order to complete the contraction procedure we should
produce the summation with respect to the index $k$. In the right hand
side of the formula the sum over $k$ can be calculated explicitly due
to the formula 
$$
\sum^3_{k=1} S^k_{\alpha}\,T^\beta_k=\delta^\beta_\alpha,
\mytag{2.6}
$$
which means $T=S^{-1}$. Due to \mythetag{2.6} upon calculating the sum
over $k$ one can calculate the sums over $\beta$ and $\alpha$. Therein
we take into account that
$$
 \sum^3_{\alpha=1}
   \tilde F^{p_1\ldots\,p_{m-1}\,\alpha\,p_m\ldots\,p_{r-1}}_{q_1
   \ldots\,q_{n-1}\,\alpha\,q_n\ldots\,q_{s-1}}=
   \tilde H^{p_1\ldots\,p_{r-1}}_{q_1\ldots\,q_{s-1}}.
$$
As a result we get the equality
$$
H^{i_1\ldots\,i_{r-1}}_{j_1\ldots\,j_{s-1}}=
\sum\Sb p_1\ldots p_{r-1}\\ q_1\ldots q_{s-1}\endSb
 S^{i_1}_{p_1}\ldots\,S^{i_{r-1}}_{p_{r-1}}\,\,
 T^{q_1}_{j_1}\ldots\,T^{q_{s-1}}_{j_{s-1}}\,
 \tilde H^{p_1\ldots\,p_{r-1}}_{q_1\ldots\,q_{s-1}},
$$
which exactly coincides with the transformation rule \mythetag{1.7}
written with respect to components of the field $\bold H$. The
correctness of the definition~\mythedefinition{2.2} is proved.\par
     The operation of contraction introduced by the
definition~\mythedefinition{2.2} implies that the positions of two indices
are specified. One of these indices should be an upper index, the other
index should be a lower index. The letter $C$ is used as a contraction
sign. The formula \mythetag{2.5} then is abbreviated as follows:
$$
\bold H=C_{m,n}(\bold F)=C(\bold F).
$$
The numbers $m$ and $n$ are often omitted since they are usually known 
from the context.\par
     A tensorial field of the type $(1,1)$ can be contracted in the unique
way. For a tensorial field $\bold F$ of the type $(2,2)$ we have two ways
of contracting. As a result of these two contractions, in general, we obtain 
two different tensorial fields of the type $(1,1)$. These tensorial fields
can be contracted again. As a result we obtain the {\it complete
contractions} of the field $\bold F$, they are scalar fields. A field of 
the type $(2,2)$ can have two complete contractions. In general case a 
field of the type $(n,n)$ has $\,n!\,$ complete contractions.\par
     The operations of tensor product and contraction often arise in a 
natural way without any special intension. For example, suppose that we
are given a vector field $\bold v$ and a covector field $\bold w$ in
the space $\Bbb E$. This means that at each point we have a vector and a
covector attached to this point. By calculating the scalar products of
these vectors and covectors we get a scalar field $f=\langle \bold w\,|\,
\bold v\rangle$. In coordinate form such a scalar field is calculated
by means of the formula 
$$
\hskip -2em
f=\sum^3_{k=1} w_i\,v^i.
\mytag{2.7}
$$
From the formula \mythetag{2.7}, it is clear that $f=C(\bold w\otimes\bold v)$.
The scalar product $f=\langle \bold w\,|\,\bold v\rangle$ is the contraction
of the tensor product of the fields $\bold w$ and $\bold v$. In a similar 
way, if an operator field $\bold F$ and a vector field $\bold v$ are given,
then applying $\bold F$ to $\bold v$ we get another vector field $\bold
u=\bold F\,\bold v$, where
$$
u^i=\sum^3_{j=1} F^i_j\,v^j.
$$
In this case we can write: $\bold u=C(\bold F\otimes\bold v)$; although this
writing cannot be uniquely interpreted. Apart from $\bold u=\bold F\,\bold
v$, it can mean the product of $\bold v$ by the trace of the operator field 
$\bold F$.
\head
\S~\mysection{3} The algebra of tensor fields.
\endhead
\rightheadtext{\S~3. The algebra of tensor fields.}
     Let $\bold v$ and $\bold w$ be two vectorial fields. Then at each 
point of the space $\Bbb E$ we have two vectors $\bold v(\bold x)$ and
$\bold w(\bold x)$. We can add them. As a result we get a new vector
field $\bold u=\bold v+\bold w$. In a similar way one can define the
addition of tensor fields. Let $\bold A$ and $\bold B$ be two tensor 
fields of the type $(r,s)$. Let's consider the sum of their components
in some Cartesian coordinate system:
$$
\hskip -2em
C^{\,i_1\ldots\,i_r}_{j_1\ldots\,j_s}=
A^{i_1\ldots\,i_r}_{j_1\ldots\,j_s}+
B^{i_1\ldots\,i_r}_{j_1\ldots\,j_s}.
\mytag{3.1}
$$
\mydefinition{3.1} The tensor field $\bold C$ of the type
$(r,s)$ whose components are calculated according to the formula 
\mythetag{3.1} is called the {\it sum} of the fields $\bold A$ and 
$\bold B$ of the type $(r,s)$.
\enddefinition
    One can easily check up the transformation rule \mythetag{1.7} 
for the components of the field $\bold C$. It is sufficient to 
write this rule \mythetag{1.7} for the components of $\bold A$ and 
$\bold B$ then add these two formulas. Therefore, the
definition~\mythedefinition{3.1} is consistent.\par
     The sum of tensor fields is commutative and associative. This
fact follows from the commutativity and associativity of the addition
of numbers due to the following obvious relationships:
$$
\align
&A^{i_1\ldots\,i_r}_{j_1\ldots\,j_s}+
B^{i_1\ldots\,i_r}_{j_1\ldots\,j_s}=
B^{i_1\ldots\,i_r}_{j_1\ldots\,j_s}+
A^{i_1\ldots\,i_r}_{j_1\ldots\,j_s},\\
\vspace{1ex}
&\left(A^{i_1\ldots\,i_r}_{j_1\ldots\,j_s}+
B^{i_1\ldots\,i_r}_{j_1\ldots\,j_s}\right)+
C^{\,i_1\ldots\,i_r}_{j_1\ldots\,j_s}=
A^{i_1\ldots\,i_r}_{j_1\ldots\,j_s}+
\left(B^{i_1\ldots\,i_r}_{j_1\ldots\,j_s}+
C^{\,i_1\ldots\,i_r}_{j_1\ldots\,j_s}\right).
\endalign
$$\par
     Let's denote by $T_{(r,s)}$ the set of tensor fields of the type
$(r,s)$. The tensor multiplication introduced by the 
definition~\mythedefinition{2.1} is the following binary operation: 
$$
\hskip -2em
T_{(r,\,s)}\times T_{(m,\,n)}\to T_{(r+m,\,s+n)}.
\mytag{3.2}
$$
The operations of tensor addition and tensor multiplication 
\mythetag{3.2} are related to each other by the distributivity laws:
$$
\hskip -2em
\aligned
&(\bold A+\bold B)\otimes\bold C=
\bold A\otimes\bold C+\bold B\otimes\bold C,\\
&\bold C\otimes(\bold A+\bold B)=
\bold C\otimes\bold A+\bold C\otimes\bold B.
\endaligned
\mytag{3.3}
$$
The distributivity laws \mythetag{3.3} follow from the distributivity
of the multiplication of numbers. Their proof is given by the following
obvious formulas:
$$
\align
&\aligned
  &\left(A^{i_1\ldots\,i_r}_{j_1\ldots\,j_s}+
  B^{i_1\ldots\,i_r}_{j_1\ldots\,j_s}\right)
  C^{\,i_{r+1}\ldots\,i_{r+m}}_{j_{s+1}\ldots\,j_{s+n}}=\\
  \vspace{1ex}
  &\qquad=A^{i_1\ldots\,i_r}_{j_1\ldots\,j_s}\,
  C^{\,i_{r+1}\ldots\,i_{r+m}}_{j_{s+1}\ldots\,j_{s+n}}+
  B^{i_1\ldots\,i_r}_{j_1\ldots\,j_s}\,
  C^{\,i_{r+1}\ldots\,i_{r+m}}_{j_{s+1}\ldots\,j_{s+n}},
 \endaligned\\
 \vspace{1.2ex}
&\aligned
 &C^{\,i_1\ldots\,i_r}_{j_1\ldots\,j_s}\left(
  A^{i_{r+1}\ldots\,i_{r+m}}_{j_{s+1}\ldots\,j_{s+n}}+
  B^{i_{r+1}\ldots\,i_{r+m}}_{j_{s+1}\ldots\,j_{s+n}}
  \right)=\\
  \vspace{1ex}
 &\qquad=C^{\,i_1\ldots\,i_r}_{j_1\ldots\,j_s}\,
  A^{i_{r+1}\ldots\,i_{r+m}}_{j_{s+1}\ldots\,j_{s+n}}+
  C^{\,i_1\ldots\,i_r}_{j_1\ldots\,j_s}\,
  B^{i_{r+1}\ldots\,i_{r+m}}_{j_{s+1}\ldots\,j_{s+n}}.
 \endaligned
\endalign
$$\par
     Due to \mythetag{3.2} the set of scalar fields $K=T_{(0,0)}$
(which is simply the set of numeric functions) is closed with respect
to tensor multiplication $\otimes$, which coincides here with the 
regular multiplication of numeric functions. The set $K$ is a 
{\it commutative ring {\rm (see \mycite{3})} with the unity}. The constant
function equal to $1$ at each point of the space $\Bbb E$ plays the role
of the unit element in this ring.\par
     Let's set $m=n=0$ in the formula \mythetag{3.2}. In this case 
it describes the multiplication of tensor fields from $T_{(r,s)}$ 
by numeric functions from the ring $K$. The tensor product 
of a field $\bold A$ and a scalar filed $\xi\in K$ is commutative:
$\bold A\otimes\xi=\xi\otimes\bold A$. Therefore, the multiplication
of tensor fields by numeric functions is denoted by standard sign of
multiplication: $\xi\otimes\bold A=\xi\cdot\bold A$. The operation 
of addition and the operation of multiplication by scalar fields in 
the set $T_{(r,s)}$ possess the following properties:
\roster
\item $\bold A+\bold B=\bold B+\bold A$;
\item $(\bold A+\bold B)+\bold C=\bold A+(\bold B+\bold C)$;
\item there exists a field $\bold 0\in T_{(r,s)}$ such that 
      $\bold A+\bold 0=\bold A$ for an arbitrary tensor field 
      $\bold A\in T_{(r,s)}$;
\item for any tensor field $\bold A\in T_{(r,s)}$ there exists
      an opposite field $\bold A'$ such that $\bold A+\bold A'
      =\bold 0$;
\item $\xi\cdot(\bold A+\bold B)=\xi\cdot\bold A+\xi\cdot\bold B$
      for any function $\xi$ from the ring $K$ and for any two
      fields $\bold A,\bold B\in T_{(r,s)}$;
\item $(\xi+\zeta)\cdot\bold A=\xi\cdot\bold A+\zeta\cdot\bold A$
      for any tensor field $\bold A\in T_{(r,s)}$ and for any two
      functions $\xi,\zeta\in K$;
\item $(\xi\,\zeta)\cdot\bold A=\xi\cdot(\zeta\cdot\bold A)$
      for any tensor field $\bold A\in T_{(r,s)}$ and for any two
      functions $\xi,\zeta\in K$;
\item $1\cdot\bold A=\bold A$ for any field $\bold A\in T_{(r,s)}$.
\endroster
     The tensor field with identically zero components plays the role 
of zero element in the property \therosteritem{3}. The field $\bold A'$ 
in the property \therosteritem{4} is defined as a field whose components
are obtained from the components of $\bold A$ by changing the sign.\par
     The properties \therosteritem{1}-\therosteritem{8} listed above 
almost literally coincide with the axioms of a linear vector space (see 
\mycite{1}). The only discrepancy is that the set of functions $K$ is a
ring, not a numeric field as it should be in the case of a linear vector
space. The sets defined by the axioms \therosteritem{1}-\therosteritem{8}
for some ring $K$ are called {\it modules over the ring $K$} or {\it
$K$-modules}. Thus, each of the sets $T_{(r,s)}$ is a module over the ring
of scalar functions $K=T_{(0,0)}$.\par
     The ring $K=T_{(0,0)}$ comprises the subset of constant functions
which is naturally identified with the set of real numbers $\Bbb R$. 
Therefore the set of tensor fields $T_{(r,s)}$ in the space $\Bbb E$ 
is a linear vector space over the field of real numbers $\Bbb R$.\par
     If $r\geqslant 1$ and $s\geqslant 1$, then in the set $T_{(r,s)}$
the operation of contraction with respect to various pairs of indices
are defined. These operations are linear, i\.\,e\. the following
relationships are fulfilled:
$$
\hskip -2em
\aligned
&C(\bold A+\bold B)=C(\bold A)+C(\bold B),\\
&C(\xi\cdot\bold A)=\xi\cdot C(\bold A).
\endaligned
\mytag{3.4}
$$
The relationships \mythetag{3.4} are proved by direct calculations 
in coordinates. For the field $\bold C=\bold A+\bold B$ from \mythetag{2.5}
we derive
$$
\align
&H^{\,i_1\ldots\,i_{r-1}}_{j_1\ldots\,j_{s-1}}=
\sum^3_{k=1}
C^{\,i_1\ldots\,i_{m-1}\,k\,i_m\ldots\,i_{r-1}}_{j_1\ldots\,
j_{n-1}\,k\,j_n\ldots\,j_{s-1}}=\\
\vspace{-0.7ex}
&=\sum^3_{k=1}
A^{i_1\ldots\,i_{m-1}\,k\,i_m\ldots\,i_{r-1}}_{j_1\ldots\,
j_{n-1}\,k\,j_n\ldots\,j_{s-1}}
+\sum^3_{k=1}
B^{i_1\ldots\,i_{m-1}\,k\,i_m\ldots\,i_{r-1}}_{j_1\ldots\,
j_{n-1}\,k\,j_n\ldots\,j_{s-1}}.
\endalign
$$
This equality proves the first relationship \mythetag{3.4}. In order to
prove the second one we take $\bold C=\xi\cdot\bold A$. Then the
second relationship \mythetag{3.4} is derived as a result of the 
following calculations:
$$
\align
&H^{\,i_1\ldots\,i_{r-1}}_{j_1\ldots\,j_{s-1}}=
\sum^3_{k=1}
C^{\,i_1\ldots\,i_{m-1}\,k\,i_m\ldots\,i_{r-1}}_{j_1\ldots\,
j_{n-1}\,k\,j_n\ldots\,j_{s-1}}=\\
\vspace{-0.7ex}
&=\sum^3_{k=1}
\xi\,A^{i_1\ldots\,i_{m-1}\,k\,i_m\ldots\,i_{r-1}}_{j_1\ldots\,
j_{n-1}\,k\,j_n\ldots\,j_{s-1}}
=\xi\,\sum^3_{k=1}
A^{i_1\ldots\,i_{m-1}\,k\,i_m\ldots\,i_{r-1}}_{j_1\ldots\,
j_{n-1}\,k\,j_n\ldots\,j_{s-1}}.
\endalign
$$\par
     The tensor product of two tensors from $T_{(r,s)}$ belongs to
$T_{(r,s)}$ only if $r=s=0$ (see formula \mythetag{3.2}). In all other
cases one cannot perform the tensor multiplication staying within one $K$-module $T_{(r,s)}$. In order to avoid this restriction the 
following direct sum is usually considered:
$$
\hskip -2em
T=\bigoplus^\infty_{r=0}\bigoplus^\infty_{s=0}
T_{(r,s)}.
\mytag{3.5}
$$
The set \mythetag{3.5} consists of finite formal sums $\bold A^{(1)}+\ldots+\bold A^{(k)}$, where each summand belongs to some
of the $K$-modules $T_{(r,s)}$. The operation of tensor product 
is extended to the $K$-module $T$ by means of the formula:
$$
(\bold A^{(1)}+\ldots+\bold A^{(k)})\otimes
(\bold A^{(1)}+\ldots+\bold A^{(q)})=
\sum^k_{i=1}\sum^q_{j=1}\bold A^{(i)}\otimes\bold A^{(j)}.
$$
This extension of the operation of tensor product is a bilinear 
binary operation in the set $T$. It possesses the following 
additional properties:
\roster
\item[9] $(\bold A+\bold B)\otimes\bold C=\bold A\otimes\bold C+
      \bold B\otimes\bold C$;
\item $(\xi\cdot\bold A)\otimes\bold C=\xi\cdot(\bold A\otimes
      \bold C)$;
\item $\bold C\otimes(\bold A+\bold B)=\bold C\otimes\bold A+
      \bold C\otimes\bold B$;
\item $\bold C\otimes(\xi\cdot\bold B)=\xi\cdot(\bold C\otimes
      \bold B)$.
\endroster
These properties of the operation of tensor product in $T$ are
easily derived from \mythetag{3.3}. Note that a $K$-module equipped
with an additional bilinear binary operation of multiplication
is called an {\it algebra over the ring $K$} or a {\it $K$-algebra}.
Therefore the set $T$ is called the {\it algebra of tensor fields}.
\par
     The algebra $T$ is a direct sum of separate $K$-modules
$T_{(r,s)}$ in \mythetag{3.5}. The operation of multiplication
is concordant with this expansion into a direct sum; this fact
is expressed by the relationship \mythetag{3.2}. Such structures
in algebras are called {\it gradings}, while algebras with
gradings are called {\it graded algebras}.
\head
\S~\mysection{4} Symmetrization and alternation.
\endhead
\rightheadtext{\S~4. Symmetrization and alternation.}
     Let $\bold A$ be a tensor filed of the type $(r,s)$ and let 
$r\geqslant 2$. The number of upper indices in the components
of the field $\bold A$ is greater than two. Therefore, we can
perform the permutation of some pair of them. Let's denote 
$$
\hskip -2em
B^{i_1\,\ldots\,i_m\ldots\,i_n\ldots\,i_r}_{j_1\,\ldots\ldots\ldots
\ldots\ldots\ldots j_s}=A^{i_1\,\ldots\,i_n\ldots\,i_m\ldots\,i_r}_{j_1
\,\ldots\ldots\ldots\ldots\ldots\ldots j_s}.
\mytag{4.1}
$$
The quantities $B^{i_1\,\ldots\,i_r}_{j_1\,\ldots\,j_s}$ in \mythetag{4.1}
are produced from the components of the tensor field $\bold A$ by the 
transposition of the pair of upper indices $i_m$ and $i_n$.
\mytheorem{4.1} The quantities $B^{\,i_1\,\ldots\,i_r}_{j_1\,
\ldots\,j_s}$ produced from the components of a tensor field $\bold A$
by the transposition of any pair of upper indices define another tensor
field $\bold B$ of the same type as the original field $\bold A$.
\endproclaim
\demo{Proof} In order to prove the theorem let's check up that the
quantities \mythetag{4.1} obey the transformation rule \mythetag{1.7} 
under a change of a coordinate system:
$$
B^{i_1\ldots\,i_r}_{j_1\ldots\,j_s}=
\sum\Sb p_1\ldots p_r\\ q_1\ldots q_s\endSb
   S^{i_1}_{p_1}\ldots\,S^{i_n}_{p_m}\ldots\,
   S^{i_m}_{p_n}\ldots\,S^{i_r}_{p_r}\,\,
   T^{q_1}_{j_1}\ldots\,T^{q_s}_{j_s}\,\,
   \tilde A^{p_1\ldots\,p_r}_{q_1\ldots\,q_s}.
$$
Let's rename the summation indices $p_m$ and $p_n$ in this formula:
let's denote $p_m$ by $p_n$ and vice versa. As a result the $S$ matrices
will be arranged in the order of increasing numbers of their upper and 
lower indices. However, the indices $p_m$ and $p_n$ in $\tilde A^{p_1
\ldots\,p_r}_{q_1\ldots\,q_s}$ will exchange their positions. It is clear
that the procedure of renaming summation indices does not change the
value of the sum:
$$
B^{i_1\ldots\,i_r}_{j_1\ldots\,j_s}=
\sum\Sb p_1\ldots p_r\\ q_1\ldots q_s\endSb
   S^{i_1}_{p_1}\ldots\,S^{i_r}_{p_r}\,\,
   T^{q_1}_{j_1}\ldots\,T^{q_s}_{j_s}\,\,
   \tilde A^{p_1\ldots\,p_n\ldots\,p_m\ldots\,p_r}_{q_1
   \ldots\,q_s}.
$$
Due to the equality $\tilde B^{p_1\ldots\,p_r}_{q_1\ldots\,q_s}=
\tilde A^{p_1\ldots\,p_n\ldots\,p_m\ldots\,p_r}_{q_1\,\ldots\ldots
\ldots\ldots\ldots\ldots\,q_s}$ the above formula is exactly the 
transformation rule \mythetag{1.7} written for the quantities
\mythetag{4.1}. Hence, they define a tensor field $\bold B$. The 
theorem is proved.
\qed\enddemo
     There is a similar theorem for transpositions of lower indices.
Let again $\bold A$ be a tensor field of the type $(r,s)$ and let
$s\geqslant 2$. Denote
$$
B^{i_1\,\ldots\ldots\ldots\ldots\ldots\ldots\,i_r}_{j_1\,\ldots\,
j_m\ldots\,j_n\ldots\,j_s}=A^{i_1\,\ldots\ldots\ldots\ldots\ldots
\ldots\,i_r}_{j_1\,\ldots\,j_n\,\ldots\,j_m\ldots\,j_s}.
\mytag{4.2}
$$
\mytheorem{4.2} The quantities $B^{\,i_1\,\ldots\,i_r}_{j_1\,
\ldots\,j_s}$ produced from the components of a tensor field $\bold A$
by the transposition of any pair of lower indices define another tensor
field $\bold B$ of the same type as the original field $\bold A$.
\endproclaim
     The proof of the theorem~\mythetheorem{4.2} is completely analogous 
to the proof of the theorem~\mythetheorem{4.1}. Therefore we do not give 
it here. Note that one cannot transpose an upper index and a lower index.
The set of quantities obtained by such a transposition does not obey the
transformation rule \mythetag{1.7}.\par
     Combining various pairwise transpositions of indices \mythetag{4.1}
and \mythetag{4.2} we can get any transposition from the symmetric group
$\goth S_r$ in upper indices and any transposition from the symmetric 
group $\goth S_s$ in lower indices. This is a well-known fact from the
algebra (see \mycite{3}). Thus the theorems~\mythetheorem{4.1} and 
\mythetheorem{4.2} define the action of the groups $\goth S_r$ and $\goth
S_s$ on the $K$-module $T_{(r,s)}$ composed of the tensor fields of the
type $(r,s)$. This is the action by linear operators, i\.\,e\.
$$
\hskip -2em
\aligned
&\sigma\compos\tau(\bold A+\bold B)=
\sigma\compos\tau(\bold A)+
\sigma\compos\tau(\bold B),\\
&\sigma\compos\tau(\xi\cdot A)=
\xi\cdot(\sigma\compos\tau(\bold A))
\endaligned
\mytag{4.3}
$$
for any two transpositions $\sigma\in\goth S_r$ and $\tau\in\goth S_s$. 
When written in coordinate form, the relationship $\bold B=\sigma\compos
\tau(\bold A)$ looks like
$$
\hskip -2em
B^{i_1\ldots\,i_r}_{j_1\ldots\,j_s}=
A^{i_{\sigma(1)}\ldots\,i_{\sigma(r)}}_{j_{\tau(1)}\ldots\,
j_{\tau(s)}},
\mytag{4.4}
$$
where the umbers $\sigma(1),\ldots,\sigma(r)$ and $\tau(1),\ldots,\tau(s)$
are obtained by applying $\sigma$ and $\tau$ to the numbers $1,\ldots,r$
and $1,\ldots,s$.\par
\mydefinition{4.1} A tensorial field $\bold A$ of the type 
$(r,s)$ is said to be {\it symmetric} in $m$-th and $n$-th upper 
(or lower) indices if $\sigma(\bold A)=\bold A$, where $\sigma$ is
the permutation of the indices given by the formula \mythetag{4.1} (or
the formula \mythetag{4.2}).
\enddefinition
\mydefinition{4.2} A tensorial field $\bold A$ of the type 
$(r,s)$ is said to be {\it skew-symmetric} in $m$-th and $n$-th upper 
(or lower) indices if $\sigma(\bold A)=-\bold A$, where $\sigma$ is
the permutation of the indices given by the formula \mythetag{4.1} (or
the formula \mythetag{4.2}).
\enddefinition
     The concepts of symmetry and skew-symmetry can be extended to
the case of arbitrary (not necessarily pairwise) transpositions. Let
$\varepsilon=\sigma\compos\tau$ be some transposition of upper and 
lower indices from \mythetag{4.4}. It is natural to treat it 
as an element of direct product of two symmetric groups: $\varepsilon\in\goth S_r \times\goth S_s$ (see \mycite{3}).
\mydefinition{4.3} A tensorial field $\bold A$ of the type
$(r,s)$ is {\it symmetric} or {\it skew-symmetric} with respect to
the transposition $\varepsilon\in\goth S_r\times\goth S_s$, if
one of the following relationships is fulfilled: $\varepsilon(\bold A)
=\bold A$ or $\varepsilon(\bold A)=(-1)^\varepsilon\cdot\bold A$.
\enddefinition
    If the field $\bold A$ is symmetric with respect to the transpositions
$\varepsilon_1$ and $\varepsilon_2$, then it is symmetric with respect to
the composite transposition $\varepsilon_1\compos\varepsilon_2$ and with
respect to the inverse transpositions $\varepsilon_1^{-1}$ and $\varepsilon_2^{-1}$. Therefore the symmetry always takes place for some subgroup $\goth G\in\goth S_r\times\goth S_s$. The same is true for the
skew-symmetry.\par
     Let $\goth G\subset\goth S_r\times\goth S_s$ be a subgroup in the 
direct product of symmetric groups and let $\bold A$ be a tensor field 
from $T_{(r,s)}$. The passage from $\bold A$ to the field 
$$
\hskip -2em
\bold B=\frac{1}{|\goth G|}\,\sum_{\varepsilon\in\goth G}
\varepsilon(\bold A)
\mytag{4.5}
$$
is called the {\it symmetrization\/} of the tensor field $\bold A$ by the
subgroup $\goth G\subset\goth S_r\times\goth S_s$. Similarly, the passage
from $\bold A$ to the field 
$$
\hskip -2em
\bold B=\frac{1}{|\goth G|}\,\sum_{\varepsilon\in\goth G}
(-1)^\varepsilon\cdot\varepsilon(\bold A)
\mytag{4.6}
$$
is called the {\it alternation\/} of the tensor field $\bold A$ by the 
subgroup $\goth G\subset\goth S_r\times\goth S_s$.\par
     The operations of symmetrization and alternation are linear operations,
this fact follows from \mythetag{4.3}. As a result of symmetrization 
\mythetag{4.5} one gets a field $\bold B$ symmetric with respect to $\goth G$.
As a result of alternation \mythetag{4.6} one gets a field skew-symmetric 
with respect to $\goth G$. If $\goth G=\goth S_r\times\goth S_s$ then
the operation \mythetag{4.5} is called the {\it complete symmetrization},
while the \mythetag{4.6} is called the {\it complete alternation}.
\head
\S~\mysection{5} Differentiation of tensor fields.
\endhead
\rightheadtext{\S~5. Differentiation of tensor fields.}
    The smoothness class of a tensor field $\bold A$ in the space $\Bbb E$
is determined by the smoothness of its components.
\mydefinition{5.1} A tensor field $\bold A$ is called an {\it
$m$-times continuously differen\-tiable field} or a {\it field of the class
$C^m$} if all its components in some Cartesian system are $m$-times
continuously differentiable functions.
\enddefinition
    Tensor fields of the class $C^1$ are often called {\it differentiable tensor fields}, while fields of the class $C^\infty$ are called {\it smooth
tensor fields}. Due to the formula \mythetag{1.7} the choice of a Cartesian 
coordinate system does not affect the smoothness class of a field $\bold A$
in the definition~\mythedefinition{5.1}. The components of a field of the
class $C^m$ are the functions of the class $C^m$ in any Cartesian coordinate
system. This fact proves that the definition~\mythedefinition{5.1} is
consistent.\par
     Let's consider a differentiable tensor field of the type $(r,s)$ and
let's consider all of the partial derivatives of its components:
$$
\hskip -2em
B^{i_1\ldots\,i_r}_{j_1\ldots\,j_s\,j_{s+1}}=
\frac{\partial A^{i_1\ldots\,i_r}_{j_1\ldots\,j_s}}
{\partial x^{j_{s+1}}}.
\mytag{5.1}
$$
The number of such partial derivatives \mythetag{5.1} is the same in all Cartesian coordinate systems. This number coincides with the number of
components of a tensor field of the type $(r,s+1)$. This coincidence
is not accidental.
\mytheorem{5.1} The partial derivatives of the components of a differentiable tensor field $\bold A$ of the type $(r,s)$ calculated
in an arbitrary Cartesian coordinate system according to the formula
\mythetag{5.1} are the components of another tensor filed $\bold B$ of
the type $(r,s+1)$.
\endproclaim
\demo{Proof} The proof consists in checking up the transformation rule 
\mythetag{1.7} for the quantities $B^{\,i_1\ldots\,i_r}_{j_1\ldots\,
j_s\,j_{s+1}}$ in \mythetag{5.1}. Let $O,\bold e_1,\bold e_2,\bold e_3$ 
and $O',\tilde\bold e_1,\tilde\bold e_2,\tilde\bold e_3$ be two
Cartesian coordinate systems. By tradition we denote by $S$ and $T$
the direct and inverse transition matrices. Let's write the first relationship \mythetag{1.7} for the field $\bold A$ and let's 
differentiate both sides of it with respect to the variable 
$x^{j_{s+1}}$:
$$
\frac{\partial A^{i_1\ldots\,i_r}_{j_1\ldots\,j_s}(\bold x)}
{\partial x^{j_{s+1}}}=
\sum\Sb p_1\ldots p_r\\ q_1\ldots q_s\endSb
   S^{i_1}_{p_1}\ldots\,S^{i_r}_{p_r}\,\,
   T^{q_1}_{j_1}\ldots\,T^{q_s}_{j_s}\,
\frac{\partial\tilde A^{p_1\ldots\,p_r}_{q_1\ldots\,q_s}
   (\tilde\bold x)}{\partial x^{j_{s+1}}}
$$
In order to calculate the derivative in the right hand side we apply
the chain rule that determines the derivatives of a composite function:
$$
\hskip -2em
\frac{\partial\tilde A^{p_1\ldots\,p_r}_{q_1\ldots\,q_s}
   (\tilde\bold x)}{\partial x^{j_{s+1}}}=
\sum^3_{q_{s+1}=1}
\frac{\partial\tilde x^{q_{s+1}}}{\partial x^{j_{s+1}}}\,
\frac{\partial\tilde A^{p_1\ldots\,p_r}_{q_1\ldots\,q_s}
   (\tilde\bold x)}{\partial\tilde x^{q_{s+1}}}.
\mytag{5.2}
$$
The variables $\bold x=(x^1,x^2,x^3)$ and $\tilde\bold x=(\tilde x^1,
\tilde x^2,\tilde x^3)$ are related as follows:
$$
\xalignat 2
&x^i=\sum^3_{j=1}S^i_j\,\tilde x^j+a^i,
&&\tilde x^i=\sum^3_{j=1}T^i_j\,x^j+\tilde a^i.
\endxalignat
$$
One of these two relationships is included into \mythetag{1.7}, the second 
being the inversion of the first one. The components of the transition
matrices $S$ and $T$ in these formulas are constants, therefore, we have 
$$
\hskip -2em
\frac{\partial\tilde x^{q_{s+1}}}{\partial x^{j_{s+1}}}\,
=T^{q_{s+1}}_{j_{s+1}}.
\mytag{5.3}
$$
Let's substitute \mythetag{5.3} into \mythetag{5.2}, then substitute
the result into the above expression for the derivatives 
$\partial A^{i_1\ldots\,i_r}_{j_1\ldots\,j_s}/\partial x^{j_{s+1}}$. 
This yields the equality
$$
B^{i_1\ldots\,i_r}_{j_1\ldots\,j_s\,j_{s+1}}=
\sum\Sb p_1\ldots p_r\\ q_1\ldots q_{s+1}\endSb
   S^{i_1}_{p_1}\ldots\,S^{i_r}_{p_r}\,\,
   T^{q_1}_{j_1}\ldots\,T^{q_{s+1}}_{j_{s+1}}\,
\tilde B^{p_1\ldots\,p_r}_{q_1\ldots\,j_{s+1}}
$$
which coincides exactly with the transformation rule \mythetag{1.7} applied
to the quantities \mythetag{5.1}. The theorem is proved.
\qed\enddemo
     The passage from $\bold A$ to $\bold B$ in \mythetag{5.1} adds one
covariant index $j_{s+1}$. This is the reason why the tensor field
$\bold B$ is called the {\it covariant differential} of the field $\bold A$.
The covariant differential is denoted as $\bold B=\nabla\bold A$. The
upside-down triangle $\nabla$ is a special symbol, it is called {\it nabla}.
In writing the components of $\bold B$ the additional covariant index
is written beside the nabla sign:
$$
\hskip -2em
B^{i_1\ldots\,i_r}_{j_1\ldots\,j_s\,k}=
\nabla_{\!k}A^{i_1\ldots\,i_r}_{j_1\ldots\,j_s}.
\mytag{5.4}
$$
Due to \mythetag{5.1} the sign $\nabla_{\!k}$ in the formula \mythetag{5.4}
replaces the differentiation operator: $\nabla_k=\partial/\partial x^k$.
However, for $\nabla_{\!k}$ the special name is reserved, it is called
the operator of {\it covariant differentiation} or the {\it covariant
derivative}. Below (in Chapter \uppercase\expandafter{\romannumeral 3})
we shall see that the concept of covariant derivative can be extended so
that it will not coincide with the partial derivative any more.\par
     Let $\bold A$ be a differentiable tensor field of the type $(r,s)$ 
and let $\bold X$ be some arbitrary vector field. Let's consider the tensor product $\nabla\bold A\otimes\bold X$. This is the tensor field of the
type $(r+1,s+1)$. The covariant differentiation adds one covariant index,
while the tensor multiplication add one contravariant index. We denote
by $\nabla_{\bold X}\bold A=C(\nabla\bold A\otimes\bold X)$ the contraction
of the field $\nabla\bold A\otimes\bold X$ with respect to these two
additional indices. The field $\bold B=\nabla_{\bold X}\bold A$ has the
same type $(r,s)$ as the original field $\bold A$. Upon choosing some
Cartesian coordinate system we can write the relationship $\bold
B=\nabla_{\bold X}\bold A$ in coordinate form:
$$
\hskip -2em
B^{i_1\ldots\,i_r}_{j_1\ldots\,j_s}=
\sum^3_{q=1}X^q\,\nabla_qA^{i_1\ldots\,i_r}_{j_1\ldots\,j_s}.
\mytag{5.5}
$$
The tensor field $\bold B=\nabla_{\bold X}\bold A$ with components
\mythetag{5.5} is called the {\it covariant derivative of the field
$\bold A$ along the vector field $\bold X$}.
\mytheorem{5.2} The operation of covariant differentiation
of tensor fields posses\-ses the following properties
\roster
\item $\nabla_{\bold X}(\bold A+\bold B)=\nabla_{\bold X}\bold A
       +\nabla_{\bold X}\bold B$;
\item $\nabla_{\bold X+\bold Y}\bold A=\nabla_{\bold X}\bold A
       +\nabla_{\bold Y}\bold A$;
\item $\nabla_{\xi\cdot\bold X}\bold A=\xi\cdot\nabla_{\bold X}
       \bold A$;
\item $\nabla_{\bold X}(\bold A\otimes\bold B)=\nabla_{\bold X}
       \bold A\otimes\bold B+\bold A\otimes\nabla_{\bold X}
       \bold B$;
\item $\nabla_{\bold X}C(\bold A)=C(\nabla_{\bold X}\bold A)$;
\endroster
where $\bold A$ and $\bold B$ are arbitrary differentiable tensor
fields, while $\bold X$ and $\bold Y$ are arbitrary vector fields
and $\xi$ is an arbitrary scalar field.
\endproclaim
\demo{Proof} It is convenient to carry out the proof of the theorem 
in some Cartesian coordinate system. Let $\bold C=\bold A+\bold B$. 
The property \therosteritem{1} follows from the relationship
$$
\sum^3_{q=1}X^q\,
\frac{\partial C^{i_1\ldots\,i_r}_{j_1\ldots\,j_s}}
{\partial x^q}=
\sum^3_{q=1}X^q\,\frac{\partial A^{i_1\ldots\,i_r}_{j_1
\ldots\,j_s}}{\partial x^q}+
\sum^3_{q=1}X^q\,\frac{\partial B^{i_1\ldots\,i_r}_{j_1
\ldots\,j_s}}{\partial x^q}.
$$
Denote $\bold Z=\bold X+\bold Y$ and then we derive the property
\therosteritem{2} from the relationship
$$
\sum^3_{q=1}Z^q\,
\nabla_qA^{i_1\ldots\,i_r}_{j_1\ldots\,j_s}=
\sum^3_{q=1}X^q\,\nabla_qA^{i_1\ldots\,i_r}_{j_1
\ldots\,j_s}+
\sum^3_{q=1}Y^q\,\nabla_qA^{i_1\ldots\,i_r}_{j_1
\ldots\,j_s}.
$$
In order to prove the property \therosteritem{3} we set 
$\bold Z=\xi\cdot\bold X$. Then
$$
\sum^3_{q=1}Z^q\,
\nabla_qA^{i_1\ldots\,i_r}_{j_1\ldots\,j_s}=
\xi\,\sum^3_{q=1}X^q\,\nabla_qA^{i_1\ldots\,i_r}_{j_1
\ldots\,j_s}.
$$
This relationship is equivalent to the property
\therosteritem{3} in the statement of the theorem.\par
     In order to prove the fourth property in the theorem one should 
carry out the following calculations with the components of $\bold A$,
$\bold B$ and $\bold X$:
$$     
\align
\sum^3_{q=1}X^q\,
&\partial/\partial x^q\left(
A^{i_1\ldots\,i_r}_{j_1\ldots\,j_s}\,
B^{i_{r+1}\ldots\,i_{r+m}}_{j_{s+1}\ldots\,j_{s+n}}
\right)=\left(\,\shave{\sum^3_{q=1}}X^q\,
\frac{\partial A^{i_1\ldots\,i_r}_{j_1\ldots\,j_s}}
{\partial x^q}\right)\times\\
&\times
B^{i_{r+1}\ldots\,i_{r+m}}_{j_{s+1}\ldots\,j_{s+n}}
+A^{i_1\ldots\,i_r}_{j_1\ldots\,j_s}\,
\left(\,\shave{\sum^3_{q=1}}X^q\,
\frac{\partial B^{i_{r+1}\ldots\,i_{r+m}}_{j_{s+1}\ldots
\,j_{s+n}}}{\partial x^q}\right).
\endalign
$$
And finally, the following series of calculations
$$
\align
\sum^3_{q=1}X^q\,
\frac{\partial}{\partial x^q}&\left(\,\shave{\sum^3_{k=1}}
A^{i_1\ldots\,i_{m-1}\,k\,i_m\ldots\,i_{r-1}}_{j_1\ldots\,
j_{n-1}\,k\,j_n\ldots\,j_{s-1}}\right)=\\
&=\sum^3_{k=1}\sum^3_{q=1}X^q\,
\frac{\partial A^{i_1\ldots\,i_{m-1}\,k\,i_m\ldots\,
i_{r-1}}_{j_1\ldots\,j_{n-1}\,k\,j_n\ldots\,j_{s-1}}}
{\partial x^q}
\endalign
$$
proves the fifth property. This completes the proof of the theorem
in whole.
\qed\enddemo
\head
\S~\mysection{6} The metric tensor and the volume pseudotensor.
\endhead
\rightheadtext{\S~6. The metric tensor and the volume pseudotensor.}
      Let $O,\bold e_1,\bold e_2,\bold e_3$ be some Cartesian
coordinate system in the space $\Bbb E$. The space $\Bbb E$ is equipped
with the scalar product. Therefore, the basis $\bold e_1,\bold e_2,\bold
e_3$ of any Cartesian coordinate system has its Gram matrix
$$
\hskip -2em
g_{ij}=(\bold e_i\,|\,\bold e_j).
\mytag{6.1}
$$
The gram matrix $\bold g$ is positive and non-degenerate:
$$
\hskip -2em
\det\bold g>0.
\mytag{6.2}
$$
The inequality \mythetag{6.2} follows from the Silvester criterion
(see \mycite{1}). Under a change of a coordinate system the quantities
\mythetag{6.1} are transformed as the components of a tensor of the
type $(0,2)$. Therefore, we can define the tensor field $\bold g$
whose components in any Cartesian coordinate system are the constant
functions coinciding with the components of the Gram matrix:
$$
g_{ij}(\bold x)=g_{ij}=\const.
$$
The tensor field $\bold g$ with such components is called the {\it
metric tensor}. The metric tensor is a special tensor field. One 
should not define it. Its existence is providentially built into 
the geometry of the space $\Bbb E$.\par
     Since the Gram matrix $\bold g$ is non-degenerate, one can
determine the inverse matrix $\hat\bold g=\bold g^{-1}$. The components
of such matrix are denoted by $g^{ij}$, the indices $i$ and $j$ are
written in the upper position. Then
$$
\hskip -2em
\sum^3_{j=1} g^{ij}\,g_{jk}=\delta^i_j.
\mytag{6.3}
$$
\mytheorem{6.1} The components of the inverse Gram matrix
$\hat\bold g$ are transformed as the components of a tensor field of the
type $(2,0)$ under a change of coordinates.
\endproclaim
\demo{Proof} Let's write the transformation rule \mythetag{1.7} for the
components of the metric tensor $\bold g$:
$$
g_{ij}=\sum^3_{p=1}\sum^3_{q=1}T^p_i\,T^q_j\,\tilde g_{pq}.
$$
In matrix form this relationship is written as
$$
\hskip -2em
\bold g=T^{\tr}\,\tilde\bold g\,T.
\mytag{6.4}
$$
Since $\bold g$, $\tilde\bold g$, and $T$ are non-degenerate, we
can pass to the inverse matrices:
$$
\hskip -2em
\bold g^{-1}=(T^{\tr}\,\tilde\bold g\,T)^{-1}=
S\,\tilde\bold g^{-1}\,S^{\tr}.
\mytag{6.5}
$$
Now we can write \mythetag{6.5} back in coordinate form. This yields
$$
\hskip -2em
g^{ij}=\sum^3_{p=1}\sum^3_{q=1}S^i_p\,S^j_q\,\tilde g^{pq}.
\mytag{6.6}
$$
The relationship \mythetag{6.6} is exactly the transformation rule
\mythetag{1.7} written for the components of a tensor field of the
type $(2,0)$. Thus, the theorem is proved.
\qed\enddemo
     The tensor field $\hat\bold g=\bold g^{-1}$ with the components
$g^{ij}$ is called the {\it inverse metric tensor\/} or the {\it dual 
metric tensor}. The existence of the inverse metric tensor also follows
from the nature of the space $\bold E$ which has the pre-built
scalar product.\par
     Both tensor fields $\bold g$ and $\hat\bold g$ are symmetric. The
symmetruy of $g_{ij}$ with respect to the indices $i$ and $j$ follows
from \mythetag{6.1} and from the properties of a scalar product. The
matrix inverse to the symmetric one is a symmetric matrix too. Therefore,
the components of the inverse metric tensor $g^{ij}$ are also symmetric
with respect to the indices $i$ and $j$.\par
     The components of the tensors $\bold g$ and $\hat\bold g$ in any 
Cartesian coordinate system are constants. Therefore, we have
$$
\xalignat 2
&\hskip -2em
\nabla\bold g=0,
&&\nabla\hat\bold g=0.
\mytag{6.7}
\endxalignat
$$
These relationships follow from the formula \mythetag{5.1} for the 
components of the covariant differential in Cartesian coordinates.
\par
     In the course of analytical geometry (see, for instance, \mycite{4})
the indexed object $\varepsilon_{ijk}$ is usually considered, which is 
called the Levi-Civita symbol. Its nonzero components are determined
by the parity of the transposition of indices:
$$
\hskip -2em
\varepsilon_{ijk}=\varepsilon^{ijk}=
\cases
\format\r&\quad\l\\
0&\text{if $i=j$, \ $i=k$, \ or \ $j=k$,}\\
1&\text{if $(ijk)$ is even, i\.\,e\. \ $\sign(ijk)=1$,}\\
-1&\text{if $(ijk)$ is odd, i\.\,e\. \ $\sign(ijk)=-1$.}
\endcases
\mytag{6.8}
$$
Recall that the Levi-Civita symbol \mythetag{6.8} \pagebreak is used for 
calculating the vectorial product\footnote{ \ It is also called 
the {\it cross product\/} of vectors.} and the mixed product\footnote{ 
\ The mixed product is defined as $(\bold X,\,\bold Y,\,\bold Z)=
(\bold X\,|\,[\bold Y,\,\bold Z])$.} through the coordinates 
of vectors in a rectangular Cartesian coordinate system with a right orthonormal basis $\bold e_1,\,\bold e_2,\,\bold e_3$:
\adjustfootnotemark{-2}
$$
\hskip -2em
\aligned
&[\bold X,\,\bold Y]=\sum^3_{i=1}\bold e_i\left(
\shave{\,\sum^3_{j=1}\sum^3_{k=1}} \varepsilon_{ijk}\,X^j\,Y^k
\right),\\
&(\bold X,\,\bold Y,\,\bold Z)=
\sum^3_{i=1}\sum^3_{j=1}\sum^3_{k=1}
\varepsilon_{ijk}\,X^i\,Y^j\,Z^k.
\endaligned
\mytag{6.9}
$$\par
     The usage of upper or lower indices in writing the components
of the Levi-Civita symbol in \mythetag{6.8} and \mythetag{6.9} makes no
difference since they do not define a tensor. However, there is a
tensorial object associated with the Levi-Civita symbol. In order
to construct such an object we apply the relationship which is\linebreak
usually proved in analytical geometry:
$$
\hskip -2em
\sum^3_{p=1}\sum^3_{q=1}\sum^3_{l=1}\varepsilon_{pql}\,
M_{ip}\,M_{jq}\,M_{kl}=\det\bold M\cdot\varepsilon_{ijk}
\mytag{6.10}
$$
(see proof in \mycite{4}). Here $\bold M$ is some square $3\times 3$
matrix. The matrix $\bold M$ can be the matrix of the components 
for some tensorial field of the type $(2,0)$, $(1,1)$, or $(0,2)$.
However, it can be a matrix without any tensorial interpretation
as well. The relationship \mythetag{6.10} is valid for any square 
$3\times 3$ matrix.\par
     Using the Levi-Civita symbol and the matrix of the metric tensor
$\bold g$ in some Cartesian coordinate system, we  construct the
following quantities:
$$
\hskip -2em
\omega_{ijk}=\sqrt{\det\bold g}\,\,\varepsilon_{ijk}.
\mytag{6.11}
$$
Then we study how the quantities $\omega_{ijk}$ and $\tilde\omega_{pql}$
constructed in two different Cartesian coordinate systems $O,\,\bold
e_1,\,\bold e_2,\,\bold e_3$ and $O',\,\tilde\bold e_1,\,\tilde\bold
e_2,\,\tilde\bold e_3$ are related to each other. From the identity
\mythetag{6.10} we derive
$$
\hskip -2em
\sum^3_{p=1}\sum^3_{q=1}\sum^3_{l=1}
T^p_i\,T^q_j\,T^l_k\,\tilde\omega_{pql}=
\sqrt{\det\tilde\bold g}\,\,\det T\,\,\varepsilon_{ijk}.
\mytag{6.12}
$$
In order to transform further the sum \mythetag{6.12} we use the
relationship \mythetag{6.4}, as an immediate consequence of it we 
obtain the formula $\det\bold g=(\det T)^2\,\det
\tilde\bold g$. Applying this formula to \mythetag{6.12}, we derive
$$
\hskip -2em
\sum^3_{p=1}\sum^3_{q=1}\sum^3_{l=1}
T^p_i\,T^q_j\,T^l_k\,\tilde\omega_{pql}=
\sign(\det T)\,\,\sqrt{\det\bold g}\,\,\varepsilon_{ijk}.
\mytag{6.13}
$$
Note that the right hand side of the relationship \mythetag{6.13} 
differs from $\omega_{ijk}$ in \mythetag{6.11} only by the sign
of the determinant: $\sign(\det T)=\sign(\det S)=\pm 1$. Therefore,
we can write the relationship \mythetag{6.13} as
$$
\hskip -2em
\omega_{ijk}=\sign(\det S)\,\,
\sum^3_{p=1}\sum^3_{q=1}\sum^3_{l=1}
T^p_i\,T^q_j\,T^l_k\,\,\tilde\omega_{pql}.
\mytag{6.14}
$$
Though the difference is only in sign, the relationship \mythetag{6.14}
differs from the transformation rule \mythetag{1.6} for the components
of a tensor of the type $(0,3)$. The formula \mythetag{6.14} gives the
cause for modifying the transformation rule \mythetag{1.6}: 
$$
\hskip -2em
F^{i_1\ldots\,i_r}_{j_1\ldots\,j_s}=
   \sum\Sb p_1\ldots\,p_r\\ q_1\ldots\,q_s\endSb
   (-1)^S\,S^{i_1}_{p_1}\ldots\,S^{i_r}_{p_r}\,
   T^{q_1}_{j_1}\ldots\,T^{q_s}_{j_s}\,\,
   \tilde F^{p_1\ldots\,p_r}_{q_1\ldots\,
   q_s}.
\mytag{6.15}
$$
Here $(-1)^S=\sign(\det S)=\pm 1$. The corresponding modification 
for the concept of a tensor is given by the following definition.
\mydefinition{6.1} A {\it pseudotensor} $\bold F$ of the type 
$(r,s)$ is a geometric object whose components in an arbitrary basis
are enumerated by $(r+s)$ indices and obey the transformation rule 
\mythetag{6.15} under a change of basis.
\enddefinition
     Once some pseudotensor of the type $(r,s)$ is given at each point
of the space $\Bbb E$, we have a pseudotensorial field of the type $(r,s)$. 
Due to the above definition~\mythedefinition{6.1} and due to 
\mythetag{6.14} the quantities $\omega_{ijk}$ from \mythetag{6.11} define
a pseudotensorial field $\boldsymbol\omega$ of the type $(0,3)$. This field
is called the {\it volume pseudotensor}. Like metric tensors $\bold g$ and
$\hat\bold g$, the volume pseudotensor is a special field pre-built
into the space $\Bbb E$. Its existence is due to the existence of the
pre-built scalar product in $\Bbb E$.\par
\head
\S~\mysection{7} The properties of pseudotensors.
\endhead
\rightheadtext{\S~7. The properties of pseudotensors.}
     Pseudotensors and pseudotensorial fields are closely relative
objects for tensors and tensorial fields. In this section we repeat
most of the results of previous sections as applied to pseudotensors.
The proofs of these results are practically the same as in purely 
tensorial case. Therefore, below we do not give the proofs.\par
      Let $\bold A$ and $\bold B$ be two pseudotensorial fields 
of the type $(r,s)$. Then the formula \mythetag{3.1} determines 
a third field $\bold C=\bold A+\bold B$ which appears to be a
pseudotensorial field of the type $(r,s)$. It is important to note 
that \mythetag{3.1} is not a correct procedure if one tries to add
a tensorial field $\bold A$ and a pseudotensorial field $\bold B$.
The sum $\bold A+\bold B$ of such fields can be understood only as 
a formal sum like in \mythetag{3.5}.\par
      The formula \mythetag{2.2} for the tensor product appears to
be more universal. It defines the product of a field $\bold A$ of 
the type $(r,s)$ and a field $\bold B$ of the type $(m,n)$. 
Therein each of the fields can be either a tensorial or a 
pseudotensorial field. The tensor product possesses the following
properties:
\roster
\item the tensor product of two tensor fields is a tensor field;
\item the tensor product of two pseudotensorial fields is a tensor 
      field;
\item the tensor product of a tensorial field and a pseudotensorial 
      field is a pseudotensorial field.
\endroster\par
     Let's denote by $P_{(r,s)}$ the set of pseudotensorial fields
of the type $(r,s)$. Due to the properties \therosteritem{1}-\therosteritem{3} and due to the distributivity
relationships \mythetag{3.3}, which remain valid for pseudotensorial 
fields too, the set $P_{(r,s)}$ is a module over the ring of scaral
fields $K=T_{(0,0)}$. As for the properties
\therosteritem{1}-\therosteritem{3}, they can be expressed in form
of the relationships 
$$
\xalignat 2
&\hskip -2em
\aligned
&T_{(r,s)}\times T_{(m,n)}\to T_{(r+m,s+n)},\\
&P_{(r,s)}\times P_{(m,n)}\to T_{(r+m,s+n)},
\endaligned
&&\aligned
&T_{(r,s)}\times P_{(m,n)}\to P_{(r+m,s+n)},\\
&P_{(r,s)}\times T_{(m,n)}\to P_{(r+m,s+n)}
\endaligned\quad
\mytag{7.1}
\endxalignat
$$
that extend the relationship \mythetag{3.2} from the section~3.\par
     The formula \mythetag{2.5} defines the operation of contraction
for a field $\bold F$ of the type $(r,s)$, where $r\geqslant 1$ 
and $s\geqslant 1$. The operation of contraction \mythetag{2.5} is
applicable to tensorial and pseudotensorial fields. It has the
following properties:
\roster
\item the contraction of a tensorial field is a tensorial field;
\item the contraction of a pseudotensorial field is a pseudotensorial field.
\endroster
The operation of contraction extended to the case of pseudotensorial fields
preserve its linearity given by the equalities \mythetag{3.4}.\par
     The covariant differentiation of pseudotensorial fields in a Cartesian
coordinate system is determined by the formula \mythetag{5.1}. The covariant
differential $\nabla\bold A$ of a tensorial field is a tensorial field;
the covariant differential of a pseudotensorial field is a pseudotensorial
field. It is convenient to express the properties of the covariant
differential through the properties of the covariant derivative
$\nabla_{\bold X}$ in the direction of a field $\bold X$. Now $\bold X$ is 
either a vectorial or a pseudovectorial field. All propositions of the
theorem~\mythetheorem{5.2} for $\nabla_{\bold X}$ remain valid.\par
\head
\S~\mysection{8} A note on the orientation.
\endhead
\rightheadtext{\S~8. A note on the orientation.}
     Pseudoscalar fields form a particular case of pseudotensorial
fields. Scalar fields can be interpreted as functions whose argument
is a point of the space $\Bbb E$. In this interpretation they do not 
depend on the choice of a coordinate system. Pseudoscalar fields even
in such interpretation preserve some dependence on a coordinate system,
though this dependence is rather weak. Let $\xi$ be a pseudoscalar 
field. In a fixed Cartesian coordinate system the field $\xi$ is
represented by a scalar function $\xi(P)$, where $P\in\Bbb E$. The value 
of this function $\xi$ at a point $P$ does not change if we pass to
another coordinate system of the same orientation, i\.\,e\. if
the determinant of the transition matrix $S$ is positive. When passing
to a coordinate system of the opposite orientation the function $\xi$
changes the sign: $\xi(P)\to -\xi(P)$. Let's consider a nonzero constant
pseudoscalar field $\xi$. In some coordinate systems $\xi=c=\const$,
in others $\xi=-c=\const$. Without loss of generality we can take $c=1$.
Then such a pseudoscalar field $\xi$ can be used to distinguish the 
coordinate systems where $\xi=1$ from those of the opposite orientation
where $\xi=-1$.
\proclaim{Proposition} Defining a unitary constant pseudoscalar field
$\xi$ is equivalent to choosing some preferable orientation in the
space $\Bbb E$.
\endproclaim
     From purely mathematical point of view the space $\Bbb E$, which is
a three-dimensional Euclidean point space (see definition in \mycite{1}),
has no preferable orientation. However, the real physical space $\Bbb E$
(where we all live) has such an orientation. Therein we can distinguish
the left hand from the right hand. This difference in the nature is not formal and purely terminological: the left hemisphere of a human brain 
is somewhat different from the right hemisphere in its functionality, in
many substances of the organic origin some isomers prevail over the
mirror symmetric isomers. The number of left-handed people and the number
of right-handed people in the mankind is not fifty-fifty as well. The asymmetry of the left and right is observed even in basic forms of the
matter:  it is reflected in modern theories of elementary particles.
Thus, we can assume the space $\Bbb E$ to be canonically equipped with 
some pseudoscalar field $\xi_E$ whose values are given by the formula
$$
\xi_E=
\cases
\format\r&\quad\l\\
1 &\text{in right-oriented coordinate systems},\\
-1&\text{in left-oriented coordinate systems}.
\endcases
$$\par
     Multiplying by $\xi_E$, we transform a tensorial field $\bold A$ into
the pseudotensorial field $\xi_E\otimes\bold A=\xi_E\cdot\bold A$. Multiplying by $\xi_E$ once again, we transform $\xi_E\cdot\bold A$ back to
$\bold A$. Therefore, in the space $\Bbb E$ equipped with the preferable
orientation in form of the field $\xi_E$ one can not to consider 
pseudotensors at all, considering only tensor fields. The components of
the volume tensor $\boldsymbol\omega$ in this space should be defined as
$$
\hskip -2em
\omega_{ijk}=\xi_E\,\sqrt{\det\bold g}\,\,\varepsilon_{ijk}.
\mytag{8.1}
$$
Let $\bold X$ and $\bold Y$ be two vectorial fields. Then we can define 
the vectorial field $\bold Z$\linebreak with the following components:
$$
\hskip -2em
Z^q=\sum^3_{i=1}\sum^3_{j=1}\sum^3_{k=1}
g^{ki}\,\omega_{ijk}\,\,X^j\,Y^k.
\mytag{8.2}
$$
From \mythetag{8.2}, it is easy to see that $\bold Z$ is derived as the
contraction of the field $\hat\bold g\otimes\boldsymbol\omega\otimes
\bold a\otimes\bold b$. In a rectangular Cartesian coordinate system 
with right-oriented orthonormal basis the formula \mythetag{8.2} takes
the form of the well-known formula for the components of the vector
product $\bold Z=[\bold X,\,\bold Y]$ (see \mycite{4} and the formula
\mythetag{6.9} above). In a space without preferable orientation, where
$\omega_{ijk}$ is given by the formula \mythetag{6.11}, the vector product
of two vectors is a pseudovector.\par
     Now let's consider three vectorial fields $\bold X$, $\bold Y$, and
$\bold Z$ and let's construct the scalar field $u$ by means of the
following formula:
$$
\hskip -2em
u=\sum^3_{i=1}\sum^3_{j=1}\sum^3_{k=1}
\omega_{ijk}\,\,X^i\,Y^j\,Z^k.
\mytag{8.3}
$$
When \mythetag{8.3} is written in a rectangular Cartesian coordinate system 
with right-oriented orthonormal basis, one easily sees that the field 
\mythetag{8.3} coincides with the mixed product $(\bold X,\,\bold Y\,\bold Z)$
of three vectors (see \mycite{4} and the formula \mythetag{6.9}
above). In a space without preferable orientation the mixed product of 
three vector fields determined by the volume pseudotensor \mythetag{6.11}
appears to be a pseudoscalar field.\par
\head
\S~\mysection{9} Raising and lowering indices.
\endhead
\rightheadtext{\S~9. Raising and lowering indices.}
     Let $\bold A$ be a tensor field or a pseudotensor field of the type
$(r,s)$ in the space $\Bbb E$ and let $r\geqslant 1$. Let's construct the
tensor product $\bold A\otimes\bold g$ of the field $A$ and the metric
tensor $\bold g$, then define the field $\bold B$ of the type $(r-1,s+1)$ 
as follows:
$$
\hskip -2em
B^{i_1\ldots\,i_{r-1}}_{j_1\ldots\,j_{s+1}}=
\sum^3_{k=1}
A^{i_1\ldots\,i_{m-1}\,k\,i_m\ldots\,i_{r-1}}_{j_1\ldots\,
j_{n-1}\,j_{n+1}\ldots\,j_{s+1}}\,\,
g_{kj_n}.
\mytag{9.1}
$$
The passage from the field $\bold A$ to the field $\bold B$ according to
the formula \mythetag{9.1} is called the {\it index lowering procedure of
the $m$-th upper index to the $n$-th lower position}.\par
     Using the inverse metric tensor, one can invert the operation 
\mythetag{9.1}. Let $\bold B$ be a tensorial or a pseudotensorial field
of the type $(r,s)$ and let $s\geqslant 1$. Then we define the field
$\bold A=C(\bold B\otimes\hat\bold g)$ of the type $(r+1,s-1)$ according
to the formula:
$$
A^{i_1\ldots\,i_{r+1}}_{j_1\ldots\,j_{s-1}}=
\sum^3_{q=1}
B^{i_1\ldots\,i_{m-1}\,i_{m+1}\ldots\,i_{r+1}}_{j_1\ldots\,
j_{n-1}\,q\,j_n\ldots\,j_{s-1}}\,\,g^{qi_m}.
\mytag{9.2}
$$
The passage from the field $\bold B$ to the field $\bold A$ according to 
the formula \mythetag{9.2} is called the {\it index raising procedure of 
the $n$-th lower index to the $m$-th upper position}.\par
     The operations of lowering and raising indices are inverse to each
other. Indeed, we can perform the following calculations:
$$
\align
C^{i_1\ldots\,i_r}_{j_1\ldots\,j_s}&=
\sum^3_{q=1}\sum^3_{k=1}
A^{i_1\ldots\,i_{m-1}\,k\,i_{m+1}\ldots\,i_r}_{j_1\ldots\,
j_{n-1}\,j_n\ldots\,j_s}\,\,
g_{kq}\,g^{qi_m}=\\
\allowdisplaybreak
&=\sum^3_{k=1}
A^{i_1\ldots\,i_{m-1}\,k\,i_{m+1}\ldots\,i_r}_{j_1\ldots\,
j_{n-1}\,j_n\ldots\,j_s}\,\,
\delta^{i_m}_k=A^{i_1\ldots\,i_r}_{j_1\ldots\,j_s}.
\endalign
$$
The above calculations show that applying the index lowering and the
index raising procedures successively $\bold A\to\bold B\to\bold C$, 
we get the field $\bold C=\bold A$. Applying the same procedures in the reverse order yields the same result. This follows from the calculations just below:
$$
\align
C^{i_1\ldots\,i_r}_{j_1\ldots\,j_s}&=
\sum^3_{k=1}\sum^3_{q=1}
A^{i_1\ldots\,i_{m-1}\,i_m\ldots\,i_r}_{j_1\ldots\,
j_{n-1}\,q\,j_{n+1}\ldots\,j_s}\,\,
g^{qk}\,g_{kj_n}=\\
\allowdisplaybreak
&=\sum^3_{q=1}
A^{i_1\ldots\,i_{m-1}\,i_m\ldots\,i_r}_{j_1\ldots\,
j_{n-1}\,q\,j_{n+1}\ldots\,j_s}\,\,
\delta^q_{i_n}=A^{i_1\ldots\,i_r}_{j_1\ldots\,j_s}.
\endalign
$$\par
     The existence of the index raising and index lowering procedures
follows from the very nature of the space $\Bbb E$ which is equipped
with the scalar product and, hence, with the metric tensors $\bold g$ 
and $\hat\bold g$. Therefore, any tensorial (or pseudotensorial) field
of the type $(r,s)$ in such a space can be understood as originated
from some purely covariant field of the type $(0,r+s)$ as a result of
raising a part of its indices. Therein a little bit different way of 
setting indices is used. Let's consider a field $\bold A$ of the type
$(0,4)$ as an example. We denote by $A_{i_1i_2i_3i_4}$ its components in
some Cartesian coordinate system. By raising one of the four indices in
$\bold A$ one can get the four fields of the type $(1,3)$. Their components
are denoted as 
$$
\hskip -2em
A^{i_1\hphantom{i_2}\hphantom{i_3}
\hphantom{i_4}}_{\hphantom{i_1}i_2i_3i_4},
\quad
A^{\hphantom{i_1}i_2\hphantom{i_3}
\hphantom{i_4}}_{i_1\hphantom{i_2}i_3i_4},
\quad
A^{\hphantom{i_1}\hphantom{i_2}i_3
\hphantom{i_4}}_{i_1i_2\hphantom{i_3}i_4},
\quad
A^{\hphantom{i_1}\hphantom{i_2}\hphantom{i_3}
i_4}_{i_1i_2i_3\hphantom{i_4}}.
\mytag{9.3}
$$
Raising one of the indices in \mythetag{9.3}, we get an empty place underneath
it in the list of lower indices, while the numbering of the indices at that
place remains unbroken. In this way of writing indices, each index has 
{\tencyr\char '074}its fixed position{\tencyr\char '076} no matter what
index it is --- a lower or an upper index. Therefore, in the writing below
we easily guess the way in which the components of tensors are derived:
$$
\pagebreak
\hskip -2em
A^{i\hphantom{j}\hphantom{k}
\hphantom{q}}_{\hphantom{i}jkq},\
A^{ij\hphantom{k}
\hphantom{q}}_{\hphantom{i}\hphantom{j}kq},\
A^{i\hphantom{j}k\hphantom{q}}_{\hphantom{i}j
\hphantom{k}q},\
A^{\hphantom{i}jk\hphantom{q}}_{i\hphantom{j}
\hphantom{k}q}.
\mytag{9.4}
$$
Despite to some advantages of the above form of index setting
in \mythetag{9.3} and \mythetag{9.4}, it is not commonly admitted.
The matter is that it has a number of disadvantages either.
For example, the writing of general formulas \mythetag{1.6}, 
\mythetag{2.2}, \mythetag{2.5}, and some others becomes huge and
inconvenient for perception. In what follows we shall not change
the previous way of index setting.\par
\head
\S~\mysection{10} Gradient, divergency, and rotor. Some identities
of the vectorial analysis.
\endhead
\rightheadtext{\S~10. Gradient, divergency and rotor.}
     Let's consider a scalar field or, in other words, a function
$f$. Then apply the operator of covariant differentiation $\nabla$, 
as defined in \S\,5, to the field $f$. The covariant differential 
$\nabla f$ is a covectorial field (a field of the type $(0,1)$). 
Applying the index raising procedure \mythetag{9.2} to the covector
field $\nabla f$, we get the vector field $\bold F$. Its components
are given by the following formula:
$$
\hskip -2em
F^i=\sum^3_{k=1}g^{ik}\,\frac{\partial f}{\partial x^k}.
\mytag{10.1}
$$
\mydefinition{10.1} The vector field $\bold F$ in the space
$\Bbb E$ whose components are calculated by the formula \mythetag{10.1}
is called the {\it gradient} of a function $f$. The gradient is denoted 
as $\bold F=\grad f$.
\enddefinition
    Let $\bold X$ be a vectorial field in $\Bbb E$. Let's consider the
scalar product of the vectorial fields $\bold X$ and $\grad f$. Due to
the formula \mythetag{10.1} such scalar product of two vectors is reduced
to the scalar product of the vector $\bold X$ and the covector
$\nabla f$:
$$
\hskip -2em
(\bold X\,|\,\grad f)=\sum^3_{k=1} X^k\,\frac{\partial f}
{\partial x^k}=\langle\nabla f\,|\,\bold X\rangle.
\mytag{10.2}
$$
The quantity \mythetag{10.2} is a scalar quantity. It does not depend on
the choice of a coordinate system where the components of $\bold X$ and
$\nabla f$ are given. Another form of writing the formula \mythetag{10.2}
is due to the covariant differentiation along the vector field $\bold X$,
it was introduced by the formula \mythetag{5.5} above:
$$
\hskip -2em
(\bold X\,|\,\grad f)=\nabla_{\bold X}f.
\mytag{10.3}
$$
By analogy with the formula \mythetag{10.3}, the covariant differential 
$\nabla\bold F$ of an arbitrary tensorial field $\bold F$ is sometimes 
called the {\it covariant gradient} of the field $\bold F$.\par
     Let $\bold F$ be a vector field. Then its covariant differential 
$\nabla\bold F$ is an operator field, i\.\,e\. a field of the type
$(1,1)$. Let's denote by $\varphi$ the contraction of the field 
$\nabla\bold F$:
$$
\hskip -2em
\varphi=C(\nabla\bold F)=\sum^3_{k=1}\frac{\partial F^k}
{\partial x^k}.
\mytag{10.4}
$$
\mydefinition{10.2} The scalar field $\varphi$ in the space
$\Bbb E$ determined by the formula \mythetag{10.4} is called the
{\it divergency} of a vector field $\bold F$. It is denoted $\varphi
=\divr\bold F$.
\enddefinition
     Apart from the scalar field $\divr\bold F$, one can use $\nabla
\bold F$ in order to build a vector field. Indeed, let's consider
the quantities
$$
\pagebreak
\hskip -2em
\rho^m=\sum^3_{i=1}\sum^3_{j=1}\sum^3_{k=1}\sum^3_{q=1}
g^{mi}\,\omega_{ijk}\,g^{jq}\,
\nabla_qF^k,
\mytag{10.5}
$$
where $\omega_{ijk}$ are the components of the volume tensor given
by the formula \mythetag{8.1}.\par
\mydefinition{10.3} The vector field $\boldsymbol\rho$ in the
space $\Bbb E$ determined by the formula \mythetag{10.5} is called the
{\it rotor\footnote{ \ The term {\eightcyr\char '074}curl{\eightcyr
\char '076} is also used for the rotor.}} of a vector field $\bold F$. 
It is denoted $\boldsymbol\rho=\rot\bold F$.
\enddefinition
\adjustfootnotemark{-1}
     Due to \mythetag{10.5} the rotor or a vector field $\bold F$ is
the contraction of the tensor field $\hat\bold g\otimes\boldsymbol
\omega\otimes\hat\bold g\otimes\nabla\bold F$ with respect to four
pairs of indices: $\rot\bold F=C(\hat\bold g\otimes\boldsymbol\omega
\otimes\hat\bold g\otimes\nabla\bold F)$.
\subhead Remark\endsubhead If $\omega_{ijk}$ in \mythetag{10.5} are
understood as components of the volume pseudotensor \mythetag{6.11}, 
then the rotor of a vector field should be understood as a 
pseudovectorial field.\par
     Suppose that $O,\,\bold e_1,\,\bold e_2,\,\bold e_3$ is a
rectangular Cartesian coordinate system in $\Bbb E$ with orthonormal 
right-oriented basis $\bold e_1,\,\bold e_2,\,\bold e_3$. The Gram 
matrix of the basis $\bold e_1,\,\bold e_2,\,\bold e_3$ is the unit 
matrix. Therefore, we have
$$
g_{ij}=g^{ij}=\delta^i_j=
\cases
1&\text{for \ $i=j$,}\\
0&\text{for \ $i\neq j$.}\\
\endcases
$$
The pseudoscalar field $\xi_E$ defining the orientation in $\Bbb E$ 
is equal to unity in a right-oriented coordinate system: $\xi_E
\equiv 1$. Due to these circumstances the formulas \mythetag{10.1} and
\mythetag{10.5} for the components of $\grad f$ and $\rot\bold F$ 
simplifies substantially:
$$
\align
&\hskip -2em
(\grad f)^i=\frac{\partial f}{\partial x^i},
\mytag{10.6}\\
&\hskip -2em
(\rot\bold F)^i=\sum^3_{j=1}\sum^3_{k=1}
\varepsilon_{ijk}\,\frac{\partial F^k}{\partial x^j}.
\mytag{10.7}
\endalign
$$
The formula \mythetag{10.4} for the divergency remains unchanged:
$$
\hskip -2em
\divr\bold F=\sum^3_{k=1}\frac{\partial F^k}
{\partial x^k}.
\mytag{10.8}
$$
The formula \mythetag{10.7} for the rotor has an elegant representation
in form of the determinant of a $3\times 3$ matrix:
$$
\hskip -2em
\rot\bold F=
\vmatrix
\bold e_1 &\bold e_2&\bold e_3\\
\vspace{1ex}
\dsize\frac{\partial}{\partial x^1}&
\dsize\frac{\partial}{\partial x^2}&
\dsize\frac{\partial}{\partial x^3}\\
\vspace{1.2ex}
F^1 &F^2 &F^3
\endvmatrix.
\mytag{10.9}
$$
The formula \mythetag{8.2} for the vector product in right-oriented 
rectangular Cartesian coordinate system takes the form of \mythetag{6.9}.
It can also be represented in the form of the formal determinant of
a $3\times 3$ matrix:
$$
\pagebreak
\hskip -2em
[\bold X,\,\bold Y]=
\vmatrix
\bold e_1 &\bold e_2&\bold e_3\\
\vspace{1ex}
X^1 &X^2 &X^3\\
\vspace{1.2ex}
Y^1 &Y^2 &Y^3
\endvmatrix.
\mytag{10.10}
$$
Due to the similarity of \mythetag{10.9} and \mythetag{10.10} one can
formally represent the operator of covariant differentiation $\nabla$
as a vector with components $\partial/\partial x^1,\,\partial
/\partial x^2,\,\partial/\partial x^3$. Then the divergency and rotor 
are represented as the scalar and vectorial products:
$$
\xalignat 2
&\divr\bold F=(\nabla\,|\,\bold F),
&&\rot\bold F=[\nabla,\,\bold F].
\endxalignat
$$\par
\mytheorem{10.1} For any scalar field $\varphi$ of the smoothness
class $C^2$ the equality \ $\rot\grad\varphi=0$ is identically fulfilled.
\endproclaim
\demo{Proof} Let's choose some right-oriented rectangular Cartesian 
coordinate system and then use the formulas \mythetag{10.6} and\mythetag{10.7}.
Let $\bold F=\rot\grad\varphi$. Then
$$
\hskip -2em
F^i=\sum^3_{j=1}\sum^3_{k=1}
\varepsilon_{ijk}\,\frac{\partial^2\varphi}{\partial x^j\,
\partial x^k}.
\mytag{10.11}
$$
Let's rename the summation indices in \mythetag{10.11}. The index $j$ is
replaced by the index $k$ and vice versa. Such a swap of indices does
not change the value of the sum in \mythetag{10.11}. Therefore, we have
$$
F^i=\sum^3_{j=1}\sum^3_{k=1}
\varepsilon_{ikj}\,\frac{\partial^2\varphi}{\partial x^k\,
\partial x^j}=-\sum^3_{j=1}\sum^3_{k=1}
\varepsilon_{ijk}\,\frac{\partial^2\varphi}{\partial x^j\,
\partial x^k}=-F^i.
$$
Here we used the skew-symmetry of the Levi-Civuta symbol with respect
to the pair of indices $j$ and $k$ and the symmetry of the second order
partial derivatives of the function $\varphi$ with respect to the same
pair of indices:
$$
\hskip -2em
\frac{\partial^2\varphi}{\partial x^j\,
\partial x^k}=
\frac{\partial^2\varphi}{\partial x^k\,
\partial x^j}.
\mytag{10.12}
$$
For $C^2$ class functions the value of second order partial derivatives
\mythetag{10.12} do not depend on the order of differentiation. The equality
$F^i=-F^i$ now immediately yields $F^i=0$. The theorem is proved.
\qed\enddemo
\mytheorem{10.2} For any vector field $\bold F$ of the smoothness
class $C^2$ the equality \ $\divr\rot\bold F=0$ is identically fulfilled.
\endproclaim
\demo{Proof} Here, as in the case of the theorem~\mythetheorem{10.1}, we
choose a right-oriented rectangular Cartesian coordinate system, then we 
use the formulas \mythetag{10.7} and \mythetag{10.8}. For the scalar field
$\varphi=\divr\rot\bold F$ from these formulas we derive
$$
\hskip -2em
\varphi=\sum^3_{i=1}\sum^3_{j=1}\sum^3_{k=1}
\varepsilon_{ijk}\,\frac{\partial^2 F^k}{\partial x^j\,
\partial x^i}.
\mytag{10.13}
$$
Using the relationship analogous to \mythetag{10.12} for the
partial derivatives
$$
\frac{\partial^2 F^k}{\partial x^j\,
\partial x^i}=
\frac{\partial^2 F^k}{\partial x^i\,
\partial x^j}
$$
and using the skew-symmetry of $\varepsilon_{ijk}$ with respect to
indices $i$ and $j$, from the formula \mythetag{10.13} we easily
derive $\varphi=-\varphi$. Hence, $\varphi=0$. \qed\enddemo
     Let $\varphi$ be a scalar field of the smoothness class $C^2$.
The quantity $\divr\grad\varphi$ in general case is nonzero. It is
denoted $\triangle\varphi=\divr\grad\varphi$. The sign $\triangle$ 
denotes the differential operator of the second order that transforms
a scalar field $\varphi$ to another scalar field $\divr\grad\varphi$.
It is called the {\it Laplace operator} or the {\it laplacian}. In a
rectangular Cartesian coordinate system it is given by the formula 
$$
\hskip -2em
\triangle=\left(\frac{\partial}{\partial x^1}\right)^{\!2}+
\left(\frac{\partial}{\partial x^2}\right)^{\!2}+
\left(\frac{\partial}{\partial x^3}\right)^{\!2}.
\mytag{10.14}
$$
Using the formulas \mythetag{10.6} and \mythetag{10.8} one can calculate the
Laplace operator in a skew-angular Cartesian coordinate system:
$$
\hskip -2em
\triangle=\sum^3_{i=1}\sum^3_{j=1}
g^{ij}\,\frac{\partial^2}{\partial x^i\,\partial x^j}.
\mytag{10.15}
$$
Using the signs of covariant derivatives $\nabla_i=\partial/\partial x^i$
we can write the Laplace operator \mythetag{10.15} as follows:
$$
\triangle=\sum^3_{i=1}\sum^3_{j=1}
g^{ij}\,\nabla_i\nabla_j.
\mytag{10.16}
$$
The equality \mythetag{10.16} differs from \mythetag{10.15} not only in
special notations for the derivatives. The Laplace operator defined
as $\triangle\varphi=\divr\grad\varphi$ can be applied only to a scalar 
field $\varphi$. The formula \mythetag{10.16} extends it, and now we can
apply the Laplace operator to any twice continuously differentiable
tensor field $\bold F$ of any type $(r,s)$. Due to this formula 
$\triangle\bold F$ is the result of contracting the tensor product
$\hat\bold g\otimes\nabla\nabla\bold F$ with respect to two pairs of indices: $\triangle\bold F=C(\hat\bold g\otimes\nabla\nabla\bold F)$. 
The resulting field $\triangle\bold F$ has the same type $(r,s)$ as
the original field $\bold F$. The laplace operator in the form of
\mythetag{10.16} is sometimes called the {\it Laplace-Beltrami} operator.
\par
\head
\S~\mysection{11} Potential and vorticular vector fields.
\endhead
\rightheadtext{\S~11. Potential and vorticular vector fields.}
\mydefinition{11.1} A differentiable vector field $\bold F$ in
the space $\Bbb E$ is called a {\it potential} field if \ $\rot\bold F=0$.
\enddefinition
\mydefinition{11.2} A differentiable vector field $\bold F$ in
the space $\Bbb E$ is called a {\it vorticular} field if \ $\divr\bold F=0$.
\enddefinition
     The theorem~\mythetheorem{10.1} yields some examples of potential
vector fields, while the theorem~\mythetheorem{10.2} yields the examples 
of vorticular fields. Indeed, any field of the form $\grad\varphi$ is a potential field, and any field of the form $\rot\bold F$ is a vorticular
one. As it appears, the theorems~\mythetheorem{10.1} and 
\mythetheorem{10.2} can be strengthened. 
\mytheorem{11.1} Any potential vector field $\bold F$ in the space
$\Bbb E$ is a gradient of some scalar field $\varphi$, i\.\,e\. $\bold
F=\grad\varphi$.
\endproclaim
\demo{Proof} Let's choose a rectangular Cartesian coordinate system
$O,\,\bold e_1,\,\bold e_2,\,\bold e_3$ with orthonormal right-oriented
basis $\bold e_1,\,\bold e_2,\,\bold e_3$. In this coordinate system
the potentiality condition $\rot\bold F=0$ for the vector field $\bold F$ 
is equivalent tho the following three relationships for its components:
$$
\align
&\hskip -2em\frac{\partial F^1(x^1,x^2,x^3)}{\partial x^2}=
\frac{\partial F^2(x^1,x^2,x^3)}{\partial x^1},
\mytag{11.1}\\
\vspace{1ex}
&\hskip -2em\frac{\partial F^2(x^1,x^2,x^3)}{\partial x^3}=
\frac{\partial F^3(x^1,x^2,x^3)}{\partial x^2},
\mytag{11.2}\\
\vspace{1ex}
&\hskip -2em\frac{\partial F^3(x^1,x^2,x^3)}{\partial x^1}=
\frac{\partial F^1(x^1,x^2,x^3)}{\partial x^3}.\mytag{11.3}
\endalign
$$
The relationships \mythetag{11.1}, \mythetag{11.2}, and \mythetag{11.3} are
easily derived from \mythetag{10.7} or from \mythetag{10.9}. Let's define
the function $\varphi(x^1,x^2,x^3)$ as the sum of three integrals: 
$$
\hskip -2em
\aligned
&\varphi(x^1,x^2,x^3)=c+\int\limits^{\ x_1}_{0\,}
F^1(x^1,0,0)\,dx^1+\\
&+\int\limits^{\ x_2}_{0\,} F^2(x^1,x^2,0)\,dx^2+
\int\limits^{\ x_3}_{0\,} F^3(x^1,x^2,x^3)\,dx^3.
\endaligned
\mytag{11.4}
$$
Here $c$ is an arbitrary constant. Now we only have to check up that
the function \mythetag{11.4} is that very scalar field for which
$\bold F=\grad\varphi$.\par
      Let's differentiate the function $\varphi$ with respect to
the variable $x^3$. The constant $c$ and the first two integrals
in \mythetag{11.4} do not depend on $x^3$. Therefore, we have
$$
\hskip -2em
\frac{\partial\varphi}{\partial x^3}=
\frac{\partial}{\partial x^3}\int\limits^{\ x_3}_{0\,}
F^3(x^1,x^2,x^3)\,dx^3=F^3(x^1,x^2,x^3).
\mytag{11.5}
$$
In deriving the relationship \mythetag{11.5} we used the rule of differentiation of an integral with variable upper limit (see
\mycite{2}).\par
     Now let's differentiate the function $\varphi$ with respect to
$x^2$. The constant $c$ and the first integral in \mythetag{11.4} does
not depend on $x^2$. Differentiating the rest two integrals, we get
the following expression:
$$
\frac{\partial\varphi}{\partial x^2}=
F^2(x^1,x^2,0)+\frac{\partial}{\partial x^2}
\int\limits^{\ x_3}_{0\,} F^3(x^1,x^2,x^3)\,dx^3.
$$
The operations of differentiation with respect to $x^2$ and 
integration with respect to $x^3$ in the above formula are
commutative (see \mycite{2}). Therefore, we have
$$
\hskip -2em
\frac{\partial\varphi}{\partial x^2}=
F^2(x^1,x^2,0)+\int\limits^{\ x_3}_{0\,}
\frac{\partial F^3(x^1,x^2,x^3)}{\partial x^2}\,dx^3.
\mytag{11.6}
$$
In order to transform the expression being integrated in \mythetag{11.6}
we use the formula \mythetag{11.2}. This leads to the following result:
$$
\hskip -2em
\aligned
&\frac{\partial\varphi}{\partial x^2}=
F^2(x^1,x^2,0)+\int\limits^{\ x_3}_{0\,}
\frac{\partial F^2(x^1,x^2,x^3)}{\partial x^3}\,dx^3=\\
\vspace{1ex}
&=F^2(x^1,x^2,0)+F^2(x^1,x^2,x)\,\vbox{\hrule width 0.5pt
height 12pt depth 8pt}^{\,\,x=x_3}_{\,\,x=0}=
F^2(x^1,x^2,x^3).
\endaligned
\mytag{11.7}
$$
In calculating the derivative $\partial\varphi/\partial x^1$ we use that
same tricks as in the case of the other two derivatives $\partial\varphi/
\partial x^3$ and $\partial\varphi/\partial x^2$:
$$
\aligned
\frac{\partial\varphi}{\partial x^1}&=
\frac{\partial}{\partial x^1}\int\limits^{\ x_1}_{0\,}
F^1(x^1,0,0)\,dx^1+\frac{\partial}{\partial x^1}
\int\limits^{\ x_2}_{0\,} F^2(x^1,x^2,0)\,dx^2+\\
&+\frac{\partial}{\partial x^1}\int\limits^{\ x_3}_{0\,}
F^3(x^1,x^2,x^3)\,dx^3=F^1(x^1,0,0)+\\
&+\int\limits^{\ x_2}_{0\,}\frac{\partial F^2(x^1,x^2,0)}
{\partial x^1}\,dx^2+\int\limits^{\ x_3}_{0\,}
\frac{\partial F^3(x^1,x^2,x^3)}{\partial x^1}\,dx^3.
\endaligned
$$
To transform the last two integrals we use the relationships
\mythetag{11.1} and \mythetag{11.3}:
$$
\hskip -2em
\aligned
\frac{\partial\varphi}{\partial x^1}&=
F^1(x^1,0,0)+
F^1(x^1,x,0)\,\vbox{\hrule width 0.5pt
height 12pt depth 8pt}^{\,\,x=x_2}_{\,\,x=0}+\\
\vspace{1ex}
&+F^1(x^1,x^2,x)\,\vbox{\hrule width 0.5pt
height 12pt depth 8pt}^{\,\,x=x_3}_{\,\,x=0}=
F^1(x^1,x^2,x^3).
\endaligned
\mytag{11.8}
$$
The relationships \mythetag{11.5}, \mythetag{11.7}, and \mythetag{11.8}
show that $\grad\varphi=\bold F$ for the function $\varphi(x^1,x^2,x^3)$
given by the formula \mythetag{11.4}. The theorem is proved.\qed\enddemo
\mytheorem{11.2} Any vorticular vector field $\bold F$ in the space
$\Bbb E$ is the rotor of\linebreak some other vector field $\bold A$, i\.\,e\.
$\bold F=\rot\bold A$.
\endproclaim
\demo{Proof} We perform the proof of this theorem in some rectangular
Cartesian coordinate system with the orthonormal basis $\bold e_1,\bold
e_2,\bold e_3$. The condition of vorticity for the field $\bold F$ in
such a coordinate system is expressed by a single equation:
$$
\hskip -2em
\frac{\partial F^1(\bold x)}{\partial x^1}+
\frac{\partial F^2(\bold x)}{\partial x^2}+
\frac{\partial F^3(\bold x)}{\partial x^3}=0.
\mytag{11.9}
$$
Let's construct the vector field $\bold A$ defining its components 
in the chosen coordinate system by the following three formulas:
$$
\align
\hskip -1em A^1&=\int\limits^{\ x_3}_{0\,}
F^2(x^1,x^2,x^3)\,dx^3-\int\limits^{\ x_2}_{0\,}
F^3(x^1,x^2,0)\,dx^2,
\hskip -5em\\
\hskip -1em A^2&=-\int\limits^{\ x_3}_{0\,}
F^1(x^1,x^2,x^3)\,dx^3,\mytag{11.10}\\
\vspace{3ex}
\hskip -1em A^3&=0.
\endalign
$$
Let's show that the field $\bold A$ with components \mythetag{11.10}
is that very field for which $\,\rot\bold A=\bold F$. We shall do
it calculating directly the components of the rotor in the chosen
coordinate system. For the first component we have
$$
\frac{\partial A^3}{\partial x^2}-
\frac{\partial A^2}{\partial x^3}=
\frac{\partial}{\partial x^3} \int\limits^{\ x_3}_{0\,}
F^1(x^1,x^2,x^3)\,dx^3=F^1(x^1,x^2,x^3).
$$
Here we used the rule of differentiation of an integral with variable 
upper limit. In calculating the second component we take into account 
that the second integral in the expression for the component $A^1$ in
\mythetag{11.10} does not depend on $x^3$:
$$
\frac{\partial A^1}{\partial x^3}-
\frac{\partial A^3}{\partial x^1}=
\frac{\partial}{\partial x^3} \int\limits^{\ x_3}_{0\,}
F^2(x^1,x^2,x^3)\,dx^3=F^2(x^1,x^2,x^3).
$$
And finally, for the third components of the rotor we derive
$$
\gather
\frac{\partial A^2}{\partial x^1}-
\frac{\partial A^1}{\partial x^2}=
-\int\limits^{\ x_3}_{0\,}\left(
\frac{\partial F^1(x^1,x^2,x^3)}{\partial x^1}+
\frac{\partial F^2(x^1,x^2,x^3)}{\partial x^2}\right)\,dx^3+\\
+\,\frac{\partial}{\partial x^2}\int\limits^{\ x_2}_{0\,}
F^3(x^1,x^2,0)\,dx^2=\int\limits^{\ x_3}_{0\,}
\frac{\partial F^3(x^1,x^2,x^3)}{\partial x^3}\,dx^3+\\
\vspace{1ex}
+\,F^3(x^1,x^2,0)=
F^3(x^1,x^2,x)\,\vbox{\hrule width 0.5pt
height 12pt depth 8pt}^{\,\,x=x_3}_{\,\,x=0}+
F^3(x^1,x^2,0)=F^3(x^1,x^2,x^3).
\endgather
$$
In these calculations we used the relationship \mythetag{11.9} in
order to replace the sum of two partial derivatives
$\partial F^1/\partial x^1+\partial F^2/\partial x^2$ by
$-\partial F^3/\partial x^3$. Now, bringing together the results of
calculating all three components of the rotor, we see that 
$\,\rot\bold A=\bold F$. Hence, the required field $\bold A$ can
indeed be chosen in the form of \mythetag{11.10}.
\qed\enddemo
\newpage
\topmatter
\title\chapter{3}
Curvilinear coordinates
\endtitle
\endtopmatter
\chapternum=3
\document
\head
\S~\mysection{1} Some examples of curvilinear coordinate systems.
\endhead
\setfirstpage
\leftheadtext{Chapter \uppercase\expandafter{\romannumeral 3}.
CURVILINEAR COORDINATES.}
\rightheadtext{\S~1. Some examples of curvilinear coordinate systems.}
     The main purpose of Cartesian coordinate systems is the numeric
representation of points: each point of the space $\Bbb E$ is represented
by some unique triple of numbers $(x^1,x^2,x^3)$. Curvilinear coordinate
systems serve for the same purpose. We begin considering such coordinate
systems with some examples.\par
     {\it Polar coordinates}. Let's consider a plane, choose some point
$O$ on it (this will be the {\it pole}) and some ray $OX$ coming out
from this point. For an arbitrary point $A\neq O$ of that plane 
\vadjust{\vskip 5pt\hbox to 0pt{\kern -20pt
\includegraphics{ris05.eps}\hss}\vskip 167pt}its
position is determined by two parameters: the length of its 
radius-vector $\rho=|\overrightarrow{OA}|$ and the value of the angle
$\varphi$ between the ray $OX$ and the radius-vector of the point $A$. Certainly, one should also choose a positive (counterclockwise) direction
to which the angle $\varphi$ is laid (this is equivalent to choosing 
a preferable orientation on the plane). Angles laid to the opposite 
direction are understood as negative angles. The numbers $\rho$ and
$\varphi$ are called the {\it polar coordinates} of the point $A$.
\par
     Let's associate some Cartesian coordinate system with the polar coordinates as shown on Fig\.~1.2. We choose the point $O$ as an origin,
then direct the abscissa axis along the ray $OX$ and get the ordinate
axis from the abscissa axis rotating it by $90^{\circ}$. Then the
Cartesian coordinates of the point $A$ are derived from its polar
coordinates by means of the formulas
$$
\pagebreak
\hskip -2em
\cases
x^1=\rho\,\cos(\varphi),\\
x^2=\rho\,\sin(\varphi).
\endcases
\mytag{1.1}
$$
Conversely, one can express $\rho$ and $\varphi$ through $x^1$ and $x^2$
as follows:
$$
\hskip -2em
\cases
\rho=\sqrt{(x^1)^{\raise 2pt\hbox{$\ssize 2$}}
+(x^2)^{\raise 2pt\hbox{$\ssize 2$}}},\\
\varphi=\arctg(x^2/x^1).
\endcases
\mytag{1.2}
$$
\par
\parshape 18 0cm 360pt 0cm 360pt
180pt 180pt 180pt 180pt 180pt 180pt 180pt 180pt 
180pt 180pt 180pt 180pt 180pt 180pt 180pt 180pt 
180pt 180pt 180pt 180pt 180pt 180pt 180pt 180pt 
180pt 180pt 180pt 180pt 180pt 180pt 
0cm 360pt
\noindent
The pair of numbers $(\rho,\varphi)$ can be treated as an element
of the two-dimensional space $\Bbb R^2$. In order to express $\Bbb R^2$
\vadjust{\vskip 5pt\hbox to 0pt{\kern -20pt
\includegraphics{ris06.eps}\hss}\vskip -5pt}
visually we represent this space as a coordinate plane. The coordinate
plane $(\rho,\varphi)$ has not its own geometric interpretation, it is
called the {\it map} of the polar coordinate system. Not all points of
the map cor\-respond to the real geometric points. The condition
$\rho\geqslant 0$ excludes the whole half-plane of the map. The sine and
cosine both are periodic functions with the period $2\pi=360^{\circ}$.
Therefore there are different points of the map that re\-present the same
geometric point. Thus, the mapping $(\rho,\varphi)\to (x^1,x^2)$ given
by the formulas \mythetag{1.1} is not injective.
     Let $U$ be the unbounded domain highlighted with the light 
blue color on Fig\.~1.3 (the points of the boundary are not included).
Denote by $V$ the image of the domain $U$ under the mapping 
\mythetag{1.1}. It is easy to understand that $V$ is the set of all
points of the $(x^1,\,x^2)$ plane except for those lying on the ray 
$OX$. If we restrict the mapping \mythetag{1.1} to the domain $U$, we 
get the bijective mapping $m:U\to V$.\par
     Note that the formula \mythetag{1.2} is not an exact expression 
for the inverse mapping $m^{-1}\!:\,V\to U$. The matter is that the values 
of $\tg(\varphi)$ at the points $(x^1,x^2)$ and $(-x^1,-x^2)$ do coincide.
In order to express $m^{-1}$ exactly it would be better to use the tangent 
of the half angle:
$$
\tg(\varphi/2)=\frac{x^2}{x^1+
\sqrt{(x^1)^{\raise 2pt\hbox{$\ssize 2$}}
+(x^2)^{\raise 2pt\hbox{$\ssize 2$}}}}.
$$
However, we prefer the not absolutely exact expression 
for $\varphi$ from \mythetag{1.2} since it is relatively simple.\par
     Let's draw the series of equidistant straight lines parallel to
the axes on the map $\Bbb R^2$ of the polar coordinate system (see
Fig\.~1.4 below). The mapping \mythetag{1.1} takes them to the series 
of rays and concentric circles on the $(x^1,\,x^2)$ plane.
The straight lines on Fig\.~1.4 and the rays and circles on Fig\.~1.5
compose the coordinate network of the polar coordinate system. By
reducing the intervals between the lines one can obtain a more dense
coordinate network. This procedure can be repeated infinitely many
times producing more and more dense networks in each step. Ultimately
(in the continuum limit), one can think the coordinate network to be
maximally dense. Such a network consist of two families of lines:
the first family is given by the condition $\varphi=\const$, the second 
one --- by the similar\linebreak\pagebreak condition $\rho=\const$.\par
     On Fig\.~1.4 exactly two coordinate lines pass through each point 
of the map: one is from the first family and the other is from the second
family. On the $(x^1,x^2)$ plane this condition is fulfilled at all points
except for the origin $O$. Here all coordinate lines of the first family
are crossed. The origin $O$ is the only singular point of the polar coordinate system.
\vadjust{\vskip 5pt\hbox to 0pt{\kern -20pt
\includegraphics{ris07.eps}\hss}\vskip 200pt}\par
\parshape 14 0cm 360pt 0cm 360pt
160pt 200pt 160pt 200pt 160pt 200pt 160pt 200pt 
160pt 200pt 160pt 200pt 160pt 200pt 160pt 200pt 
160pt 200pt 160pt 200pt 160pt 200pt 
0cm 360pt
     The {\it cylindrical coordinate system} in the space $\Bbb E$
is obtained from the polar coordinates on a plane by adding the
third coordinate $h$. As in the case of polar coordinate system,
\vadjust{\vskip 5pt\hbox to 0pt{\kern -25pt
\includegraphics{ris08.eps}\hss}\vskip -5pt}we 
associate some Cartesian coordinate system with the cylin\-drical
coordinate system (see Fig\.~1.6). Then 
$$
\hskip -2em
\cases
x^1=\rho\,\cos(\varphi),\\
x^2=\rho\,\sin(\varphi),\\
x^3=h.
\endcases
\mytag{1.3}
$$
Conversely, one can pass from Cartesian to cylindrical coordinates 
by means of the formu\-la analogous to \mythetag{1.2}:
$$
\cases
\rho=\sqrt{(x^1)^{\raise 2pt\hbox{$\ssize 2$}}
+(x^2)^{\raise 2pt\hbox{$\ssize 2$}}},\\
\varphi=\arctg(x^2/x^1),\\
h=x^3.
\endcases
\mytag{1.4}
$$
The coordinate network of the cylindrical coordinate system
consists of three families of lines. These are the horizontal
rays coming out from the points of the vertical axis $Ox^3$,
the horizontal circles with the centers at the points of the
axis $Ox^3$, and the vertical straight lines parallel to the
axis $Ox^3$. The singular points of the cylindrical coordinate
system fill the axis $Ox^3$. Exactly three coordinate lines 
(one from each family) pass through each regular point of the 
space $\Bbb E$, i\.\,e\. through each point that does not lie 
\pagebreak on the axis $Ox^3$.\par
\parshape 15 
160pt 200pt 160pt 200pt 160pt 200pt 160pt 200pt 
160pt 200pt 160pt 200pt 160pt 200pt 160pt 200pt 
160pt 200pt 160pt 200pt 160pt 200pt 160pt 200pt 
160pt 200pt 160pt 200pt 
0cm 360pt
     The {\it spherical coordinate system} in the spa\-ce $\Bbb E$
\vadjust{\vskip 5pt\hbox to 0pt{\kern -25pt
\includegraphics{ris09.eps}\hss}\vskip -5pt}is
obtained by slight modification of the cylindrical coordinates.
The coordinate $h$ is replaced by the angular coordinate 
$\vartheta$, while the quantity $\rho$ in spherical coordinates
denotes the length of the radius-vector of the point $A$ (see
Fig\.~1.7). Then
$$
\hskip -2em
\cases
x^1=\rho\,\sin(\vartheta)\,\cos(\varphi),\\
x^2=\rho\,\sin(\vartheta)\,\sin(\varphi),\\
x^3=\rho\,\cos(\vartheta).
\endcases
\mytag{1.5}
$$
The spherical coordinates of a point are usual\-ly written in 
the following order: $\rho$ is the first coordinate, $\vartheta$ is 
the second one, and $\varphi$ is the third coordinate. The converse 
transition from Cartesian coor\-dinates to these quantities is 
given by the formula:
$$
\hskip -2em
\cases
\rho=\sqrt{(x^1)^{\raise 2pt\hbox{$\ssize 2$}}
+(x^2)^{\raise 2pt\hbox{$\ssize 2$}}
+(x^3)^{\raise 2pt\hbox{$\ssize 2$}}},\\
\vartheta=\arccos\left(x^3/\sqrt{(
(x^1)^{\raise 2pt\hbox{$\ssize 2$}}
+(x^2)^{\raise 2pt\hbox{$\ssize 2$}}
+(x^3)^{\raise 2pt\hbox{$\ssize 2$}}}\,\right),\\
\varphi=\arctg(x^2/x^1).\\
\endcases
\mytag{1.6}
$$
Coordinate lines of spherical coordinates form three families.
The first family is composed of the rays coming out from the 
point $O$; the second family is formed by circles that lie in
various vertical planes passing through the axix $Ox^3$; and
the third family consists of horizontal circles whose centers 
are on the axis $Ox^3$. Exactly three coordinate lines pass
through each regular point of the space $\Bbb E$, one line from 
each family.\par
     The condition $\rho=\const$ specifies the sphere of the radius
$\rho$ in the space $\Bbb E$. The coordinate lines of the second and 
third families define the network of meridians and parallels on this 
sphere exactly the same as used in geography to define the coordinates 
on the Earth surface.\par
\head
\S~\mysection{2} Moving frame of a curvilinear coordinate system.
\endhead
\rightheadtext{\S~2. Moving frame of a curvilinear coordinate system.}
     Let $D$ be some domain in the space $\Bbb E$. Saying domain, one 
usually understand a connected open set. An open set $D$ means that
along with each its point $A\in D$ the set $D$ comprises some spherical
neighborhood $O(A)$ of this point. A connected set $D$ means that any 
two points of this set can be connected by a smooth curve lying within
$D$. See more details in \mycite{2}. Let's consider three numeric functions
$u^1(\bold x)$, $u^2(\bold x)$, and $u^3(\bold x)$ defined in the domain
$D$. Generally speaking, their domains could be wider, but we need them
only within $D$. The values of three functions $u^1,\,u^2,\,u^3$ at each
point form a triple of numbers, they can be interpreted as a point of
the space $\Bbb R^3$. Then the triple of functions $u^1,\,u^2,\,u^3$ 
define a mapping $\bold u\!:\,D\to\Bbb R^3$.
\mydefinition{2.1} A triple of differentiable functions 
$u^1,\,u^2,\,u^3$ is called {\it regular at a point $A$} of the space
$\Bbb E$ if the gradients of these functions $\grad u^1$, $\grad u^2$, 
and $\grad u^3$ are linearly independent at the point $A$.
\enddefinition
     Let's choose some Cartesian coordinate system in $\Bbb E$, in this
coordinate system the above functions $u^1,\,u^2,\,u^3$ are represented
by the functions $u^i=u^i(x^1,x^2,x^3)$ of Cartesian coordinates of a
point. The gradients of the differentiable functions $u^1,\,u^2,\,u^3$
form the triple of covectorial fields whose components are given by the
partial derivatives of $u^1,\,u^2,\,u^3$ with respect to $x^1,\,x^2,\,x^3$:
$$
\hskip -2em
\grad u^i=\left(\frac{\partial u^i}{\partial x^1},\
\frac{\partial u^i}{\partial x^2},\
\frac{\partial u^i}{\partial x^3}\right).
\mytag{2.1}
$$
Let's compose a matrix of the gradients \mythetag{2.1}:
$$
\hskip -2em
J=\left\Vert
\vphantom{\vrule height 39pt depth 39pt}
\matrix
\dsize\frac{\partial u^1}{\partial x^1}&
\dsize\frac{\partial u^1}{\partial x^2}&
\dsize\frac{\partial u^1}{\partial x^3}\\
\vspace{1.2ex}
\dsize\frac{\partial u^2}{\partial x^1}&
\dsize\frac{\partial u^2}{\partial x^2}&
\dsize\frac{\partial u^2}{\partial x^3}\\
\vspace{1.2ex}
\dsize\frac{\partial u^3}{\partial x^1}&
\dsize\frac{\partial u^3}{\partial x^2}&
\dsize\frac{\partial u^3}{\partial x^3}
\endmatrix\right\Vert.
\mytag{2.2}
$$
The matrix $J$ of the form \mythetag{2.2} is called the {\it Jacobi
matrix} of the mapping $\bold u\!:\,D\to R^3$ given by the triple
of the differentiable functions $u^1,\,u^2,\,u^3$ in the domain $D$.
It is obvious that the regularity of the functions $u^1,\,u^2,\,u^3$ 
at a point is equivalent to the non-degeneracy of the Jacobi matrix 
at that point: $\det J\neq 0$.\par
\mytheorem{2.1} If continuously differentiable functions
$u^1,\,u^2,\,u^3$ with the do\-main $D$ are regular at a point $A$, 
then there exists some neighborhood $O(A)$ of the point $A$ and
some neighborhood $O(\bold u(A))$ of the point $\bold u(A)$ in
the space $\Bbb R^3$ such that the following conditions are fulfilled:
\roster
\item the mapping $\bold u\!:\,O(A)\to O(\bold u(A))$ is bijective;
\item the inverse mapping $\bold u^{-1}\!:\,O(\bold u(A))\to O(A)$ 
is continuously differentiable.
\endroster
\endproclaim
The theorem~\mythetheorem{2.1} or propositions equivalent to it are 
usually proved in the course of mathematical analysis (see \mycite{2}).
They are known as the theorems on implicit functions.\par
\mydefinition{2.2} Say that an ordered triple of continuously differentiable functions $u^1,\,u^2,\,u^3$ with the domain $D\subset
\Bbb E$ define a {\it curvilinear coordinate system} in $D$ if it is
regular at all points of $D$ and if the mapping $\bold u$ determined 
by them is a bijective mapping from $D$ to some domain  $U\subset\Bbb R^3$.
\enddefinition
     The cylindrical coordinate system is given by three functions
$u^1=\rho(\bold x)$, $u^2=\varphi(\bold x)$, and $u^3=h(\bold x)$ from
\mythetag{1.4}, while the spherical coordinate system is given by the
functions \mythetag{1.6}. However, the triples of functions \mythetag{1.4} 
and \mythetag{1.6} satisfy the conditions from the 
definition~\mythedefinition{2.2} only after reducing somewhat their
domains. Upon proper choice of a domain $D$ for \mythetag{1.4} and
\mythetag{1.6} the inverse mappings $\bold u^{-1}$ are given by the
formulas \mythetag{1.3} and \mythetag{1.5}.\par
     Suppose that in a domain $D\subset\Bbb E$ a curvilinear coordinate
system $u^1,\,u^2,\,u^3$ is given. Let's choose an auxiliary Cartesian 
coordinate system in $\Bbb E$. Then $u^1,\,u^2,\,u^3$ is a triple of
functions defining a map $\bold u$ from $D$ onto some domain $U\subset
\Bbb R^3$:
$$
\hskip -2em
\cases
u^1=u^1(x^1,x^2,x^3),\\
u^2=u^2(x^1,x^2,x^3),\\
u^3=u^3(x^1,x^2,x^3).
\endcases
\mytag{2.3}
$$
The domain $D$ is called the {\it domain being mapped}, the domain
$U\subset\Bbb R^3$ is called the {\it map} or the {\it chart}, while 
$\bold u^{-1}\!:\,U\to D$ is called the {\it chart mapping}.
The chart mapping is given by the following three functions:
$$
\hskip -2em
\cases
x^1=x^1(u^1,u^2,u^3),\\
x^2=x^2(u^1,u^2,u^3),\\
x^3=x^3(u^1,u^2,u^3).
\endcases
\mytag{2.4}
$$
Denote by $\bold r$ the radius-vector $\bold r$ of the point with Cartesian
coordinates $x^1,\,x^2,\,x^3$. Then instead of three scalar functions
\mythetag{2.4} we can use one vectorial function 
$$
\hskip -2em
\bold r(u^1,u^2,u^3)=\sum^3_{q=1}x^q(u^1,u^2,u^3)\cdot
\bold e_q.
\mytag{2.5}
$$
Let's fix some two of three coordinates $u^1,\,u^2,\,u^3$ and let's vary
the third of them. Thus we get three families of straight lines within the
domain $U\subset\Bbb R^3$:
$$
\xalignat 3
&\hskip -2em
\cases
u^1=t,\\
u^2=c^2,\\
u^3=c^3,
\endcases
&&\cases
u^1=c^1,\\
u^2=t,\\
u^3=c^3,
\endcases
&&\cases
u^1=c^1,\\
u^2=c^2,\\
u^3=t.
\endcases
\mytag{2.6}
\endxalignat
$$
Here $c^1,\,c^2,\,c^3$ are constants. The straight lines \mythetag{2.6}
form a rectangular coordinate network within the chart $U$. Exactly one 
straight line from each of the families \mythetag{2.6} passes through each
point of the chart. Substituting \mythetag{2.6} into \mythetag{2.5} we map
the rectangular network from $U$ onto a curvilinear network in the domain
$D\subset E$. Such a network is called the {\it coordinate network} of a
curvilinear coordinate system.\par
     The coordinate network of a curvilinear coordinate system on the
domain $D$ consists of three families of lines. Due to the bijectivity
of the mapping $\bold u\!:\,D\to U$ exactly three coordinate lines pass
through each point of the domain $D$ --- one line from each family. Each coordinate line has its canonical parametrization: $t=u^1$ is the parameter
for the lines of the first family, $t=u^2$ is the parameter for 
the lines of the second family, and finally, $t=u^3$ is the parameter for
the lines of the third family. At each point of the domain $D$ we have 
three tangent vectors, they are tangent to the coordinate lines of the
three families passing through that point. Let's denote them $\bold E_1$,
$\bold E_2$, $\bold E_3$. The vectors $\bold E_1,\,\bold E_2,\,\bold E_3$
are obtained by differentiating the radius-vector $\bold r(u^1,u^2,u^3)$ 
with respect to the parameters $u^1$, $u^2$, $u^3$ of coordinate lines.
Therefore, we can write
$$
\hskip -2em
\bold E_j(u^1,u^2,u^3)=\frac{\partial\,\bold r(u^1,u^2,u^3)}
{\partial u^j}.
\mytag{2.7}
$$
Let's substitute \mythetag{2.5} into \mythetag{2.7}. The basis vectors 
$\bold e_1,\,\bold e_2,\,\bold e_3$ do not depend on the variables
$u^1,\,u^2,\,u^3$, hence, we get
$$
\bold E_j(u^1,u^2,u^3)=\sum^3_{q=1}
\frac{\partial x^q(u^1,u^2,u^3)}{\partial u^j}
\cdot\bold e_q.
\mytag{2.8}
$$
The formula \mythetag{2.8} determines the expansion of the vectors
$\bold E_1,\,\bold E_2,\,\bold E_3$ in the basis $\bold e_1,\,
\bold e_2,\,\bold e_3$. The column-vectors of the coordinates of
$\bold E_1$, $\bold E_2$, and $\bold E_3$ can be concatenated
into the following matrix:
$$
\hskip -2em
I=\left\Vert
\vphantom{\vrule height 39pt depth 39pt}
\matrix
\dsize\frac{\partial x^1}{\partial u^1}&
\dsize\frac{\partial x^1}{\partial u^2}&
\dsize\frac{\partial x^1}{\partial u^3}\\
\vspace{1.2ex}
\dsize\frac{\partial x^2}{\partial u^1}&
\dsize\frac{\partial x^2}{\partial u^2}&
\dsize\frac{\partial x^2}{\partial u^3}\\
\vspace{1.2ex}
\dsize\frac{\partial x^3}{\partial u^1}&
\dsize\frac{\partial x^3}{\partial u^2}&
\dsize\frac{\partial x^3}{\partial u^3}
\endmatrix\right\Vert.
\mytag{2.9}
$$
Comparing \mythetag{2.9} and \mythetag{2.2}, we see that \mythetag{2.9} 
is the Jacobi matrix for the mapping $\bold u^{-1}\!:\,U\to D$ 
given by the functions \mythetag{2.4}. Let's substitute \mythetag{2.4} 
into \mythetag{2.3}:
$$
\hskip -2em
u^i(x^1(u^1,u^2,u^3),x^2(u^1,u^2,u^3),x^3(u^1,u^2,u^3))=u^i.
\mytag{2.10}
$$
The identity \mythetag{2.10} follows from the fact that the functions
\mythetag{2.3} and \mythetag{2.4} define two mutually inverse mappings 
$\bold u$ and $\bold u^{-1}$. Let's differentiate the identity \mythetag{2.10} with respect to the variable $u^j$:
$$
\hskip -2em
\sum^3_{q=1}\frac{\partial u^i(x^1,x^2,x^3)}{\partial x^q}\,
\frac{\partial x^q(u^1,u^2,u^3)}{\partial u^j}=\delta^i_j.
\mytag{2.11}
$$
Here we used the chain rule for differentiating the composite function
in \mythetag{2.10}. The relationship \mythetag{2.11} shows that the matrices
\mythetag{2.2} and \mythetag{2.9} are inverse to each other. More precisely,
we have the following relationship
$$
\hskip -2em
I(u^1,u^2,u^3)=J(x^1,x^2,x^3)^{-1},
\mytag{2.12}
$$
where $x^1,\,x^2,\,x^3$ should be expressed through $u^1,\,u^2,\,u^3$ by means of \mythetag{2.4}, or conversely, $u^1,\,u^2,\,u^3$ should be expressed
through $x^1,\,x^2,\,x^3$ by means of \mythetag{2.3}. The arguments shown
in the relationship \mythetag{2.12} are the natural arguments for the
components of the Jacobi matrices $I$ and $J$. However, one can pass to 
any required set of variables by means of \mythetag{2.3} or \mythetag{2.4}
whenever it is necessary.\par
     The regularity of the triple of functions \mythetag{2.3} defining a
curvilinear coordinate system in the domain $D$ means that the matrix
\mythetag{2.2} is non-degenerate. Then, due to \mythetag{2.12}, the inverse
matrix \mythetag{2.9} is also non-degenerate. Therefore, the vectors 
$\bold E_1,\,\bold E_2,\,\bold E_3$ given by the formula \mythetag{2.8}
are linearly independent at any point of the domain $D$. Due to the
linear independence of the {\it coordinate tangent vectors} $\bold E_1,\,
\bold E_2,\,\bold E_3$ they form a {\it moving frame} which is usually
called the {\it coordinate frame} of the curvilinear coordinate system.
The formula \mythetag{2.8} now can be interpreted as the transition formula
for passing from the basis of the auxiliary Cartesian coordinate system 
to the basis formed by the vectors of the frame $\bold E_1,\,\bold E_2,
\,\bold E_3$:
$$
\hskip -2em
\bold E_j=\sum^3_{q=1} S^q_j(u^1,u^2,u^3)\cdot\bold e_q.
\mytag{2.13}
$$
The transition matrix $S$ in the formula \mythetag{2.13} coincides with
the Jacobi matrix \mythetag{2.9}, therefore its components depend on
$u^1,\,u^2,\,u^3$. These are the natural variables for the components
of $S$.\par
     The inverse transition from the basis $\bold E_1,\,\bold E_2,\,
\bold E_3$ to the basis of the Cartesian coordinate system is given
by the inverse matrix $T=S^{-1}$. Due to \mythetag{2.12} the inverse
transition matrix coincides with the Jacobi matrix \mythetag{2.2}. 
Therefore, $x^1,\,x^2,\,x^3$ are the natural variables for the 
components of the matrix $T$:
$$
\bold e_q=\sum^3_{i=1} T^i_q(x^1,x^2,x^3)\cdot\bold E_i.
\mytag{2.14}
$$
The vectors $\bold E_1,\,\bold E_2,\,\bold E_3$ of the moving frame 
depend on a point of the domain $D\subset E$. Since in a curvilinear
coordinate system such a point is represented by its coordinates
$u^1,\,u^2,\,u^3$, these variables are that very arguments which are
natural for the vectors of the moving frame: $\bold E_i=\bold E_i(u^1,
u^2,u^3)$.\par
\head
\S~\mysection{3} Change of curvilinear coordinates.
\endhead
\rightheadtext{\S~3. Change of curvilinear coordinates.}
     Let $u^1,\,u^2,\,u^3$ be some curvilinear coordinates in some
domain $D_1$ and let $\tilde u^1,\,\tilde u^2,\,\tilde u^3$ be some
other curvilinear coordinates in some other domain $D_2$. If the
domains $D_1$ and $D_2$ do intersect, then in the domain $D=D_1\cap
\tilde D_2$ we have two coordinate systems.
\vadjust{\vskip 5pt\hbox to 0pt{\kern 5pt
\includegraphics{ris10.eps}\hss}\vskip 160pt}We denote 
by $U$ and $\tilde U$ the preimages of the domain $D$ in the maps 
$U_1$ and $U_2$, i\.\,e\. we denote $U=\bold u(D_1\cap\tilde D_2)$ 
and we denote $\tilde U=\tilde\bold u(D_1\cap\tilde D_2)$.
Due to the chart mappings the points of the domain $D$ are in one-to-one
correspondence with the points of the domains $U$ and $\tilde U$. As for
the chart mappings $\bold u^{-1}$ and $\tilde\bold u^{-1}$, they are
given by the following functions:
$$
\xalignat 2
&\hskip -2em
\cases
x^1=x^1(u^1,u^2,u^3),\\
x^2=x^2(u^1,u^2,u^3),\\
x^3=x^3(u^1,u^2,u^3),
\endcases
&&\cases
x^1=x^1(\tilde u^1,\tilde u^2,\tilde u^3),\\
x^2=x^2(\tilde u^1,\tilde u^2,\tilde u^3),\\
x^3=x^3(\tilde u^1,\tilde u^2,\tilde u^3).
\endcases
\mytag{3.1}
\endxalignat
$$
The mappings $\bold u$ and $\tilde\bold u$ inverse to the chart mappings 
are given similarly:
$$
\xalignat 2
&\hskip -2em
\cases
u^1=u^1(x^1,x^2,x^3),\\
u^2=u^2(x^1,x^2,x^3),\\
u^3=u^3(x^1,x^2,x^3),
\endcases
&&\cases
\tilde u^1=\tilde u^1(x^1,x^2,x^3),\\
\tilde u^2=\tilde u^2(x^1,x^2,x^3),\\
\tilde u^3=\tilde u^3(x^1,x^2,x^3).
\endcases
\mytag{3.2}
\endxalignat
$$
Let's substitute the first set of functions \mythetag{3.1} into the
\pagebreak 
arguments of the second set of functions \mythetag{3.2}. Similarly, 
we substitute the second set of functions \mythetag{3.1} into the arguments
of the first set of functions in \mythetag{3.2}. As a result we get the
functions
$$
\align
&\hskip -2em
\cases
\tilde u^1(x^1(u^1,u^2,u^3),x^2(u^1,u^2,u^3),
x^3(u^1,u^2,u^3)),\\
\tilde u^2(x^1(u^1,u^2,u^3),x^2(u^1,u^2,u^3),
x^3(u^1,u^2,u^3)),\\
\tilde u^3(x^1(u^1,u^2,u^3),x^2(u^1,u^2,u^3),
x^3(u^1,u^2,u^3)),
\endcases
\mytag{3.3}\\
\vspace{1ex}
&\hskip -2em
\cases
u^1(x^1(\tilde u^1,\tilde u^2,\tilde u^3),
        x^2(\tilde u^1,\tilde u^2,\tilde u^3),
        x^3(\tilde u^1,\tilde u^2,\tilde u^3)),\\
u^2(x^1(\tilde u^1,\tilde u^2,\tilde u^3),
        x^2(\tilde u^1,\tilde u^2,\tilde u^3),
        x^3(\tilde u^1,\tilde u^2,\tilde u^3)),\\
u^3(x^1(\tilde u^1,\tilde u^2,\tilde u^3),
        x^2(\tilde u^1,\tilde u^2,\tilde u^3),
        x^3(\tilde u^1,\tilde u^2,\tilde u^3))
\endcases
\mytag{3.4}
\endalign
$$
which define the pair of mutually inverse mappings $\tilde\bold u
\compos\bold u^{-1}$ and $\bold u\compos\tilde\bold u^{-1}$. For
the sake of brevity we write these sets of functions as follows:
$$
\xalignat 2
&\hskip -2em
\cases
\tilde u^1=\tilde u^1(u^1,u^2,u^3),\\
\tilde u^2=\tilde u^2(u^1,u^2,u^3),\\
\tilde u^3=\tilde u^3(u^1,u^2,u^3),
\endcases
&&\cases
u^1=u^1(\tilde u^1,\tilde u^2,\tilde u^3),\\
u^2=u^2(\tilde u^1,\tilde u^2,\tilde u^3),\\
u^3=u^3(\tilde u^1,\tilde u^2,\tilde u^3).
\endcases
\mytag{3.5}
\endxalignat
$$
The formulas \mythetag{3.5} express the coordinates of a point from the
domain $D$ in some curvilinear coordinate system through its coordinates
in some other coordinate system. These formulas are called the 
{\it transformation formulas} or the formulas for {\it changing the 
curvilinear coordinates}.\par
     Each of the two curvilinear coordinate systems has its own moving
frame within the domain $D=D_1\cap D_2$. Let's denote by 
$S$ and $T$ the transition matrices relating these two moving frames.
Then we can write
$$
\xalignat 2
&\hskip -2em
\tilde\bold E_j=\sum^3_{i=1} S^i_j\cdot\bold E_i,
&&\bold E_i=\sum^3_{k=1} T^k_i\cdot\tilde\bold E_k.
\mytag{3.6}
\endxalignat
$$
\mytheorem{3.1} The components of the transition matrices 
$S$ and $T$ for the moving frames of two curvilinear coordinate system
in \mythetag{3.6} are determined by the partial derivatives of the
functions \mythetag{3.5}:
$$
\xalignat 2
&\hskip -2em
S^i_j(\tilde u^1,\tilde u^2,\tilde u^3)=\frac{\partial u^i}
{\partial\tilde u^j},
&&T^k_i(u^1,u^2,u^3)=\frac{\partial\tilde u^k}{\partial u^i}.
\mytag{3.7}
\endxalignat
$$
\endproclaim
\demo{Proof} We shall prove only the first formula in \mythetag{3.7}. The
proof of the second formula is absolutely analogous to the proof of the
first one. Let's choose some auxiliary Cartesian coordinate system and
then write the formula \mythetag{2.8} applied to the frame vectors of the 
second curvilinear coordinate system:
$$
\hskip -2em
\tilde\bold E_j(\tilde u^1,\tilde u^2,\tilde u^3)=\sum^3_{q=1}
\frac{\partial x^q(\tilde u^1,\tilde u^2,\tilde u^3)}{\partial
\tilde u^j}\cdot\bold e_q.
\mytag{3.8}
$$
Applying the formula \mythetag{2.14}, we express $\bold e_1,\,\bold e_2,
\,\bold e_3$ through $\bold E_1,\,\bold E_2,\,\bold E_3$. Remember that 
the matrix $T$ in \mythetag{2.14} coincides with the Jacobi \pagebreak
matrix $J(x^1,x^2,x^3)$ from \mythetag{2.2}. Therefore, we can write the
following formula:
$$
\hskip -2em
\bold e_q=\sum^3_{i=1}
\frac{\partial u^i(x^1,x^2,x^3)}
{\partial x^q}\cdot\bold E_i.
\mytag{3.9}
$$
Now let's substitute \mythetag{3.9} into \mythetag{3.8}. As a result we 
get the formula relating the frame vectors of two curvilinear 
coordinate systems:
$$
\tilde\bold E_j=\sum^3_{i=1}\left(\shave{\,\sum^3_{q=1}}
\frac{\partial u^i(x^1,x^2,x^3)}
{\partial x^q}\,\frac{\partial x^q(\tilde u^1,\tilde u^2,
\tilde u^3)}{\partial\tilde u^j}\right)\cdot\bold E_i.
\mytag{3.10}
$$
Comparing \mythetag{3.10} and \mythetag{3.6}, from this comparison for
the components of $S$ we get
$$
S^i_j=\sum^3_{q=1}\frac{\partial u^i(x^1,x^2,x^3)}
{\partial x^q}\,\frac{\partial x^q(\tilde u^1,\tilde u^2,
\tilde u^3)}{\partial\tilde u^j}.
\mytag{3.11}
$$
Remember that the Cartesian coordinates $x^1,\,x^2,\,x^3$ in the above
formula \mythetag{3.11} are related to the curvilinear coordinates
$\tilde u^1,\,\tilde u^2,\,\tilde u^3$ by means of \mythetag{3.1}. 
Hence, the sum in right hand side of \mythetag{3.11} can be transformed
to the partial derivative of the composite function
$u^i((x^1(\tilde u^1,\tilde u^2,\tilde u^3),x^2(\tilde u^1,\tilde u^2,\tilde u^3),x^3(\tilde u^1,\tilde u^2,\tilde u^3))$ from \mythetag{3.4}:
$$
S^i_j=\frac{\partial u^i}{\partial\tilde u^j}.
$$
Note that the functions \mythetag{3.4} written in the form of 
\mythetag{3.5} are that very functions relating $\tilde u^1,\,
\tilde u^2,\,\tilde u^3$ and $u^1,\,u^2,\,u^3$, and their 
derivatives are in formula \mythetag{3.7}. The theorem 
is proved.\qed\enddemo
\subhead A remark on the orientation\endsubhead From the 
definition~\mythedefinition{2.2} we derive that the functions
\mythetag{2.3} are continuously differentiable. Due to the 
theorem~\mythetheorem{2.1} the functions \mythetag{2.4} representing 
the inverse mappings are also continuously differentiable. Then the
components of the matrix $S$ in the formula \mythetag{2.13} coinciding 
with the components of the Jacobi matrix \mythetag{2.9} are 
continuous functions within the domain $U$. The same is true
for the determinant of the matrix $S$: the determinant $\det
S(u^1,u^2,u^3)$ is a continuous function in the domain $U$
which is nonzero at all points of this domain. A nonzero 
continuous real function in a connected set $U$ cannot take 
the values of different signs in $U$. This means that
$\det S>0$ or $\det S<0$. This means that the orientation 
of the triple of vectors forming the moving frame of a 
curvilinear coordinate system is the same for all points of
a domain where it is defined. Since the space $\Bbb E$ is 
equipped with the preferable orientation, we can subdivide
all curvilinear coordinates in $\Bbb E$ into right-oriented 
and left-oriented coordinate systems.\par
\subhead  A remark on the smoothness\endsubhead The 
definition~\mythedefinition{2.2} yields the concept of a continuously
differentiable curvilinear coordinate system. However, the functions
\mythetag{2.3} could belong to a higher smoothness class $C^m$. In this
case we say that we have a {\it curvilinear coordinate system of the
smoothness class $C^m$}. The components of the Jacobi matrix \mythetag{2.2}
for such a coordinate system are the functions of the class $C^{m-1}$. Due
to the relationship \mythetag{2.12} the components of the Jacobi matrix
\mythetag{2.9} belong to the same smoothness class $C^{m-1}$. Hence, the
functions \mythetag{2.4} belong to the smoothness class $C^m$.\par
     If we have two curvilinear coordinate systems of the smoothness 
class $C^m$, then, according to the above considerations, the 
transformation functions \mythetag{3.5} belong to the class $C^m$, while
the components of the transition matrices $S$ and $T$ given by the
formulas \mythetag{3.7} belong to the smoothness class $C^{m-1}$.\par
\head
\S~\mysection{4} Vectorial and tensorial fields\\
 in curvilinear coordinates.
\endhead
\rightheadtext{\S~4. Vectorial and tensorial fields \dots}
     Let $u^1,\,u^2,\,u^3$ be some curvilinear coordinate system
in some domain $D\subset\Bbb E$ and let $\bold F$ be some vector 
field defined at the points of the domain $D$. Then at a point
with coordinates $u^1,\,u^2,\,u^3$ we have the field vector 
$\bold F(u^1,u^2,u^3)$ and the triple of the frame vectors 
$\bold E_1(u^1,u^2,u^3),\,\bold E_2(u^1,u^2,u^3),\,\bold E_3(u^1,
u^2,u^3)$. Let's expand the field vector $\bold F$ in the
basis formed by the frame vectors:
$$
\hskip -2em
\bold F(u^1,u^2,u^3)=\sum^3_{i=1} F^i(u^1,u^2,u^3)\cdot
\bold E_i(u^1,u^2,u^3).
\mytag{4.1}
$$
The quantities $F^i(u^1,u^2,u^3)$ in such expansion are naturally called
the {\it components} of the vector field $\bold F$ in the curvilinear
coordinates $u^1,\,u^2,\,u^3$. If we have another curvilinear coordinate 
system $\tilde u^1,\,\tilde u^2,\,\tilde u^3$ in the domain $D$, then we
have the other expansion of the form \mythetag{4.1}:
$$
\hskip -2em
\bold F(\tilde u^1,\tilde u^2,\tilde u^3)=\sum^3_{i=1}\tilde
F^i(\tilde u^1,\tilde u^2,\tilde u^3)\cdot
\tilde\bold E_i(\tilde u^1,\tilde u^2,\tilde u^3).
\mytag{4.2}
$$
By means of the formulas \mythetag{3.6} one can easily derive the
relationships binding the components of the field $\bold F$ in the
expansions \mythetag{4.1} and \mythetag{4.2}:
$$
\hskip -2em
\aligned
&F^i(\bold u)=\sum^3_{j=1} S^i_j(\tilde\bold u)\,\,
\tilde F^j(\tilde\bold u),\\
\vspace{0.5ex}
&u^i=u^i(\tilde u^1,\tilde u^2,\tilde u^3).
\endaligned
\mytag{4.3}
$$
The relationships \mythetag{4.3} are naturally interpreted as the
generalizations for the relationships \mythetagchapter{1.2}{2} from  
Chapter \uppercase\expandafter{\romannumeral 2} for the case of
curvilinear coordinates.\par
     Note that Cartesian coordinate systems can be treated as
a special case of curvilinear coordinates. The transition functions
$u^i=u^i(\tilde u^1,\tilde u^2,\tilde u^3)$ in the case of a pair
of Cartesian coordinate systems are linear, therefore the matrix
$S$ calculated according to the theorem~\mythetheorem{3.1} in this 
case is a constant matrix.\par
     Now let $\bold F$ be either a field of covectors, a field of
linear operators, or a field of bilinear forms. In any case the
components of the field $F$ at some point are determined by fixing
some basis attached to that point. The vectors of the moving frame
of a curvilinear coordinate system at a point with coordinates 
$u^1,\,u^2,\,u^3$ provide the required basis. The components of the
field $\bold F$ determined by this basis are called the {\it components}
of the field $\bold F$ in that curvilinear coordinates. The transformation
rules for the components of the fields listed above under a change of
curvilinear coordinates generalize the formulas \mythetagchapter{1.3}{2},
\mythetagchapter{1.4}{2}, and \mythetagchapter{1.5}{2} from
Chapter~\uppercase\expandafter{\romannumeral 2}. For a covectorial field $\bold F$ the transformation rule for its components under a change of coordinates looks like
$$
\hskip -2em
\aligned
&F_i(\bold u)=\sum^3_{j=1} T^j_i(\bold u)\,\,
\tilde F_j(\tilde\bold u),\\
\vspace{0.5ex}
&u^i=u^i(\tilde u^1,\tilde u^2,\tilde u^3).
\endaligned
\mytag{4.4}
$$
The transformation rule for the components of an operator field 
$\bold F$ is written as
$$
\hskip -2em
\aligned
&F^i_j(\bold x)=\sum^3_{p=1}\sum^3_{q=1} S^i_p(\tilde\bold u)\,
T^q_j(\bold u)\,\,\tilde F^p_q(\tilde\bold x),\\
\vspace{0.5ex}
&u^i=u^i(\tilde u^1,\tilde u^2,\tilde u^3).
\endaligned
\mytag{4.5}
$$
In the case of a field of bilinear (quadratic) forms the generalization
of the formula \mythetagchapter{1.5}{2} from Chapter
\uppercase\expandafter{\romannumeral 2} looks like
$$
\hskip -2em
\aligned
&F_{ij}(\bold u)=\sum^3_{p=1}\sum^3_{q=1} T^p_i(\bold u)\,
T^q_j(\bold u)\,\,\tilde F_{pq}(\tilde\bold u),\\
\vspace{0.5ex}
&u^i=u^i(\tilde u^1,\tilde u^2,\tilde u^3).
\endaligned
\mytag{4.6}
$$\par
     Let $\bold F$ be a tensor field of the type $(r,s)$. In contrast
to a vectorial field, the value of such a tensorial field at a point
have no visual embodiment in form of an arrowhead segment. Moreover,
in general case there is no visually explicit way of finding the
numerical values for the components of such a field in a given basis.
However, according to the definition~\mythedefinitionchapter{1.1}{2}
from Chapter \uppercase\expandafter{\romannumeral 2}, a tensor is a
geometric object that for each basis has an array of components associated
with this basis. Let's denote by $\bold F(u^1,u^2,u^3)$ the value of the
field $\bold F$ at the point with coordinates $u^1,\,u^2,\,u^3$. This is
a tensor whose components in the basis $\bold E_1(u^1,u^2,u^3),\,\bold
E_2(u^1,u^2,u^3),\,\bold E_3(u^1,u^2,u^3)$ are called components of the
field $\bold F$ in a given curvilinear coordinate system. The transformation
rules for the components of a tensor field under a change of a coordinate
system follow from the formula \mythetagchapter{1.6}{2} in
Chapter \uppercase\expandafter{\romannumeral 2}. For a tensorial field of
the type $(r,s)$ it looks like
$$
\hskip -2em
\aligned
&\aligned F^{i_1\ldots\,i_r}_{j_1\ldots\,j_s}(\bold u)=
  \sum\Sb p_1\ldots\,p_r\\ q_1\ldots\,q_s\endSb
    &S^{i_1}_{p_1}(\tilde\bold u)\ldots\,S^{i_r}_{p_r}
    (\tilde\bold u)\times\\
    \vspace{-3ex}
    &\ \times T^{q_1}_{j_1}(\bold u)\ldots\,T^{q_s}_{j_s}(\bold u)
    \,\,\tilde F^{p_1\ldots\,p_r}_{q_1\ldots\,q_s}(\tilde\bold u),
    \hskip -1cm
 \endaligned\\
\vspace{0.9ex}
&u^i=u^i(\tilde u^1,\tilde u^2,\tilde u^3).
\endaligned
\mytag{4.7}
$$
The formula \mythetag{4.7} has two important differences as compared 
to the corresponding formula \mythetagchapter{1.7}{2} in 
Chapter \uppercase\expandafter{\romannumeral 2}. In the case of
curvilinear coordinates 
\roster
\rosteritemwd=10pt
\item the transition functions $u^i(\tilde u^1,\tilde u^2,\tilde u^3)$
      should not be linear functions;
\item the transition matrices $S(\tilde\bold u)$ and $T(\bold u)$ are not
      necessarily constant matrices.
\endroster
Note that these differences do not affect the algebraic operations
with tensorial fields. The operations of addition, tensor product,
contraction, index permutation, symmetrization, and alternation are
implemented by the same formulas as in Cartesian coordinates. The
differences \therosteritem{1} and \therosteritem{2} reveal only in
the operation of covariant differentiation of tensor fields.\par
     Any curvilinear coordinate system is naturally equipped with the
the metric tensor $\bold g$. This is a tensor whose components are
given by mutual scalar products of the frame vectors for a given
coordinate system:
$$
\hskip -2em
g_{ij}=(\bold E_i(\bold u)\,|\,\bold E_j(\bold u)).
\mytag{4.8}
$$
The components of the inverse metric tensor $\hat\bold g$ are obtained
by inverting the matrix $\bold g$. In a curvilinear coordinates the
quantities $g_{ij}$ and $g^{ij}$ are not necessarily constants any more.
\par
     We already know that the metric tensor $\bold g$ defines the volume
pseudotensor $\boldsymbol\omega$. As before, in curvilinear coordinates
its components are given by the formula \mythetagchapter{6.11}{2} from 
Chapter \uppercase\expandafter{\romannumeral 2}. Since the space $\Bbb E$
has the preferable orientation, the volume pseudotensor can be transformed
to the volume tensor $\boldsymbol\omega$. The formula
\mythetagchapter{8.1}{2} from Chapter \uppercase\expandafter{\romannumeral
2} for the components of this tensor remains valid in a curvilinear
coordinate system either.\par
\head
\S~\mysection{5} Differentiation of tensor fields\\
in curvilinear coordinates.
\endhead
\rightheadtext{\S~5. Differentiation of tensor fields \dots}
      Let $\bold A$ be a differentiable tensor field of the type $(r,s)$. 
In \S\,5 of Chapter \uppercase\expandafter{\romannumeral 2} we have
defined the concept of covariant differential. The covariant differential
$\nabla\bold A$ of a field $\bold A$ is a tensorial field of the type
$(r,s+1)$. In an arbitrary Cartesian coordinate system the components
of the field $\nabla\bold A$ are obtained by differentiating the components
of the original field $\bold A$ with respect to $x^1$, $x^2$, and $x^3$.
The use of curvilinear coordinates does not annul the operation of 
covariant differentiation. However, the procedure of deriving the
components  of the field $\nabla\bold A$ from the components of
$\bold A$  in curvilinear coordinates is more complicated.\par
      Let $u^1,\,u^2,\,u^3$ be some curvilinear coordinate system in a
domain $D\subset\Bbb E$. Let's derive the rule for covariant differentiation of tensor fields in a curvilinear coordinate system.
We consider a vectorial field $\bold A$ to begin with. This is a
field whose components are specified by one upper index: $A^i(u^1,u^2,
u^3)$. In order to calculate the components of the field $\bold B=\nabla
\bold A$ we choose some auxiliary Cartesian coordinate system $\tilde x^1,
\,\tilde x^2,\,\tilde x^3$. Then we need to do the following maneuver:
first we transform the components of $\bold A$ from curvilinear 
coordinates to Cartesian ones, then calculate the components of the
field $\bold B=\nabla\bold A$ by means of the formula
\mythetagchapter{5.1}{2} from
Chapter \uppercase\expandafter{\romannumeral 2}, and finally, we
transform the components of $\nabla\bold A$ from Cartesian coordinates
back to the original curvilinear coordinates.\par
     The Cartesian coordinates $\tilde x^1,\,\tilde x^2,\,\tilde x^3$ 
and the curvilinear coordinates $u^1,\,u^2,\,u^3$ are related by the
following transition functions:
$$
\xalignat 2
&\hskip -2em
\cases
\tilde x^1=\tilde x^1(u^1,u^2,u^3),\\
\tilde x^2=\tilde x^2(u^1,u^2,u^3),\\
\tilde x^3=\tilde x^3(u^1,u^2,u^3),
\endcases
&&\cases
u^1=u^1(\tilde x^1,\tilde x^2,\tilde x^3),\\
u^2=u^2(\tilde x^1,\tilde x^2,\tilde x^3),\\
u^3=u^3(\tilde x^1,\tilde x^2,\tilde x^3).
\endcases
\mytag{5.1}
\endxalignat
$$
The components of the corresponding transition matrices are calculated
according to the formula \mythetag{3.7}. When applied to \mythetag{5.1},
this formula yields
$$
\xalignat 2
&\hskip -2em
S^i_j(\tilde\bold x)=\frac{\partial u^i}
{\partial\tilde x^j},
&&T^k_i(\bold u)=\frac{\partial\tilde x^k}{\partial u^i}.
\mytag{5.2}
\endxalignat
$$
Denote by $\tilde A^k(\tilde x^1,\tilde x^2,\tilde x^3)$ the components
of the vector field $\bold A$ in the Cartesian coordinate system 
$\tilde x^1,\,\tilde x^2,\,\tilde x^3$. Then we get
$$
\tilde A^k=\sum^3_{p=1} T^k_p(\bold u)\,A^p(\bold u).
$$
For the components of the field $\bold B=\nabla\bold A$ in these
Cartesian coordinates, applying the formula \mythetagchapter{5.1}{2}
from Chapter \uppercase\expandafter{\romannumeral 2}, we get
$$
\hskip -2em
\tilde B^k_q=\frac{\partial\tilde A^k}{\partial\tilde x^q}
=\sum^3_{p=1}\frac{\partial}{\partial\tilde x^q}\left(
T^k_p(\bold u)\,A^p(\bold u)\right).
\mytag{5.3}
$$
Now we perform the inverse transformation of the components of 
$\bold B$ from the Cartesian coordinates $\tilde x^1,\,\tilde x^2,
\,\tilde x^3$ back to the curvilinear coordinates $u^1,\,u^2,\,u^3$:
$$
\hskip -2em
\nabla_jA^i=B^i_j(\bold u)=
\sum^3_{k=1}\sum^3_{q=1}S^i_k(\tilde\bold x)\,T^q_j(\bold u)
\,\tilde B^k_q.
\mytag{5.4}
$$
Let's apply the Leibniz rule for calculating the partial derivative 
in \mythetag{5.3}. As a result we get two sums. Then, substituting these
sums into \mythetag{5.4}, we obtain
$$
\align
\nabla_jA^i&=\sum^3_{q=1}\sum^3_{p=1}\left(
\shave{\,\sum^3_{k=1}}S^i_k(\tilde\bold x)\,T^k_p(\bold u)
\right)\,T^q_j(\bold u)\,\frac{\partial A^p(\bold u)}
{\partial\tilde x^q}+\\
&+\sum^3_{p=1}\left(\shave{\sum^3_{q=1}\sum^3_{p=1}}
S^i_k(\tilde\bold x)\,T^q_j(\bold u)\,
\frac{\partial T^k_p(\bold u)}{\partial\tilde x^q}
\right)\,A^p(\bold u).
\endalign
$$
Note that the matrices $S$ and $T$ are inverse to each other.
Therefore, we can calculate the sums over $k$ and $p$ in the first
summand. Moreover, we replace $T^q_j(\bold u)$ by the derivatives
$\partial\tilde x^q/\partial u^j$ due to the formula \mythetag{5.2}, 
and we get
$$
\sum^3_{q=1}T^q_j(\bold u)\,
\frac{\partial}{\partial\tilde x^q}
=\sum^3_{q=1}\frac{\partial\tilde x^q}{\partial u^j}\,
\frac{\partial}{\partial\tilde x^q}=
\frac{\partial}{\partial u^j}.
$$
Taking into account all the above arguments, we transform the
formula for the covariant derivative $\nabla_{\!j}A^i$ into the 
following one:
$$
\nabla_jA^i(\bold u)=\frac{\partial A^i(\bold u)}
{\partial u^j}+\sum^3_{p=1}\left(\,\shave{\sum^3_{k=1}}
S^i_k(\tilde\bold x)\,\frac{\partial T^k_p(\bold u)}
{\partial u^j}\right)\,A^p(\bold u).
$$
We introduce the special notation for the sum enclosed into the
round brackets in the above formula, we denote it by $\Gamma^i_{jp}$:
$$
\pagebreak
\hskip -2em
\Gamma^i_{jp}(\bold u)=\sum^3_{k=1}
S^i_k(\tilde\bold x)\,\frac{\partial T^k_p(\bold u)}
{\partial u^j}.
\mytag{5.5}
$$
Taking into account the notations \mythetag{5.5}, 
now we can write the rule of covariant differentiation of a 
vector field in curvilinear coordinates as follows:
$$
\hskip -2em
\nabla_jA^i=\frac{\partial A^i}
{\partial u^j}+\sum^3_{p=1}\Gamma^i_{jp}\,A^p.
\mytag{5.6}
$$
The quantities $\Gamma^i_{jp}$ calculated according to \mythetag{5.5}
are called the {\it connection compo\-nents\/} or the {\it Christoffel
symbols}. These quantities are some inner characteristics of a 
curvilinear coordinate system. This fact is supported by the following
lemma.
\mylemma{5.1} The connection components $\Gamma^i_{jp}$ of a
curvilinear coordinate system $u^1,\,u^2,\,u^3$ given by the formula
\mythetag{5.5} do not depend on the choice of an auxiliary Cartesian
coordinate system $\tilde x^1,\,\tilde x^2,\,\tilde x^3$.
\endproclaim
\demo{Proof} Let's multiply both sides of the equality \mythetag{5.5} by
the frame vector $\bold E_i$ and then sum over the index $i$:
$$
\hskip -2em
\sum^3_{i=1}\Gamma^i_{jp}(\bold u)\,\bold E_i(\bold u)=
\sum^3_{i=1}\sum^3_{k=1}
\frac{\partial T^k_p(\bold u)}{\partial u^j}
S^i_k(\tilde\bold x)\,\bold E_i(\bold u).
\mytag{5.7}
$$
The sum over $i$ in right hand side of the equality \mythetag{5.7} can be
calculated explicitly due to the first of the following two formulas:
$$
\xalignat 2
&\hskip -2em
\tilde\bold e_k=\sum^3_{i=1}S^i_k\,\bold E_i,
&&\bold E_p=\sum^3_{k=1}T^k_p\,\tilde\bold e_k.
\mytag{5.8}
\endxalignat
$$
These formulas \mythetag{5.8} relate the frame vectors $\bold E_1,\,
\bold E_2,\,\bold E_3$ and the basis vectors $\tilde\bold e_1,\,
\tilde\bold e_2,\,\tilde\bold e_3$ of the auxiliary Cartesian 
coordinate system. Now \mythetag{5.7} is written as: 
$$
\sum^3_{i=1}\Gamma^i_{jp}\,\bold E_i=
\sum^3_{k=1}\frac{\partial T^k_p(\bold u)}{\partial u^j}\,
\tilde\bold e_k=\sum^3_{k=1}\frac{\partial}{\partial u^j}
\left(T^k_p(\bold u)\,\tilde\bold e_k\right).
$$
The basis vector $\tilde\bold e_k$ does not depend on $u^1,\,u^2,\,
u^3$. Therefore, it is brought into the brackets under the differentiation
with respect to $u^j$. The sum over $k$ in right hand side of the above formula is calculated explicitly due to the second formula \mythetag{5.8}. As a result the relationship \mythetag{5.7} is transformed to the following one:
$$
\hskip -2em
\frac{\partial\bold E_p}{\partial u^j}=
\sum^3_{i=1}\Gamma^i_{jp}\cdot\bold E_i.
\mytag{5.9}
$$
The formula \mythetag{5.9} expresses the partial derivatives of the frame
vectors back through these vectors. It can be understood as another one 
way for calculating the connection components $\Gamma^i_{jp}$. This formula
comprises nothing related to the auxiliary Cartesian coordinates 
$\tilde x^1,\,\tilde x^2,\,\tilde x^3$. The vector $\bold E_p(u^1,u^2,
u^3)$ is determined by the choice of curvilinear coordinates $u^1,\,
u^2,\,u^3$ in the domain $D$. It is sufficient to differentiate this
vector with respect to $u^j$ and expand the resulting vector in the 
basis of the frame vectors $\bold E_1,\,\bold E_2,\,\bold E_3$. Then
the coefficients of this expansion yield the required values for
$\Gamma^i_{jp}$. It is obvious that these values do not depend on the
choice of the auxiliary Cartesian coordinates $\tilde x^1,\,\tilde x^2,
\,\tilde x^3$ above.\qed\enddemo
     Now let's proceed with deriving the rule for covariant differentiation
of an arbitrary tensor field $\bold A$ of the type $(r,s)$ in curvilinear
coordinates. For this purpose we need another one expression for the
connection components. It is derived from \mythetag{5.5}. Let's transform the formula \mythetag{5.5} as follows:
$$
\Gamma^i_{jp}(\bold u)=\sum^3_{k=1}\frac{\partial}{\partial u^j}
\left(S^i_k(\tilde\bold x)\,T^k_p(\bold u)\right)-
\sum^3_{k=1}T^k_p(\bold u)\,\frac{\partial S^i_k(\tilde\bold x)}
{\partial u^j}.
$$
The matrices $S$ and $T$ are inverse to each other. Therefore, upon
performing the summation over $k$ in the first term we find that it vanishes. Hence, we get
$$
\Gamma^i_{jp}(\bold u)=-\sum^3_{k=1}T^k_p(\bold u)\,
\frac{\partial S^i_k(\tilde\bold x)}
{\partial u^j}.
\mytag{5.10}
$$
Let $A^{i_1\ldots\,i_r}_{j_1\ldots\,j_s}$ be the components of a
tensor field $\bold A$ of the type $(r,s)$ in curvilinear coordinates.
In order to calculate the components of $\bold B=\nabla\bold A$ we do
the same maneuver as above. First of all we transform the components 
of $\bold A$ to some auxiliary Cartesian coordinate system:
$$
\tilde A^{p_1\ldots\,p_r}_{q_1\ldots\,q_s}=
\sum\Sb v_1\ldots\,v_r\\ w_1\ldots\,w_s\endSb
T^{p_1}_{v_1}\ldots\,T^{p_r}_{v_r}\,\,
S^{w_1}_{q_1}\ldots\,S^{w_s}_{q_s}\,\,
A^{v_1\ldots\,v_r}_{w_1\ldots\,w_s}.
$$
Then we calculate the components of the field $\bold B$ in this auxiliary
Cartesian coordinate system simply by differentiating:
$$
\tilde B^{p_1\ldots\,p_r}_{q_1\ldots\,q_{s+1}}=
\sum\Sb v_1\ldots\,v_r\\ w_1\ldots\,w_s\endSb
\frac{\partial\left(T^{p_1}_{v_1}\ldots\,T^{p_r}_{v_r}\,
\,S^{w_1}_{q_1}\ldots\,S^{w_s}_{q_s}\,\,
A^{v_1\ldots\,v_r}_{w_1\ldots\,w_s}\right)}
{\partial\tilde x^{q_{s+1}}}.
$$
Then we perform the inverse transformations of the components of
$\bold B$ from the Cartesian coordinates back to the original
curvilinear coordinate system:
$$
\gather
\hskip -2em
B^{i_1\ldots\,i_r}_{j_1\ldots\,j_{s+1}}=
\sum\Sb p_1\ldots\,p_r\\ q_1\ldots\,q_{s+1}\endSb
S^{i_1}_{p_1}\ldots\,S^{i_r}_{p_r}\,\,
T^{q_1}_{j_1}\ldots\,T^{q_{s+1}}_{j_{s+1}}\,\,
\tilde B^{p_1\ldots\,p_r}_{q_1\ldots\,q_{s+1}}=\\
\hskip -2em
\aligned
 =\sum\Sb p_1\ldots\,p_r\\ q_1\ldots\,q_{s+1}\endSb
 \ &\sum\Sb v_1\ldots\,v_r\\ w_1\ldots\,w_s\endSb
 S^{i_1}_{p_1}\ldots\,S^{i_r}_{p_r}\,\,
 T^{q_1}_{j_1}\ldots\,T^{q_{s+1}}_{j_{s+1}}\times\\
 &\times
 \frac{\partial\left(T^{p_1}_{v_1}\ldots\,T^{p_r}_{v_r}\,
 \,S^{w_1}_{q_1}\ldots\,S^{w_s}_{q_s}\,\,
 A^{v_1\ldots\,v_r}_{w_1\ldots\,w_s}\right)}
 {\partial\tilde x^{q_{s+1}}}.
 \endaligned
 \mytag{5.11}
\endgather
$$
Applying the Leibniz rule for differentiating in \mythetag{5.11}, as
a result we get three groups of summands. The summands of the first 
group correspond to differentiating the components of the matrix 
$T$, the summands of the second group arise when we differentiate 
the components of the matrix $S$ in \mythetag{5.11}, and finally, 
the unique summand in the third group is produced by differentiating
$A^{v_1\ldots\,v_r}_{w_1\ldots\,w_s}$. In any one of these summands 
if the term $T^{p_m}_{v_m}$ or the term $S^{w_n}_{q_n}$ is not
differentiated, then this term is built into a sum that can be evaluated
explicitly:
$$
\xalignat 2
&\sum^3_{p_m=1} S^{i_m}_{p_m}\,T^{p_m}_{v_m}=\delta^{i_m}_{v_m},
&&\sum^3_{q_n=1} T^{q_n}_{j_n}\,S^{w_n}_{q_n}=\delta^{w_n}_{j_n}.
\endxalignat
$$
Therefore, one can evaluate explicitly the most part of the sums
in the formula \mythetag{5.11}. Moreover, we have the following equality:
$$
\sum^3_{q_{s+1}=1}T^{q_{s+1}}_{j_{s+1}}\,
\frac{\partial}{\partial\tilde x^{q_{s+1}}}=
\sum^3_{q_{s+1}=1}\frac{\partial\tilde x^{q_{s+1}}}
{\partial u^{j_{s+1}}}\,\frac{\partial}
{\partial\tilde x^{q_{s+1}}}=
\frac{\partial}{\partial u^{j_{s+1}}}.
$$
Taking into account all the above facts, we can bring \mythetag{5.11} to 
$$
\align
\nabla_{j_{s+1}}&A^{i_1\ldots\,i_r}_{j_1\ldots\,j_s}=
 \sum^r_{m=1}\sum^3_{v_m=1}\left(\,\shave{\sum^3_{p_m=1}}
 S^{i_m}_{p_m}\,\frac{\partial T^{p_m}_{v_m}}
 {\partial u^{j_{s+1}}}\right)\,
 A^{i_1\ldots\,v_m\ldots\,i_r}_{j_1\ldots\,j_s}+\\
&+\sum^s_{n=1}\sum^3_{w_n=1}\left(\,\shave{\sum^3_{q_n=1}}
 T^{q_n}_{j_n}\,\frac{\partial S^{w_n}_{q_n}}
 {\partial u^{j_{s+1}}}\right)\,
 A^{i_1\ldots\,i_r}_{j_1\ldots\,w_n\ldots\,j_s}+
 \frac{\partial A^{i_1\ldots\,i_r}_{j_1\ldots\,j_s}}
 {\partial u^{j_{s+1}}}.
\endalign
$$
Due to the formulas \mythetag{5.5} and \mythetag{5.10} one can express the 
sums enclosed into round brackets in the above equality through the Christoffel symbols. Ultimately, the formula \mythetag{5.11} is brought 
to the following form:
$$
\hskip -2em
\gathered
\nabla_{\!j_{s+1}}A^{i_1\ldots\,i_r}_{j_1\ldots\,j_s}=
 \frac{\partial A^{i_1\ldots\,i_r}_{j_1\ldots\,j_s}}
 {\partial u^{j_{s+1}}}\,+\\
 +\sum^r_{m=1}\sum^3_{v_m=1}\Gamma^{i_m}_{j_{s+1}\,v_m}\,
 A^{i_1\ldots\,v_m\ldots\,i_r}_{j_1\ldots\,j_s}
 -\sum^s_{n=1}\sum^3_{w_n=1}
  \Gamma^{w_n}_{j_{s+1}\,j_n}\,
 A^{i_1\ldots\,i_r}_{j_1\ldots\,w_n\ldots\,j_s}.
\endgathered
\mytag{5.12}
$$
The formula \mythetag{5.12} is the rule for covariant differentiation 
of a tensorial field $\bold A$ of the type $(r,s)$ in an arbitrary
curvilinear coordinate system. This formula can be commented as follows:
the covariant derivative $\nabla_{\!j_{s+1}}$ is obtained from the partial 
derivative $\partial/\partial u_{j_{s+1}}$ by adding $r+s$ terms --- one
per each index in the components of the field $\bold A$. The terms
associated with the upper indices enter with the positive sign, the other terms associated with the lower indices enter with the negative sign.
In such additional terms each of the upper indices $i_m$ and each of the
lower indices $j_n$ are sequentially moved to the Christoffel symbol,
while in its place we write the summation index $v_m$ or $w_n$. The lower
index $j_{s+1}$ added as a result of covariant differentiation is always
written as the first lower index in Christoffel symbols. The position of
the summation indices $v_m$ and $w_n$ in Christoffel symbols is always 
complementary to their positions in the components of the field $\bold A$
so that they always form a pair of upper and lower indices. Though the  formula \mythetag{5.12} is rather huge, we hope that due to the above
comments one can easily remember it and reproduce it in any particular 
case.\par
\head
\S~\mysection{6} Transformation of the connection components
under a change of a coordinate system.
\endhead
\rightheadtext{\S~6. Transformation of connection components \dots}
     In deriving the formula for covariant differentiation of tensorial
fields in curvilinear coordinates we discovered a new type of indexed
objects --- these are Christoffel symbols. The quantities $\Gamma^k_{ij}$
are enumerated  by one upper index and two lower indices, and their values 
are determined by the choice of a coordinate system. However, they are not 
the components of a tensorial fields of the type $(1,2)$. Indeed, the
values of all $\Gamma^k_{ij}$ in a Cartesian coordinate system are 
identically zero (this follows from the comparison of \mythetag{5.12} with
the formula \mythetagchapter{5.1}{2} in Chapter
\uppercase\expandafter{\romannumeral 2}). But a tensorial field with purely
zero components in some coordinate system cannot have nonzero components in
any other coordinate system. Therefore, Christoffel symbols are the
components of a non-tensorial geometric object which is called a {\it
connection field} or simply a  {\it connection}.
\mytheorem{6.1} Let $u^1,\,u^2,\,u^3$ and $\tilde u^1,\,
\tilde u^2,\,\tilde u^3$ be two coordinate systems in a domain 
$D\subset\Bbb E$. Then the connection components in these two 
coordinate systems are related to each other by means of the 
following equality:
$$
\hskip -2em
\Gamma^k_{ij}=\sum^3_{m=1}\sum^3_{p=1}\sum^3_{q=1}
S^k_m\,T^p_i\,T^q_j\,\tilde\Gamma^m_{pq}
+\sum^3_{m=1}S^k_m\,\frac{\partial T^m_i}{\partial u^j}.
\mytag{6.1}
$$
Here $S$ and $T$ are the transition matrices given by the 
formulas \mythetag{3.7}.
\endproclaim
\subhead A remark on the smoothness\endsubhead The derivatives 
of the components of $T$ in  \mythetag{6.1} and the formulas 
\mythetag{3.7}, where the components of $T$ are defined as the 
partial derivatives of the transition functions \mythetag{3.5}, 
show that the connection components can be correctly defined
only for coordinate systems of the smoothness class not lower
than $C^2$. The same conclusion follows from the formula 
\mythetag{5.5} for $\Gamma^i_{jp}$.
\demo{Proof} In order to prove the theorem~\mythetheorem{6.1} we 
apply the formula \mythetag{5.9}. Let's write it for the frame 
vectors $\bold E_1,\,\bold E_2,\,\bold E_3$, then apply the formula
\mythetag{3.6} for to express $\bold E_j$ through the vectors $\tilde
\bold E_1$, $\tilde\bold E_2$, and $\tilde\bold E_3$:
$$
\hskip -2em
\sum^3_{k=1}\Gamma^k_{ij}\,\bold E_k=
\frac{\partial\tilde\bold E_j}{\partial u^i}=
\sum^3_{m=1}\frac{\partial}{\partial u^i}\left(
T^m_j\,\tilde\bold E_m\right).
\mytag{6.2}
$$
Applying the Leibniz rule to the right hand side of \mythetag{6.2}, we get
two terms:
$$
\hskip -2em
\sum^3_{k=1}\Gamma^k_{ij}\,\bold E_k=
\sum^3_{m=1}\frac{\partial T^m_j}{\partial u^i}\,
\tilde\bold E_m
+\sum^3_{q=1}T^q_j\,\frac{\partial\tilde\bold E_q}
{\partial u^i}.
\mytag{6.3}
$$
In the first term in the right hand side of \mythetag{6.3} we express
$\tilde\bold E_m$ through the vectors $\bold E_1$, $\bold E_2$, and
$\bold E_3$. In the second term we apply the chain rule and express
the derivative with respect to $u^i$ through the derivatives 
with respect to $\tilde u^1,\,\tilde u^2,\,\tilde u^3$:
$$
\sum^3_{k=1}\Gamma^k_{ij}\,\bold E_k=
\sum^3_{k=1}\sum^3_{m=1}S^k_m\,\frac{\partial T^m_j}
{\partial u^i}\,\bold E_k+
\sum^3_{q=1}\sum^3_{p=1}T^q_j\,\frac{\partial\tilde u^p}
{\partial u^i}\,\frac{\partial\tilde\bold E_q}{\partial\tilde u^p}.
$$
Now let's replace $\partial\tilde u^p/\partial u^i$ by $T^p_i$ relying
upon the formulas \mythetag{3.7} and then apply the relationship \mythetag{5.9} once more in the form of 
$$
\frac{\partial\tilde\bold E_q}{\partial\tilde u^p}=
\sum^3_{m=1}\tilde\Gamma^m_{pq}\,\tilde\bold E_m.
$$
As a result of the above transformations we can write the equality \mythetag{6.3} as follows:
$$
\sum^3_{k=1}\Gamma^k_{ij}\,\bold E_k=
\sum^3_{k=1}\sum^3_{m=1}S^k_m\,\frac{\partial T^m_j}
{\partial u^i}\,\bold E_k+
\sum^3_{q=1}\sum^3_{p=1}\sum^3_{m=1}T^p_i\,T^q_j\,
\tilde\Gamma^m_{pq}\,\tilde\bold E_m.
$$
Now we need only to express $\tilde\bold E_m$ through the frame vectors
$\bold E_1,\,\bold E_2,\,\bold E_3$ and collect the similar terms in the
above formula:
$$
\sum^3_{k=1}\left(\Gamma^k_{ij}-
\shave{\sum^3_{m=1}\sum^3_{q=1}\sum^3_{p=1}}S^k_m\,T^p_i\,
T^q_j\,\tilde\Gamma^m_{pq}-\shave{\sum^3_{m=1}}
S^k_m\,\frac{\partial T^m_j}{\partial u^i}\right)\bold E_k=0.
$$
Since the frame vectors $\bold E_1,\,\bold E_2,\,\bold E_3$ are
linearly independent, the expression enclosed into round brackets 
should vanish. As a result we get the equality exactly equivalent 
to the relationship \mythetag{6.1} that we needed to prove.
\qed\enddemo
\head
\S~\mysection{7} Concordance of metric and connection. 
Another formula for Christoffel symbols.
\endhead
\rightheadtext{\S~7. Concordance of metric and connection.}
     Let's consider the metric tensor $\bold g$. The covariant 
differential $\nabla\bold g$ of the field $\bold g$ is equal to 
zero (see formulas \mythetagchapter{6.7}{2} in Chapter
\uppercase\expandafter{\romannumeral 2}). This is 
because in any Cartesian coordinates $x^1,\,x^2,\,x^3$ in
$\Bbb E$ the components $g_{ij}$ of the metric tensor do not depend
on $x^1,\,x^2,\,x^3$. In a curvilinear coordinate system the components
of the metric tensor $g_{ij}(u^1,u^2,u^3)$ usually are not constants.
However, being equal to zero in Cartesian coordinates, the tensor
$\nabla\bold g$ remains zero in any other coordinates:
$$
\hskip -2em
\nabla_kg_{ij}=0.
\mytag{7.1}
$$
The relationship \mythetag{7.1} is known as the {\it concordance condition}
for a metric and a connection. Taking into account \mythetag{5.12}, we can
rewrite this condition as 
$$
\hskip -2em
\frac{g_{ij}}{\partial u^k}-
\sum^3_{r=1}\Gamma^r_{ki}\,g_{rj}-
\sum^3_{r=1}\Gamma^r_{kj}\,g_{ir}=0.
\mytag{7.2}
$$
The formula \mythetag{7.2} relates the connection components $\Gamma^k_{ij}$
and the components of the metric tensor $g_{ij}$. Due to this relationship
we can express $\Gamma^k_{ij}$ through the components of the metric tensor
provided we remember the following very important property of the connection
components \mythetag{5.5}.
\mytheorem{7.1} The connection given by the formula \mythetag{5.5} is
a symmetric con\-nection, i\.\,e\. $\Gamma^k_{ij}=\Gamma^k_{j\,i}$.
\endproclaim
\demo{Proof} From \mythetag{5.2} and \mythetag{5.5} for $\Gamma^k_{ij}$ we derive the following expression:
$$
\hskip -2em
\Gamma^k_{ij}(\bold u)=\sum^3_{q=1}
S^k_q\,\frac{\partial T^q_j(\bold u)}
{\partial u^i}=\sum^3_{q=1} S^k_q\,
\frac{\partial^2\tilde x^q}
{\partial u^j\,\partial u^i}.
\mytag{7.3}
$$
For the functions of the smoothness class $C^2$ the mixed second order partial derivatives do not depend on the order of differentiation:
$$
\frac{\partial^2\tilde x^q}{\partial u^j\,\partial u^i}=
\frac{\partial^2\tilde x^q}{\partial u^i\,\partial u^j}.
$$
This fact immediately proves the symmetry of the Christoffel symbols
given by the formula \mythetag{7.3}. Thus, the proof is over.
\qed\enddemo
     Now, returning back to the formula \mythetag{7.2} relating 
$\Gamma^k_{ij}$ and $g_{ij}$, we introduce the following notations that
simplify the further calculations:
$$
\hskip -2em
\Gamma_{ijk}=\sum^3_{r=1}\Gamma^r_{ij}\,g_{kr}.
\mytag{7.4}
$$
It is clear that the quantities $\Gamma_{ijk}$ in \mythetag{7.4} are produced from $\Gamma_{ij}^k$ by means of index lowering procedure
described in Chapter \uppercase\expandafter{\romannumeral 2}. Therefore,
conversely, $\Gamma^k_{ij}$ are obtained from $\Gamma_{ijk}$ according to
the following formula:
$$
\hskip -2em
\Gamma^k_{ij}=\sum^3_{r=1}g^{kr}\,\Gamma_{ijr}.
\mytag{7.5}
$$
From the symmetry of $\Gamma^k_{ij}$ it follows that the quantities
$\Gamma_{ijk}$ in \mythetag{7.4} are also symmetric with respect to
the indices $i$ and $j$, i\.\,e\. $\Gamma_{ijk}=\Gamma_{jik}$. Using 
the notations \mythetag{7.4} and the symmetry of the metric tensor, the
relationship \mythetag{7.2} can be rewritten in the following way:
$$
\hskip -2em
\frac{\partial g_{ij}}{\partial u^k}-\Gamma_{kij}-\Gamma_{kji}=0.
\mytag{7.6}
$$
Let's complete \mythetag{7.6} with two similar relationships applying 
two cyclic transposi\-tions of the indices $i\to j\to k\to i$ to the
formula \mythetag{7.6}. As a result we obtain
$$
\hskip -2em
\aligned
&\frac{\partial g_{ij}}{\partial u^k}-\Gamma_{kij}-\Gamma_{kji}=0,\\
\vspace{1ex}
&\frac{\partial g_{jk}}{\partial u^i}-\Gamma_{ijk}-\Gamma_{ikj}=0,\\
\vspace{1ex}
&\frac{\partial g_{ki}}{\partial u^j}-\Gamma_{jki}-\Gamma_{jik}=0.
\endaligned
\mytag{7.7}
$$
Let's add the last two relationships \mythetag{7.7} and subtract the first
one from the sum. Taking into account the symmetry of $\Gamma_{ijk}$ with
respect to $i$ and $j$, we get
$$
\frac{\partial g_{jk}}{\partial u^i}+\frac{\partial g_{ki}}
{\partial u^j}-\frac{\partial g_{ij}}{\partial u^k}-
2\,\Gamma_{ijk}=0.
$$
Using this equality, one can easily express $\Gamma_{ijk}$ through
the components of the metric tensor. Then one can substitute this
expression into \mythetag{7.5} and derive
$$
\hskip -2em
\Gamma^k_{ij}=\frac{1}{2}\sum^3_{r=1}g^{kr}\left(
\frac{\partial g_{rj}}{\partial u^i}+\frac{\partial g_{ir}}
{\partial u^j}-\frac{\partial g_{ij}}{\partial u^r}\right).
\mytag{7.8}
$$
The relationship \mythetag{7.8} is another formula for the Christoffel
symbols $\Gamma^k_{ij}$, it follows from the symmetry of $\Gamma^k_{ij}$
and from the concordance condition for the metric and connection. It is
different from \mythetag{5.5} and \mythetag{5.10}. The relationship 
\mythetag{7.8} has the important advantage as compared to \mythetag{5.5}:
one should not use an auxiliary Cartesian coordinate system for to
apply it. As compared to \mythetag{5.9}, in \mythetag{7.8} one should not 
deal with vector-functions $\bold E_i(u^1,u^2,u^3)$. All calculations
in \mythetag{7.8} are performed within a fixed curvilinear coordinate
system provided the components of the metric tensor in this coordinate
system are known.\par
\head
\S~\mysection{8} Parallel translation. The equation\\
of a straight line in curvilinear coordinates.
\endhead
\rightheadtext{\S~8. Parallel translation.}
\vskip 1.2cm
\special{em:graph ris-1-11.bmp}
\vskip -1.2cm
\parshape 16 0pt 360pt 0pt 360pt
180pt 180pt 180pt 180pt 180pt 180pt 180pt 180pt 
180pt 180pt 180pt 180pt 180pt 180pt 180pt 180pt 
180pt 180pt 180pt 180pt 180pt 180pt 180pt 180pt 
180pt 180pt 0pt 360pt
     Let $\bold a$ be a nonzero vector attached to some point $A$ in the
space $\Bbb E$. In a Euclidean space there is a procedure of parallel
translation; applying this procedure one can bring the vector $\bold a$
\vadjust{\vskip 5pt\hbox to 0pt{\kern -10pt
\includegraphics{ris11.eps}\hss}\vskip -5pt}
from the point $A$ to some other point $B$. This procedure does change
neither the modulus nor the direction of the vector $\bold a$ being
translated. In a Cartesian coordinate system the procedure of parallel
translation is des\-cribed in the most simple way: the ori\-ginal vector
$\bold a$ at the point $A$ and the translated vector $\bold a$ at the
point $B$ have the equal coordinates. In a curvilinear coordinate system
the frame vectors at the point $A$ and the frame vector at the point
$B$ form two different bases. Therefore, the components of the vector
$\bold a$ in the following two expansions
$$
\hskip -2em
\aligned
&\bold a=a^1(A)\cdot\bold E_1(A)+
a^2(A)\cdot\bold E_2(A)+a^3(A)\cdot\bold E_3(A),\\
&\bold a=a^1(B)\cdot\bold E_1(B)+
a^2(B)\cdot\bold E_2(B)+a^3(B)\cdot\bold E_3(B)
\endaligned
\mytag{8.1}
$$
in general case are different. If the points $A$ and $B$ are closed to
each other, then the triples of vectors $\bold E_1(A),\,\bold E_2(A),\,
\bold E_3(A)$ and $\bold E_1(B),\,\bold E_2(B),\,\bold E_3(B)$ are
approximately the same. Hence, in this case the components of the vector
$\bold a$ in the expansions \mythetag{8.1} are slightly different from each other. This consideration shows that in curvilinear coordinates the
parallel translation should be performed gradually: one should first
converge the point $B$ with the point $A$, then slowly move the point
$B$ toward its ultimate position and record the coordinates of the
vector $\bold a$ in the second expansion \mythetag{8.1} at each intermediate
position of the point $B$. The most simple way to implement this plan is
to link $A$ and $B$ with some smooth parametric curve $\bold r=\bold r(t)$,
where $t\in [0,\,1]$. In a curvilinear coordinate system a parametric
curve is given by three functions $u^1(t),\,u^2(t),\,u^3(t)$ that for each
$t\in [0,\,1]$ yield the coordinates of the corresponding point on the
curve.
\mytheorem{8.1} For a parametric curve given by three functions
$u^1(t)$, $u^2(t)$, and $u^3(t)$ in some curvilinear coordinate system
the components of the tangent vector $\boldsymbol\tau(t)$ in the moving
frame of that coordinate system are determined by the derivatives 
$\dot u^1(t),\,\dot u^2(t),\,\dot u^3(t)$.
\endproclaim
\demo{Proof} Curvilinear coordinates $u^1,\,u^2,\,u^3$ determine the
position of a point in the space by means of the vector-function
$\bold r=\bold r(u^1,\,u^2,\,u^3)$, where $\bold r$ is the radius-vector
of that point in some auxiliary Cartesian coordinate system (see formulas
\mythetag{2.4} and \mythetag{2.5}). Therefore, the vectorial-parametric
equation of the curve is represented in the following way:
$$
\hskip -2em
\bold r=\bold r(u^1(t),u^2(t),u^3(t)).
\mytag{8.1}
$$
Applying the chain rule to the function $\bold r(t)$ in \mythetag{8.1},
we get
$$
\hskip -2em
\boldsymbol\tau(t)=\frac{d\bold r}{dt}=
\sum^3_{j=1}\frac{\partial\bold r}{\partial u^j}\cdot\dot u^j(t).
\mytag{8.2}
$$
Remember that due to the formula \mythetag{2.7} the partial derivatives
in \mythetag{8.2} coincide with the frame vectors of the curvilinear
coordinate system. Therefore the formula \mythetag{8.2} itself can be
rewritten as follows:
$$
\boldsymbol\tau(t)=
\sum^3_{j=1}\dot u^j(t)\cdot\bold E_j(u^1(t),u^2(t),u^3(t)).
\mytag{8.3}
$$
It is easy to see that \mythetag{8.3} is the expansion of the tangent
vector $\boldsymbol\tau(t)$ in the basis formed by the frame vectors
of the curvilinear coordinate system. The components of the vector
$\boldsymbol\tau(t)$ in the expansion \mythetag{8.3} are the derivatives
$\dot u^1(t),\,\dot u^2(t),\,\dot u^3(t)$. The theorem is proved.
\qed\enddemo
     Let's apply the procedure of parallel translation to the vector
$\bold a$ and translate this vector to all points of the curve linking
the points $A$ and $B$ (see Fig\.~8.1). Then we can write the following 
expansion for this vector
$$
\hskip -2em
\bold a=\sum^3_{i=1}a^i(t)\cdot\bold E_i(u^1(t),u^2(t),u^3(t)).
\mytag{8.4}
$$
This expansion is analogous to \mythetag{8.3}. Let's differentiate the relationship \mythetag{8.4} with respect to the parameter $t$ and take 
into account that $\bold a=\const$:
$$
0=\frac{d\bold a}{dt}=\sum^3_{i=1}\dot a^i\cdot\bold E_i+
\sum^3_{i=1}\sum^3_{j=1}a^i\,\frac{\partial\bold E_i}
{\partial u^j}
\,\dot u^j.
$$
Now let's use the formula \mythetag{5.9} in order to differentiate 
the frame vectors of the curvilinear coordinate system. As a result
we derive
$$
\sum^3_{i=1}\left(\dot a^i+\shave{\sum^3_{j=1}\sum^3_{k=1}}
\Gamma^i_{jk}\,\dot u^j\,a^k\right)\cdot\bold E_i=0.
$$
Since the frame vectors $\bold E_1,\,\bold E_2,\,\bold E_3$ are 
linearly independent, we obtain
$$
\hskip -2em
\dot a^i+\sum^3_{j=1}\sum^3_{k=1}
\Gamma^i_{jk}\,\dot u^j\,a^k=0.
\mytag{8.5}
$$
The equation \mythetag{8.5} is called the {\it differential 
equation of the parallel translation of a vector along a curve}.
This is the system of three linear differential equations of the
first order with respect to the components of the vector 
$\bold a$. Actually, in order to perform the parallel translation
of a vector $\bold a$ from the point $A$ to the point $B$ in
curvilinear coordinates one should set the initial data for the
components of the vector $\bold a$ at the point $A$ (i\.\,e\. for 
$t=0$) and then solve the Cauchy problem for the equations 
\mythetag{8.5}.\par
     The procedure of the parallel translation of vectors along
curves leads us to the situation where at each point of a curve 
in $\Bbb E$ we have some vector attached to that point. The same
situation arises in considering the vectors $\boldsymbol\tau$, 
$\bold n$, and $\bold b$ that form the Frenet frame of a curve
in $\Bbb E$ (see Chapter \uppercase\expandafter{\romannumeral 1}).
Generalizing this situation one can consider the set of tensors
of the type $(r,s)$ attached to the points of some curve. Defining
such a set of tensors differs from defining a tensorial field in
$\Bbb E$ since in order to define a tensor field in $\Bbb E$ one
should attach a tensor to each point of the space, not only to
the points of a curve. In the case, where the tensors of the type
$(r,s)$ are defined only at the points of a curve, we say that 
a {\it tensor field of the type $(r,s)$ on a curve\/} is given.
In order to write the components of such a tensor field $\bold A$
we can use the moving frame $\bold E_1,\,\bold E_2,\,\bold E_3$
of some curvilinear coordinate system in some  neighborhood of
the curve. These components form a set of functions of the scalar
parameter $t$ specifying the points of the curve:
$$
\hskip -2em
A^{i_1\ldots\,i_r}_{j_1\ldots\,j_s}=
A^{i_1\ldots\,i_r}_{j_1\ldots\,j_s}(t).
\mytag{8.6}
$$
Under a change of curvilinear coordinate system the quantities
\mythetag{8.6} are transfor\-med according to the standard rule 
$$
\hskip -2em
\aligned
A^{i_1\ldots\,i_r}_{j_1\ldots\,j_s}(t)=
\sum\Sb p_1\ldots\,p_r\\ q_1\ldots\,q_s\endSb
 &S^{i_1}_{p_1}(t)\ldots\,S^{i_r}_{p_r}(t)\times\\
 \vspace{-1ex}
 &\ \times T^{q_1}_{j_1}(t)\ldots\,T^{q_s}_{j_s}(t)
  \,\,\tilde A^{p_1\ldots\,p_r}_{q_1\ldots\,q_s}(t),
\endaligned
\mytag{8.7}
$$
where $S(t)$ and $T(t)$ are the values of the transition matrices
at the points of the curve. They are given by the following
formulas:
$$
\hskip -2em
\aligned
&S(t)=S(\tilde u^1(t),\tilde u^2(t),\tilde u^3(t)),\\
&T(t)=T(u^1(t),u^2(t),u^3(t)).
\endaligned
\mytag{8.8}
$$
We cannot use the formula \mythetag{5.12} for differentiating the
field $\bold A$ on the curve since the only argument, which the
functions \mythetag{8.6} depend on, is the parameter $t$. Therefore,
we need to modify the formula \mythetag{5.12} as follows:
$$
\gathered
\nabla_{\!t}A^{i_1\ldots\,i_r}_{j_1\ldots\,j_s}=
 \frac{dA^{i_1\ldots\,i_r}_{j_1\ldots\,j_s}}{dt}\,+\\
 +\sum^r_{m=1}\sum^3_{q=1}\sum^3_{v_m=1}
 \Gamma^{i_m}_{q\,v_m}\,\dot u^q\,
 A^{i_1\ldots\,v_m\ldots\,i_r}_{j_1\ldots\,j_s}-
 \sum^s_{n=1}\sum^3_{q=1}\sum^3_{w_n=1}
  \Gamma^{w_n}_{q\,j_n}\,\dot u^q\,
 A^{i_1\ldots\,i_r}_{j_1\ldots\,w_n\ldots\,j_s}.
\endgathered\quad
\mytag{8.9}
$$
The formula \mythetag{8.9} expresses the {\it rule for covariant
differentiation of a tensor field $\bold A$ with respect to the
parameter $t$ along a parametric curve} in curvilinear \pagebreak
coordinates
$u^1,u^2,u^3$. Unlike \mythetag{5.12}, the index $t$ beside the nabla
sign is not an additional index. It is set only for to denote the 
variable $t$ with respect to which the differentiation in the formula
\mythetag{8.9} is performed.
\mytheorem{8.2} Under a change of coordinates $u^1,\,u^2,\,u^3$ 
for other coordinates $\tilde u^1,\,\tilde u^2,\,\tilde u^3$ the quantities 
$B^{i_1\ldots\,i_r}_{j_1\ldots\,j_s}=\nabla_{\!t}
A^{i_1\ldots\,i_r}_{j_1\ldots\,j_s}$ calculated by means of the formula
\mythetag{8.9} are transformed according to the rule \mythetag{8.7} and define
a tensor field \ $\bold B=\nabla_{\!t}\bold A$ of the type $(r,s)$ which
is called the {\it covariant derivative of the field $\bold A$ with
respect to the parameter $t$ along a curve}.
\endproclaim
\demo{Proof} The proof of this theorem is pure calculations. Let's begin
with the first term in \mythetag{8.9}. Let's express $A^{i_1\ldots\,i_r}_{j_1
\ldots\,j_s}$ through the components of the field $\bold A$ in the other
coordinates $\tilde u^1,\,\tilde u^2,\,\tilde u^3$ by means of \mythetag{8.7}.
In calculating $dA^{i_1\ldots\,i_r}_{j_1\ldots\,j_s}/dt$ this is equivalent
to differentiating both sides of \mythetag{8.7} with respect to $t$:
$$
\aligned
&\frac{dA^{i_1\ldots\,i_r}_{j_1\ldots\,j_s}}{dt}=
\sum\Sb p_1\ldots\,p_r\\ q_1\ldots\,q_s\endSb
 S^{i_1}_{p_1}\ldots\,S^{i_r}_{p_r}\,T^{q_1}_{j_1}
 \ldots\,T^{q_s}_{j_s}\,\frac{d\tilde A^{p_1\ldots\,
 p_r}_{q_1\ldots\,q_s}}{dt}+\\
&+\sum^r_{m=1}\sum\Sb p_1\ldots\,p_r\\ q_1\ldots\,q_s\endSb
 S^{i_1}_{p_1}\ldots\,\dot S^{i_m}_{p_m}\ldots\,S^{i_r}_{p_r}\,
 T^{q_1}_{j_1}\ldots\,T^{q_s}_{j_s}\,
 \tilde A^{p_1\ldots\,p_r}_{q_1\ldots\,q_s}+\\
&+\sum^s_{n=1}\sum\Sb p_1\ldots\,p_r\\ q_1\ldots\,q_s\endSb
 S^{i_1}_{p_1}\ldots\,S^{i_r}_{p_r}\,
 T^{q_1}_{j_1}\ldots\,\dot T^{q_n}_{j_n}\ldots\,T^{q_s}_{j_s}\,
 \tilde A^{p_1\ldots\,p_r}_{q_1\ldots\,q_s}.
\endaligned
\mytag{8.10}
$$
For to calculate the derivatives $\dot S^{i_m}_{p_m}$ and $\dot
T^{q_n}_{j_n}$ in \mythetag{8.10} we use the fact that the transition
matrices $S$ and $T$ are inverse to each other:
$$
\align
\dot S^{i_m}_{p_m}&=\sum^3_{k=1}\sum^3_{v_m=1}\dot S^{i_m}_k
\,T^k_{v_m}\,S^{v_m}_{p_m}=
\sum^3_{v_m=1}\sum^3_{k=1}\frac{d\left(S^{i_m}_k\,
T^k_{v_m}\right)}{dt}\,S^{v_m}_{p_m}-\\
&-\sum^3_{v_m=1}\sum^3_{k=1}S^{i_m}_k\,\frac{dT^k_{v_m}}{dt}
\,S^{v_m}_{p_m}=-\sum^3_{v_m=1}\left(\,\shave{\sum^3_{k=1}}
S^{i_m}_k\,\frac{dT^k_{v_m}}{dt}\right)\,S^{v_m}_{p_m},\\
\vspace{1ex}
\dot T^{q_n}_{j_n}&=\sum^3_{k=1}\sum^3_{w_n=1}
\dot T^k_{j_n}\,S^{w_n}_k\,T^{q_n}_{w_n}=\sum^3_{w_n=1}
\left(\,\shave{\sum^3_{k=1}}\frac{dT^k_{j_n}}{dt}\,S^{w_n}_k
\right)\,T^{q_n}_{w_n}.
\endalign
$$
In order to transform further the above formulas for the 
derivatives $\dot S^{i_m}_{p_m}$ and $\dot T^{q_n}_{j_n}$ 
we use the second formula in \mythetag{8.8}:
$$
\align
&\hskip -2em
\dot S^{i_m}_{p_m}=-\sum^3_{v_m=1}\left(\,\shave{\sum^3_{k=1}
\sum^3_{q=1}} S^{i_m}_k\,\frac{\partial T^k_{v_m}}{\partial u^q}
\,\dot u^q\right)\,S^{v_m}_{p_m},\hskip -1em
\mytag{8.11}\\
&\hskip -2em
\dot T^{q_n}_{j_n}=\sum^3_{w_n=1}\left(\,\shave{\sum^3_{k=1}
\sum^3_{q=1}} S^{w_n}_k\,\frac{\partial T^k_{j_n}}{\partial u^q}
\,\dot u^q\right)\,T^{q_n}_{w_n}.
\mytag{8.12}
\endalign
$$
Let's substitute \mythetag{8.11} and \mythetag{8.12} into \mythetag{8.10}.
Then, taking into account the relationship \mythetag{8.7}, we can perform
the summation over $p_1,\ldots,\,p_r$ and $q_1,\ldots,\,q_s$ in the second
and the third terms in \mythetag{8.10} thus transforming \mythetag{8.10} to
$$
\aligned
&\frac{dA^{i_1\ldots\,i_r}_{j_1\ldots\,j_s}}{dt}=
\sum\Sb p_1\ldots\,p_r\\ q_1\ldots\,q_s\endSb
 S^{i_1}_{p_1}\ldots\,S^{i_r}_{p_r}\,T^{q_1}_{j_1}
 \ldots\,T^{q_s}_{j_s}\,\frac{d\tilde A^{p_1\ldots\,
 p_r}_{q_1\ldots\,q_s}}{dt}\,-\\
&-\sum^r_{m=1}\sum^3_{q=1}\sum^3_{v_m=1}
 \left(\,\shave{\sum^3_{k=1}}
 S^{i_m}_k\,\frac{\partial T^k_{v_m}}{\partial u^q}
 \right)\,\dot u^q\,
 A^{i_1\ldots\,v_m\ldots\,i_r}_{j_1\ldots\,j_s}\,+\\
&+\sum^s_{n=1}\sum^3_{q=1}\sum^3_{w_n=1}
 \left(\,\shave{\sum^3_{k=1}} S^{w_n}_k\,
 \frac{\partial T^k_{j_n}}{\partial u^q}
 \right)\,\dot u^q\,
 A^{i_1\ldots\,i_r}_{j_1\ldots\,w_n\ldots\,j_s}.
\endaligned
\mytag{8.13}
$$
The second and the third terms in \mythetag{8.9} and \mythetag{8.13} 
are similar in their structure. Therefore, one can collect the similar
terms upon substituting \mythetag{8.13} into \mythetag{8.9}. Collecting
these similar terms, we get the following two expressions
$$
\xalignat 2
&\Gamma^{i_m}_{q\,v_m}-\sum^3_{k=1} S^{i_m}_k\,
 \frac{\partial T^k_{v_m}}{\partial u^q},
&&\Gamma^{w_n}_{q\,j_n}-\sum^3_{k=1} S^{w_n}_k\,
 \frac{\partial T^k_{j_n}}{\partial u^q}
\endxalignat
$$
as the coefficients. Let's apply \mythetag{6.1} to these expressions:
$$
\hskip -2em
\aligned
&\Gamma^{i_m}_{q\,v_m}-\sum^3_{k=1} S^{i_m}_k\,
 \frac{\partial T^k_{v_m}}{\partial u^q}
=\sum^3_{r=1}\sum^3_{p=1}\sum^3_{k=1}
S^{i_m}_r\,\tilde\Gamma^r_{pk}\,T^p_q\,T^k_{v_m},\\
&\Gamma^{w_n}_{q\,j_n}-\sum^3_{k=1} S^{w_n}_k\,
 \frac{\partial T^k_{j_n}}{\partial u^q}
=\sum^3_{r=1}\sum^3_{p=1}\sum^3_{k=1}
S^{w_n}_r\,\tilde\Gamma^r_{pk}\,T^p_q\,T^k_{j_n}.
\endaligned
\mytag{8.14}
$$
If we take into account \mythetag{8.14} when substituting 
\mythetag{8.13} into \mythetag{8.9}, then the equality \mythetag{8.9} 
is written in the following form:
$$
\aligned
&\nabla_{\!t}A^{i_1\ldots\,i_r}_{j_1\ldots\,j_s}=
\sum\Sb p_1\ldots\,p_r\\ q_1\ldots\,q_s\endSb
 S^{i_1}_{p_1}\ldots\,S^{i_r}_{p_r}\,T^{q_1}_{j_1}
 \ldots\,T^{q_s}_{j_s}\,\frac{d\tilde A^{p_1\ldots\,
 p_r}_{q_1\ldots\,q_s}}{dt}+\\
&+\sum^r_{m=1}\sum^3_{q=1}\sum^3_{v_m=1}\sum^3_{p_m=1}\sum^3_{p=1}
\sum^3_{k=1}
S^{i_m}_{p_m}\,\tilde\Gamma^{p_m}_{pk}\,T^p_q\,T^k_{v_m}
 \,\dot u^q\,
 A^{i_1\ldots\,v_m\ldots\,i_r}_{j_1\ldots\,j_s}-\\
&-\sum^s_{n=1}\sum^3_{q=1}\sum^3_{w_n=1}\sum^3_{q_n=1}\sum^3_{p=1}
\sum^3_{k=1}
S^{w_n}_k\,\tilde\Gamma^k_{p\,q_n}\,T^p_q\,T^{q_n}_{j_n}
 \,\dot u^q\,
 A^{i_1\ldots\,i_r}_{j_1\ldots\,w_n\ldots\,j_s}.
\endaligned
$$
In order to transform it further we express
$A^{i_1\ldots\,v_m\ldots\,i_r}_{j_1\ldots\,j_s}$ and
$A^{i_1\ldots\,i_r}_{j_1\ldots\,w_n\ldots\,j_s}$ through the 
components of the field $\bold A$ in the other coordinate system
by means of the formula \mythetag{8.7}. Moreover, we take into
account that $T^p_q\,\dot u^q$ upon summing up over $q$ yields
$\dot{\tilde u}\vphantom{u}^p$. As a result we obtain:
$$
\hskip -0.1em
\gathered
\nabla_{\!t}A^{i_1\ldots\,i_r}_{j_1\ldots\,j_s}=
\sum\Sb p_1\ldots\,p_r\\ q_1\ldots\,q_s\endSb
 S^{i_1}_{p_1}\ldots\,S^{i_r}_{p_r}\,T^{q_1}_{j_1}
 \ldots\,T^{q_s}_{j_s}\left(\frac{d\tilde A^{p_1\ldots\,
 p_r}_{q_1\ldots\,q_s}}{dt}\right.+\\
 +\sum^r_{m=1}\sum^3_{p=1}\sum^3_{v_m=1}
 \tilde\Gamma^{p_m}_{p\,v_m}\,\dot{\tilde u}\vphantom{u}^p\,
 \tilde A^{p_1\ldots\,v_m\ldots\,p_r}_{q_1\ldots\,q_s}-
 \left.\sum^s_{n=1}\sum^3_{p=1}\sum^3_{w_n=1}
  \tilde\Gamma^{w_n}_{p\,q_n}\,\dot{\tilde u}\vphantom{u}^p\,
 \tilde A^{p_1\ldots\,p_r}_{q_1\ldots\,w_n\ldots\,q_s}\right)\!.
\endgathered
\mytag{8.15}
$$
Note that the expression enclosed into round brackets in 
\mythetag{8.15} is $\nabla_{\!t}\tilde A^{p_1\ldots\,p_r}_{q_1\ldots\,q_s}$
exactly. Therefore, the formula \mythetag{8.15} means that the components
of the field\linebreak $\nabla_{\!t}\bold A$ on a curve calculated according 
to the formula \mythetag{8.9} obey the transformation rule \mythetag{8.7}. Thus, the theorem~\mythetheorem{8.2} is proved.
\qed\enddemo
     Now let's return to the formula \mythetag{8.5}. The left hand side
of this formula coincides with the expression \mythetag{8.9} for the 
covariant derivative of the vector field $\bold a$ with respect to the
parameter $t$. Therefore, the equation of parallel translation can be
written as $\nabla_{\!t}\bold a=0$. In this form, the equation of parallel
translation can be easily generalized for the case of an arbitrary 
tensor $\bold A$:
$$
\hskip -2em
\nabla_{\!t}\bold A=0.
\mytag{8.16}
$$
The equation \mythetag{8.16} cannot be derived directly since the procedure
of parallel translation for arbitrary tensors has no visual representation
like Fig\.~8.1.\par
     Let's consider a segment of a straight line given parametrically by 
the functions $u^1(t),\,u^2(t),\,u^3(t)$ in a curvilinear coordinates.
Let $t=s$ be the natural parameter on this straight line. Then the tangent
vector $\boldsymbol\tau(t)$ is a vector of the unit length at all points
of the line. Its direction is also unchanged. Therefore, its components
$\dot u^i$ satisfy the equation of parallel translation. Substituting
$a^i=\dot u^i$ into \mythetag{8.5}, we get
$$
\hskip -2em
\ddot u^i+\sum^3_{j=1}\sum^3_{k=1}
\Gamma^i_{jk}\,\dot u^j\,\dot u^k=0.
\mytag{8.17}
$$
The equation \mythetag{8.17} is the differential equation of a straight line
in curvilinear coordinates (written for the natural parametrization $t=s$).
\par
\head
\S~\mysection{9} Some calculations in polar,\\
cylindrical, and spherical coordinates.
\endhead
\rightheadtext{\S~9. Some calculations \dots}
     Let's consider the polar coordinate system on a plane. It is given
by formulas \mythetag{1.1}. Differentiating the expressions \mythetag{1.1}, 
we find the components of the frame vectors for the polar coordinate 
system:
$$
\xalignat 2
&\hskip -2em
\bold E_1=
\Vmatrix
\cos(\varphi)\\
\vspace{1ex}
\sin(\varphi)
\endVmatrix,
&&\bold E_2=
\Vmatrix
-\rho\sin(\varphi)\\
\vspace{1ex}
\rho\cos(\varphi)
\endVmatrix.
\mytag{9.1}
\endxalignat
$$
The column-vectors \mythetag{9.1} are composed by the coordinates 
of the vectors $\bold E_1$ and $\bold E_2$ in the orthonormal 
basis. Therefore, we can calculate their scalar products and
thus find the components of direct and inverse metric tensors
$\bold g$ and $\hat\bold g$:
$$
\xalignat 2
&\hskip -2em
g_{ij}=
\Vmatrix
1 &0\\
\vspace{1ex}
0 &\rho^2
\endVmatrix,
&
&g^{ij}=
\Vmatrix
\format\l&\quad\l\\
1 &0\\
\vspace{1ex}
0 &\rho^{-2}
\endVmatrix.
\mytag{9.2}
\endxalignat
$$
Once the components of $\bold g$ and $\hat\bold g$ are known,
we can calculate the Christoffel symbols. For this purpose
we apply the formula \mythetag{7.8}:
$$
\xalignat 3
&\hskip -2em
\Gamma^1_{11}=0,
&&\Gamma^1_{12}=\Gamma^1_{21}=0,
&&\Gamma^1_{22}=-\rho,\\
\vspace{-1ex}
&&&&&\mytag{9.3}\\
\vspace{-1ex}
&\hskip -2em
\Gamma^2_{11}=0,
&&\Gamma^2_{12}=\Gamma^2_{21}=\rho^{-1},
&&\Gamma^2_{22}=0.
\endxalignat
$$
Let's apply the connection components \mythetag{9.3} in order to calculate
the Laplace operator $\triangle$ in polar coordinates. Let $\psi$ be some
scalar field: $\psi=\psi(\rho,\varphi)$. Then
$$
\hskip -2em
\triangle\psi=\sum^2_{i=1}\sum^2_{j=1}
g^{ij}\left(\frac{\partial^2\psi}{\partial u^i\,
\partial u^j}-\shave{\sum^2_{k=1}}\Gamma^k_{ij}
\frac{\partial\psi}{\partial u^k}\right).
\mytag{9.4}
$$
The formula \mythetag{9.4} is a two-dimensional version of the formula
\mythetagchapter{10.15}{2} from
Chapter~\uppercase\expandafter{\romannumeral 2} applied to a scalar field.
Substituting \mythetag{9.3} into \mythetag{9.4}, we get
$$
\hskip -2em
\triangle\psi=\frac{\partial^2\psi}{\partial\rho^2}+
\frac{1}{\rho}\,\frac{\partial\psi}{\partial\rho}+
\frac{1}{\rho^2}\,\frac{\partial^2\psi}{\partial
\varphi^2}.
\mytag{9.5}
$$\par
     Now let's consider the cylindrical coordinate system. For the
components of metric tensors $\bold g$ and $\hat\bold g$ in this 
case we have
$$
\xalignat 2
&\hskip -2em
g_{ij}=
\Vmatrix
1 &0 &0\\
\vspace{1ex}
0 &\rho^2 &0\\
\vspace{1ex}
0 &0 &1
\endVmatrix,
&&g^{ij}=
\Vmatrix
\format\l&\quad\l&\quad\l\\
1 &0 &0\\
\vspace{1ex}
0 &\rho^{-2} &0\\
\vspace{1ex}
0 &0 &1
\endVmatrix.
\mytag{9.6}
\endxalignat
$$
From \mythetag{9.6} by means of \mythetag{7.8} we derive the connection
components:
$$
\xalignat 3
&\hskip -2em
\Gamma^1_{11}=0,
&&\Gamma^1_{12}=0,
&&\Gamma^1_{21}=0,\\
\vspace{0.5ex}
&\hskip -2em
\Gamma^1_{13}=0,
&&\Gamma^1_{31}=0,
&&\Gamma^1_{22}=-\rho,
\mytag{9.7}\\
\vspace{0.5ex}
&\hskip -2em
\Gamma^1_{23}=0,
&&\Gamma^1_{32}=0,
&&\Gamma^1_{33}=0,\\
&\vbox{\hsize=7cm\hrule width 7cm}\hskip -7cm &&&&\\
\vspace{1ex}
&\hskip -2em
\Gamma^2_{11}=0,
&&\Gamma^2_{12}=\rho^{-1},
&&\Gamma^2_{21}=\rho^{-1},\\
\vspace{0.5ex}
&\hskip -2em
\Gamma^2_{13}=0,
&&\Gamma^2_{31}=0,
&&\Gamma^2_{22}=0,
\mytag{9.8}\\
\vspace{0.5ex}
&\hskip -2em
\Gamma^2_{23}=0,
&&\Gamma^2_{32}=0,
&&\Gamma^2_{33}=0,\\
&\vbox{\hsize=7cm\hrule width 7cm}\hskip -7cm &&&&\\
\vspace{1ex}
&\hskip -2em
\Gamma^3_{11}=0,
&&\Gamma^3_{12}=0,
&&\Gamma^3_{21}=0,\\
\vspace{0.5ex}
&\hskip -2em
\Gamma^3_{13}=0,
&&\Gamma^3_{31}=0,
&&\Gamma^3_{22}=0,
\mytag{9.9}\\
\vspace{0.5ex}
&\hskip -2em
\Gamma^3_{23}=0,
&&\Gamma^3_{32}=0,
&&\Gamma^3_{33}=0.
\endxalignat
$$
Let's rewrite in the dimension $3$ the relationship \mythetag{9.4} for
the Laplace operator applied to a scalar field $\psi$:
$$
\hskip -2em
\triangle\psi=\sum^3_{i=1}\sum^3_{j=1}
g^{ij}\left(\frac{\partial^2\psi}{\partial u^i\,
\partial u^j}-\shave{\sum^3_{k=1}}\Gamma^k_{ij}
\frac{\partial\psi}{\partial u^k}\right).
\mytag{9.10}
$$
Substituting \mythetag{9.7}, \mythetag{9.8}, and \mythetag{9.9} into the
formula \mythetag{9.10}, we get
$$
\hskip -2em
\triangle\psi=\frac{\partial^2\psi}{\partial\rho^2}+
\frac{1}{\rho}\,\frac{\partial\psi}{\partial\rho}+
\frac{1}{\rho^2}\,\frac{\partial^2\psi}{\partial
\varphi^2}+\frac{\partial^2\psi}{\partial h^2}.
\mytag{9.11}
$$\par
Now we derive the formula for the components of rotor in cylindrical
coordina\-tes. Let $\bold A$ be a vector field and let 
$A^1,\,A^2,\,A^3$ be its components in cylindrical coordinates. In order
to calculate the components of the field $\bold F=\rot\bold A$ we use the
formula \mythetagchapter{10.5}{2} from Chapter
\uppercase\expandafter{\romannumeral 2}. This formula comprises the volume
tensor whose components are calculated by formula \mythetagchapter{8.1}{2}
from Chapter \uppercase\expandafter{\romannumeral 2}. The sign factor
$\xi_E$ in this formula is determined by the orientation of a coordinate
system. The cylindrical coordinate system can be either right-oriented or
left-oriented. It depends on the orientation of the auxiliary Cartesian
coordinate system $x^1,x^2,x^3$ which is related to the cylindrical
coordinates by means of the relationships \mythetag{1.3}. For the sake of
certainty we assume that the right-oriented cylindrical coordinates are
chosen. Then $\xi_E=1$ and for the components of the rotor $\bold
F=\rot\bold A$ we derive
$$
\hskip -2em
F^m=\sqrt{\det\bold g}\
\sum^3_{i=1}\sum^3_{j=1}\sum^3_{k=1}\sum^3_{q=1}
g^{mi}\,\varepsilon_{ijk}\,g^{jq}\,
\nabla_qF^k.
\mytag{9.12}
$$
Taking into account \mythetag{9.7}, \mythetag{9.8}, \mythetag{9.9}, \mythetag{9.6}
and using \mythetag{9.12}, we get
$$
\aligned
F^1&=\frac{1}{\rho}\,\frac{\partial A^3}{\partial\varphi}-
\rho\,\frac{\partial A^2}{\partial h},\\
\vspace{1ex}
F^2&=\frac{1}{\rho}\,\frac{\partial A^1}{\partial h}-
\frac{1}{\rho}\,\frac{\partial A^3}{\partial\rho},\\
\vspace{1ex}
F^3&=\rho\,\frac{\partial A^2}{\partial\rho}-
     \frac{1}{\rho}\,\frac{\partial A^1}{\partial\varphi}+
     2\,A^2.
\endaligned
\mytag{9.13}
$$
The relationships \mythetag{9.13} can be written in form of the
determinant:
$$
\hskip -2em
\rot\bold A=
\frac{1}{\rho}\,\left|\matrix
\bold E_1 &\bold E_2 &\bold E_3\\
\vspace{2ex}
\dsize\frac{\partial}{\partial\rho}
&\dsize\frac{\partial}{\partial\varphi}
&\dsize\frac{\partial}{\partial h}\\
\vspace{2ex}
A^1 &\rho^2A^2 &A^3
\endmatrix\right|,
\mytag{9.14}
$$
Here $\bold E_1,\,\bold E_2,\,\bold E_3$ are the frame vectors of the
cylindrical coordinates.\par
     In the case of spherical coordinates, we begin the calculations by
deriving the formula for the components of the metric tensor $\bold g$:
$$
\hskip -2em
g_{ij}=
\Vmatrix
\format
\c\qquad\quad &\c\quad &\c\\
1 &0 &0\\
\vspace{2ex}
0 &\rho^2 &0\\
\vspace{2ex}
0 &0 &\rho^2\sin^2(\vartheta)
\endVmatrix.
\mytag{9.15}
$$
Then we calculate the connection components and write then in form of
the array:
$$
\xalignat 3
&\hskip -2em
\Gamma^1_{11}=0,
&&\Gamma^1_{12}=0,
&&\Gamma^1_{21}=0,\\
\vspace{0.5ex}
&\hskip -2em
\Gamma^1_{13}=0,
&&\Gamma^1_{31}=0,
&&\Gamma^1_{22}=-\rho,\quad\qquad
\mytag{9.16}\\
\vspace{0.5ex}
&\hskip -2em
\Gamma^1_{23}=0,
&&\Gamma^1_{32}=0,
&&\Gamma^1_{33}=-\rho\sin^2(\vartheta),\hskip -6em\\
\displaybreak
&\hskip -2em
\Gamma^2_{11}=0,
&&\Gamma^2_{12}=\rho^{-1},
&&\Gamma^2_{21}=\rho^{-1},\\
\vspace{0.5ex}
&\hskip -2em
\Gamma^2_{13}=0,
&&\Gamma^2_{31}=0,
&&\Gamma^2_{22}=0,\quad\qquad
\mytag{9.17}\\
\vspace{0.5ex}
&\hskip -2em
\Gamma^2_{23}=0,
&&\Gamma^2_{32}=0,
&&\Gamma^2_{33}=-\frac{\sin(2\vartheta)}{2},\hskip -6em\\
&\hskip -2em
\vbox{\hsize=7cm\hrule width 8cm}\hskip -8cm &&&&\\
\vspace{1ex}
&\hskip -2em
\Gamma^3_{11}=0,
&&\Gamma^3_{12}=0,
&&\Gamma^3_{21}=0,\\
\vspace{0.5ex}
&\hskip -2em
\Gamma^3_{13}=\rho^{-1},
&&\Gamma^3_{31}=\rho^{-1},
&&\Gamma^3_{22}=0,\quad\qquad
\mytag{9.18}\\
\vspace{0.5ex}
&\hskip -2em
\Gamma^3_{23}=\ctg(\vartheta),\hskip -6em
&&\Gamma^3_{32}=\ctg(\vartheta),\hskip -6em
&&\Gamma^3_{33}=0.
\endxalignat
$$
Substituting \mythetag{9.16}, \mythetag{9.17}, and \mythetag{9.18}
into the relationship \mythetag{9.10}, we get
$$
\hskip -2em
\triangle\psi=\frac{\partial^2\psi}{\partial\rho^2}+
\frac{2}{\rho}\,\frac{\partial\psi}{\partial\rho}+
\frac{1}{\rho^2}\,\frac{\partial^2\psi}{\partial
\vartheta^2}+\frac{\ctg(\vartheta)}{\rho^2}\,\frac{\partial\psi}
{\partial\vartheta}+\frac{1}{\rho^2\sin^2(\vartheta)}\,
\frac{\partial^2\psi}{\partial\varphi^2}.
\mytag{9.19}
$$\par
     Let $\bold A$ be a vector field with the components 
$A^1,\,A^2,\,A^3$ in the right-oriented spherical coordinates.
Denote $\bold F=\rot\bold A$. Then for the components of
$\bold F$ we get
$$
\hskip -2em
\aligned
F^1&=\sin(\vartheta)\,\frac{\partial A^3}{\partial\vartheta}-
\frac{1}{\sin(\vartheta)}\,\frac{\partial A^2}{\partial\varphi}+
2\,\cos(\vartheta)\,A^3,\\
\vspace{1ex}
F^2&=\frac{1}{\rho^2\sin(\vartheta)}\,
\frac{\partial A^1}{\partial\varphi}-
\sin(\vartheta)\,\frac{\partial A^3}{\partial\rho}-
\frac{2\sin(\vartheta)}{\rho}\,A^3,\\
\vspace{1ex}
F^3&=\frac{1}{\sin(\vartheta)}\,\frac{\partial A^2}{\partial\rho}-
     \frac{1}{\rho^2\sin(\vartheta)}\,\frac{\partial A^1}
     {\partial\vartheta}+\frac{2}{\rho\sin(\vartheta)}\,A^2.
\endaligned
\mytag{9.20}
$$
Like \mythetag{9.13}, the formulas \mythetag{9.20} can be written in form 
of the determinant:
$$
\hskip -2em
\rot\bold A=
\frac{\rho^{-2}}{\sin(\vartheta)}\,\left|\matrix
\format \c\qquad\quad &\c\quad &\c\\
\bold E_1 &\bold E_2 &\bold E_3\\
\vspace{2ex}
\dsize\frac{\partial}{\partial\rho}
&\dsize\frac{\partial}{\partial\vartheta}
&\dsize\frac{\partial}{\partial\varphi}\\
\vspace{2ex}
A^1 &\rho^2A^2 &\rho^2\sin^2(\vartheta)A^3
\endmatrix\right|.
\mytag{9.21}
$$
The formulas \mythetag{9.5}, \mythetag{9.11}, and \mythetag{9.19} for
the Laplace operator and the formulas \mythetag{9.14} and \mythetag{9.21}
for the rotor is the main goal of the calculations performed just
above in this section. They are often used in applications and can be
found in some reference books for engineering computations.\par
     The matrices $\bold g$ in all of the above coordinate systems are
diagonal. Such coordinate systems are called {\it orthogonal}, while
the quantities $H_i=\sqrt{g_{ii}}$ are called the {\it Lame coefficients}
of orthogonal coordinates. Note that there is no orthonormal curvilinear
coordinate system. All such systems are necessarily Cartesian, this fact
follows from \mythetag{7.8} and \mythetag{5.9}.\par
\newpage
\topmatter
\title\chapter{4}
Geometry of surfaces.
\endtitle
\endtopmatter
\chapternum=4
\document
\head
\S~\mysection{1} Parametric surfaces.\\
Curvilinear coordinates on a surface.
\endhead
\rightheadtext{\S~1. Parametric surfaces.}
\setfirstpage
\leftheadtext{Chapter \uppercase\expandafter{\romannumeral 4}.
GEOMETRY OF SURFACES.}
     A surface is a two-dimensional spatially extended geometric object.
There are several ways for expressing quantitatively (mathematically) 
this fact of two-dimensionality of surfaces. In the three-dimensional
Euclidean space $\Bbb E$ the choice of an arbitrary point implies three
degrees of freedom: a point is determined by three coordinates. In order
to decrease this extent of arbitrariness we can bind three coordinates 
of a point by an equation:
$$
\hskip -2em
F(x^1,x^2,x^3)=0.
\mytag{1.1}
$$
Then the choice of two coordinates determines the third coordinate of a
point. This means that we can define a surface by means of an equation
in some coordinate system (for the sake of simplicity we can choose a
Cartesian coordinate system). We have already used this method of defining
surfaces (see formula \mythetagchapter{1.2}{1} in
Chapter \uppercase\expandafter{\romannumeral 1}) when we defined a curve as
an intersection of two surfaces.\par
     Another way of defining a surface is the {\it parametric method}. Unlike curves, surfaces are parameterized by two parameters. Let's
denote them $u^1$ and $u^2$:
$$
\hskip -2em
\bold r=\bold r(u^1,u^2)=\Vmatrix x^1(u^1,u^2)\\
\vspace{2.5ex}
x^2(u^1,u^2)\\
\vspace{2.5ex}
x^3(u^1,u^2)\endVmatrix.
\mytag{1.2}
$$
The formula \mythetag{1.2} expresses the radius-vector of the points of
a surface in some Cartesian coordinate system as a function of two
parameters $u^1,\,u^2$. Usually, only a part of a surface is represented
in parametric form. Therefore, considering the pair of numbers $(u^1,\,
u^2)$ as a point of $\Bbb R^2$, we can assume that the point $(u^1,\,
u^2)$ runs over some domain $U\subset\Bbb R^2$. Let's denote by $D$ the
image of the domain $U$ under the mapping \mythetag{1.2}. Then $D$ is
the {\it domain being mapped}, \ $U$ is the {\it map} or the {\it chart},
and \mythetag{1.2} is the {\it chart mapping}: it maps $U$ onto $D$.\par
     The smoothness class of the surface $D$ is determined by the smoothness class of the functions $x^1(u^1,u^2)$, $x^2(u^1,u^2)$, and $x^3(u^1,u^2)$ in formula \mythetag{1.2}. In what fallows we shall consider 
only those surfaces for which these functions are at least continuously
differentiable. Then, differentiating these functions, we \pagebreak 
can arrange their derivatives into the Jacobi matrix:
$$
\hskip -2em
I=\left\Vert
\vphantom{\vrule height 39pt depth 39pt}
\matrix
\dsize\frac{\partial x^1}{\partial u^1}&
\dsize\frac{\partial x^1}{\partial u^2}\\
\vspace{1.2ex}
\dsize\frac{\partial x^2}{\partial u^1}&
\dsize\frac{\partial x^2}{\partial u^2}\\
\vspace{1.2ex}
\dsize\frac{\partial x^3}{\partial u^1}&
\dsize\frac{\partial x^3}{\partial u^2}
\endmatrix\right\Vert.
\mytag{1.3}
$$
\mydefinition{1.1} \ A continuously differentiable mapping
\mythetag{1.2} is called {\it regular\/} at a point $(u^1,u^2)$ if the
rank of the Jacobi matrix \mythetag{1.3} calculated at that point\linebreak
is equal to $2$.
\enddefinition
\mydefinition{1.2} A set $D$ is called a {\it regular fragment\/}
of a continuously differentiable surface if there is a mapping 
$\bold u\!:\,D\to U$ from $D$ to some domain $U\subset\Bbb R^2$ and the
following conditions are fulfilled:
\roster
\rosteritemwd=5pt
\item the mapping $\bold u: D\to U$ is bijective;
\item the inverse mapping $\bold u^{-1}\!:\,U\to D$ given by three
continuously differentiable functions \mythetag{1.2} is regular at 
all points of the domain $U$.
\endroster
\enddefinition
     The Jacobi matrix \mythetag{1.3} has three minors of the order $2$.
These are the determinants of the following $2\times 2$ matrices:
$$
\xalignat 3
&\hskip -2em\left\vert
\vphantom{\vrule height 24pt depth 24pt}
\matrix
\dsize\frac{\partial x^1}{\partial u^1}&
\dsize\frac{\partial x^1}{\partial u^2}\\
\vspace{1.2ex}
\dsize\frac{\partial x^2}{\partial u^1}&
\dsize\frac{\partial x^2}{\partial u^2}
\endmatrix\right\vert,
&&\left\vert
\vphantom{\vrule height 24pt depth 24pt}
\matrix
\dsize\frac{\partial x^2}{\partial u^1}&
\dsize\frac{\partial x^2}{\partial u^2}\\
\vspace{1.2ex}
\dsize\frac{\partial x^3}{\partial u^1}&
\dsize\frac{\partial x^3}{\partial u^2}
\endmatrix\right\vert,
&&\left\vert
\vphantom{\vrule height 24pt depth 24pt}
\matrix
\dsize\frac{\partial x^3}{\partial u^1}&
\dsize\frac{\partial x^3}{\partial u^2}\\
\vspace{1.2ex}
\dsize\frac{\partial x^1}{\partial u^1}&
\dsize\frac{\partial x^1}{\partial u^2}
\endmatrix\right\vert.
\mytag{1.4}
\endxalignat
$$
In the case of regularity of the mapping \mythetag{1.2} at least one
of the determinants \mythetag{1.4} is nonzero. At the expense of renaming
the variables $x^1,\,x^2,\,x^3$ we always can do so that the first
determinant will be nonzero:
$$
\hskip -2em
\left\vert
\vphantom{\vrule height 24pt depth 24pt}
\matrix
\dsize\frac{\partial x^1}{\partial u^1}&
\dsize\frac{\partial x^1}{\partial u^2}\\
\vspace{1.2ex}
\dsize\frac{\partial x^2}{\partial u^1}&
\dsize\frac{\partial x^2}{\partial u^2}
\endmatrix\right\vert\neq 0.
\mytag{1.5}
$$
In this case we consider the first two functions $x^1(u^2,u^2)$ 
and  $x^2(u^2,u^2)$ in \mythetag{1.2} as a mapping and write them
as follows:
$$
\hskip -2em
\cases
x^1=x^1(u^1,u^2),\\
x^2=x^2(u^1,u^2).
\endcases
\mytag{1.6}
$$
Due to \mythetag{1.5} the mapping \mythetag{1.6} is locally invertible. 
Upon restricting \mythetag{1.6} to some sufficiently small neighborhood
of an arbitrary preliminarily chosen point one can construct two
continuously differentiable functions
$$
\hskip -2em
\cases
u^1=u^1(x^1,x^2),\\
u^2=u^2(x^1,x^2)
\endcases
\mytag{1.7}
$$
that implement the inverse mapping for \mythetag{1.6}. This fact is
well-known, it is a version of the theorem on implicit functions
(see \mycite{2}, see also the theorem~\mythetheoremchapter{2.1}{3} in 
Chapter \uppercase\expandafter{\romannumeral 3}). Let's substitute
$u^1$ and $u^2$ from \mythetag{1.7} into the arguments of the third
function $x^3(u^1,u^2)$ in the formula \mythetag{1.2}. As a result we 
obtain the function $F(x^1,x^2)=x^3(u^1(x^2,x^2),u^2(x^2,x^2))$ such 
that each regular fragment of a surface can locally (i\.\,e\. in some 
neighborhood of each its point) be presented as a graph of a 
continuously differentiable function of two variables:
$$
\hskip -2em
x^3=F(x^1,x^2).
\mytag{1.8}
$$
\subhead A remark on singular points\endsubhead If we give up the
regularity condition from the definition~\mythedefinition{1.2}, this may
cause the appearance of singular points on a surface. As an example we
consider two surfaces given by smooth functions:
$$
\xalignat 2
&\hskip -2em
\cases
x^1=(u^1)^3,\\
x^2=(u^2)^3,\\
x^3=(u^1)^2+(u^2)^2,
\endcases
&&\cases
x^1=(u^1)^3,\\
x^2=(u^2)^3,\\
x^3=(u^1)^4+(u^2)^4.
\endcases
\mytag{1.9}
\endxalignat
$$
In both cases the regularity condition breaks at the point $u^1=u^2=0$. 
As a result the first surface \mythetag{1.9} gains the singularity at the
origin. The second surface is non-singular despite to the breakage of the
regularity condition.\par
     Marking a regular fragment $D$ on a surface and defining a chart
mapping $\bold u^{-1}\!:\,U\to D$ can be treated
\vadjust{\vskip 5pt\hbox to 0pt{\kern -15pt
\includegraphics{ris12.eps}\hss}\vskip 210pt}as introducing
a curvilinear coordinate system on the surface. The conditions $u^1=\const$
and $u^2=\const$ determine two families of coordinate lines on the plane
of parameters $u^1,\,u^2$. They form the coordinate network in $U$. The
mapping \mythetag{1.2} maps it onto the coordinate network on the surface 
$D$ (see Fig\.~1.1 and Fig\.~1.2). Let's consider the vectors $\bold E_1$
and $\bold E_2$ tangent to the lines of the coordinate network on the surface
$D$:
$$
\hskip -2em
\bold E_i(u^1,u^2)=\frac{\partial\bold r(u^1,u^2)}
{\partial u^i}.
\mytag{1.10}
$$
The formula \mythetag{1.10} defines a pair of tangent vectors 
$\bold E_1$ and $\bold E_2$ attached to each point of the 
surface $D$.\par
     The vector-function $\bold r(u^1,u^2)$ which defines the mapping 
\mythetag{1.2} can be written in form of the expansion in the basis of 
the auxiliary Cartesian coordinate system:
$$
\hskip -2em
\bold r(u^1,u^2)=\sum^3_{q=1} x^q(u^1,u^2)\cdot\bold e_q.
\mytag{1.11}
$$
Substituting the expansion \mythetag{1.11} into \mythetag{1.10} we can 
express the tangent vectors $\bold E_1$ and $\bold E_2$ through the 
basis vectors $\bold e_1,\,\bold e_2,\,\bold e_3$:
$$
\hskip -2em
\bold E_i(u^1,u^2)=\sum^3_{q=1}\frac{\partial x^q(u^1,u^2)}
{\partial u^i}\cdot\bold e_q.
\mytag{1.12}
$$
Let's consider the column-vectors composed by the Cartesian coordinates
of the tangent vectors $\bold E_1$ and $\bold E_2$:
$$
\xalignat 2
&\hskip -2em
\bold E_1=\left\Vert
\vphantom{\vrule height 39pt depth 39pt}
\matrix
\dsize\frac{\partial x^1}{\partial u^1}\\
\vspace{1.2ex}
\dsize\frac{\partial x^2}{\partial u^1}\\
\vspace{1.2ex}
\dsize\frac{\partial x^3}{\partial u^1}
\endmatrix\right\Vert,
&&\bold E_2=\left\Vert
\vphantom{\vrule height 39pt depth 39pt}
\matrix
\dsize\frac{\partial x^1}{\partial u^2}\\
\vspace{1.2ex}
\dsize\frac{\partial x^2}{\partial u^2}\\
\vspace{1.2ex}
\dsize\frac{\partial x^3}{\partial u^2}
\endmatrix\right\Vert.
\mytag{1.13}
\endxalignat
$$
Note that the column-vectors \mythetag{1.13} coincide with the 
columns in the Jacobi matrix \mythetag{1.3}. However, from the
regularity condition (see the definition~\mythedefinition{1.1}) 
it follows that the column of the Jacobi matrix \mythetag{1.3} 
are linearly independent. This consideration proves the following
proposition.
\mytheorem{1.1} The tangent vectors $\bold E_1$ and 
$\bold E_2$ are linearly independent at each point of a surface. 
Therefore, they form the frame of the tangent vector fields in
\nolinebreak $D$.
\endproclaim
     The frame vectors $\bold E_1$ and $\bold E_2$ attached to some
point of a surface $D$ define the tangent plane at this point. Any
vector tangent to the surface at this point lies in the tangent plane,
it can be expanded in the basis formed by the vectors $\bold E_1$ and
$\bold E_2$. Let's consider some arbitrary curve $\gamma$ lying 
completely on the surface (see Fig\.~1.1 and Fig\.~1.2). In parametric 
form such a curve is given by two functions of a parameter $t$. They
define the curve as follows:
$$
\hskip -2em
\cases
u^1=u^1(t),\\
u^2=u^2(t).
\endcases
\mytag{1.14}
$$
By substituting \mythetag{1.14} into \mythetag{1.11} or into \mythetag{1.2}
we find the radius-vector of a point of the curve in the auxiliary
Cartesian coordinate system $\bold r(t)=\bold r(u^1(t),u^2(t))$. Let's
differentiate $\bold r(t)$ with respect to $t$ and find the tangent 
vector of the curve given by the above two functions \mythetag{1.14}:
$$
\boldsymbol\tau(t)=\frac{d\bold r}{dt}=
\sum^2_{i=1}\frac{\partial\bold r}{\partial u^i}\cdot
\frac{du^i}{dt}.
$$
Comparing this expression with \mythetag{1.10}, we find that $\boldsymbol
\tau$ is expressed as follows:
$$
\hskip -2em
\boldsymbol\tau(t)=\sum^2_{i=1}\dot u^i\cdot\bold E_i.
\mytag{1.15}
$$
Due to \mythetag{1.15} the vector $\boldsymbol\tau$ is a linear combination
of the vectors $\bold E_1$ and $\bold E_2$ forming the tangent frame.
Hence, if a curve $\gamma$ lies completely on the surface, its tangent
vector lies on the tangent plane to this surface, while the derivatives
of the functions \mythetag{1.14} are the components of the vector
$\boldsymbol\tau$ expanded in the basis of the frame vectors $\bold E_1$
and $\bold E_2$.\par
\head
\S~\mysection{2} Change of curvilinear coordinates on a surface.
\endhead
\rightheadtext{\S~2. Change of curvilinear coordinates on a surface.}
     Let's consider two regular fragments $D_1$ and $D_2$ on some surface,
each equipped with with its own curvilinear coordinate system. Assume that
their intersection $D=D_1\cap D_2$ is not empty. Then in $D$ we have two
curvilinear coordinate systems $u^1,\,u^2$ and $\tilde u^1,\,\tilde u^2$.
Denote by $U$ and $\tilde U$ the preimages of $D$ under the corresponding
chart mappings (see Fig\.~3.1 in 
Chapter \uppercase\expandafter{\romannumeral 3}). Due to the bijectivity
of the chart mappings (see definition~\mythedefinition{1.2}) we can
construct two mappings
$$
\xalignat 2
&\hskip -2em
\tilde\bold u\compos\bold u^{-1}\!:\,U\to\tilde U,
&&\bold u\compos\tilde\bold u^{-1}\!:\,\tilde U\to U.
\mytag{2.1}
\endxalignat
$$
The mappings $\tilde\bold u\compos\bold u^{-1}$
and $\bold u\compos\tilde\bold u^{-1}$ in \mythetag{2.1} are also
bijective, they can be represented by the following pairs of functions:
$$
\xalignat 2
&\hskip -2em
\cases
\tilde u^1=\tilde u^1(u^1,u^2),\\
\tilde u^2=\tilde u^2(u^1,u^2).
\endcases
&&\cases
u^1=u^1(\tilde u^1,\tilde u^2),\\
u^2=u^2(\tilde u^1,\tilde u^2).
\endcases
\mytag{2.2}
\endxalignat
$$
\mytheorem{2.1} The functions \mythetag{2.2} representing
$\tilde\bold u\compos\bold u^{-1}$ and $\bold u\compos\tilde\bold u^{-1}$
are conti\-nuously differentiable.
\endproclaim
\demo{Proof} We shall prove the continuous differentiability of the
second pair of functions \mythetag{2.2}. For the first pair the proof 
is analogous. Let's choose some point on the chart $U$ and map it to
$D$. Then we choose a suitable Cartesian coordinate system in $\Bbb E$
such that the condition \mythetag{1.5} is fulfilled and in some neighborhood
of the mapped point there exists the mapping \mythetag{1.7} inverse for
\mythetag{1.6}. The mapping \mythetag{1.7} is continuously differentiable.
\par
     The other curvilinear coordinate system in $D$ induces the other
pair of functions that plays the same role as the functions \mythetag{1.6}:
$$
\hskip -2em
\cases
x^1=x^1(\tilde u^1,\tilde u^2),\\
x^2=x^2(\tilde u^1,\tilde u^2).
\endcases
\mytag{2.3}
$$
These are two of three functions that determine the mapping 
$\tilde\bold u^{-1}$ in form of \mythetag{1.2}. The functions  $u^1=u^1(\tilde u^1,\tilde u^2)$ and $u^2=u^2(\tilde u^1,\tilde u^2)$
that determine the mapping $\bold u\compos\tilde\bold u^{-1}$ in
\mythetag{2.2} are obtained by substituting \mythetag{2.3} into the 
arguments of \mythetag{1.7}:
$$
\pagebreak 
\hskip -2em
\cases
u^1=u^1(x^1(\tilde u^1,\tilde u^2),x^2(\tilde u^1,\tilde u^2)),\\
u^2=u^2(x^1(\tilde u^1,\tilde u^2),x^2(\tilde u^1,\tilde u^2)).
\endcases
\mytag{2.4}
$$
The compositions of continuously differentiable functions in \mythetag{2.4}
are continuously differentiable functions. This fact completes the proof
of the theorem.\qed\enddemo
      The functions \mythetag{2.2}, whose continuous differentiability was
proved just above, perform  the {\it transformation} or the {\it change} 
of curvilinear coordinates on a surface. They are analogous to the 
functions \mythetagchapter{3.5}{3}
in Chapter \uppercase\expandafter{\romannumeral 3}.
\par
\subhead A remark on the smoothness\endsubhead If the functions 
\mythetag{1.2} of both coordinate systems $u^1,u^2$ and $\tilde u^1,
\tilde u^2$ belong to the smoothness class $C^m$, then the
transition functions \mythetag{2.2} also belong to the smoothness 
class $C^m$.\par
     Let $\bold r(u^1,u^2)$ and $\bold r(\tilde u^1,\tilde u^2)$
be two vector-functions of the form \mythetag{1.2} for two curvilinear
coordinate systems in $D$. They define the mappings $\bold u^{-1}$ and
$\tilde\bold u^{-1}$ acting from the charts $U$ and $\tilde U$ to
$D$. Due to the identity $\tilde\bold u^{-1}=\bold u^{-1}\compos
(\bold u\compos\tilde\bold u^{-1})$ the function $\bold r(\tilde u^1,\tilde
u^2)$ is obtained by substituting the corresponding transition functions
\mythetag{2.2} into the arguments of $\bold r(u^1,u^2)$:
$$
\hskip -2em
\bold r(\tilde u^1,\tilde u^2)=\bold r(u^1(\tilde u^1,\tilde u^2),
u^2(\tilde u^1,\tilde u^2)).
\mytag{2.5}
$$
Let's differentiate \mythetag{2.5} with respect to $\tilde u^j$ and take
into account the chain rule and the formula \mythetag{1.10} for the vectors
of the tangent frame:
$$
\tilde\bold E_j=\frac{\partial\bold r}{\partial\tilde u^j}=
\sum^2_{i=1}\frac{\partial\bold r}{\partial u^i}\cdot
\frac{\partial u^i}{\partial\tilde u^j}=
\sum^2_{i=1}\frac{\partial u^i}{\partial\tilde u^j}
\cdot\bold E_i.
$$
Differentiating the identity $\bold r(u^1,u^2)=\bold r(\tilde u^1(u^1,
u^2),\tilde u^2(u^1,u^2))$, we derive the analogous relationship inverse
to the previous one:
$$
\bold E_i=\frac{\partial\bold r}{\partial u^i}=
\sum^2_{k=1}\frac{\partial\bold r}{\partial\tilde u^k}\cdot
\frac{\partial\tilde u^k}{\partial u^i}=
\sum^2_{i=1}\frac{\partial\tilde u^k}{\partial u^i}
\cdot\tilde\bold E_k.
$$
It is clear that the above relationships describe the direct and inverse
transitions from some tangent frame to another. Let's write them as
$$
\xalignat 2
&\hskip -2em
\tilde\bold E_j=\sum^2_{i=1} S^i_j\cdot\bold E_i,
&&\bold E_i=\sum^2_{k=1} T^k_i\cdot\tilde\bold E_k,
\mytag{2.6}
\endxalignat
$$
where the components of the matrices $S$ and $T$ are given by the
formulas
$$
\xalignat 2
&\hskip -2em
S^i_j(\tilde u^1,\tilde u^2)=\frac{\partial u^i}
{\partial\tilde u^j},
&&T^k_i(u^1,u^2)=\frac{\partial\tilde u^k}{\partial u^i}.
\mytag{2.7}
\endxalignat
$$
From \mythetag{2.7}, we see that the transition matrices $S$ and $T$
are the Jacobi matrices for the mappings given by the transition
functions \mythetag{2.2}. They are non-degenerate and are inverse to
each other.\par
     The transformations \mythetag{2.2} and the transition matrices
$S$ and $T$ related to them are used in order to construct the theory
of tensors and tensor fields analogous to that which we considered in
Chapter \uppercase\expandafter{\romannumeral 2} and
Chapter \uppercase\expandafter{\romannumeral 3}. Tensors and tensor 
fields defined through the transformations \mythetag{2.2} and transition
matrices \mythetag{2.7} are called {\it inner} tensors and {\it inner}
tensor fields on a surface:
$$
F^{i_1\ldots\,i_r}_{j_1\ldots\,j_s}=
\sum\Sb p_1\ldots\,p_r\\ q_1\ldots\,q_s\endSb
S^{i_1}_{p_1}\ldots\,S^{i_r}_{p_r}\,\,
T^{q_1}_{j_1}\ldots\,T^{q_s}_{j_s}\,\,
\tilde F^{p_1\ldots\,p_r}_{q_1\ldots\,q_s}.
\mytag{2.8}
$$
\mydefinition{2.1} An {\it inner tensor\/} of the type $(r,s)$
on a surface is a geometric object $\bold F$ whose components in an
arbitrary curvilinear coordinate system on that surface are enumerated
by $(r+s)$ indices and under a change of coordinate system are transformed
according to the rule \mythetag{2.8}.
\enddefinition
     The formula \mythetag{2.8} differs from the formula
\mythetagchapter{1.6}{2} in Chapter \uppercase\expandafter{\romannumeral 2}
only in the range of indices. Each index here runs over the range of two
values 1 and 2. By setting the sign factor $(-1)^S=\sign(\det S)=\pm 1$
into the formula \mythetag{2.8} we get the definition of an inner
pseudotensor 
$$
\hskip -2em
F^{i_1\ldots\,i_r}_{j_1\ldots\,j_s}=
   \sum\Sb p_1\ldots\,p_r\\ q_1\ldots\,q_s\endSb
   (-1)^S\,S^{i_1}_{p_1}\ldots\,S^{i_r}_{p_r}\,
   T^{q_1}_{j_1}\ldots\,T^{q_s}_{j_s}\,\,
   \tilde F^{p_1\ldots\,p_r}_{q_1\ldots\,
   q_s}.
\mytag{2.9}
$$
\mydefinition{2.3} An {\it inner pseudotensor\/} of the type
$(r,s)$ on a surface is a geometric object $\bold F$ whose components 
in an arbitrary curvilinear coordinate system on that surface are 
enumerated by $(r+s)$ indices and under a change of coordinate system 
are transformed according to the rule \mythetag{2.9}.
\enddefinition
    Inner tensorial and pseudotensorial fields are obtained by defining
an inner tensor or pseudotensor at each point of a surface. The operations
of addition, tensor product, contraction, transposition of indices,
symmetrization and alter\-nation for such fields are defined in a way similar to that of the case of the fields in the space $\Bbb E$ (see 
Chapter \uppercase\expandafter{\romannumeral 2}). All properties of these
operations are preserved.\par
\subhead A remark on the differentiation\endsubhead The operation of
covariant differentiation of tensor fields in the space $\Bbb E$ was
first introduced for Cartesian coordinate systems. Then it was extended
to the case of curvilinear coordinates. On surfaces, as a rule, there is no
Cartesian coordinate system at all. Therefore, the operation of covariant
differentiation for inner tensor fields on a surface should be defined 
in a somewhat different way.\par
\head
\S~\mysection{3} The metric tensor and the area tensor.
\endhead
\rightheadtext{\S~3. The metric tensor and the area tensor.}
      The choice of parameters $u^1,\,u^2$ on a surface determines
the tangent frame $\bold E_1,\,\bold E_2$ on that surface. Let's
consider the scalar products of the vectors $\bold E_1,\,\bold E_2$
forming the tangent frame of the surface:
$$
\hskip -2em
g_{ij}=(\bold E_i\,|\,\bold E_j).
\mytag{3.1}
$$
They compose the $2\times 2$ Gram matrix $\bold g$ which is symmetric,
non-degenerate, and positive. Therefore, we have the inequality
$$
\hskip -2em
\det\bold g>0.
\mytag{3.2}
$$
Substituting \mythetag{2.6} into \mythetag{3.1}, we find that under a change of a coordinate system the quantities \mythetag{3.1} are transformed
as the components of an inner tensorial field of the type $(0,2)$. The
tensor $\bold g$ with the components \mythetag{3.1} is called the {\it metric tensor} of the surface. Note that the components of the metric
tensor are determined by means of the scalar product in the outer space $\Bbb E$. Therefore, we say that the tensor field $\bold g$ {\it is
induced\/} by the outer scalar product. For this reason the tensor $\bold
g$ is called the metric tensor of the {\it induced metric}.\par
     Symmetric tensors of the type $(0,2)$ are related to quadratic forms.
This fact yields another title for the tensor $\bold g$. It is called
the {\it first quadratic form} of a surface. Sometimes, for the components
of the first quadratic form the special notations are used: $g_{11}=E$,
$g_{12}=g_{21}=F$, $g_{22}=G$. These notations are especially popular in
the earlier publications on the differential geometry:
$$
\hskip -2em
g_{ij}=\left\Vert
\matrix
E &F\\ F &G
\endmatrix
\right\Vert
\mytag{3.3}
$$\par
     Since the Gram matrix $\bold g$ is non-degenerate, we can define the
inverse matrix $\hat\bold g=\bold g^{-1}$. The components of such inverse matrix are denoted by $g^{ij}$, setting the indices $i$ and $j$ to the
upper position:
$$
\hskip -2em
\sum^3_{j=1} g^{ij}\,g_{jk}=\delta^i_j.
\mytag{3.4}
$$
For the matrix $\hat\bold g$ the proposition analogous to the 
theorem~\mythetheoremchapter{6.1}{2} from Chapter
\uppercase\expandafter{\romannumeral 2} is valid. The
components of this matrix define an inner tensor field of the type 
$(2,0)$ on a surface, this field is called the {\it inverse metric
tensor} or the {\it dual metric tensor}. The proof of this proposition
is completely analogous to the proof of the
theorem~\mythetheoremchapter{6.1}{2} in Chapter
\uppercase\expandafter{\romannumeral 2}. Therefore, here we do not give
this proof.\par
     From the symmetry of the matrix $\bold g$ and from the relationships
\mythetag{3.4} it follows that the components of the inverse matrix 
$\hat\bold g$ are symmetric. The direct and inverse metric tensors are 
used in order to lower and raise indices of tensor fields. These operations
are defined by the formulas analogous to \mythetagchapter{9.1}{2} and
\mythetagchapter{9.2}{2} in Chapter \uppercase\expandafter{\romannumeral 2}:
$$
\hskip -2em
\aligned
B^{i_1\ldots\,i_{r-1}}_{j_1\ldots\,j_{s+1}}=
\sum^2_{k=1}
A^{i_1\ldots\,i_{m-1}\,k\,i_m\ldots\,i_{r-1}}_{j_1\ldots\,
j_{n-1}\,j_{n+1}\ldots\,j_{s+1}}\,\,
g_{kj_n},\\
A^{i_1\ldots\,i_{r+1}}_{j_1\ldots\,j_{s-1}}=
\sum^2_{k=1}
B^{i_1\ldots\,i_{m-1}\,i_{m+1}\ldots\,i_{r+1}}_{j_1\ldots\,
j_{n-1}\,q\,j_n\ldots\,j_{s-1}}\,\,g^{qi_m}.
\endaligned
\mytag{3.5}
$$
The only difference of the formulas \mythetag{3.5} here is that the
summation index $k$ runs over the range of two numbers 1 and 2. 
Due to \mythetag{3.4} the operations of raising and lowering indices 
\mythetag{3.5} are inverse to each other.\par
     In order to define the area tensor (or the area pseudotensor)
we need the following skew-symmetric $2\times 2$ matrix:
$$
\hskip -2em
d_{ij}=d^{\,ij}=\left\Vert
\matrix
\format \r\quad &\l\\
0 &1\\ -1 &0
\endmatrix
\right\Vert.
\mytag{3.6}
$$
The quantities \mythetag{3.6} form the two-dimensional analog of the
Levi-Civita symbol (see formula \mythetagchapter{6.8}{2} in Chapter
\uppercase\expandafter{\romannumeral 2}). These quantities satisfy the
relationship 
$$
\hskip -2em
\sum^2_{p=1}\sum^2_{q=1}d_{pq}\,M_{ip}\,M_{jq}=
\det\bold M\,\,d_{ij},
\mytag{3.7}
$$
where $\bold M$ is some arbitrary square $2\times 2$ matrix. The formula
\mythetag{3.7} is an analog of the formula \mythetagchapter{6.10}{2}
from Chapter \uppercase\expandafter{\romannumeral 2} (see proof in
\mycite{4}).\par
     Using the quantities $d_{ij}$ and the matrix of the metric tensor
$\bold g$ in some curvilinear coordinate system, we construct the following
quantities:
$$
\hskip -2em
\omega_{ij}=\sqrt{\det\bold g}\,\,d_{ij}.
\mytag{3.8}
$$
From \mythetag{3.7} one can derive the following relationship linking the
quantities $\omega_{ij}$ and $\tilde\omega_{pq}$ defined according to the
formula \mythetag{3.8} in two different coordinate systems:
$$
\hskip -2em
\omega_{ij}=\sign(\det S)\,\,
\sum^3_{p=1}\sum^3_{q=1} T^p_i\,T^q_j\,\,\tilde\omega_{pq}.
\mytag{3.9}
$$
Due to \mythetag{3.9} the quantities \mythetag{3.8} define a skew-symmetric
inner pseudotensorial field of the type $(0,2)$. It is called the 
{\it area pseudotensor}. If on a surface $D$ one of the two possible
orientations is marked, then the formula
$$
\hskip -2em
\omega_{ij}=\xi_D\,\sqrt{\det\bold g}\,\,d_{ij}
\mytag{3.10}
$$
defines a tensorial field of the type $(0,2)$. It is called the
{\it area tensor}. The formula \mythetag{3.10} differs from \mythetag{3.8} 
only in sign factor $\xi_D$ which is the unitary pseudoscalar field 
defining the orientation (compare with the formula \mythetagchapter{8.1}{2}
in Chapter \uppercase\expandafter{\romannumeral 2}). Here one should 
note that not any surface admits some preferable orientation globally.
The M\"obius strip is a well-known example of a non-orientable surface. 
\par
\head
\S~\mysection{4} Moving frame of a surface.\\
Veingarten's derivational formulas.
\endhead
\rightheadtext{\S~4. Veingarten's derivational formulas.}
     Each choice of a curvilinear coordinate system on a surface
determines some frame of two tangent vectorial fields $\bold E_1,\,
\bold E_2$ on it. The vectors of such a frame define the tangent 
plane at each point of the surface. However, they are insufficient
for to expand an arbitrary vector of the space $\Bbb E$ at that
point of the surface. Therefore, they are usually completed by 
a vector that does not belong to the tangent plane.
\mydefinition{4.1} A {\it unit normal vector } $\bold n$ 
to a surface $D$ at a point $A$ is a vector of the unit length
attached to the point $A$ and perpendicular to all vectors of the
tangent plane to $D$ at that point.
\enddefinition
     The definition~\mythedefinition{4.1} fixes the unit normal vector
$\bold n$ only up to the sign: at each point there are two opposite unit
vectors perpendicular to the tangent plane. One of the ways to fix $\bold
n$ uniquely is due to the vector product:
$$
\hskip -2em
\bold n=\frac{[\bold E_1,\,\bold E_2]}
{\left|[\bold E_1,\,\bold E_2]\right|}.
\mytag{4.1}
$$
The vector $\bold n$ determined by the formula \mythetag{4.1} depends on
the choice of a curvilinear coordinate system. Therefore, under a change
of coordinate system it can change its direction. Indeed, the relation 
of the frame vectors $\bold E_1,\bold E_2$ and $\tilde\bold E_1,\tilde
\bold E_2$ is given by the formula \mythetag{2.6}. Therefore, we write
$$
\xalignat 2
&\bold E_1=T^1_1\cdot\tilde\bold E_1+T^2_1\cdot\tilde\bold E_2,
&&\bold E_2=T^1_2\cdot\tilde\bold E_1+T^2_2\cdot\tilde\bold E_2.
\endxalignat
$$
Substituting these expressions into the vector product $[\bold E_1,\,
\bold E_2]$, we obtain
$$
[\bold E_1,\,\bold E_2]=(T^1_1\,T^2_2-T^2_1\,T^1_2)\cdot
[\tilde\bold E_1,\,\tilde\bold E_2]=\det T\cdot
[\tilde\bold E_1,\,\tilde\bold E_2].
$$
Now we easily derive the transformation rule for the normal vector 
$\bold n$:
$$
\hskip -2em
\bold n=(-1)^S\cdot\tilde\bold n.
\mytag{4.2}
$$
The sign factor $(-1)^S=\sign(\det S)=\pm 1$ here is the same as in the
formula \mythetag{2.8}.\par
     Another way of choosing the normal vector is possible if there is
a preferable orientation on a surface. Suppose that this orientation 
on $D$ is given by the unitary pseudoscalar field $\xi_D$. Then $\bold n$ is given by the formula
$$
\hskip -2em
\bold n=\xi_D\cdot\frac{[\bold E_1,\,\bold E_2]}
{\left|[\bold E_1,\,\bold E_2]\right|}.
\mytag{4.3}
$$
In this case the transformation rule for the normal vector simplifies
substantially:
$$
\bold n=\tilde\bold n.
\mytag{4.4}
$$
\mydefinition{4.2} The tangent frame $\bold E_1,\,\bold E_2$
of a curvilinear coordinate system $u^1,\,u^2$ on a surface completed
by the unit normal vector $\bold n$ is called the {\it moving frame}
or the {\it escort frame} of this surface.
\enddefinition
     If the normal vector is chosen according to the formula \mythetag{4.1},
the escort frame $\bold E_1,\,\bold E_2,\,\bold n$ is always right-oriented.
Therefore, in this case if we change the orientation of the tangent frame
$\bold E_1,\,\bold E_2$, the direction of the normal vector $\bold n$ is
changed immediately. In the other case, if $\bold n$ is determined by the
formula \mythetag{4.3}, then its direction does not depend on the choice of
the tangent frame $\bold E_1,\,\bold E_2$. This fact means that the choice
of the orientation on a surface is equivalent to choosing the normal vector
independent on the choice tangent vectors $\bold E_1,\bold E_2$.\par
     There is a special case, when such an independent choice of the normal
vector does exist. Let $D$ be the boundary of a three-dimensional domain.
Then one of two opposite normal vectors is the inner normal, the other is
the outer normal vector. Thus, we conclude that the boundary of a
three-dimensional domain in the space $\Bbb E$ is always orientable.\par
     Let $D$ be some fragment of a surface of the smoothness class $C^2$.
The vectors of the escort frame of such a surface are continuously
differentiable vector-functions of curvilinear coordinates: $\bold E_1(u^1,
u^2)$, $\bold E_2(u^1,u^2)$, and $\bold n(u^1,u^2)$. The derivatives of
such vectors are associated with the same point on the surface. Hence,
we can write the following expansions for them:
$$
\hskip -2em
\frac{\partial\bold E_j}{\partial u^i}=
\sum^2_{k=1}\Gamma^k_{ij}\cdot\bold E_k+b_{ij}\cdot\bold n.
\mytag{4.5}
$$
The derivatives of the unit vector $\bold n$ are perpendicular to this
vector (see lemma~\mythelemmachapter{3.1}{1} in Chapter
\uppercase\expandafter{\romannumeral 1}). Hence, we have the equality 
$$
\hskip -2em
\frac{\partial\bold n}{\partial u^i}=\sum^2_{k=1} c^k_i\cdot
\bold E_k.
\mytag{4.6}
$$
Let's consider the scalar product of \mythetag{4.5} and the vector 
$\bold n$. We also consider the scalar product of \mythetag{4.6} and 
the vector $\bold E_j$. Due to $(\bold E_k|\,\bold n)=0$ we get
$$
\align
&\hskip -2em
(\partial\bold E_j/\partial u^i\,|\,\bold n)=b_{ij}\,
(\bold n\,|\,\bold n)=b_{ij},
\mytag{4.7}\\
&\hskip -2em
(\bold E_j\,|\,\partial\bold n/\partial u^i)=
\sum^2_{k=1} c^k_i\,g_{kj}.
\mytag{4.8}
\endalign
$$
Let's add the left hand sides of the above formulas \mythetag{4.7} and
\mythetag{4.8}. Upon rather easy calculations we find that the sum is
equal to zero:
$$
(\partial\bold E_j/\partial u^i\,|\,\bold n)+
(\bold E_j\,|\,\partial\bold n/\partial u^i)=
\partial(\bold E_j\,|\,\bold n)/\partial u^i=0.
$$
From this equality we derive the relations of $b_{ij}$ and $c^k_j$ in \mythetag{4.5} and \mythetag{4.6}:
$$
b_{ij}=-\sum^2_{k=1} c^k_i\,g_{kj}.
$$
By means of the matrix of the inverse metric tensor $\hat\bold g$ we
can invert this relationship. Let's introduce the following quite 
natural notation:
$$
\hskip -2em
b^k_i=\sum^2_{j=1} b_{ij}\,g^{jk}.
\mytag{4.9}
$$
Then the coefficients $c^k_i$ in \mythetag{4.6} can be expressed through 
the coefficients $b_{ij}$ in \mythetag{4.5} by means of the following
formula:
$$
\hskip -2em
c^k_i=-b^k_i.
\mytag{4.10}
$$
Taking into account \mythetag{4.10}, we can rewrite \mythetag{4.5} and
\mythetag{4.6} as follows:
$$
\hskip -2em
\aligned
&\frac{\partial\bold E_j}{\partial u^i}=
\sum^2_{k=1}\Gamma^k_{ij}\cdot\bold E_k+b_{ij}\cdot\bold n,\\
\vspace{1ex}
&\frac{\partial\bold n}{\partial u^i}=-\sum^2_{k=1} b^k_i\cdot
\bold E_k.
\endaligned
\mytag{4.11}
$$
The expansions \mythetag{4.5} and \mythetag{4.6} written in form of
\mythetag{4.11} are called the {\it Veingarten's derivational 
formulas}. They determine the dynamics of the moving frame and
play the central role in the theory of surfaces.\par
\head
\S~\mysection{5} Christoffel symbols\\
and the second quadratic form.
\endhead
\rightheadtext{\S~5. Christoffel symbols and the second quadratic form.}
     Let's study the first Veingarten's derivational formula in two
different coordinate systems $u^1,\,u^2$ and $\tilde u^1,\,\tilde u^2$
on a surface. In the coordinates $u^1,\,u^2$ it is written as
$$
\hskip -2em
\frac{\partial\bold E_j}{\partial u^i}=
\sum^2_{k=1}\Gamma^k_{ij}\cdot\bold E_k+b_{ij}\cdot\bold n.
\mytag{5.1}
$$
In the other coordinates $\tilde u^1,\,\tilde u^2$ this formula is 
rewritten as
$$
\hskip -2em
\frac{\partial\tilde\bold E_q}{\partial\tilde u^p}=
\sum^2_{m=1}\tilde\Gamma^m_{pq}\cdot\tilde\bold E_m+
\tilde b_{pq}\cdot\tilde\bold n.
\mytag{5.2}
$$
Let's express the vector $\bold E_j$ in the left hand side of the formula
\mythetag{5.1} through the frame vectors of the second coordinate system.
For this purpose we use \mythetag{2.6}:
$$
\frac{\partial\bold E_j}{\partial u^i}=
\sum^2_{q=1}\frac{\partial (T^q_j\cdot\tilde\bold E_q)}
{\partial u^i}=\sum^2_{q=1}\frac{\partial T^q_j}
{\partial u^i}\cdot\tilde\bold E_q+\sum^2_{q=1}
T^q_j\cdot\frac{\partial\tilde\bold E_q}{\partial u^i}.
$$
For the further transformation of the above expression we use
the chain rule for differentiating the composite function:
$$
\frac{\partial\bold E_j}{\partial u^i}=
\sum^2_{m=1}\frac{\partial T^m_j}{\partial u^i}
\cdot\tilde\bold E_m+\sum^2_{q=1}\sum^2_{p=1}
\left(T^q_j\,\frac{\partial\tilde u^p}{\partial u^i}
\right)\cdot\frac{\partial\tilde\bold E_q}{\partial
\tilde u^p}.
$$
The values of the partial derivatives $\partial\tilde\bold E_q/
\partial\tilde u^p$ are determined by the formula \mythetag{5.2}.
Moreover, we should take into account \mythetag{2.7} in form of the
equality $\partial\tilde u^p/\partial u^i=T^p_i$:
$$
\align
&\quad\frac{\partial\bold E_j}{\partial u^i}=
\sum^2_{m=1}\frac{\partial T^m_j}{\partial u^i}
\cdot\tilde\bold E_m+\sum^2_{q=1}\sum^2_{p=1}\sum^2_{m=1}
(T^q_j\,T^p_i\,\tilde\Gamma^m_{pq})\cdot\tilde\bold E_m+\\
&\quad+\sum^2_{q=1}\sum^2_{p=1}(T^q_j\,T^p_i\,\tilde b_{pq})
\cdot\tilde\bold n=
\sum^2_{m=1}\sum^2_{k=1}\left(\frac{\partial T^m_j}{\partial u^i}
\,S^k_m\right)\cdot\bold E_k+\\
&+\sum^2_{q=1}\sum^2_{p=1}\sum^2_{m=1}\sum^2_{k=1}
(T^q_j\,T^p_i\,\tilde\Gamma^m_{pq}\,S^k_m)\cdot\bold E_k+
\sum^2_{q=1}\sum^2_{p=1}(T^q_j\,T^p_i\,\tilde b_{pq})
\cdot\tilde\bold n.
\endalign
$$
The unit normal vectors $\bold n$ and $\tilde\bold n$ can differ only in
sign: $\bold n=\pm\tilde\bold n$. Hence, the above expansion
for $\partial\bold E_j/\partial u^i$ and the expansion \mythetag{5.1} both
are the expansions in the same basis $\bold E_1,\,\bold E_2,\,\bold n$.
Therefore, we have
$$
\align
&\hskip -2em
\Gamma^k_{ij}=\sum^2_{m=1}S^k_m\,\frac{\partial T^m_j}
{\partial u^i}+\sum^2_{m=1}\sum^2_{p=1}\sum^2_{q=1}
S^k_m\,T^p_i\,T^q_j\,\,\tilde\Gamma^m_{pq},
\mytag{5.3}\\
&\hskip -2em
b_{ij}=\pm\sum^2_{p=1}\sum^2_{q=1}T^p_i\,T^q_j\,\,\tilde b_{pq}.
\mytag{5.4}\\
\endalign
$$
The formulas \mythetag{5.3} and \mythetag{5.4} express the transformation
rules for the coefficients $\Gamma^k_{ij}$ and $b_{ij}$ in Veingarten's
derivational formulas under a change of curvilinear coordinates on a 
surface.\par
     In order to make certain the choice of the sign in \mythetag{5.4}
one should fix some rule for choosing the unit normal vector. If we
choose $\bold n$ according to the formula \mythetag{4.1}, then under a 
change of coordinate system it obeys the transformation formula 
\mythetag{4.2}. In this case the equality \mythetag{5.4} is written as
$$
b_{ij}=\sum^2_{p=1}\sum^2_{q=1}(-1)^S\,T^p_i\,T^q_j\,\,
\tilde b_{pq}.
\mytag{5.5}
$$
It is easy to see that in this case $b_{ij}$ are the components of 
an inner pseudotensorial field of the type $(0,2)$ on a surface.\par
     Otherwise, if we use the formula \mythetag{4.3} for choosing the
normal vector $\bold n$, then $\bold n$ does not depend on the choice
of a curvilinear coordinate system on a surface (see formula \mythetag{4.4}).
In this case $b_{ij}$ are transformed as the components of a tensorial
field of the type $(0,2)$. The formula \mythetag{5.4} then takes the
form
$$
\hskip -2em
b_{ij}=\sum^2_{p=1}\sum^2_{q=1}T^p_i\,T^q_j\,\,\tilde b_{pq}.
\mytag{5.6}
$$\par
    Tensors of the type $(0,2)$ correspond to bilinear and quadratic
forms. Pseudo\-tensors of the type $(0,2)$ have no such interpretation.
Despite to this fact the quantities $b_{ij}$ in Veingarten's derivational 
formulas are called the components of the {\it second quadratic form
$\bold b$} of a surface. The following theorem supports this 
interpretation.
\mytheorem{5.1} The quantities $\Gamma^k_{ij}$ and $b_{ij}$
in Veingarten's derivational formulas \mythetag{4.11} are symmetric
with respect to the lower indices $i$ and $j$.
\endproclaim
\demo{Proof} In order to prove the theorem we apply the formula
\mythetag{1.10}. Let's write this formula in the following way:
$$
\hskip -2em
\bold E_j(u^1,u^2)=\frac{\partial\bold r(u^1,u^2)}
{\partial u^j}.
\mytag{5.7}
$$
Then let's substitute \mythetag{5.7} into the first formula \mythetag{4.11}:
$$
\hskip -2em
\frac{\partial^2\bold r(u^1,u^2)}
{\partial u^i\,\partial u^j}=
\sum^2_{k=1}\Gamma^k_{ij}\cdot\bold E_k+b_{ij}\cdot\bold n.
\mytag{5.8}
$$
The values of the mixed partial derivatives do not depend on the order
of differentiation. Therefore, the left hand side of \mythetag{5.8} is a
vector that does not change if we transpose indices $i$ and $j$. Hence,
the coefficients $\Gamma^k_{ij}$ and $b_{ij}$ of its expansion in the
basis $\bold E_1,\,\bold E_2,\,\bold n$ do not change under the 
transposition of the indices $i$ and $j$. The theorem is proved.
\qed\enddemo
     Sometimes, for the components of the matrix of the second quadratic
form the notations similar to \mythetag{3.3} are used:
$$
\hskip -2em
b_{ij}=\left\Vert
\matrix
L &M\\ M &N
\endmatrix
\right\Vert.
\mytag{5.9}
$$
These notations are especially popular in the earlier publications.\par
     The tensor fields $\bold g$ and $\bold b$ define a pair of quadratic
forms at each point of a surface. This fact explains in part their titles
--- the {\it first} and the {\it second} quadratic forms. The first 
quadratic form is non-degenerate and positive. This situation is well-known
in linear algebra (see \mycite{1}). Two forms, one of which is positive, can
be brought to the diagonal form simultaneously, the matrix of the positive
form being brought to the unit matrix. The diagonal elements of the second
quadratic form upon such diagonalization are called the {\it invariants of
a pair of forms}. In order to calculate these invariants we consider the
following contraction:
$$
\hskip -2em
b^k_i=\sum^2_{j=1} b_{ij}\,g^{jk}.
\mytag{5.10}
$$
The quantities $b^k_i$ enter the second Veingarten's derivational formula
\mythetag{4.11}. They define a tensor field (or a pseudotensorial field) of
the type $(1,1)$, i\.\,e\. an operator field. The operator with the matrix 
\mythetag{5.10} is called the {\it Veingarten operator}. The matrix of this
operator is diagonalized simultaneously with the matrices of the first and
the second quadratic forms, and its eigenvalues are exactly the invariants 
of that pair of forms. Let's denote them by $k_1$ and $k_2$.
\mydefinition{5.1} The eigenvalues $k_1(u^1,u^2)$ and
$k_2(u^1,u^2)$ for the matrix of the Veingarten operator are called the
{\it principal curvatures\/} of a surface at its point with the coordinates
$u^1,\,u^2$.
\enddefinition
From the computational point of view the other two invariants are more
convenient. These are the following ones:
$$
\xalignat 2
&\hskip -2em
H=\frac{k_1+k_2}{2},
&&K=k_1\,k_2.
\mytag{5.11}
\endxalignat
$$
The invariants \mythetag{5.11} can be calculated without knowing the
eigenvalues of the matrix $b^k_i$. It is sufficient to find the
trace for the matrix of the Veingarten operator and the determinant 
of this matrix:
$$
\xalignat 2
&\hskip -2em
H=\frac{1}{2}\tr(b^k_i),
&&K=\det(b^k_i).
\mytag{5.12}
\endxalignat
$$
The quantity $H$ in the formulas \mythetag{5.11} and \mythetag{5.12}
is called the {\it mean curvature}, while the quantity $K$ is called
the {\it Gaussian curvature}. There are formulas, expressing the invariants 
$H$ and $K$ through the components of the first and the second quadratic
forms \mythetag{3.3} and \mythetag{5.9}:
$$
\xalignat 2
&\qquad H=\frac{1}{2}\,\frac{EN+GL-2\,FM}{EG-F^2},
&
&K=\frac{LN-M^2}{EG-F^2}.
\mytag{5.13}
\endxalignat
$$\par
     Let $\bold v(u^1,u^2)$ and $\bold w(u^1,u^2)$ be the vectors of the
basis in which the matrix of the first quadratic form is equal to the unit
matrix, while the matrix of the second quadratic form is a diagonal matrix:
$$
\xalignat 2
&\hskip -2em
\bold v=\Vmatrix v^1\\ v^2\endVmatrix,
&&\bold w=\Vmatrix w^2\\ w^2\endVmatrix.
\mytag{5.14}
\endxalignat
$$
Then $\bold v$ and $\bold w$ are the eigenvectors of the Veingarten
operator. The vectors \mythetag{5.14} have their three-dimensional
realization in the space $\Bbb E$:
$$
\xalignat 2
&\hskip -2em
\bold v=v^1\cdot\bold E_1+v^2\cdot\bold E_2,
&&\bold w=w^1\cdot\bold E_1+w^2\cdot\bold E_2.
\mytag{5.15}
\endxalignat
$$
This is the pair of the unit vectors lying on the tangent plane and
being perpendicular to each other. The directions \pagebreak 
given by the vectors
\mythetag{5.15} are called the {\it principal directions} on a surface
at the point with coordinates $u^1,\,u^2$. If the principal curvatures
at this point are not equal to each other: $k_1\neq k_2$, then the
principal directions are determined uniquely. Otherwise, if $k_1=k_2$, 
then any two mutually perpendicular directions on the tangent plane 
can be taken for the principal directions. A point of a surface
where the principal curvatures are equal to each other ($k_1=k_2$) 
is called an {\it umbilical point}.\par
\subhead A remark on the sign\endsubhead Remember that depending on
the way how we choose the normal vector the second quadratic form is
either a tensorial field or a pseudotensorial field. The same is true
for the Veingarten operator. Therefore, in general, the principal
curvatures $k_1$ and $k_2$ are determined only up to the sign. The
mean curvature $H$ is also determined only up to the sign. As for the
Gaussian curvature, it has no uncertainty in sign. Moreover, the sign
of the Gaussian curvature divides the points of a surface into three
subsets: for any point of a surface if $K>0$, the point is called an 
{\it elliptic point}; if $K<0$, the point is called a {\it hyperbolic
point}; and finally, if $K=0$, the point is called a {\it parabolic 
point}.
\head
\S~\mysection{6} Covariant differentiation\\
of inner tensorial fields of a surface.
\endhead
\rightheadtext{\S~6. Covariant differentiation \dots}
     Let's consider the formula \mythetag{5.3} and compare it with the
formula \mythetagchapter{6.1}{3} in Chapter
\uppercase\expandafter{\romannumeral 3}. These two formulas differ only
in the ranges over which the indices run. Therefore the quantities
$\Gamma^k_{ij}$, which appear as coefficients in the Veingarten's
derivational formula, define a geometric object on a surface that is called
a {\it connection}. The connection components $\Gamma^k_{ij}$ are called
the {\it Christoffel symbols}.\par
     The main purpose of the Christoffel symbols $\Gamma^k_{ij}$ is their
usage for the covariant differentiation of tensor fields. Let's reproduce
here the formula \mythetagchapter{5.12}{3} from Chapter
\uppercase\expandafter{\romannumeral 3} for the covariant derivative 
modifying it for the two-dimensional case:
$$
\hskip -2em
\gathered
\nabla_{j_{s+1}}A^{i_1\ldots\,i_r}_{j_1\ldots\,j_s}=
 \frac{\partial A^{i_1\ldots\,i_r}_{j_1\ldots\,j_s}}
 {\partial u^{j_{s+1}}}+\\
 +\sum^r_{m=1}\sum^2_{v_m=1}\Gamma^{i_m}_{j_{s+1}\,v_m}\,
 A^{i_1\ldots\,v_m\ldots\,i_r}_{j_1\ldots\,j_s}
 -\sum^s_{n=1}\sum^2_{w_n=1}
  \Gamma^{w_n}_{j_{s+1}\,j_n}\,
 A^{i_1\ldots\,i_r}_{j_1\ldots\,w_n\ldots\,j_s}.
\endgathered
\mytag{6.1}
$$
\mytheorem{6.1} The formula \mythetag{6.1} correctly defines 
the covariant differentiation of inner tensor fields on a surface
that transforms a field of the type $(r,s)$ into a field of the
type $(r,s+1)$ if and only if the quantities $\Gamma^k_{ij}$ obey
the transformation rule \mythetag{5.3} under a change of curvilinear
coordinates on a surface.
\endproclaim
\demo{Proof} Let's begin with proving the necessity. For this purpose
we choose some arbitrary vector field $\bold A$ and produce the tensor
field $\bold B=\nabla\bold A$ of the type $(1,1)$ by means of the
formula \mythetag{6.1}. The correctness of the formula \mythetag{6.1} means
that the components of the field $\bold B$ are transformed according
to the formula \mythetag{2.8}. From this condition we should derive the
transformation formula \mythetag{5.3} for the quantities $\Gamma^k_{ij}$ 
in \mythetag{6.1}. Let's write the formula \mythetag{2.8} for the field
$\bold B=\nabla\bold A$:
$$
\frac{\partial A^k}{\partial u^i}+\sum^2_{j=1}
\Gamma^k_{ij}\,A^j=\sum^2_{m=1}\sum^2_{p=1}
S^k_m\,T^p_i\left(\frac{\partial\tilde A^m}{\partial
\tilde u^p}+\shave{\sum^2_{q=1}}\tilde\Gamma^m_{pq}\,
\tilde A^q\right).
$$
Then we expand the brackets in the right hand side of this relationship.
In the first summand we replace $T^p_i$ by $\partial\tilde u^p
/\partial u^i$ according to the formula \mythetag{2.7} and we express
$\tilde A^m$ through $A^j$ according to the transformation rule for a
vector field:
$$
\hskip -2em
\gathered
\sum^2_{p=1}T^p_i\,\frac{\partial\tilde A^m}{\partial
\tilde u^p}=\sum^2_{p=1}\frac{\partial\tilde u^p}
{\partial u^i}\,\frac{\partial\tilde A^m}{\partial
\tilde u^p}=\frac{\partial\tilde A^m}{\partial u^i}=\\
=\frac{\partial}{\partial u^i}\left(\shave{\sum^2_{k=1}}
T^m_k\,A^k\right)=\sum^2_{j=1}\frac{\partial T^m_j}
{\partial u^i}\,A^j+\sum^2_{j=1} T^m_j\,
\frac{\partial A^j}{\partial u^i}.
\endgathered
\mytag{6.2}
$$
Taking into account \mythetag{6.2}, we can cancel the partial derivatives 
in the previous equality and bring it to the following form:
$$
\sum^2_{j=1}\Gamma^k_{ij}\,A^j=\sum^2_{j=1}\sum^2_{m=1}
S^k_m\,\frac{\partial T^m_j}{\partial u^i}\,A^j+
\sum^2_{m=1}\sum^2_{p=1}\sum^2_{q=1}S^k_m\,T^p_i\,
\tilde\Gamma^m_{pq}\,\tilde A^q.
$$
Now we need only to express $\tilde A^q$ through $A^j$ applying the
transformation rule for the components of a vectorial field and then
extract the common factor $A^j$:
$$
\sum^2_{j=1}\left(\Gamma^k_{ij}-\shave{\sum^2_{m=1}}
S^k_m\,\frac{\partial T^m_j}{\partial u^i}-
\shave{\sum^2_{m=1}\sum^2_{p=1}\sum^2_{q=1}}
S^k_m\,T^p_i\,T^q_j\,\tilde\Gamma^m_{pq}\right)
A^j=0.
$$
Since $\bold A$ is an arbitrary vector field, each coefficient
enclosed into round brackets in the above sum should vanish separately.
This condition coincides exactly with the transformation rule 
\mythetag{5.3}. Thus, the necessity is proved.\par
     Let's proceed with proving the sufficiency. Suppose that the
condition \mythetag{5.3} is fulfilled. Let's choose some tensorial field
$\bold A$ of the type $(r,s)$ and prove that the quantities 
$\nabla_{j_{s+1}}A^{i_1\ldots\,i_r}_{j_1\ldots\,j_s}$ determined by the
formula \mythetag{6.1} are transformed as the components of a tensorial
field of the type $(r,s+1)$. Looking at the formula \mythetag{6.1} we see
that it contains the partial derivative $\partial A^{i_1\ldots\,i_r}_{j_1
\ldots\,j_s}/\partial u^{j_{s+1}}$ and other $r+s$ terms. Let's denote
these terms by $A^{i_1\ldots\,i_r}_{j_1\ldots\,j_{s+1}}
\tsize\thickfracwithdelims[]\thickness0{m}{0}$ and
$A^{i_1\ldots\,i_r}_{j_1\ldots\,j_{s+1}}
\tsize\thickfracwithdelims[]\thickness0{0}{n}$. Then
$$
\hskip -2em
\gathered
\nabla_{j_{s+1}}A^{i_1\ldots\,i_r}_{j_1\ldots\,j_s}=
 \partial A^{i_1\ldots\,i_r}_{j_1\ldots\,j_s}/\partial
 u^{j_{s+1}}\,+\\
  +\sum^r_{m=1}A^{i_1\ldots\,i_r}_{j_1\ldots\,j_{s+1}}
  {\tsize\thickfracwithdelims[]\thickness0{m}{0}}-
 \sum^s_{n=1}A^{i_1\ldots\,i_r}_{j_1\ldots\,j_{s+1}}
 {\tsize\thickfracwithdelims[]\thickness0{0}{n}}.
\endgathered
\mytag{6.3}
$$
The tensorial nature of $\bold A$ means that its components are transformed
according to the formula \mythetag{2.8}. Therefore, in the first term of 
\mythetag{6.3} we get:
$$
\allowdisplaybreaks
\align
&\frac{\partial A^{i_1\ldots\,i_r}_{j_1\ldots\,j_s}}
{\partial u^{j_{s+1}}}=
\frac{\partial}{\partial u^{j_{s+1}}}\left(\,
\shave{\sum\Sb p_1\ldots\,p_r\\ q_1\ldots\,q_s\endSb}
S^{i_1}_{p_1}\ldots\,S^{i_r}_{p_r}\,\,
T^{q_1}_{j_1}\ldots\,T^{q_s}_{j_s}\,\,
\tilde A^{p_1\ldots\,p_r}_{q_1\ldots\,q_s}\right)=\\
&\quad=\sum\Sb p_1\ldots\,p_r\\ q_1\ldots\,q_s\endSb
\sum^2_{q_{s+1}=1} S^{i_1}_{p_1}\ldots\,S^{i_r}_{p_r}\,\,
T^{q_1}_{j_1}\ldots\,T^{q_s}_{j_s}\,\,
\frac{\partial\tilde u^{q_{s+1}}}{\partial u^{j_{s+1}}}\,
\frac{\partial\tilde A^{p_1\ldots\,p_r}_{q_1\ldots\,q_s}}
{\partial\tilde u^{q_{s+1}}}\,+\\
&\quad+\sum^r_{m=1}\sum\Sb p_1\ldots\,p_r\\ q_1\ldots\,q_s\endSb
S^{i_1}_{p_1}\ldots\,\frac{\partial S^{i_m}_{p_m}}
{\partial u^{j_{s+1}}}\ldots\,S^{i_r}_{p_r}\,\,
T^{q_1}_{j_1}\ldots\,T^{q_s}_{j_s}\,\,\tilde A^{p_1\ldots\,
p_r}_{q_1\ldots\,q_s}\,+\\
&\quad+\sum^s_{n=1}\sum\Sb p_1\ldots\,p_r\\ q_1\ldots\,q_s\endSb
S^{i_1}_{p_1}\ldots\,S^{i_r}_{p_r}\,\,
T^{q_1}_{j_1}\ldots\,\frac{\partial T^{q_n}_{j_n}}
{\partial u^{j_{s+1}}}\ldots\,T^{q_s}_{j_s}\,\,
\tilde A^{p_1\ldots\,p_r}_{q_1\ldots\,q_s}.
\endalign
$$
Here we used the Leibniz rule for differentiating the product of
multiple functions and the chain rule in order to express the 
derivatives with respect to $u^{j_{s+1}}$ through the derivatives 
with respect to $\tilde u^{q_{s+1}}$. For the further transformation
of the above formula we replace $\partial\tilde u^{q_{s+1}}/
\partial u^{j_{s+1}}$ by $T^{q_{s+1}}_{j_{s+1}}$ according to
\mythetag{2.7} and we use the following identities based on the fact
that $S$ and $T$ are mutually inverse matrices:
$$
\hskip -2em
\aligned
&\frac{\partial S^{i_m}_{p_m}}{\partial u^{j_{s+1}}}=
-\sum^2_{v_m=1}\sum^2_{k=1}S^{i_m}_{v_m}\,
\frac{\partial T^{v_m}_k}{\partial u^{j_{s+1}}}\,
S^k_{p_m},\\
&\frac{\partial T^{q_n}_{j_n}}{\partial u^{j_{s+1}}}=
\sum^2_{w_n=1}\sum^2_{k=1}T^{w_n}_{j_n}\,S^k_{w_n}\,
\frac{\partial T^{q_n}_k}{\partial u^{j_{s+1}}}.
\endaligned
\mytag{6.4}
$$
Upon substituting \mythetag{6.4} into the preceding equality it is
convenient to transpose the summation indices: $p_m$ with $v_m$ 
and $q_n$ with $w_n$. Then we get
$$
\hskip -2em
\aligned
&\frac{\partial A^{i_1\ldots\,i_r}_{j_1\ldots\,j_s}}
{\partial u^{j_{s+1}}}=
\sum\Sb p_1\ldots\,p_r\\ q_1\ldots\,q_s\endSb
S^{i_1}_{p_1}\ldots\,S^{i_r}_{p_r}\,\,
T^{q_1}_{j_1}\ldots\,T^{q_s}_{j_s}\times\\
&\times\left(\,\shave{\sum^2_{q_{s+1}=1}}
T^{q_{s+1}}_{j_{s+1}}\,\frac{\tilde A^{p_1\ldots\,p_r}_{q_1
\ldots\,q_s}}{\partial\tilde u^{q_{s+1}}}
-\sum^r_{m=1}V{\tsize\thickfracwithdelims[]\thickness0{m}{0}}
+\sum^s_{n=1}W{\tsize\thickfracwithdelims[]\thickness0{0}{n}}
\right),
\endaligned
\mytag{6.5}
$$
where the following notations are used for the sake of simplicity:
$$
\aligned
&V{\tsize\thickfracwithdelims[]\thickness0{m}{0}}=
\sum^2_{v_m=1}\sum^2_{k=1} S^k_{v_m}\,
\frac{\partial T^{p_m}_k}{\partial u^{j_{s+1}}}\,
\tilde A^{p_1\ldots\,v_m\ldots\,p_r}_{q_1\ldots\,q_s},\\
&W{\tsize\thickfracwithdelims[]\thickness0{0}{n}}=
\sum^2_{w_n=1}\sum^2_{k=1} S^k_{q_n}\,
\frac{\partial T^{w_n}_k}{\partial u^{j_{s+1}}}\,
\tilde A^{p_1\ldots\,p_r}_{q_1\ldots\,w_n\ldots\,q_s}.
\endaligned
\mytag{6.6}
$$
No let's consider $A^{i_1\ldots\,i_r}_{j_1\ldots\,
j_{s+1}}\thickfracwithdelims[]\thickness0{m}{0}$ in \mythetag{6.3}. 
They are obviously defined as follows:
$$
A^{i_1\ldots\,i_r}_{j_1\ldots\,j_{s+1}}
{\tsize\thickfracwithdelims[]\thickness0{m}{0}}=
\sum^2_{v_m=1}\Gamma^{i_m}_{j_{s+1}\,v_m}\,
 A^{i_1\ldots\,v_m\ldots\,i_r}_{j_1\ldots\,j_s}.
$$
Applying the transformation rule \mythetag{2.8} to the components of the
field $\bold A$, we get:
$$
\aligned
A^{i_1\ldots\,i_r}_{j_1\ldots\,j_{s+1}}
&{\tsize\thickfracwithdelims[]\thickness0{m}{0}}=
\sum\Sb p_1\ldots\,p_r\\ q_1\ldots\,q_s\endSb\sum^2_{v_m=1}
S^{i_1}_{p_1}\ldots\,S^{v_m}_{p_m}\ldots\,S^{i_r}_{p_r}\times
\\
\vspace{1.5ex}
&\times T^{q_1}_{j_1}\ldots\,T^{q_s}_{j_s}\,\,
\Gamma^{i_m}_{j_{s+1}\,v_m}\,
\tilde A^{p_1\ldots\,p_r}_{q_1\ldots\,q_s}.
\endaligned
\mytag{6.7}
$$
For the further transformation of the above expression we use
\mythetag{5.3} written as
$$
\Gamma^{i_m}_{j_{s+1}\,v_m}=
\sum^2_{k=1}S^{i_m}_k\,\frac{\partial T^k_{v_m}}
{\partial u^{j_{s+1}}}+\sum^2_{k=1}\sum^2_{q=1}
\sum^2_{q_{s+1}=1} S^{i_m}_k\,T^{q_{s+1}}_{j_{s+1}}
\,T^q_{v_m}\,\,\tilde\Gamma^k_{{q_{s+1}}\,q}.
$$
Immediately after substituting this expression into \mythetag{6.7} we
perform the cyclic transposition of the summation indices: $r\to
p_m\to v_m\to r$. Some sums in the resulting expression are 
evaluated explicitly if we take into account the fact that the
transition matrices $S$ and $T$ are inverse to each other:
$$
\hskip -2em
\aligned
A^{i_1\ldots\,i_r}_{j_1\ldots\,j_{s+1}}
&{\tsize\thickfracwithdelims[]\thickness0{m}{0}}=
\sum\Sb p_1\ldots\,p_r\\ q_1\ldots\,q_s\endSb
S^{i_1}_{p_1}\ldots\,S^{i_r}_{p_r}\,\,
T^{q_1}_{j_1}\ldots\,T^{q_s}_{j_s}\times\\
&\times\left(\,\shave{\sum^2_{q_{s+1}=1}}
T^{q_{s+1}}_{j_{s+1}}\,\tilde A^{p_1\ldots\,
p_r}_{q_1\ldots\,q_{s+1}}
{\tsize\thickfracwithdelims[]\thickness0{m}{0}}
+V{\tsize\thickfracwithdelims[]\thickness0{m}{0}}
\right).
\endaligned
\mytag{6.8}
$$
By means of the analogous calculations one can derive the following
formula:
$$
\hskip -2em
\aligned
A^{i_1\ldots\,i_r}_{j_1\ldots\,j_{s+1}}
&{\tsize\thickfracwithdelims[]\thickness0{0}{n}}=
\sum\Sb p_1\ldots\,p_r\\ q_1\ldots\,q_s\endSb
S^{i_1}_{p_1}\ldots\,S^{i_r}_{p_r}\,\,
T^{q_1}_{j_1}\ldots\,T^{q_s}_{j_s}\times\\
&\times\left(\,\shave{\sum^2_{q_{s+1}=1}}
T^{q_{s+1}}_{j_{s+1}}\,\tilde A^{p_1\ldots\,
p_r}_{q_1\ldots\,q_{s+1}}
{\tsize\thickfracwithdelims[]\thickness0{0}{n}}
+W{\tsize\thickfracwithdelims[]\thickness0{0}{n}}
\right).
\endaligned
\mytag{6.9}
$$
Now we substitute \mythetag{6.5}, \mythetag{6.8}, and \mythetag{6.9} into
the formula \mythetag{6.3}. Then the entries of
$V\thickfracwithdelims[]\thickness0{m}{0}$ and 
$W\thickfracwithdelims[]\thickness0{0}{n}$ do cancel each other.
A a residue, upon collecting the similar terms and cancellations,
we get the formula expressing the transformation rule \mythetag{2.8}
applied to the components of the field $\nabla\bold A$. The 
theorem~\mythetheorem{6.1} is proved.
\qed\enddemo
     The theorem~\mythetheorem{6.1} yields a universal mechanism for
constructing the covariant differentiation. It is sufficient to have a
connection whose components are transformed according to the formula
\mythetag{5.3}. We can compare two connections: the Euclidean connection 
in the space $\Bbb E$ constructed by means of the Cartesian coordinates 
and a connection on a surface whose components are given by the
Veingarten's derivational formulas. Despite to the different origin of
these two connections, the covariant derivatives defined by them have many
common properties. It is convenient to formulate these properties using
covariant deriva\-tives along vector fields. Let $\bold X$ be a vector
field on a surface. For any tensor field $\bold A$ of the type $(r,s)$ we
define the tensor field $\bold B=\nabla_{\bold X}\bold A$ of the same type
$(r,s)$ given by the following formula:
$$
\hskip -2em
B^{i_1\ldots\,i_r}_{j_1\ldots\,j_s}=
\sum^2_{q=1}X^q\,\nabla_qA^{i_1\ldots\,i_r}_{j_1\ldots\,j_s}.
\mytag{6.10}
$$
\mytheorem{6.2} The operation of covariant differentiation 
of tensor fields posses\-ses the following properties: 
\roster
\item $\nabla_{\bold X}(\bold A+\bold B)=\nabla_{\bold X}\bold A
       +\nabla_{\bold X}\bold B$;
\item $\nabla_{\bold X+\bold Y}\bold A=\nabla_{\bold X}\bold A
       +\nabla_{\bold Y}\bold A$;
\item $\nabla_{\xi\cdot\bold X}\bold A=\xi\cdot\nabla_{\bold X}
       \bold A$;
\item $\nabla_{\bold X}(\bold A\otimes\bold B)=\nabla_{\bold X}
       \bold A\otimes\bold B+\bold A\otimes\nabla_{\bold X}
       \bold B$;
\item $\nabla_{\bold X}C(\bold A)=C(\nabla_{\bold X}\bold A)$;
\endroster
where $\bold A$ and $\bold B$ are arbitrary differentiable tensor fields,
while $\bold X$ and $\bold Y$ are arbitrary vector fields and $\xi$ is an
arbitrary scalar field.
\endproclaim
    Looking attentively at the theorem~\mythetheorem{6.2} and at the 
formula \mythetag{6.10}, we see that the theorem~\mythetheorem{6.2}
is a copy of the theorem~\mythetheoremchapter{5.2}{2} from Chapter
\uppercase\expandafter{\romannumeral 2}, while the formula
\mythetag{6.10} is a two-dimensional analog of the formula 
\mythetagchapter{5.5}{2} from
the same Chapter \uppercase\expandafter{\romannumeral 2}. However, the
proof there is for the case of the Euclidean connection in the space
$\Bbb E$. Therefore we need to give another proof.
\demo{Proof}  Let's choose some arbitrary curvilinear coordinate system
on a surface and prove the theorem by means of direct calculations in
coordinates. Denote $\bold C=\bold A+\bold B$, where $\bold A$ and $\bold B$
are two tensorial fields of the type $(r,s)$. Then
$$
C^{\,i_1\ldots\,i_r}_{j_1\ldots\,j_s}=A^{i_1\ldots\,i_r}_{j_1\ldots
\,j_s}+B^{i_1\ldots\,i_r}_{j_1\ldots\,j_s}.
$$
Substituting $C^{\,i_1\ldots\,i_r}_{j_1\ldots\,j_s}$ into \mythetag{6.1},
upon rather simple calculations we get
$$
\nabla_{j_{s+1}}C^{\,i_1\ldots\,i_r}_{j_1\ldots\,j_s}=
\nabla_{j_{s+1}}A^{i_1\ldots\,i_r}_{j_1\ldots\,j_s}+
\nabla_{j_{s+1}}B^{i_1\ldots\,i_r}_{j_1\ldots\,j_s}.
$$
The rest is to multiply both sides of the above equality by
$X^{j_{s+1}}$ and perform summation over the index $j_{s+1}$. 
Applying \mythetag{6.10}, we derive the formula of the item 
\therosteritem{1} in the theorem.\par
     Note that the quantities $B^{i_1\ldots\,i_r}_{j_1\ldots\,
j_s}$ in the formula \mythetag{6.10} are obtained as the linear
combinations of the components of $\bold X$. The items \therosteritem{2} 
and \therosteritem{3} of the theorem follow immediately from this 
fact.\par
     Let's proceed with the item \therosteritem{4}. Denote $\bold C=
\bold A\otimes\bold B$. Then for the components of the tensor field 
$\bold C$ we have the equality 
$$
\hskip -2em
C^{\,i_1\ldots\,i_{r+p}}_{j_1\ldots\,j_{s+q}}=
A^{i_1\ldots\,i_r}_{j_1\ldots\,j_s}
B^{i_{r+1}\ldots\,i_{r+p}}_{j_{s+1}\ldots\,j_{s+q}}.
\mytag{6.11}
$$
Let's substitute the quantities $C^{\,i_1\ldots\,i_{r+p}}_{j_1\ldots\,
j_{s+q}}$ from \mythetag{6.11} into the formula \mythetag{6.1} defining
the covariant derivative. As a result for $\bold D=\nabla\bold C$ we
derive
$$
\align
&D^{\,i_1\ldots\,i_{r+p}}_{j_1\ldots\,j_{s+q+1}}=
\bigl(\partial A^{\,i_1\ldots\,i_r}_{j_1\ldots\,j_s}/
\partial u^{j_{s+q+1}}\bigr)\,\,B^{i_{r+1}\ldots\,
i_{r+p}}_{j_{s+1}\ldots\,j_{s+q}}+\\
\vspace{1.5ex}
&\qquad\quad+A^{\,i_1\ldots\,i_r}_{j_1\ldots\,j_s}\,\,
\bigl(\partial B^{i_{r+1}\ldots\,i_{r+p}}_{j_{s+1}\ldots\,
j_{s+q}}/\partial u^{j_{s+q+1}}\bigr)+
\\
&+\sum^r_{m=1}\sum^2_{v_m=1}
\Gamma^{i_m}_{j_{s+q+1}\,v_m}\,
A^{i_1\ldots\,v_m\ldots\,i_r}_{j_1\ldots\,j_s}\,
B^{i_{r+1}\ldots\,i_{r+p}}_{j_{s+1}\ldots\,j_{s+q}}+
\\
&+\sum^{r+p}_{m=r+1}\sum^2_{v_m=1}
A^{\,i_1\ldots\,i_r}_{j_1\ldots\,j_s}\,
\Gamma^{i_m}_{j_{s+q+1}\,v_m}\,
B^{i_{r+1}\ldots\,v_m\ldots\,i_{r+p}}_{j_{s+1}\ldots\,j_{s+q}}-
\\
&-\sum^s_{n=1}\sum^2_{w_n=1}
\Gamma^{w_n}_{j_{s+q+1}\,j_n}\,
A^{i_1\ldots\,i_r}_{j_1\ldots\,w_n\ldots\,j_s}\,
B^{i_{r+1}\ldots\,i_{r+p}}_{j_{s+1}\ldots\,j_{s+q}}-
\\
&-\sum^{s+q}_{n=s+1}\sum^2_{w_n=1}
A^{\,i_1\ldots\,i_r}_{j_1\ldots\,j_s}\,
\Gamma^{w_n}_{j_{s+q+1}\,j_n}\,
B^{i_{r+1}\ldots\,i_{r+p}}_{j_{s+1}\ldots\,w_n\ldots\,j_{s+q}}.
\endalign
$$
Note that upon collecting the similar terms the above huge formula
can be transformed to the following shorter one:
$$
\hskip -2em
\aligned
\nabla_{j_{s+q+1}}&\bigl(A^{i_1\ldots\,i_r}_{j_1\ldots\,j_s}\,
B^{i_{r+1}\ldots\,i_{r+p}}_{j_{s+1}\ldots\,j_{s+q}}\bigr)=
\bigl(\nabla_{j_{s+q+1}}A^{i_1\ldots\,i_r}_{j_1\ldots\,j_s}
\bigr)\times
\\
\vspace{2ex}
&\times B^{i_{r+1}\ldots\,i_{r+p}}_{j_{s+1}\ldots\,j_{s+q}}
+A^{i_1\ldots\,i_r}_{j_1\ldots\,j_s}
\,\bigl(\nabla_{j_{s+q+1}}B^{i_{r+1}\ldots\,i_{r+p}}_{j_{s+1}
\ldots\,j_{s+q}}\bigr).
\endaligned
\mytag{6.12}
$$
Now in order to prove the fourth item of the theorem it is sufficient
to multiply \mythetag{6.12} by $X^{j_{s+q+1}}$ and sum up over the
index $j_{s+q+1}$.\par
     Proceeding with the last fifth item of the theorem, we consider
two tensorial fields $\bold A$ and $\bold B$ one of which is the
contraction of another:
$$
\hskip -2em
B^{i_1\ldots\,i_r}_{j_1\ldots\,j_s}=
\sum^2_{k=1}A^{i_1\ldots\,i_{p-1}\,k\,i_p\ldots\,i_r}_{j_1
\ldots\,j_{q-1}\,k\,j_q\ldots\,j_s}.
\mytag{6.13}
$$
Substituting \mythetag{6.13} into \mythetag{6.1}, for 
$\nabla_{j_{s+1}}B^{i_1\ldots\,i_r}_{j_1\ldots\,j_s}$ we obtain
$$
\align
&\hskip -2em
\nabla_{\!j_{s+1}}B^{i_1\ldots\,i_r}_{j_1\ldots\,j_s}=
\sum^2_{k=1}
\frac{\partial A^{i_1\ldots\,i_{p-1}\,k\,i_p\ldots\,
i_r}_{j_1\ldots\,j_{q-1}\,k\,j_q\ldots\,j_s}}
{\partial u^{j_{s+1}}}\,+\\
&\hskip -2em
+\sum^r_{m=1}\sum^2_{k=1}\sum^2_{v_m=1}
\Gamma^{i_m}_{j_{s+1}\,v_m}\,A^{i_1\ldots\,v_m\ldots\,k\,
\ldots\,i_r}_{j_1\ldots\,j_{q-1}\,k\,j_q\ldots\,j_s}\,-
\mytag{6.14}\\
&\hskip -2em
-\sum^s_{n=1}\sum^2_{k=1}\sum^2_{w_n=1}
\Gamma^{w_n}_{j_{s+1}\,j_n}\,A^{i_1\ldots\,i_{p-1}\,k\,i_p
\ldots\,i_r}_{j_1\ldots\,w_n\ldots\,k\,\ldots\,j_s}.
\endalign
$$
The index $v_m$ in \mythetag{6.14} can be either to the left of the index 
$k$ or to the right of it. The same is true for the index $w_n$. However,
the formula \mythetag{6.14} does not comprise the terms where $v_m$ or $w_n$
replaces the index $k$. Such terms would have the form:
$$
\gather
\hskip -2em
\sum^2_{k=1}\sum^2_{v=1}
\Gamma^k_{j_{s+1}\,v}\,\,A^{i_1\ldots\,i_{p-1}\,v\,i_p\ldots
\,i_r}_{j_1\ldots\,j_{q-1}\,k\,j_q\ldots\,j_s},
\mytag{6.15}\\
\hskip -2em
-\sum^2_{k=1}\sum^2_{w=1}
\Gamma^w_{j_{s+1}\,k}\,\,A^{i_1\ldots\,i_{p-1}\,k\,i_p\ldots
\,i_r}_{j_1\ldots\,j_{q-1}\,w\,j_q\ldots\,j_s}.
\mytag{6.16}
\endgather
$$
It is easy to note that \mythetag{6.15} and \mythetag{6.16} differ only
in sign. It is sufficient to rename $k$ to $v$ and $w$ to $k$ in the
formula \mythetag{6.16}. Adding both \mythetag{6.15} and \mythetag{6.16} to
\mythetag{6.14} would not break the equality. But upon adding them one 
can rewrite the equality \mythetag{6.14} in the following form:
$$
\hskip -2em
\nabla_{j_{s+1}}B^{i_1\ldots\,i_r}_{j_1\ldots\,j_s}=
\sum^2_{k=1}\nabla_{j_{s+1}}
A^{i_1\ldots\,i_{p-1}\,k\,i_p\ldots\,i_r}_{j_1\ldots\,j_{q-1}
\,k\,j_q\ldots\,j_s}.
\mytag{6.17}
$$
No in order to complete the proof of the item \therosteritem{5},
and thus prove the theorem in whole, it is sufficient to multiply
the equality \mythetag{6.17} by $X^{j_{s+1}}$ and sum up over the 
index $j_{s+1}$.
\qed\enddemo
     Among the five properties of the covariant derivative listed 
in the theorem~\mythetheorem{6.2} the fourth property written as
\mythetag{6.12} is most often used in calculations. Let's rewrite the
equality \mythetag{6.12} as follows:
$$
\nabla_{\!k}\bigl(A^{i_1\ldots\,i_r}_{j_1\ldots\,j_s}\,
B^{v_1\ldots\,v_p}_{w_1\ldots\,w_q}\bigr)=
\bigl(\nabla_{\!k}A^{i_1\ldots\,i_r}_{j_1\ldots\,j_s}
\bigr)\,B^{v_1\ldots\,v_p}_{w_1\ldots\,w_q}
+A^{i_1\ldots\,i_r}_{j_1\ldots\,j_s}
\,\bigl(\nabla_{\!k}B^{v_1\ldots\,v_p}_{w_1\ldots\,w_q}
\bigr).\quad
\mytag{6.18}
$$
The formula \mythetag{6.18} is produced from \mythetag{6.12} simply by renaming
the indices; however, it is more convenient for reception.\par
\head
\S~\mysection{7} Concordance of metric and connection on a surface.
\endhead
\rightheadtext{\S~7. Concordance of metric and connection on a surface.}
     Earlier we have noted that the covariant differential of the metric tensor $\bold g$ in the Euclidean space $\Bbb E$ is equal to zero (see
formula \mythetagchapter{6.7}{2}
in Chapter \uppercase\expandafter{\romannumeral 2}).
This property was called the concordance of the metric and connection.
Upon passing to curvilinear coordinates we used this property in order to express the components of the Euclidean connection $\Gamma^k_{ij}$ through
the components of the metric tensor (see formula \mythetagchapter{7.8}{3}
in Chapter \uppercase\expandafter{\romannumeral 3}). Let's study whether
the metric and the connection on surfaces are concordant or not. The answer
here is also positive. It is given by the following theorem.
\mytheorem{7.1} The components of the metric tensor $g^{ij}$ and 
the connection components $\Gamma^k_{ij}$ in arbitrary coordinates on a
surface are related by the equality 
$$
\hskip -2em
\nabla_kg_{ij}=\frac{\partial g_{ij}}{\partial u^k}-
\sum^2_{q=1}\Gamma^q_{ki}\,g_{qj}-
\sum^2_{q=1}\Gamma^q_{kj}\,g_{iq}=0
\mytag{7.1}
$$
which expresses the concordance condition for the metric and connection.
\endproclaim
\demo{Proof} Let's consider the first Veingarten's derivational formula 
in \mythetag{4.11} and let's rewrite it renaming some indices:
$$
\hskip -2em
\frac{\partial\bold E_i}{\partial u^k}=
\sum^2_{q=1}\Gamma^q_{ki}\cdot\bold E_q+b_{ki}\cdot\bold n.
\mytag{7.2}
$$
Let's take the scalar products of both sides of \mythetag{7.2} by $\bold E_j$
and remember that the vectors $\bold E_j$ and $\bold n$ are perpendicular.
The scalar product of $\bold E_q$ and $\bold E_j$ in the right hand side 
yields the element $g_{qj}$ of the Gram matrix:
$$
\hskip -2em
(\partial\bold E_i/\partial u^k\,|\,\bold E_j)=
\sum^2_{q=1}\Gamma^q_{ki}\,g_{qj}.
\mytag{7.3}
$$
Now let's transpose the indices $i$ and $j$ in \mythetag{7.3} and take into
account the symmetry of the Gram matrix. As a result we obtain
$$
\hskip -2em
(\bold E_i\,|\,\partial\bold E_j/\partial u^k)=
\sum^2_{q=1}\Gamma^q_{kj}\,g_{iq}.
\mytag{7.4}
$$
Then let's add \mythetag{7.3} with \mythetag{7.4} and remember the Leibniz
rule as applied to the differentiation of the scalar product in the space
$\Bbb E$:
$$
\align
&(\partial\bold E_i/\partial u^k\,|\,\bold E_j)+
(\bold E_i\,|\,\partial\bold E_j/\partial u^k)=
\partial(\bold E_i\,|\,\bold E_j)/\partial u^k=\\
\vspace{1ex}
&\qquad=\partial g_{ij}/\partial u^k=
\sum^2_{q=1}\Gamma^q_{ki}\,g_{qj}+
\sum^2_{q=1}\Gamma^q_{kj}\,g_{iq}.
\endalign
$$
Now it is easy to see that the equality just obtained coincides in
essential with \mythetag{7.1}. The theorem is proved.
\qed\enddemo
    As an immediate consequence of the theorems~\mythetheorem{7.1} 
and \mythetheorem{5.1} we get the following formula for the connection components:
$$
\hskip -2em
\Gamma^k_{ij}=\frac{1}{2}\sum^2_{r=1}g^{kr}\left(
\frac{g_{rj}}{\partial u^i}+\frac{g_{ir}}{\partial u^j}-
\frac{g_{ij}}{\partial u^r}\right).
\mytag{7.5}
$$
We do not give its prove here since, in essential, it is the same as in
the case of the formula \mythetagchapter{7.8}{3} in Chapter
\uppercase\expandafter{\romannumeral 3}.\par
     From the condition $\nabla_qg_{ij}=0$ one can easily derive that
the covariant derivatives of the inverse metric tensor are also equal 
to zero. For this purpose one should apply the formula \mythetag{3.4}.
The covariant derivatives of the identical operator field with the
components $\delta^i_k$ are equal to zero. Indeed, we have
$$
\hskip -2em
\nabla_q\delta^i_k=\frac{\partial(\delta^i_k)}
{\partial u^q}+\sum^2_{r=1}\Gamma^i_{qr}\delta^r_k-
\sum^2_{r=1}\Gamma^r_{qk}\delta^i_r=0.
\mytag{7.6}
$$
Let's differentiate both sides of \mythetag{3.4} and take into account
\mythetag{7.6}:
$$
\align
\hskip -2em
\nabla_q&\left(\,\shave{\sum^2_{j=1}} g^{ij}\,g_{jk}\right)=
\sum^2_{j=1}(\nabla_qg^{ij}\,g_{jk}+g^{ij}\,\nabla_qg_{jk})=\\
\hskip -2em
&\qquad=\sum^2_{j=1}\nabla_qg^{ij}\,g_{jk}=\nabla_q\delta^i_k=0.
\mytag{7.7}
\endalign
$$
In deriving \mythetag{7.7} we used the items \therosteritem{4} and
\therosteritem{5} from the theorem~\mythetheorem{6.2}. The procedure 
of lowering $j$ by means of the contraction with the metric tensor 
$g_{jk}$ is an invertible operation. Therefore, \mythetag{7.7} yields
$\nabla_qg^{ij}=0$. Now the concordance condition for the metric
and connection is written as a pair of two relationships 
$$
\xalignat 2
&\hskip -2em
\nabla\bold g=0,
&&\nabla\hat\bold g=0,
\mytag{7.8}
\endxalignat
$$
which look exactly like the relationships \mythetagchapter{6.7}{2} in
Chapter \uppercase\expandafter{\romannumeral 2} for the case of metric
tensors in the space $\Bbb E$.\par
     Another consequence of the theorem~\mythetheorem{7.1} is that the
index raising and the index lowering operations performed by means of
contractions with the metric tensors $\hat\bold g$ and $\bold g$ commute
with the operations of covariant differentiations. This fact is presented 
by the following two formulas:
$$
\hskip -2em
\aligned
&\nabla_{\!q}\!\left(\,\shave{\sum^2_{k=1}}g_{ik}\,
 A^{\ldots\,k\,\ldots}_{\ldots\,\hphantom{k}\,\ldots}\right)=
 \sum^2_{k=1}g_{ik}\,
 \nabla_{\!q}A^{\ldots\,k\,\ldots}_{\ldots\,\hphantom{k}\,\ldots},\\
&\nabla_{\!q}\!\left(\,\shave{\sum^2_{k=1}}g^{ik}\,
 A^{\ldots\,\hphantom{k}\,\ldots}_{\ldots\,k\,\ldots}\right)=
 \sum^2_{k=1}g^{ik}\,
 \nabla_{\!q}A^{\ldots\,\hphantom{k}\,\ldots}_{\ldots\,k\,\ldots}.
\endaligned
\mytag{7.9}
$$
The relationship \mythetag{7.9} is easily derived from \mythetag{7.8} 
using the items \therosteritem{4} and \therosteritem{5} in the
theorem~\mythetheorem{6.2}.\par
\mytheorem{7.2} The covariant differential of the area
pseudotensor \mythetag{3.8} on any surface is equal to zero:
$\nabla\boldsymbol\omega=0$.
\endproclaim
In order to prove this theorem we need two auxiliary propositions
which are formulated as the following lemmas.
\mylemma{7.1} For any square matrix $M$ whose components are
differentiable functions of some parameter $x$ the equality
$$
\hskip -2em
\frac{d(\ln\det M)}{dx}=\tr(M^{-1}\,M')
\mytag{7.10}
$$
is fulfilled, where $M'$ is the square matrix composed by the derivatives
of the corresponding components of the matrix $M$.
\endproclaim
\mylemma{7.2} For any square $2\times 2$ matrix $M$ the equality
$$
\hskip -2em
\sum^2_{q=1}\bigl(M_{iq}\,d_{qj}+M_{jq}\,d_{iq}\bigr)=
\tr M\,d_{ij}
\mytag{7.11}
$$
is fulfilled, where $d_{ij}$ are the components of the skew-symmetric
matrix determined by the relationship \mythetag{3.6}.
\endproclaim
    The proof of these two lemmas~\mythelemma{7.1} and \mythelemma{7.2}
as well as the proof of the above formula \mythetag{3.7} from \S\,3 can be
found in \mycite{4}.\par
     Let's apply the lemma~\mythelemma{7.1} to the matrix of the metric
tensor. Let $x=u^k$. Then we rewrite the relationship \mythetag{7.10} as
follows:
$$
\hskip -2em
\frac{1}{\sqrt{\det\bold g}}\,
\frac{\partial\sqrt{\det\bold g}}{\partial u^k}=
\frac{1}{2}\sum^2_{q=1}\sum^2_{p=1}g^{qp}\,
\frac{\partial g_{qp}}{\partial u^k}.
\mytag{7.12}
$$
Note that in \mythetag{7.11} any array of four numbers enumerated with 
two indices can play the role of the matrix $M$. Having fixed the index
$k$, one can use the connection components $\Gamma^j_{ki}$ as such an
array. Then we obtain
$$
\hskip -2em
\sum^2_{q=1}\bigl(\Gamma^q_{ki}\,d_{qj}+\Gamma^q_{kj}
\,d_{iq}\bigr)=\sum^2_{q=1}\Gamma^q_{kq}\,d_{ij}.
\mytag{7.13}
$$
\demo{Proof for the theorem~\mythetheorem{7.2}} The components of the area pseudotensor
$\boldsymbol\omega$ are determined by the formula \mythetag{3.8}. In order
to find the components of the pseudotensor $\nabla\boldsymbol\omega$ we
apply the formula \mythetag{6.1}. It yields
$$
\aligned
&\nabla_{\!k}\omega_{ij}=\frac{\partial\sqrt{\det\bold g}}
{\partial u^k}\,\,d_{ij}-\sum^2_{q=1}\sqrt{\det\bold g}\,
\bigl(\Gamma^q_{ki}\,d_{qj}+\Gamma^q_{kj}\,d_{iq}\bigr)=\\
&=\sqrt{\det\bold g}\left(\frac{1}{\sqrt{\det\bold g}}\,
\frac{\partial\sqrt{\det\bold g}}{\partial u^k}\,d_{ij}-
\shave{\sum^2_{q=1}}\bigl(\Gamma^q_{ki}\,d_{qj}+
\Gamma^q_{kj}\,d_{iq}\bigr)\right).
\endaligned
$$
For the further transformation of this expression we apply 
\mythetag{7.12} and \mythetag{7.13}:
$$
\nabla_k\omega_{ij}=\sqrt{\det\bold g}\left(
\frac{1}{2}\shave{\sum^2_{q=1}\sum^2_{p=1}}g^{qp}\,
\frac{\partial g_{qp}}{\partial u^k}-
\shave{\sum^2_{q=1}}\Gamma^q_{kq}\right)d_{ij}.
\hskip -1.5em
\mytag{7.14}
$$
Now let's express $\Gamma^q_{kq}$ through the components of the
metric tensor by means of the formula \mythetag{7.5}. Taking into
account the symmetry of $g^{pq}$, we get
$$
\sum^2_{q=1}\Gamma^q_{kq}=
\frac{1}{2}\sum^2_{q=1}\sum^2_{p=1}g^{qp}
\left(\frac{g_{pq}}{\partial u^k}+\frac{g_{kp}}
{\partial u^q}-\frac{g_{kq}}{\partial u^p}\right)=
\frac{1}{2}\sum^2_{q=1}\sum^2_{p=1}g^{qp}\,\frac{g_{qp}}
{\partial u^k}.
$$
Substituting this expression into the formula \mythetag{7.14}, we find 
that it vanishes. Hence, $\nabla_k\omega_{ij}=0$. The theorem is
proved.\qed\enddemo
\subhead A remark on the sign\endsubhead The area tensor differs from
the area pseudotensor only by the scalar sign factor $\xi_D$. Therefore,
the proposition of the theorem~\mythetheorem{7.2} for the area tensor of an arbitrary
surface is also valid. 
\subhead A remark on the dimension\endsubhead For the volume tensor 
(and for the volume pseudotensor) in the Euclidean space $\Bbb E$ 
we have the analogous proposition: it states that $\nabla\boldsymbol
\omega=0$. The proof of this proposition is even more simple than the
proof of the theorem~\mythetheorem{7.2}. The components of the field
$\boldsymbol\omega$ in any Cartesian coordinate system in $\Bbb E$ are
constants. Hence, their derivatives are zero.
\head
\S~\mysection{8} Curvature tensor.
\endhead
\rightheadtext{\S~8. Curvature tensor.}
     The covariant derivatives in the Euclidean space $\Bbb E$ are reduced
to the partial derivatives in any Cartesian coordinates. Is there such a
coordinate system for covariant derivatives on a surface\,? The answer to
this question is related to the commutators. Let's choose a vector field
$\bold X$ and calculate the tensor field $\bold Y$ of the type $(1,2)$
with the following components:
$$
\hskip -2em
Y^k_{ij}=\nabla_i\nabla_jX^k-\nabla_j\nabla_iX^k.
\mytag{8.1}
$$
In order to calculate the components of the field $\bold Y$ we apply the
formula \mythetag{6.1}:
$$
\hskip -2em
\aligned
&\nabla_i\nabla_jX^k=\frac{\partial(\nabla_jX^k)}
{\partial u^i}+\sum^2_{q=1}\Gamma^k_{iq}\nabla_jX^q-
\sum^2_{q=1}\Gamma^q_{ij}\nabla_qX^k,\\
&\nabla_j\nabla_iX^k=\frac{\partial(\nabla_iX^k)}
{\partial u^j}+\sum^2_{q=1}\Gamma^k_{jq}\nabla_iX^q-
\sum^2_{q=1}\Gamma^q_{ji}\nabla_qX^k.
\endaligned
\mytag{8.2}
$$
Let's subtract the second relationship \mythetag{8.2} from the first 
one. Then the last terms in them do cancel each other due to the
symmetry of $\Gamma^k_{ij}$:
$$
\align
&Y^k_{ij}=
\frac{\partial}{\partial u^i}
\left(\frac{\partial X^k}{\partial u^j}+
\shave{\sum^2_{r=1}}\Gamma^k_{jr}\,X^r\right)
-\frac{\partial}{\partial u^j}
\left(\frac{\partial X^k}{\partial u^i}+
\shave{\sum^2_{r=1}}\Gamma^k_{ir}\,X^r\right)+\\
&+\sum^2_{q=1}\Gamma^k_{iq}\left(\frac{\partial X^q}
{\partial u^j}+\shave{\sum^2_{r=1}}\Gamma^q_{jr}\,X^r\right)
-\sum^2_{q=1}\Gamma^k_{jq}\left(\frac{\partial X^q}
{\partial u^i}+\shave{\sum^2_{r=1}}\Gamma^q_{ir}\,
X^r\right).
\endalign
$$
Upon expanding the brackets and some cancellations we get
$$
\hskip -2em
Y^k_{ij}=\sum^2_{r=1}\left(
\frac{\partial\Gamma^k_{jr}}{\partial u^i}
-\frac{\partial\Gamma^k_{ir}}{\partial u^j}+
\shave{\sum^2_{q=1}}\bigl(\Gamma^k_{iq}\Gamma^q_{jr}-
\Gamma^k_{jq}\Gamma^q_{ir}\bigr)
\right)X^r.\hskip -2ex
\mytag{8.3}
$$
It is important to note that the formula \mythetag{8.3} does not contain
the derivatives of the components of $\bold X$ --- they are canceled.
Let's denote
$$
\hskip -2em
R^k_{rij}=
\frac{\partial\Gamma^k_{jr}}{\partial u^i}
-\frac{\partial\Gamma^k_{ir}}{\partial u^j}+
\sum^2_{q=1}\Gamma^k_{iq}\Gamma^q_{jr}-
\sum^2_{q=1}\Gamma^k_{jq}\Gamma^q_{ir}.
\mytag{8.4}
$$
The formula \mythetag{8.3} for the components of the field
\mythetag{8.1} then can be written as 
$$
\hskip -2em
(\nabla_i\nabla_j-\nabla_j\nabla_i)X^k=
\sum^2_{r=1}R^k_{rij}\,X^r.
\mytag{8.5}
$$
Let's replace the vector field $\bold X$ by a covector field.
Performing the similar calculations, in this case we obtain
$$
(\nabla_i\nabla_j-\nabla_j\nabla_i)X_k=
-\sum^2_{r=1}R^r_{kij}\,X_r.
\mytag{8.6}
$$
The formulas \mythetag{8.5} and\mythetag{8.6} can be generalized 
for the case of an arbitrary tensor field $\bold X$ of the
type $(r,s)$:
$$
\hskip -2em
\aligned
(\nabla_i\nabla_j-\nabla_j\nabla_i)&X^{i_1\ldots\,i_r}_{j_1
\ldots\,j_s}=\sum^r_{m=1}\sum^2_{v_m=1}
R^{i_m}_{v_m\,ij}\,X^{i_1\ldots\,v_m\ldots\,i_r}_{j_1\ldots\,j_s}-\\
&-\sum^s_{n=1}\sum^2_{w_n=1} R^{w_n}_{j_n\,ij}\,
X^{i_1\ldots\,i_r}_{j_1\ldots\,w_n\ldots\,j_s}.
\endaligned
\mytag{8.7}
$$
Comparing \mythetag{8.5}, \mythetag{8.6}, and \mythetag{8.7}, we see that
all of them contain the quantities $R^k_{rij}$ given by the formula
\mythetag{8.4}.
\mytheorem{8.1} The quantities $R^k_{rij}$ introduced by the 
formula \mythetag{8.4} define a tensor field of the type $(1,3)$. This
tensor field is called the {\it curvature tensor} or the {\it Riemann
tensor}.
\endproclaim
     The theorem~\mythetheorem{8.1} can be proved directly on the base
of the formula \mythetag{5.3}. However, we give another proof which is 
more simple. 
\mylemma{8.1} Let $\bold R$ be a geometric object which is
presented by a four-dimensio\-nal array $R^k_{rij}$ in coordinates. 
If the contraction of $\bold R$ with an arbitrary vector $\bold X$
$$
Y^k_{ij}=\sum^2_{q=1}R^k_{qij}\,X^q
\mytag{8.8}
$$
is a tensor of the type $(1,2)$, then the object $\bold R$ itself is 
a tensor of the type $(1,3)$.
\endproclaim
\demo{Proof of the lemma} Let $u^1,\,u^2$ and $\tilde u^1,\,\tilde u^2$
be two curvilinear coordinate systems on a surface. Let's fix some
numeric value of the index $r$ ($r=1$ or $r=2$). Since $\bold X$ is an
arbitrary vector, we choose this vector so that its $r$-th component
in the coordinate system $u^1,\,u^2$ is equal to unity, while all other
components are equal to zero. Then for $Y^k_{ij}$ in this coordinate 
system we get 
$$
\hskip -2em
Y^k_{ij}=\sum^2_{q=1}R^k_{qij}\,X^q=R^k_{rij}.
\mytag{8.9}
$$
For the components of the vector $\bold X$ in the other coordinate
system we derive
$$
\tilde X^m=\sum^2_{q=1}T^m_q\,X^q=T^m_r,
$$
then we apply the formula \mythetag{8.8} on order to calculate the components
of the tensor $\bold Y$ in the second coordinate system:
$$
\hskip -2em
\tilde Y^n_{pq}=\sum^2_{m=1}\tilde R^n_{mpq}\,\tilde X^m=
\sum^2_{m=1}\tilde R^n_{mpq}\,T^m_r.
\mytag{8.10}
$$
The rest is to relate the quantities $Y^k_{ij}$ from \mythetag{8.9} and
the quantities $\tilde Y^n_{pq}$ from \mythetag{8.10}. From the statement
of the theorem we know that these quantities are the components of the
same tensor in two different coordinate systems. Hence, we get
$$
\hskip -2em
Y^k_{ij}=\sum^2_{n=1}\sum^2_{p=1}\sum^2_{q=1}
S^k_n\,T^p_i\,T^q_j\,\,\tilde Y^n_{pq}.
\mytag{8.11}
$$
Substituting \mythetag{8.9} and \mythetag{8.10} into the formula
\mythetag{8.11}, we find
$$
R^k_{rij}=\sum^2_{n=1}\sum^2_{m=1}\sum^2_{p=1}\sum^2_{q=1}
S^k_n\,T^m_rT^p_i\,T^q_j\,\,\tilde R^n_{mpq}.
$$
This formula exactly coincides with the transformation rule for the
components of a tensorial field of the type $(1,3)$ under a change
of coordinates. Thus, the lemma is proved.
\qed\enddemo
     The theorem~\mythetheorem{8.1} is an immediate consequence of the
lemma~\mythelemma{8.1}. Indeed, the left hand side of the formula
\mythetag{8.6} defines a tensor of the type $(1,2)$ for any choice of 
the vector field $\bold X$, while the right hand side is the contraction 
of $\bold R$ and $\bold X$.\par
     The components of the curvature tensor given by the formula 
\mythetag{8.4} are enumerated by three lower indices and one upper index.
Upon lowering by means of the metric tensor the upper index is usually
written in the first position:
$$
\hskip -2em
R_{qrij}=\sum^2_{k=1}R^k_{rij}\,g_{kq}.
\mytag{8.12}
$$
The tensor of the type $(0,4)$ given by the formula \mythetag{8.12} is
denoted by the same letter $\bold R$. Another tensor is derived
from \mythetag{8.4} by raising the first lower index:
$$
\hskip -2em
R^{kq}_{ij}=\sum^2_{r=1}R^k_{rij}\,g^{rq}.
\mytag{8.13}
$$
The raised lower index is usually written as the second upper index.
The tensors of the type $(0,4)$ and $(2,2)$ with the components 
\mythetag{8.12} and \mythetag{8.13} are denoted by the same letter 
$\bold R$ and called the curvature tensors.
\mytheorem{8.2} The components of the curvature tensor 
$\bold R$ determined by the connection \mythetag{7.5} according
to the formula \mythetag{8.4} satisfy the following relationships:
\roster
\item $R^k_{rij}=-R^k_{rji}$;\vphantom{\vrule height 6pt depth 7pt}
\item $R_{qrij}=-R_{rqij}$;\vphantom{\vrule height 6pt depth 7pt}
\item $R_{qrij}=R_{ijqr}$;\vphantom{\vrule height 6pt depth 7pt}
\item $R^k_{rij}+R^k_{ijr}+R^k_{jri}=0$.
      \vphantom{\vrule height 6pt depth 7pt}
\endroster
\endproclaim
\demo{Proof} The first relationship is an immediate consequence of the
formula \mythetag{8.4} itself. When transposing the indices $i$ and $j$
the right hand side of \mythetag{8.4} changes the sign. Hence, we get
the identity \therosteritem{1} which means that the curvature tensor is
skew-symmetric with respect to the last pair of its lower indices.\par
     In order to prove the identity in the item \therosteritem{2}
we apply \mythetag{8.7} to the metric tensor. As a result we get the
following equality:
$$
(\nabla_i\nabla_j-\nabla_j\nabla_i)\,g_{qr}=
\sum^2_{k=1}\bigl(R^k_{qij}\,g_{kr}+R^k_{rij}\,g_{qk}\bigr).
$$
Taking into account \mythetag{8.12}, this equality can be rewritten as

$$
(\nabla_i\nabla_j-\nabla_j\nabla_i)\,g_{qr}=
R_{rqij}+R_{qrij}.
\mytag{8.14}
$$
Remember that due to the concordance of the metric and connection
the covariant derivatives of the metric tensor are equal to zero
(see formula \mythetag{7.1}). Hence, the left hand side of the equality
\mythetag{8.14} is equal to zero, and as a consequence we get the
identity from the item \therosteritem{2} of the theorem.\par
     Let's drop for a while the third item of the theorem and prove
the fourth item by means of the direct calculations on the base of 
the formula \mythetag{8.4}. Let's write the relationship \mythetag{8.4} 
and perform twice the cyclic transposition of the indices in it:
$i\to j\to r\to i$. As a result we get the following three equalities:
$$
\align
&R^k_{rij}=
\frac{\partial\Gamma^k_{jr}}{\partial u^i}
-\frac{\partial\Gamma^k_{ir}}{\partial u^j}+
\sum^2_{q=1}\Gamma^k_{iq}\Gamma^q_{jr}-
\sum^2_{q=1}\Gamma^k_{jq}\Gamma^q_{ir},\\
&R^k_{ijr}=
\frac{\partial\Gamma^k_{ri}}{\partial u^j}
-\frac{\partial\Gamma^k_{ji}}{\partial u^r}+
\sum^2_{q=1}\Gamma^k_{jq}\Gamma^q_{ri}-
\sum^2_{q=1}\Gamma^k_{rq}\Gamma^q_{ji},\\
&R^k_{jri}=
\frac{\partial\Gamma^k_{ij}}{\partial u^i}
-\frac{\partial\Gamma^k_{rj}}{\partial u^i}+
\sum^2_{q=1}\Gamma^k_{rq}\Gamma^q_{ij}-
\sum^2_{q=1}\Gamma^k_{iq}\Gamma^q_{rj}.
\endalign
$$
Let's add all the three above equalities and take into account
the symmetry of the Christoffer symbols with respect to their 
lower indices. It is easy to see that the sum in the right hand
side will be zero. This proves the item \therosteritem{4} of the 
theorem.\par
     The third item of the theorem follows from the first, the 
second, and the third items. In the left hand side of the equality
that we need to prove we have $R_{qrij}$. The simultaneous transposition
of the indices $q\leftrightarrow r$ and $i\leftrightarrow j$ does not
change this quantity, i\.\,e\. we have the equality
$$
\hskip -2em
R_{qrij}=R_{rqji}.
\mytag{8.15}
$$
This equality follows from the item \therosteritem{1} and the item
\therosteritem{2}. Let's apply the item \therosteritem{4} to the
quantities in both sides of the equality \mythetag{8.15}:
$$
\hskip -2em
\aligned
&R_{qrij}=-R_{qijr}-R_{qjri},\\
&R_{rqji}=-R_{rjiq}-R_{riqj}.
\endaligned
\mytag{8.16}
$$
Now let's perform the analogous manipulations with the quantity $R_{ijqr}$:
$$
\hskip -2em
R_{ijqr}=R_{jirq}.
\mytag{8.17}
$$
To each quantity in \mythetag{8.17} we apply the item \therosteritem{4}
of the theorem:
$$
\hskip -2em
\aligned
&R_{ijqr}=-R_{iqrj}-R_{irjq},\\
&R_{jirq}=-R_{jrqi}-R_{jqir}.
\endaligned
\mytag{8.18}
$$
Let's add the equalities \mythetag{8.16} and subtract from the sum the
equalities \mythetag{8.18}. It is easy to verify that due to the items
\therosteritem{1} and \therosteritem{2} of the theorem the right hand 
side of the resulting equality is zero. Then, using \mythetag{8.15} and
\mythetag{8.17}, we get
$$
2\,R_{qrij}-2\,R_{ijqr}=0.
$$
Dividing by $2$, now we get the identity that we needed to prove. Thus,
the theorem is completely proved.
\qed\enddemo
     The curvature tensor $\bold R$ given by its components \mythetag{8.4}
has the indices on both levels.  Therefore, we can consider the
contraction:
$$
\hskip -2em
R_{rj}=\sum^2_{k=1}R^k_{rkj}.
\mytag{8.19}
$$
The formula \mythetag{8.19} for $R_{rj}$ can be transformed as follows:
$$
R_{rj}=\sum^2_{i=1}\sum^2_{k=1}g^{ik}\,R_{irkj}.
$$
From this equality due to the symmetry $g^{ik}$ and due to the item
\therosteritem{4} of the theorem~\mythetheorem{8.2} we derive the 
symmetry of the tensor $R_{rj}$:
$$
R_{rj}=R_{jr}.
\mytag{8.20}
$$
The symmetric tensor of the type $(0,2)$ with the components
\mythetag{8.19} is called the {\it Ricci tensor}. It is denoted by the 
same letter $\bold R$ as the curvature tensor.\par
     Note that there are other two contractions of the curvature tensor.
However, these contractions do not produce new tensors:
$$
\xalignat 2
&\sum^2_{k=1}R^k_{krj}=0,
&&\sum^2_{k=1}R^k_{rik}=-R_{ri}.
\endxalignat
$$
Using the Ricci tensor, one can construct a scalar field $R$ by means
of the formula
$$
R=\sum^2_{r=1}\sum^2_{j=1}R_{rj}\,g^{rj}.
\mytag{8.21}
$$
The scalar $R(u^1,u^2)$ defined by the formula \mythetag{8.21} is called
the {\it scalar curvature} of a surface at the point with the coordinats
$u^1,u^2$. The scalar curvature is a result of total contraction of the
curvature tensor $\bold R$ given by the formula \mythetag{8.13}:
$$
R=\sum^2_{i=1}\sum^2_{j=1}R^{ij}_{ij}.
\mytag{8.22}
$$
The formula \mythetag{8.22} is easily derived from \mythetag{8.21}. Any other ways of contracting the curvature tensor do not give other scalars
essentially different from \mythetag{8.21}.\par
     In general, passing from the components of the curvature tensor
$R^{kr}_{ij}$ to the scalar curvature, we should lose a substantial part of
the information contained in the tensor $\bold R$: this means that we
replace $16$ quantities by the only one. However, due to the 
theorem~\mythetheorem{8.2} in two-dimensional case we do not lose the
information at all. Indeed, due to the theorem~\mythetheorem{8.2} the
components of the curvature tensor $R^{kr}_{ij}$ are skew-symmetric both
with respect to upper and lower indices. If $k=r$ or $i=j$, they do vanish.
Therefore, the only nonzero components are $R^{12}_{12}$, $R^{21}_{12}$,
$R^{12}_{21}$, $R^{21}_{21}$, and they satisfy the equalities
$R^{12}_{12}=R^{21}_{21}=-R^{21}_{12}=-R^{12}_{21}$. Hence, we get
$$
R=R^{12}_{12}+R^{21}_{21}=2\,R^{12}_{12}.
$$
Let's consider the tensor $\bold D$ with the following components:
$$
D^{kr}_{ij}=\frac{R}{2}\left(\delta^k_i\,\delta^r_j-
\delta^k_j\,\delta^r_i\right).
$$
The tensor $\bold D$ is also skew-symmetric with respect to upper and 
lower indices and $D^{12}_{12}=R^{12}_{12}$. Hence, these tensors do
coincide: $\bold D=\bold R$. In coordinates we have
$$
\hskip -2em
R^{kr}_{ij}=\frac{R}{2}\left(\delta^k_i\,\delta^r_j-
\delta^k_j\,\delta^r_i\right).
\mytag{8.23}
$$
By lowering the upper index $r$, from \mythetag{8.23} we derive
$$
\hskip -2em
R^k_{rij}=\frac{R}{2}\left(\delta^k_i\,g_{rj}-
\delta^k_j\,g_{ri}\right).
\mytag{8.24}
$$
The formula \mythetag{8.24} determines the components of the curvature tensor
on an arbitrary surface. For the Ricci tensor this formula yields
$$
\hskip -2em
R_{ij}=\frac{R}{2}\,g_{ij}.
\mytag{8.25}
$$
The Ricci tensor of an arbitrary surface is proportional to the metric
tensor.\par
\subhead A remark\endsubhead The curvature tensor determined by the
symmetric connection \mythetag{7.5} possesses another one (fifth) property
expressed by the identity
$$
\hskip -2em
\nabla_kR^q_{rij}+\nabla_iR^q_{rjk}+\nabla_jR^q_{rki}=0.
\mytag{8.26}
$$
The relationship \mythetag{8.23} is known as the {\it Bianchi identity}.
However, in the case of surfaces (in the dimension $2$) it appears to
be a trivial consequence from the item \therosteritem{1} of the
theorem~\mythetheorem{8.2}. Therefore, we do not give it here.
\head
\S~\mysection{9} Gauss equation and Peterson-Codazzi equation.
\endhead
\rightheadtext{\S~9. Gauss equation and Peterson-Codazzi equation.}
     Let's consider the Veingarten's derivational formulas \mythetag{4.11}.
They can be treated as a system of nine vectorial equations with respect 
to three vector-functions $\bold E_1(u^1,u^2)$, $\bold E_2(u^1,u^2)$, and
$\bold n(u^1,u^2)$. So, the number of the equations is greater than the
number functions. Such systems are said to be {\it overdetermined}. 
Overdetermined systems are somewhat superfluous. One usually can derive
new equations of the same or lower order from them. Such equations are
called {\it differential consequences\/} or {\it compatibility
conditions\/} of the original equations.\par
     As an example we consider the following system of two partial
differential equations with respect to the function $f=f(x,y)$:
$$
\xalignat 2
&\hskip -2em
\frac{\partial f}{\partial x}=a(x,y),
&&\frac{\partial f}{\partial y}=b(x,y).
\mytag{9.1}
\endxalignat
$$
Let's differentiate the first equation \mythetag{9.1} with respect to $y$
and the second equation with respect to $x$. Then we subtract one from
another:
$$
\hskip -2em
\frac{\partial a}{\partial y}-
\frac{\partial b}{\partial x}=0.
\mytag{9.2}
$$
The equation \mythetag{9.2} is a compatibility condition for the equations
\mythetag{9.1}. It is a necessary condition for the existence of the
function satisfying the equations \mythetag{9.1}.\par
     Similarly, one can derive the compatibility conditions for the system
of Veingarten's derivational equations \mythetag{4.11}. Let's write the first
of them as
$$
\frac{\partial\bold E_k}{\partial u^j}=
\sum^2_{q=1}\Gamma^q_{jk}\cdot\bold E_q+b_{jk}\cdot\bold n.
\mytag{9.3}
$$
Then we differentiate \mythetag{9.3} with respect to $u^i$ and express 
the derivatives $\partial\bold E_k/\partial u^i$ and $\partial\bold
n/\partial u^i$ arising therein by means of the derivational formulas
\mythetag{4.11}:
$$
\aligned
&\frac{\partial\bold E_k}{\partial u^i\,\partial u^j}=
\left(\frac{\partial b_{jk}}{\partial u^i}+
\shave{\sum^2_{q=1}}\Gamma^q_{jk}\,b_{iq}\right)\cdot\bold n+\\
&\qquad+\sum^2_{q=1}\left(\frac{\partial\Gamma^q_{jk}}{\partial u^i}+
\shave{\sum^2_{s=1}}\Gamma^s_{jk}\,\Gamma^q_{is}-b_{jk}\,b^q_i
\right)\cdot\bold E_q.\hskip -2em
\endaligned
\mytag{9.4}
$$
Let's transpose indices $i$ and $j$ in the formula \mythetag{9.4}. The value
of the second order mixed partial derivative does not depend on the order
of differentiation. Therefore, the value of the left hand side of 
\mythetag{9.4} does not change under the transposition of indices $i$ and
$j$. Let's subtract from \mythetag{9.4} the relationship obtained by
transposing the indices. As a result we get
$$
\gather
\sum^2_{q=1}\left(\frac{\partial\Gamma^q_{jk}}{\partial u^i}
-\frac{\partial\Gamma^q_{ik}}{\partial u^j}+
\shave{\sum^2_{s=1}}\Gamma^s_{jk}\,\Gamma^q_{is}-
\shave{\sum^2_{s=1}}\Gamma^s_{ik}\,\Gamma^q_{js}+
b_{ik}\,b^q_j-b_{jk}\,b^q_i\right)\cdot\bold E_q+\\
\vspace{1.5ex}
+\left(\frac{\partial b_{jk}}{\partial u^i}+
\shave{\sum^2_{q=1}}\Gamma^q_{jk}\,b_{iq}-
\frac{\partial b_{ik}}{\partial u^j}-
\shave{\sum^2_{q=1}}\Gamma^q_{ik}\,b_{jq}\right)\cdot\bold n=0.
\endgather
$$
The vectors $\bold E_1$, $\bold E_2$, and $\bold n$ composing the moving
frame are linearly independent.  Therefore the above equality can be
broken into two separate equalities
$$
\aligned
 &\frac{\partial\Gamma^q_{jk}}{\partial u^i}
  -\frac{\partial\Gamma^q_{ik}}{\partial u^j}+
   \sum^2_{s=1}\Gamma^s_{jk}\,\Gamma^q_{is}-
   \sum^2_{s=1}\Gamma^s_{ik}\,\Gamma^q_{js}=
 b_{jk}\,b^q_i-b_{ik}\,b^q_j,\\
\vspace{1.5ex}
&\frac{\partial b_{jk}}{\partial u^i}-
 \shave{\sum^2_{q=1}}\Gamma^q_{ik}\,b_{jq}=
 \frac{\partial b_{ik}}{\partial u^j}-
 \shave{\sum^2_{q=1}}\Gamma^q_{jk}\,b_{iq}.
\endaligned
$$
Note that the left hand side of the first of these relationships 
coincides with the formula for the components of the curvature
tensor (see \mythetag{8.4}). Therefore, we can rewrite the first
relationship as follows:
$$
\hskip -2em
R^q_{kij}=b_{jk}\,b^q_i-b_{ik}\,b^q_j.
\mytag{9.5}
$$
The second relationship can also be simplified:
$$
\hskip -2em
\nabla_ib_{jk}=\nabla_jb_{ik}.
\mytag{9.6}
$$
It is easy to verify this fact immediately by transforming \mythetag{9.6}
back to the initial form applying the formula \mythetag{6.1}.\par
     The equations \mythetag{9.5} and \mythetag{9.6} are differential 
consequences of the Veingarten's derivational formulas \mythetag{4.11}.
The first of them is known as the {\it Gauss equation\/} and the second one
is known as the {\it Peterson-Codazzi equation}.\par
     The tensorial Gauss equation \mythetag{9.5} contains 16 separate
equalities. However, due to the relationship \mythetag{8.24} not all of them
are independent. In order to simplify \mythetag{9.5} let's raise the index $k$ in it. As a result we get
$$
\hskip -2em
\frac{R}{2}\left(\delta^q_i\,\delta^k_j-
\delta^q_j\,\delta^k_i\right)=
b^q_i\,b^k_j-b^q_j\,b^k_i.
\mytag{9.7}
$$
The expression in right hand side of \mythetag{9.7} is skew-symmetric 
both with respect to upper and lower pairs of indices and each index 
in \mythetag{9.7} runs over only two values. Therefore the right hand
side of the equation \mythetag{9.7} can be transformed as
$$
\hskip -2em
b^q_i\,b^k_j-b^q_j\,b^k_i=B\,\left(\delta^q_i\,\delta^k_j-
\delta^q_j\,\delta^k_i\right).
\mytag{9.8}
$$
Substituting $q=1$, $k=2$, $i=1$, $j=2$ into \mythetag{9.8}, for $B$ in \mythetag{9.8} we get
$$
B=b^1_1\,b^2_2-
b^1_2\,b^2_1=\det(b^k_i)=K,
$$
where $K$ is the Gaussian curvature of a surface (see formula
\mythetag{5.12}). The above considerations show that the Gauss equation
\mythetag{9.5} is equivalent to exactly one scalar equation which is written
as follows:
$$
\hskip -2em
R=2\,K.
\mytag{9.9}
$$
This equation relates the scalar and Gaussian curvatures of a surface. It
is also called the {\it Gauss equation}.\par
\newpage
\topmatter
\title\chapter{5}
Curves on surfaces
\endtitle
\endtopmatter
\chapternum=5
\document
\setfirstpage
\head
\S~\mysection{1} Parametric equations of a curve on a surface.
\endhead
\leftheadtext{CHAPTER \uppercase\expandafter{\romannumeral 5}.
CURVES ON SURFACES.}
\rightheadtext{\S~1. Parametric equations of a curve on a surface.}
     Let $\bold r(t)$ be the vectorial-parametric equation 
of a differentiable curve all points of which lie on some 
differentiable surface. Suppose that a fragment $D$ containing 
the points of the curve is charted, i\.\,e\. it is equipped with
curvilinear coordinates $u^1,\,u^2$. This means that there is a
bijective mapping $\bold u\!:\,D\to U$ that maps the points of 
the curve to some domain $U\subset\Bbb R^2$. The curve in the 
chart $U$ is represented not by three, but by two functions of
the parameter $t$:
$$
\hskip -2em
\cases
u^1=u^1(t),\\
u^2=u^2(t)
\endcases
\mytag{1.1}
$$
(compare with the formulas \mythetagchapter{1.14}{4} from Chapter
\uppercase\expandafter{\romannumeral 4}). The inverse mapping
$\bold u^{-1}$ is represented by the vector-function
$$
\hskip -2em
\bold r=\bold r(u^1,u^2).
\mytag{1.2}
$$
It determines the radius-vectors of the points of the surface. 
Therefore, we have
$$
\hskip -2em
\bold r(t)=\bold r(u^1(t),u^2(t)).
\mytag{1.3}
$$
Differentiating \mythetag{1.3} with respect to $t$, we obtain the
vector $\boldsymbol\tau$:
$$
\hskip -2em
\boldsymbol\tau(t)=\sum^2_{i=1}
\dot u^i\cdot\bold E_i(u^1(t),u^2(t)).
\mytag{1.4}
$$
This is the tangent vector of the curve (compare with the formulas
\mythetagchapter{1.15}{4}
in Chapter \uppercase\expandafter{\romannumeral 4}). The formula
\mythetag{1.4} shows that the vector $\boldsymbol\tau$ lies in the 
tangent plane of the surface. This is the consequence of the fact 
that the curve in whole lies on the surface.\par
     Under a change of curvilinear coordinates on the surface the
derivatives $\dot u^i$ are transformed as the components of a tensor
of the type $(1,0)$. They determine the inner (two-dimensional)
representation of the vector $\tau$ in the chart. The formula
\mythetag{1.4} is used to pass from inner to outer (tree-dimensional)
representation of this vector. {\bf Our main goal} in this chapter
is to describe the geometry of curves lying on a surface in terms
of its two-dimensional representation in the chart.\par
     The length integral is an important object in the theory of
curves, see formula \mythetagchapter{2.3}{1} in Chapter
\uppercase\expandafter{\romannumeral 1}. Substituting \mythetag{1.4}
into this formula, we get
$$
\hskip -2em
L=\int\limits^{\,b}_{a\,}\sqrt{\shave{\sum^2_{i=1}\sum^2_{j=1}}
g_{ij}\,\dot u^i\,\dot u^j\,\,}\,dt.
\mytag{1.5}
$$
The expression under integration in \mythetag{1.5} is the length of the
vector $\boldsymbol\tau$ in its inner representation. If $s=s(t)$ is
the natural parameter of the curve, then, denoting $du^i=\dot u^i\,dt$, 
we can write the following formula:
$$
\hskip -2em
ds^2=\sum^2_{i=1}\sum^2_{j=1}g_{ij}\,du^i\,du^j.
\mytag{1.6}
$$
The formula \mythetag{1.6} approves the title {\tencyr\char '074}the
first quadratic form{\tencyr\char '076} for the metric tensor. Indeed,
the square of the length differential $ds^2$ is a quadratic form of
differentials of the coordinate functions in the chart. If $t=s$ is
the natural parameter of the curve, then there is the equality
$$
\hskip -2em
\sum^2_{i=1}\sum^2_{j=1}g_{ij}\,\dot u^i\,\dot u^j=1
\mytag{1.7}
$$
that expresses the fact that the length of the tangent vector $\boldsymbol
\tau$ of a curve in the natural parametrization is equal to unity (see
\S\,2 in Chapter \uppercase\expandafter{\romannumeral 1}).
\head
\S~\mysection{2} Geodesic and normal curvatures of a curve.
\endhead
\rightheadtext{\S~2. Geodesic and normal curvatures of a curve.}
     Let $t=s$ be the natural parameter of a parametric curve given by
the equations \mythetag{1.1} in curvilinear coordinates on some surface.
Let's differentiate the tangent vector $\boldsymbol\tau(s)$ of this
curve \mythetag{1.4} with respect to the parameter $s$. The derivative 
of the left hand side of \mythetag{1.4} is given by the Frenet formulas
\mythetagchapter{3.8}{1} from Chapter 
\uppercase\expandafter{\romannumeral 1}:
$$
k\cdot\bold n_{\,\text{curv}}=
\sum^2_{k=1}\ddot u^k\cdot\bold E_k+
\sum^2_{i=1}\sum^2_{j=1}\dot u^i\cdot
\frac{\partial\bold E_i}{\partial u^j}\cdot\dot u^j.
\mytag{2.1}
$$
By $\bold n_{\,\text{curv}}$ we denote the unit normal vector of the curve
in order to distinguish it from the unit normal vector $\bold n$ of the
surface. For to calculate the derivatives $\partial\bold E_i/\partial u^j$
we apply the Veingarten's derivational formulas \mythetag{4.11}:
$$
k\cdot\bold n_{\,\text{curv}}=
\sum^2_{k=1}\left(\ddot u^k+\shave{\sum^2_{i=1}\sum^2_{j=1}}
\Gamma^k_{ji}\,\dot u^i\,\dot u^j\right)\cdot\bold E_k+
\left(\,\shave{\sum^2_{i=1}\sum^2_{j=1}}b_{ij}\,\dot u^i\,
\dot u^j\right)\cdot\bold n.\quad
\mytag{2.2}
$$
Let's denote by $k_{\,\text{norm}}$ the coefficient of the vector $\bold n$
in the formula \mythetag{2.2}. This quantity is called the {\it normal
curvature} of a curve:
$$
\hskip -2em
k_{\,\text{norm}}=\sum^2_{i=1}\sum^2_{j=1}b_{ij}\,\dot u^i\,\dot u^j.
\mytag{2.3}
$$
In contrast to the curvature $k$, which is always a non-negative quantity,
the normal curvature of a curve \mythetag{2.3} can be either positive, or
zero, or negative. Taking into account \mythetag{2.3}, we can rewrite the
relationship \mythetag{2.2} itself as follows:
$$
\hskip -2em
k\cdot\bold n_{\,\text{curv}}-k_{\,\text{norm}}\cdot\bold n=
\sum^2_{k=1}\left(\ddot u^k+\shave{\sum^2_{i=1}\sum^2_{j=1}}
\Gamma^k_{ji}\,\dot u^i\,\dot u^j\right)\cdot\bold E_k.
\hskip -2em
\mytag{2.4}
$$
The vector in the right hand side of \mythetag{2.4} is a linear
combination of the vectors $\bold E_1$ and $\bold E_2$ that compose
a basis in the tangent plane. Therefore, this vector lies in the 
tangent plane. Its length is called the {\it geodesic curvature\/} 
of a curve:
$$
\hskip -2em
k_{\,\text{geod}}=\left|\,
\shave{\sum^2_{k=1}}\left(\ddot u^k+\shave{\sum^2_{i=1}\sum^2_{j=1}}
\Gamma^k_{ji}\,\dot u^i\,\dot u^j\right)\cdot\bold E_k
\right|.
\mytag{2.5}
$$
Due to the formula \mythetag{2.5} the geodesic curvature of a curve 
is always non-negative. If $k_{\,\text{geod}}\neq 0$, then, taking into
account the relationship \mythetag{2.5}, one can define the unit vector 
$\bold n_{\,\text{inner}}$ and rewrite the formula \mythetag{2.4} as follows:
$$
\hskip -2em
k\cdot\bold n_{\,\text{curv}}-k_{\,\text{norm}}\cdot\bold n=
k_{\,\text{geod}}\cdot\bold n_{\,\text{inner}}.
\mytag{2.6}
$$
The unit vector $\bold n_{\,\text{inner}}$ in the formula \mythetag{2.6}
is called the {\it inner normal vector\/} of a curve on a surface.\par
\parshape 18 
180pt 180pt 180pt 180pt 180pt 180pt 180pt 180pt 
180pt 180pt 180pt 180pt 180pt 180pt 180pt 180pt 
180pt 180pt 180pt 180pt 180pt 180pt 180pt 180pt 
180pt 180pt 180pt 180pt 180pt 180pt 180pt 180pt 
180pt 180pt 0pt 360pt
     Due to \mythetag{2.6} the vector $\bold n_{\,\text{inner}}$ is a 
linear combination of the vectors $\bold n_{\,\text{curv}}$ and $\bold n$
which \vadjust{\vskip 5pt\hbox to 0pt{\kern -5pt
\includegraphics{ris13.eps}\hss}\vskip -5pt}are perpendicular
to the unit vector $\boldsymbol\tau$ lying in the tangent plane. Hence, 
$\bold n_{\,\text{inner}}\perp\boldsymbol\tau$. On the other hand, being
a linear combination of the vectors $\bold E_1$ and $\bold E_2$, the
vector $\bold n_{\,\text{inner}}$ itself lies in the tangent plane. 
Therefore, it is determined up to the sign:
$$
\hskip -2em
\bold n_{\,\text{inner}}=\pm[\boldsymbol\tau,\,\bold n].
\mytag{2.7}
$$
Using the relationship \mythetag{2.7}, sometimes one can extend the definition of the vec\-tor $\bold n_{\,\text{inner}}$ even to those 
points of the curve where $k_{\,\text{geod}}=0$.\par
     Let's move the term $k_{\,\text{norm}}\cdot\bold n$ to the right
hand side of the formula \mythetag{2.6}. Then this formula is rewritten
as follows:
$$
\hskip -2em
k\cdot\bold n_{\,\text{curv}}=k_{\,\text{geod}}\cdot\bold n_{\,\text{inner}}
+k_{\,\text{norm}}\cdot\bold n.
\mytag{2.8}
$$
The relationship \mythetag{2.8} can be treated as an expansion of the vector
$k\cdot\bold n_{\,\text{curv}}$ as a sum of two mutually prpendicular 
components. Hence, we have
$$
\hskip -2em
k^2=(k_{\,\text{geod}})^2+(k_{\,\text{norm}})^2.
\mytag{2.9}
$$\par
     The formula \mythetag{2.3} determines the value of the normal 
curvature of a curve in the natural parametrization $t=s$. Let's 
rewrite it as  follows:
$$
k_{\,\text{norm}}=\frac{\dsize\sum^2_{i=1}\sum^2_{j=1}b_{ij}\,
\dot u^i\,\dot u^j}{\dsize\sum^2_{i=1}\sum^2_{j=1}g_{ij}\,
\dot u^i\,\dot u^j}.
\mytag{2.10}
$$
In the natural parametrization the formula \mythetag{2.10} coincides with
\mythetag{2.3} because of \mythetag{1.7}. When passing to an arbitrary 
parametrization all derivatives $\dot u^i$ are multiplied by the same
factor. Indeed, we have
$$
\hskip -2em
\frac{du^i}{dt}=\frac{du^i}{ds}\,\frac{ds}{dt}.
\mytag{2.11}
$$
But the right hand side of \mythetag{2.10} is insensitive to such a change
of $\dot u^i$. Therefore, \mythetag{2.10} is a valid formula for the normal
curvature in any parametrization.\par
     The formula \mythetag{2.10} shows that the normal curvature is a very
rough characte\-ristic of a curve. It is determined only by the direction
of its tangent vector $\boldsymbol\tau$ in the tangent plane. The components of the matrices $g_{ij}$ and $b_{ij}$ characterize not the
curve, but the point of the surface through which this curve passes.
\par
     Let $\bold a$ be some vector tangent to the surface. In curvilinear coordinates $u^1,\,u^2$ it is given by two numbers $a^1,\,a^2$ --- they 
are the coefficients in its expansion in the basis of two frame vectors
$\bold E_1$ and $\bold E_2$. Let's consider the value of the second
quadratic form of the surface on this vector:
$$
\hskip -2em
\bold b(\bold a,\bold a)=
\sum^2_{i=1}\sum^2_{j=1}b_{ij}\,a^i\,a^j.
\mytag{2.12}
$$
\mydefinition{2.1} The direction given by a nonzero vector
$\bold a$ in the tangent plane is called an {\it asymptotic direction\/}
if the value of the second quadratic form \mythetag{2.12} on such vector
is equal to zero.
\enddefinition
    Note that asymptotic directions do exist only at those points of a surface where the second quadratic form is indefinite in sign or 
degenerate. In the first case the Gaussian curvature is negative:
$K<0$, in the second case it is equal to zero: $K=0$. At those points
where $K>0$ there are no asymptotic directions.\par
\mydefinition{2.2} A curve on a surface whose tangent vector 
$\boldsymbol\tau$ at all its points lies in an asymptotic direction
is called an {\it asymptotic line}.
\enddefinition
Due to \mythetag{2.12} the equation of an asymptotic line has the following form:
$$
\hskip -2em
\bold b(\boldsymbol\tau,\boldsymbol\tau)=
\sum^2_{i=1}\sum^2_{j=1}b_{ij}\,\dot u^i\,\dot u^j=0.
\mytag{2.13}
$$
Comparing \mythetag{2.13} and \mythetag{2.10}, we see that asymptotic lines
are the lines with zero normal curvature: $k_{\,\text{norm}}=0$. On the
surfaces with negative Gaussian curvature $K<0$ at each point there are
two asymptotic directions. Therefore, on such surfaces always there are
two families of asymptotic lines, they compose the {\it asymptotic
network} of such a surface. On any surface of the negative Gaussian
curvature there exists a curvilinear coordinate system $u^1,\,u^2$
whose coordinate network coincides with the asymptotic network of this
surface. However, we shall not prove this fact here.\par
     The {\it curvature lines\/} are defined by analogy with the 
asymptotic lines. These are the curves on a surface whose tangent vector
lies in a principal direction at each point (see formulas 
\mythetagchapter{5.14}{4} and \mythetagchapter{5.15}{4} in \S\,5 of Chapter
\uppercase\expandafter{\romannumeral 4}). The curvature lines do exist on any surface, there are no restrictions for the Gaussian curvature of a
surface in this case.\par
\mydefinition{2.3} A {\it geodesic line} on a surface is a curve
whose geodesic curvature is identically equal to zero:
$k_{\,\text{geod}}=0$.
\enddefinition
    From the Frenet formula $\dot{\boldsymbol\tau}=k\cdot\bold
n_{\,\text{curv}}$ and from the relationship \mythetag{2.8} for the geodesic lines 
we derive the following equality:
$$
\hskip -2em
\frac{d\boldsymbol\tau}{ds}=k_{\,\text{norm}}\cdot\bold n.
\mytag{2.14}
$$
In the other words, the derivative of the unit normal vector on a geodesic
line is directed along the unit normal vector of a surface. This is the
external description of geodesic lines. The inner description is derived
from the formula \mythetag{2.5}:
$$
\hskip -2em
\ddot u^k+\sum^2_{i=1}\sum^2_{j=1}\Gamma^k_{ji}\,
\dot u^i\,\dot u^j=0.
\mytag{2.15}
$$
The equations \mythetag{2.15} are the {\it differential equations of geodesic
lines in natural parametrization}. One can pass from the natural
parametrization to an arbitrary one by means of the formula \mythetag{2.11}.\par
\head
\S~\mysection{3} Extremal property of geodesic lines.
\endhead
\rightheadtext{\S~3. Extremal property of geodesic lines.}
     Let's compare the equations of geodesic lines \mythetag{2.15} 
with the equations of straight lines in curvilinear coordinates
\mythetagchapter{8.17}{3} which we have derived in 
Chapter~\uppercase\expandafter{\romannumeral 3}. These equations have 
the similar structure. They differ only in the ranges of indices: in the
case of geodesic lines on a surface they run over two numbers instead of
three numbers in the case of straight lines. Therefore, geodesic lines are
natural analogs of the straight lines in the inner geometry of surfaces.
The following theorem strengthens this analogy.
\mytheorem{3.1} A geodesic line connecting two given points $A$ and
$B$ on a surface has the {\it extremal} length in the class of curves
connecting these two points.
\endproclaim
    It is known that in the Euclidean space $\Bbb E$ the shortest path
from a point $A$ to another point $B$ is the segment of straight line
connecting these points. The theorem~\mythetheorem{3.1} proclaims a very
similar proposition for geodesic lines on a surface. Remember that real
functions can have local maxima and minima --- they are called {\it
extrema}. Apart from maxima and minima, there are also {\it conditional
extrema\/} (saddle points), for example, the point $x=0$ for the function
$y=x^3$. All those points are united by the common property --- the linear
part of the function increment at any such point is equal to zero:
$f(x_0+h)=f(x_0)+O(h^2)$.\par
     In the case of a geodesic line connecting the points $A$ and $B$ on a 
surface we should slightly deform (variate) this line keeping it to be
a line on the surface connecting the same two points $A$ and $B$. The
deformed line is not a geodesic line. Its length differs from the length
of the original line. The condition of the {\it extremal length\/} in the
theorem~\mythetheorem{3.1} means that the linear part of the length
increment is equal to zero.\par
     Let's specify the method of deforming the curve. For the sake of
simplicity assume that the points $A$ and $B$, and the geodesic line
connecting these points in whole lie within some charted fragment $D$
of the surface. Then this geodesic line is given by two functions
\mythetag{1.1}. Let's increment by one the number of arguments in these
functions. Then we shall assume that these functions are sufficiently 
many times differentiable with respect to all their arguments:
$$
\hskip -2em
\cases
u^1=u^1(t,h),\\
u^2=u^2(t,h).
\endcases
\mytag{3.1}
$$
For each fixed $h$ in \mythetag{3.1} we have the functions of the parameter
$t$, they define a curve on the surface. Changing the parameter $h$, we
deform the curve so that in the process of this deformation its point are
always on the surface. The differentiability of the functions \mythetag{3.1}
guarantees that small deformations of the curve correspond to small changes of the parameter $h$.\par
     Let's impose to the functions \mythetag{3.1} a series of restrictions
which are easy to satisfy. Assume that the length of the initial geodesic
line is equal to $a$ and let the parameter $t$ run over the segment
$[0,\,a]$. Let
$$
\xalignat 2
&\hskip -2em
u^k(0,h)=u^k(0,0),
&&u^k(a,h)=u^k(a,0).
\mytag{3.2}
\endxalignat
$$
The condition \mythetag{3.2} means that under a change of the parameter $h$
the initial point $A$ and the ending point $B$ of the curve do not move.\par
     For the sake of brevity let's denote the partial derivatives of the functions $u^i(t,h)$ with respect to $t$ by setting the dot. Then the
quantities $\dot u^i=\partial u^i/\partial t$ determine the inner
representation of the tangent vector to the curve.\par
     Assume that the initial line correspond to the value $h=0$ of the
parameter $h$. Assume also that for $h=0$ the parameter $t$ coincides
with the natural parameter of the geodesic line. Then for $h=0$ the
functions \mythetag{3.1} satisfy the equations \mythetag{1.7} and \mythetag{2.15}
simultaneously. For $h\neq 0$ the parameter $t$ should not coincide with
the natural parameter on the deformed curve, and the deformed curve itself
should not be a geodesic line in this case.\par
    Let's calculate the lengths of the deformed curves. It is the function
of the parameter $h$ determined by the length integral of the form
\mythetag{1.5}:
$$
\hskip -2em
L(h)=\int\limits^{\,a}_{0\,}\sqrt{\shave{\sum^2_{i=1}\sum^2_{j=1}}
g_{ij}\,\dot u^i\,\dot u^j\,\,}\,dt.
\mytag{3.3}
$$
For $h=0$ we have $L(0)=a$. The proposition of the
theorem~\mythetheorem{3.1} on the extremity of the length now is formulated
as $L(h)=a+O(h^2)$ or, equivalently, as 
$$
\hskip -2em
\frac{dL(h)}{dh}\,
\vrule height 14pt depth 10pt
\lower 10 pt \hbox{$\ssize\,h=0$}
=0.
\mytag{3.4}
$$
\demo{Proof of the theorem~\mythetheorem{3.1}} Let's prove the 
equality \mythetag{3.4} for the length integral \mythetag{3.3} under 
the deformations of the curve described just above. Denote by
$\lambda(t,h)$ the expression under the square root in the formula
\mythetag{3.3}. Then by direct differentiation of \mythetag{3.3} we 
obtain
$$
\hskip -2em
\frac{dL(h)}{dh}=\int\limits^{\,a}_{0\,}
\frac{\partial\lambda/\partial h}
{2\,\sqrt{\lambda}}\,dt.
\mytag{3.5}
$$
Let's calculate the derivative in the numerator of the fraction
\mythetag{3.5}:
$$
\align
&\frac{\partial\lambda}{\partial h}=
\frac{\partial}{\partial h}\left(\,\shave{\sum^2_{i=1}
\sum^2_{j=1}} g_{ij}\,\dot u^i\,\dot u^j\right)=
\sum^2_{i=1}\sum^2_{j=1}\sum^2_{k=1}
\frac{\partial g_{ij}}{\partial u^k}
\frac{\partial u^k}{\partial h}\,\dot u^i\,\dot u^j+\\
\vspace{1ex}
&+\sum^2_{k=1}\sum^2_{i=1}\sum^2_{j=1}g_{ij}\,
\frac{\partial(\dot u^i\,\dot u^j)}{\partial\dot u^k}\,
\frac{\partial\dot u^k}{\partial h}=
\sum^2_{i=1}\sum^2_{j=1}\sum^2_{k=1}
\frac{\partial g_{ij}}{\partial u^k}
\frac{\partial u^k}{\partial h}\,\dot u^i\,\dot u^j+\\
\vspace{1ex}
&+\sum^2_{k=1}\sum^2_{i=1}\sum^2_{j=1}g_{ij}\,
\delta^i_k\,\dot u^j\,\frac{\partial\dot u^k}{\partial h}+
\sum^2_{k=1}\sum^2_{i=1}\sum^2_{j=1}g_{ij}\,
\dot u^i\,\delta^j_k\,\frac{\partial\dot u^k}{\partial h}.
\endalign
$$
Due to the Kronecker symbols $\delta^i_k$ and $\delta^j_k$ in the
above expression we can perform explicitly the summation over $k$ 
in the last two terms.  Moreover, due to the symmetry of $g_{ij}$ 
they are equal to each other:
$$
\frac{\partial\lambda}{\partial h}=
\sum^2_{i=1}\sum^2_{j=1}\sum^2_{k=1}
\frac{\partial g_{ij}}{\partial u^k}
\frac{\partial u^k}{\partial h}\,\dot u^i\,\dot u^j+
2\,\sum^2_{i=1}\sum^2_{j=1}g_{ij}\,
\dot u^i\,\frac{\partial\dot u^j}{\partial h}.
$$
Substituting this expression into \mythetag{3.5}, we get two integrals:
$$
\align
&\hskip -2em
I_1=\sum^2_{i=1}\sum^2_{j=1}\sum^2_{k=1}
\int\limits^{\,a}_{0\,}
\frac{\partial g_{ij}}{\partial u^k}\,
\frac{\dot u^i\,\dot u^j}{2\,\sqrt{\lambda}}\,
\frac{\partial u^k}{\partial h}\,
\,dt,
\mytag{3.6}\\
&\hskip -2em
I_2=\sum^2_{i=1}\sum^2_{j=1}
\int\limits^{\,a}_{0\,}
\frac{g_{ik}\,\dot u^i}{\sqrt{\lambda}}\,
\frac{\partial\dot u^k}{\partial h}\,dt.
\mytag{3.7}
\endalign
$$
The integral \mythetag{3.7} contain the second order mixed partial
derivatives of \mythetag{3.1}:
$$
\frac{\partial\dot u^k}{\partial h}=
\frac{\partial^2 u^k}{\partial t\,\partial h}.
$$
In order to exclude such  derivatives we integrate \mythetag{3.7} by
parts:
$$
\int\limits^{\,a}_{0\,}
\frac{g_{ik}\,\dot u^i}{\sqrt{\lambda}}\,
\frac{\partial\dot u^k}{\partial h}\,dt=
\frac{g_{ik}\,\dot u^i}{\sqrt{\lambda}}\,
\frac{\partial u^k}{\partial h}\,
{\vrule height 18.5pt depth 14pt}^{\,a}_{\,0}
-\int\limits^{\,a}_{0\,}
\frac{\partial}{\partial t}\left(
\frac{g_{ik}\,\dot u^i}{\sqrt{\lambda}}\right)\,
\frac{\partial u^k}{\partial h}\,dt.
$$
Let's differentiate the equalities \mythetag{3.2} with respect to $h$. As a 
result we find that the derivatives $\partial u^k/\partial h$ vanish at
the ends of the integration segment over $t$. This means that non-integral
terms in the above formula do vanish. Hence, for the integral $I_2$ in
\mythetag{3.7} we obtain
$$
I_2=-\sum^2_{i=1}\sum^2_{k=1}
\int\limits^{\,a}_{0\,}
\frac{\partial}{\partial t}\left(
\frac{g_{ik}\,\dot u^i}{\sqrt{\lambda}}\right)\,
\frac{\partial u^k}{\partial h}\,dt.
\mytag{3.8}
$$
Now let's add the integrals $I_1$ and $I_2$ from \mythetag{3.6} and
\mythetag{3.8}. As a result for the derivative $dL/dh$ in \mythetag{3.5}
we derive the following equality:
$$
\frac{dL(h)}{dh}=
\sum^2_{i=1}\sum^2_{k=1}
\int\limits^{\,a}_{0\,}
\left(\,\shave{\sum^2_{j=1}}\frac{\partial g_{ij}}
{\partial u^k}\,\frac{\dot u^i\,\dot u^j}{2\,\sqrt{\lambda}}-
\frac{\partial}{\partial t}\left(
\frac{g_{ik}\,\dot u^i}{\sqrt{\lambda}}\right)\right)
\frac{\partial u^k}{\partial h}\,dt.
$$
In this equality the only derivatives with respect to the parameter 
$h$ are $\partial u^k/\partial h$. For their values at $h=0$ we introduce
the following notations:
$$
\hskip -2em
\delta u^k=\frac{\partial u^k}{\partial h}\,
\vrule height 14pt depth 10pt
\lower 10 pt \hbox{$\ssize\,h=0$}
\mytag{3.9}
$$
The quantities $\delta u^k=\delta u^k(t)$ in \mythetag{3.9} are called the
{\it variations of the coordinates} on the initial curve. Note that under
a change of curvilinear coordinates these quantities are transformed as 
the components of a vector (although this fact does not matter for proving
the theorem).\par
     Let's substitute $h=0$ into the above formula for the derivative
$dL/dh$. When substituted, the quantity $\lambda$ in the denominators 
of the fractions becomes equal to unity: $\lambda(t,0)=1$. This fact
follows from \mythetag{1.7} since $t$ coincides with the natural parameter
on the initial geodesic line. Then
$$
\frac{dL(h)}{dh}\,
\vrule height 14pt depth 10pt
\lower 10 pt \hbox{$\ssize\,h=0$}=
\sum^2_{i=1}\sum^2_{k=1}
\int\limits^{\,a}_{0\,}
\left(\,\shave{\sum^2_{j=1}}\frac{\partial g_{ij}}
{\partial u^k}\,\frac{\dot u^i\,\dot u^j}{2}-
\frac{d(g_{ik}\,\dot u^i)}{dt}\right)
\delta u^k\,dt.
$$
Since the above equality does not depend on $h$ any more, we replace the
partial derivative with respect to $t$ by $d/dt$. All of the further calculations in the right hand side are for the geodesic line where
$t$ is the natural parameter.\par
    Let's move the sums over $i$ and $k$ under the integration and let's
calculate the coefficients of $\delta u^k$ denoting these coefficients 
by $U_k$:
$$
\hskip -2em
\aligned
U_k&=\sum^2_{i=1}\left(\,\shave{\sum^2_{j=1}}
\frac{\partial g_{ij}}{\partial u^k}\,
\frac{\dot u^i\,\dot u^j}{2}-
\frac{d(g_{ik}\,\dot u^i)}{dt}\right)=\\
&=\sum^2_{i=1}\sum^2_{j=1}\left(\frac{1}{2}
\frac{\partial g_{ij}}{\partial u^k}\,
-\frac{\partial g_{ik}}{\partial u^j}\right)
\dot u^i\,\dot u^j-\sum^2_{i=1}g_{ik}\,\ddot u^i.
\endaligned
\mytag{3.10}
$$
Due to the symmetry of $\dot u^i\,\dot u^j$ the second term within
round brackets in the formula \mythetag{3.10} can be broken into two
terms. This yields
$$
U_k=\sum^2_{i=1}\sum^2_{j=1}\frac{1}{2}
\left(
\frac{\partial g_{ij}}{\partial u^k}\,
-\frac{\partial g_{ik}}{\partial u^j}
-\frac{\partial g_{jk}}{\partial u^i}
\right)
\dot u^i\,\dot u^j-\sum^2_{i=1}g_{ik}\,\ddot u^i.
$$
Let's raise the index $k$ in $U_k$, i\.\,e\. consider the
quantities $U^q$ given by the formula
$$
U^q=\sum^2_{k=1}g^{qk}\,U_k.
$$
For these quantities from the previously derived formula we obtain
$$
-U^q=\ddot u^q+\sum^2_{i=1}\sum^2_{j=1}
\sum^2_{k=1}\frac{g^{qk}}{2}
\left(
 \frac{\partial g_{kj}}{\partial u^i}
+\frac{\partial g_{ik}}{\partial u^j}
-\frac{\partial g_{ij}}{\partial u^k}
\right)\dot u^i\,\dot u^j.
$$
Let's compare this formula with the formula \mythetagchapter{7.5}{4}
in Chapter \uppercase\expandafter{\romannumeral 4} that determines the
connection components. As a result we get:
$$
-U^q=\ddot u^q+\sum^2_{i=1}\sum^2_{j=1}
\Gamma^q_{ij}\,\dot u^i\,\dot u^j.
\mytag{3.11}
$$
Now it is sufficient to compare \mythetag{3.11} with the equation of 
geodesic lines \mythetag{2.15} and derive $U^q=0$. The quantities $U_k$ 
are obtained from $U^q$ by lowering the index:
$$
U_k=\sum^2_{q=1}g_{kq}\,U^q.
$$
Therefore, the quantities $U_k$ are also equal to zero. From this fact we immediately derive the equality \mythetag{3.4} which means exactly that the
extremity condition for the geodesic lines is fulfilled. The theorem is
proved.\qed\enddemo
\head
\S~\mysection{4} Inner parallel translation on a surface.
\endhead
\rightheadtext{\S~4. Inner parallel translation on a surface.}
     The equation of geodesic lines in the Euclidean space $\Bbb E$ in
form of \mythetagchapter{8.17}{3} was derived in Chapter
\uppercase\expandafter{\romannumeral 3} when considering the parallel
translation of vectors in curvilinear coordinates. The differential equation
of the parallel translation \mythetagchapter{8.5}{3} can be rewritten now
in the two-dimensional case:
$$
\hskip -2em
\dot a^i+\sum^2_{j=1}\sum^2_{k=1}
\Gamma^i_{jk}\,\dot u^j\,a^k=0.
\mytag{4.1}
$$
The equation \mythetag{4.1} is called the {\it equation of the inner
parallel translation\/} of vectors along curves on a surface.\par 
     Suppose that we have a surface on some fragment of which the
curvilinear coordinates $u^1,\,u^2$ and a parametric curve \mythetag{1.1}
are given. Let's consider some tangent vector $\bold a$ to the surface
at the initial point of the curve, i\.\,e\. at $t=0$. The vector 
$\bold a$ has the inner representation in form of two numbers $a^1,\,
a^2$, they are its components. Let's set the Cauchy problem for the
differential equations \mythetag{4.1} given by the following initial data 
at $t=0$:
$$
\xalignat 2
&\hskip -2em
a^1(t)\,\vrule height 14pt depth 10pt
\lower 10 pt \hbox{$\ssize\,t=0$}=a^1,
&&a^2(t)\,\vrule height 14pt depth 10pt
\lower 10 pt \hbox{$\ssize\,t=0$}=a^2.
\mytag{4.2}
\endxalignat
$$
Solving the Cauchy problem \mythetag{4.2}, we get two functions $a^1(t)$ 
and $a^2(t)$ which determine the vectors $\bold a(t)$ at all points of
the curve. The procedure described just above is called the {\it inner
parallel translation of the vector $\bold a$ along a curve on a surface}.
\par
     Let's consider the inner parallel translation of the vector $\bold a$
from the outer point of view, i\.\,e\. as a process in outer
(three-dimensional) geometry of the space $\Bbb E$ where the surface 
under consideration is embedded. The relation of inner and outer 
representations of tangent vectors of the surface is given by the formula:
$$
\hskip -2em
\bold a=\sum^2_{i=1}\,a^i\cdot\bold E_i.
\mytag{4.3}
$$
Let's differentiate the equality \mythetag{4.3} with respect to $t$ 
assuming that $a^1$ and $a^2$ depend on $t$ as solutions of the
differential equations \mythetag{4.1}:
$$
\hskip -2em
\frac{d\bold a}{dt}=\sum^2_{i=1}\dot a^i\cdot\bold E_i
+\sum^2_{i=1} \sum^2_{j=1}a^i\cdot\frac{\partial\bold E_i}
{\partial u^j}\cdot\dot u^j.
\mytag{4.4}
$$
The derivatives $\partial\bold E_i/\partial u^j$ are calculated 
according to Veingarten's derivational formulas (see formulas
\mythetagchapter{4.11}{4} in Chapter 
\uppercase\expandafter{\romannumeral 4}). Then
$$
\frac{d\bold a}{dt}=\sum^2_{i=1}
\left(\dot a^i+\shave{\sum^2_{j=1}\sum^2_{k=1}}
\Gamma^i_{jk}\,\dot u^j\,a^k\right)\cdot\bold E_i+
\left(\shave{\sum^2_{j=1}\sum^2_{k=1}}b_{jk}\,
\dot u^j\,a^k\right)\cdot\bold n.
$$
Since the functions $a^i(t)$ satisfy the differential equations
\mythetag{4.1}, the coefficients at the vectors $\bold E_i$ in this 
formula do vanish:
$$
\hskip -2em
\frac{d\bold a}{dt}=
\left(\shave{\sum^2_{j=1}\sum^2_{k=1}}b_{jk}\,
\dot u^j\,a^k\right)\cdot\bold n.
\mytag{4.5}
$$
The coefficient at the normal vector $\bold n$ in the above formula
\mythetag{4.5} is determined by the second quadratic form of the surface.
This is the value of the corresponding symmetric bilinear form on the pair
of vectors $\bold a$ and $\boldsymbol\tau$. Therefore, the formula
\mythetag{4.5} is rewritten in a vectorial form as follows:
$$
\hskip -2em
\frac{d\bold a}{dt}=\bold b(\boldsymbol\tau,\bold a)
\cdot\bold n.
\mytag{4.6}
$$
The vectorial equation \mythetag{4.6} is called the {\it outer equation of
the inner parallel translation on surfaces}.\par
     The operation of parallel translation can be generalized to the case
of inner tensors of the arbitrary type $(r,s)$. For this purpose we have 
introduced the operation of covariant differentiation of tensorial function
with respect to the parameter $t$ on curves (see formula 
\mythetagchapter{8.9}{3} in Chapter 
\uppercase\expandafter{\romannumeral 3}). Here is the two-dimensional
version of this formula:
$$
\pagebreak
\gathered
\nabla_{\!t}A^{i_1\ldots\,i_r}_{j_1\ldots\,j_s}=
 \frac{dA^{i_1\ldots\,i_r}_{j_1\ldots\,j_s}}{dt}+\\
 +\sum^r_{m=1}\sum^2_{q=1}\sum^2_{v_m=1}
 \Gamma^{i_m}_{q\,v_m}\,\dot u^q\,
 A^{i_1\ldots\,v_m\ldots\,i_r}_{j_1\ldots\,j_s}
 -\sum^s_{n=1}\sum^2_{q=1}\sum^2_{w_n=1}
  \Gamma^{w_n}_{q\,j_n}\,\dot u^q\,
 A^{i_1\ldots\,i_r}_{j_1\ldots\,w_n\ldots\,j_s}.
\endgathered
\mytag{4.7}
$$
In terms of the covariant derivative \mythetag{4.7} the equation of the
inner parallel translation for the tensorial field $\bold A$ is written
as
$$
\hskip -2em
\nabla_{\!t}\bold A=0.
\mytag{4.8}
$$
The consistence of defining the inner parallel translation by means
of the equation \mythetag{4.8} follows from the two-dimensional analog
of the theorem~\mythetheoremchapter{8.2}{3} from Chapter \uppercase\expandafter{\romannumeral 3}.
\mytheorem{4.1} For any inner tensorial function $\bold A(t)$
determined at the points of a parametric curve on some surface the
quantities $B^{i_1\ldots\,i_r}_{j_1\ldots\,j_s}=\nabla_{\!t}A^{i_1\ldots\,
i_r}_{j_1\ldots\,j_s}$ calculated according to the formula \mythetag{4.7}
define a tensorial function $\bold B(t)=\nabla_{\!t}\bold A$ of the same 
type $(r,s)$ as the original function $\bold A(t)$.
\endproclaim
    The proof of this theorem almost literally coincides with the 
proof of the theorem~\mythetheoremchapter{8.2}{3} in Chapter
\uppercase\expandafter{\romannumeral 3}. Therefore, we do not give it
here.\par
    The covariant differentiation  $\nabla_{\!t}$ defined by the formula
\mythetag{4.7} possesses a series of properties similar to those of the
covariant differentiation along a vector field $\nabla_{\bold X}$ (see formula \mythetagchapter{6.10}{4} and theorem~\mythetheoremchapter{6.2}{4}
in Chapter \uppercase\expandafter{\romannumeral 4}).\par
\mytheorem{4.2} The operation of covariant differentiation of 
tensor-valued func\-tions with respect to the parameter $t$ along a curve
on a surface possesses the following properties:
\roster
\item $\nabla_{\!t}(\bold A+\bold B)=\nabla_{\!t}\bold A
       +\nabla_{\!t}\bold B$;
\item $\nabla_{\!t}(\bold A\otimes\bold B)=\nabla_{\!t}
       \bold A\otimes\bold B+\bold A\otimes\nabla_{\!t}\bold B$;
\item $\nabla_{\!t}C(\bold A)=C(\nabla_{\!t}\bold A)$.
\endroster
\endproclaim
\demo{Proof} Let's choose some curvilinear coordinate system and prove the
theorem by means of direct calculations in coordinates. Let $\bold C=\bold
A+\bold B$. Then for the components of the tensor-valued function
$\bold C(t)$ we have
$$
C^{\,i_1\ldots\,i_r}_{j_1\ldots\,j_s}=A^{i_1\ldots\,i_r}_{j_1\ldots
\,j_s}+B^{i_1\ldots\,i_r}_{j_1\ldots\,j_s}.
$$
Substituting $C^{\,i_1\ldots\,i_r}_{j_1\ldots\,j_s}$ into \mythetag{4.7},
for the covariant derivative $\nabla_{\!t}\bold C$ we get
$$
\nabla_{\!t}C^{\,i_1\ldots\,i_r}_{j_1\ldots\,j_s}=
\nabla_{\!t}A^{i_1\ldots\,i_r}_{j_1\ldots\,j_s}+
\nabla_{\!t}B^{i_1\ldots\,i_r}_{j_1\ldots\,j_s}.
$$
This equality proves the first item of the theorem.\par
     Let's proceed with the item \therosteritem{2}. Denote $\bold C=
\bold A\otimes\bold B$. Then for the components of the tensor-valued
function $\bold C(t)$ we have
$$
\hskip -2em
C^{\,i_1\ldots\,i_{r+p}}_{j_1\ldots\,j_{s+q}}=
A^{i_1\ldots\,i_r}_{j_1\ldots\,j_s}
B^{i_{r+1}\ldots\,i_{r+p}}_{j_{s+1}\ldots\,j_{s+q}}.
\mytag{4.9}
$$
Let's substitute the quantities $C^{\,i_1\ldots\,i_{r+p}}_{j_1\ldots\,
j_{s+q}}$ from \mythetag{4.9} into the formula \mythetag{4.8} for the 
covariant derivative. As a result for the components of $\nabla_{\!t}
\bold C$ we derive 
$$
\gather
\nabla_{\!t}C^{\,i_1\ldots\,i_{r+p}}_{j_1\ldots\,j_{s+q}}=
\bigl(dA^{\,i_1\ldots\,i_r}_{j_1\ldots\,j_s}/dt\bigr)\,
\,B^{i_{r+1}\ldots\,i_{r+p}}_{j_{s+1}\ldots\,j_{s+q}}\,+\\
\vspace{1.5ex}
\qquad\qquad+A^{\,i_1\ldots\,i_r}_{j_1\ldots\,j_s}
\,\,\bigl(dB^{i_{r+1}\ldots\,i_{r+p}}_{j_{s+1}\ldots\,j_{s+q}}/
dt\bigr)\,+\\
\displaybreak
+\sum^r_{m=1}\sum^2_{q=1}\sum^2_{v_m=1}
\Gamma^{i_m}_{\,v_m}\,\dot u^q\,
A^{i_1\ldots\,v_m\ldots\,i_r}_{j_1\ldots\,j_s}\,
B^{i_{r+1}\ldots\,i_{r+p}}_{j_{s+1}\ldots\,j_{s+q}}+
\\
+\sum^{r+p}_{m=r+1}\sum^2_{q=1}\sum^2_{v_m=1}
A^{\,i_1\ldots\,i_r}_{j_1\ldots\,j_s}\,
\Gamma^{i_m}_{\,v_m}\,\dot u^q\,
B^{i_{r+1}\ldots\,v_m\ldots\,i_{r+p}}_{j_{s+1}\ldots\,j_{s+q}}-
\\
-\sum^s_{n=1}\sum^2_{q=1}\sum^2_{w_n=1}
\Gamma^{w_n}_{\,j_n}\,\dot u^q\,
A^{i_1\ldots\,i_r}_{j_1\ldots\,w_n\ldots\,j_s}\,
B^{i_{r+1}\ldots\,i_{r+p}}_{j_{s+1}\ldots\,j_{s+q}}-\\
-\sum^{s+q}_{n=s+1}\sum^2_{q=1}\sum^2_{w_n=1}
A^{\,i_1\ldots\,i_r}_{j_1\ldots\,j_s}\,
\Gamma^{w_n}_{\,j_n}\,\dot u^q\,
B^{i_{r+1}\ldots\,i_{r+p}}_{j_{s+1}\ldots\,w_n\ldots\,j_{s+q}}.
\endgather
$$
Note that upon collecting the similar terms the above huge formula
can be transformed to the following one:
$$
\hskip -2em
\aligned
\nabla_{\!t}&\bigl(A^{i_1\ldots\,i_r}_{j_1\ldots\,j_s}\,
B^{i_{r+1}\ldots\,i_{r+p}}_{j_{s+1}\ldots\,j_{s+q}}\bigr)=
\bigl(\nabla_{\!t}A^{i_1\ldots\,i_r}_{j_1\ldots\,j_s}
\bigr)\times
\\
\vspace{2ex}
&\times B^{i_{r+1}\ldots\,i_{r+p}}_{j_{s+1}\ldots\,j_{s+q}}
+A^{i_1\ldots\,i_r}_{j_1\ldots\,j_s}
\,\bigl(\nabla_{\!t}B^{i_{r+1}\ldots\,i_{r+p}}_{j_{s+1}
\ldots\,j_{s+q}}\bigr).
\endaligned
\mytag{4.10}
$$
Now it is easy to see that the formula \mythetag{4.10} proves the
second item of the theorem.\par
     Let's choose two tensor-valued functions $\bold A(t)$ and
$\bold B(t)$ one of which is the contraction of another. In
coordinates this fact looks like
$$
\hskip -2em
B^{i_1\ldots\,i_r}_{j_1\ldots\,j_s}=
\sum^2_{k=1}A^{i_1\ldots\,i_{p-1}\,k\,i_p\ldots\,i_r}_{j_1
\ldots\,j_{q-1}\,k\,j_q\ldots\,j_s}.
\mytag{4.11}
$$
Let's substitute \mythetag{4.11} into the formula \mythetag{4.7}.
For $\nabla_{\!t}B^{i_1\ldots\,i_r}_{j_1\ldots\,j_s}$ we derive
$$
\align
&\hskip -2em
\nabla_{\!t}B^{i_1\ldots\,i_r}_{j_1\ldots\,j_s}=
\sum^2_{k=1}
\frac{dA^{i_1\ldots\,i_{p-1}\,k\,i_p\ldots\,
i_r}_{j_1\ldots\,j_{q-1}\,k\,j_q\ldots\,j_s}}{dt}\ +\\
&\hskip -2em
+\sum^r_{m=1}\sum^2_{k=1}\sum^2_{q=1}\sum^2_{v_m=1}
\Gamma^{i_m}_{q\,v_m}\,\dot u^q\,A^{i_1\ldots\,v_m\ldots\,
k\,\ldots\,i_r}_{j_1\ldots\,j_{q-1}\,k\,j_q\ldots\,j_s}-
\mytag{4.12}\\
&\hskip -2em
-\sum^s_{n=1}\sum^2_{k=1}\sum^2_{q=1}\sum^2_{w_n=1}
\Gamma^{w_n}_{q\,j_n}\,\dot u^q\,A^{i_1\ldots\,i_{p-1}\,
k\,i_p\ldots\,i_r}_{j_1\ldots\,w_n\ldots\,k\,\ldots\,j_s}.
\endalign
$$
In the formula \mythetag{4.12} the index $v_m$ sequentially occupies the
positions to the left of the index $k$ and to the right of it. The same 
is true for the index $w_n$. However, the formula \mythetag{4.12} has no
terms where $v_m$ or $w_n$ replaces the index $k$. Such terms, provided they would be present, according to \mythetag{4.7}, would have the form
$$
\gather
\hskip -2em
\sum^2_{k=1}\sum^2_{q=1}\sum^2_{v=1}
\Gamma^k_{q\,v}\,\dot u^q\,\,A^{i_1\ldots\,i_{p-1}\,v\,i_p\ldots
\,i_r}_{j_1\ldots\,j_{q-1}\,k\,j_q\ldots\,j_s},
\mytag{4.13}\\
\hskip -2em
-\sum^2_{k=1}\sum^2_{q=1}\sum^2_{w=1}
\Gamma^w_{q\,k}\,\dot u^q\,\,A^{i_1\ldots\,i_{p-1}\,k\,
i_p\ldots \,i_r}_{j_1\ldots\,j_{q-1}\,w\,j_q\ldots\,j_s}.
\mytag{4.14}
\endgather
$$
It is easy to note that \mythetag{4.13} and \mythetag{4.14} differ only in
sign. Indeed, it is sufficient to rename $k$ to $v$ and $w$ to $k$ in the
formula \mythetag{4.14}. If we add simultaneously \mythetag{4.13} and
\mythetag{4.14} to \mythetag{4.12}, their contributions cancel each other
thus keeping the equality valid. Therefore, \mythetag{4.12} can be written 
as
$$
\nabla_{j_{s+1}}B^{i_1\ldots\,i_r}_{j_1\ldots\,j_s}=
\sum^2_{k=1}\nabla_{j_{s+1}}
A^{i_1\ldots\,i_{p-1}\,k\,i_p\ldots\,i_r}_{j_1\ldots\,j_{q-1}
\,k\,j_q\ldots\,j_s}.
\mytag{4.15}
$$
The relationship \mythetag{4.15} proves the third item of the theorem and completes the proof in whole.
\qed\enddemo
     Under a reparametrization of a curve a new parameter $\tilde t$ should
be a strictly monotonic function of the old parameter $t$ (see details in
\S\,2 of Chapter \uppercase\expandafter{\romannumeral 1}). Under such a
reparametrization $\nabla_{\!\tilde t}$ and $\nabla_{\!t}$ are related to each other by the formula
$$
\hskip -2em
\nabla_{\!t}\bold A=\frac{d\tilde t(t)}{dt}\cdot
\nabla_{\tilde t}\bold A
\mytag{4.16}
$$
for any tensor-valued function $\bold A$ on a curve. This relationship 
is a simple consequence from \mythetag{4.7} and from the chain rule for
differentiating a composite function. It is an analog of the item 
\therosteritem{3} in the theorem~\mythetheoremchapter{6.2}{4} of Chapter
\uppercase\expandafter{\romannumeral 4}.\par
     Let $\bold A$ be a tensor field of the type  $(r,s)$ on a surface.
This means that at each point of the surface some tensor of the type 
$(r,s)$ is given. If we mark only those points of the surface which 
belong to some curve, we get a tensor-valued function $\bold A(t)$ on
that curve. In coordinates this is written as 
$$
\hskip -2em
A^{i_1\ldots\,i_r}_{j_1\ldots\,j_s}(t)=
A^{i_1\ldots\,i_r}_{j_1\ldots\,j_s}(u^1(t),u^2(t)).
\mytag{4.17}
$$
The function $\bold A(t)$ constructed in this way is called the 
{\it restriction of a tensor field $\bold A$} to a curve. The specific
feature of the restrictions of tensor fields on curves expressed by the
formula \mythetag{4.17} reveals in differentiating them:
$$
\hskip -2em
\frac{dA^{i_1\ldots\,i_r}_{j_1\ldots\,j_s}}{dt}=
\sum^2_{q=1}\frac{\partial A^{i_1\ldots\,i_r}_{j_1\ldots\,j_s}}
{\partial u^q}\,\dot u^q.
\mytag{4.18}
$$
Substituting \mythetag{4.18} into the formula \mythetag{4.7}, we can 
extract the common factor $\dot u^q$ in the sum over $q$. Upon extracting
this common factor we find 
$$
\hskip -2em
\nabla_{\!t}A^{i_1\ldots\,i_r}_{j_1\ldots\,j_s}=
\sum^2_{q=1}\dot u^q\,\nabla_qA^{i_1\ldots\,
i_r}_{j_1\ldots\,j_s}.
\mytag{4.19}
$$
The formula \mythetag{4.19} means that the covariant derivative of the
restriction of a tensor field $\bold A$ to a curve is the contraction 
of the covariant differential $\nabla\bold A$ with the tangent vector 
of the curve.\par
     Assume that $\nabla\bold A=0$. Then due to \mythetag{4.19} the
restriction of the field $\bold A$ to any curve is a tensor-valued
function satisfying the equation of the parallel translation \mythetag{4.8}.
The values of such a field $\bold A$ at various points are related to each other by parallel translation along any curve connecting these points.\par
\mydefinition{4.1} A tensor field $\bold A$ is called an {\it autoparallel}
field or a {\it covariantly constant} field if its covariant differential
is equal to zero identically: $\nabla\bold A=0$.
\enddefinition
    Some of the well-known tensor fields have identically zero covariant
differentials: this is the metric tensor $\bold g$, the inverse metric
tensor $\hat\bold g$, and the area tensor (pseudotensor) $\boldsymbol\tau$.
The autoparallelism of these fields plays the important role for describing
the inner parallel translation.\par
    Let $\bold a$ and $\bold b$ be two tangent vectors of the surface at
the initial point of some curve. Their scalar product is calculated through
their components:
$$
(\bold a\,|\,\bold b)=\sum^2_{i=1}\sum^2_{j=1}
g_{ij}\,a^i\,b^j.
$$
Let's perform the parallel translation of the vectors $\bold a$ and 
$\bold b$ along the curve solving the equation \mythetag{4.8} and using 
the components of $\bold a$ and $\bold b$ as initial data in Cauchy
problems. As a result we get two vector-valued functions $\bold a(t)$
and $\bold b(t)$ on the curve. Let's consider the function
$\psi(t)$ equal to their scalar product:
$$
\hskip -2em
\psi(t)=(\bold a(t)\,|\,\bold b(t))=\sum^2_{i=1}\sum^2_{j=1}
g_{ij}(t)\,\,a^i(t)\,b^j(t).
\mytag{4.20}
$$
According to the formula \mythetag{4.7} the covariant derivative
$\nabla_{\!t}\psi$ coincides with the regular derivative. Therefore, we
have
$$
\frac{d\psi}{dt}=\nabla_{\!t}\psi=\sum^2_{i=1}\sum^2_{j=1}
\bigl(\nabla_{\!t}g_{ij}\,\,a^i\,b^j+g_{ij}\,\nabla_{\!t}a^i\,\,b^j+
g_{ij}\,a^i\,\,\nabla_{\!t}b^j\bigr).
$$
Here we used the items \therosteritem{2} and \therosteritem{3} of the
theorem~\mythetheorem{4.2}. But $\nabla_{\!t}a^i=0$ and $\nabla_{\!t}b^j=0$
since we $\bold a(t)$ and $\bold b(t)$ are obtained as a result of parallel
translation of the vectors $\bold a$ and $\bold b$. Moreover,
$\nabla_{\!t}g_{ij}=0$ due to autoparallelism of the metric tensor.
For the scalar function $\psi(t)$ defined by \mythetag{4.20} this yields
$d\psi/dt=0$ and $\psi(t)=(\bold a\,|\,\bold b)=\const$. As a result of
these considerations we have proved the following theorem.
\mytheorem{4.3} The operation of inner parallel translation of
vectors along cur\-ves preserves the scalar product of vectors.
\endproclaim
    Preserving the scalar product, the operation of inner parallel
translation preserves the length of vectors and the angles between
them.\par
    From the autoparallelism of metric tensors $\bold g$ and $\hat\bold g$
we derive the following formulas analogous to the formulas 
\mythetagchapter{7.9}{4} in
Chapter \uppercase\expandafter{\romannumeral 4}:
$$
\hskip -2em
\aligned
&\nabla_{\!t}\left(\,\shave{\sum^2_{k=1}}g_{ik}\,
 A^{\ldots\,k\,\ldots}_{\ldots\,\hphantom{k}\,\ldots}\right)=
 \sum^2_{k=1}g_{ik}\,
 \nabla_{\!t}A^{\ldots\,k\,\ldots}_{\ldots\,\hphantom{k}\,\ldots},\\
&\nabla_{\!t}\left(\,\shave{\sum^2_{k=1}}g^{ik}\,
 A^{\ldots\,\hphantom{k}\,\ldots}_{\ldots\,k\,\ldots}\right)=
 \sum^2_{k=1}g^{ik}\,
 \nabla_{\!t}A^{\ldots\,\hphantom{k}\,\ldots}_{\ldots\,k\,\ldots}.
\endaligned
\mytag{4.21}
$$
Then from the formulas \mythetag{4.21} we derive the following fact.
\mytheorem{4.4} The operation of inner parallel translation
of tensors commutate with the operations of index raising and index 
lowering.
\endproclaim
\head
\S~\mysection{5} Integration on surfaces. Green's formula.
\endhead
\rightheadtext{\S~5. Integration on surfaces. Green's formula.}
\parshape 3 
0pt 360pt 0pt 360pt 180pt 180pt 
     Let's consider the two-dimensional space $\Bbb R^2$. Let's draw it as
a coordinate plane $u^1,\,u^2$. Let's choose some simply connected domain $\Omega$ \vadjust{\vskip 5pt\hbox to 0pt{\kern -5pt
\includegraphics{ris14.eps}\hss}\vskip -5pt}outlined by a
closed piecewise continuously differentiable con\-tour $\gamma$ on the
coordinate plane $u^1,\,u^2$. Then we mark the direction (orientation) on
the contour $\gamma$ so that when moving in this direction the domain
$\Omega$ lies to the left of the contour. On Fig\.~5.1 this direction is
marked by the arrow. In other words, we choose the orientation on $\gamma$
induced from the orientation of $\Omega$,\linebreak i\.\,e\. $\gamma=\partial\Omega$.
\par
\parshape 6 180pt 180pt 180pt 180pt 180pt 180pt 180pt 180pt
180pt 180pt 0pt 360pt 
    Let's consider a pair of continuously differentiable functions on the
coordinate plane: $P(u^1,u^2)$ and $Q(u^1,u^2)$. Then, if all the above
conditions are fulfilled, there is the following integral identity:
$$
\hskip -2em
\oint\limits_{\gamma\,\,}\bigl(P\,du^1+Q\,du^2\bigr)=\iint
\limits_{\Omega}\left(\frac{\partial Q}{\partial u^1}
-\frac{\partial P}{\partial u^2}\right)\,du^1du^2.
\mytag{5.1}
$$
The identity \mythetag{5.1} is known as Green's formula (see \mycite{2}). 
The equality \mythetag{5.1} is an equality for a plane. We need its generalization for the case of an arbitrary surface in the space
$\Bbb E$. In such generalization the coordinate plane $u^1,\,u^2$ 
or some its part plays the role of a chart, while the real geometric
domain and its boundary contour should be placed on a surface. 
Therefore, the integrals in both parts of Green's formula should be
transformed so that one can easily write them for any curvilinear
coordinates on a surface and their values should not depend on a
particular choice of such coordinate system.\par
     Let's begin with the integral in the left hand side of \mythetag{5.1}.
Such integrals are called {\it path integrals of the second kind}. Let's
rename $P$ to $v_1$ and $Q$ to $v_2$. Then the integral in the left hand
side of \mythetag{5.1} is written as 
$$
\hskip -2em
I=\oint\limits_{\gamma\,\,}\sum^2_{i=1}
v_i(u^1,u^2)\,du^i.
\mytag{5.2}
$$
In order to calculate the integral \mythetag{5.2} practically the contour
$\gamma$ should be parametrized, i\.\,e\. it should be represented as a
parametric curve \mythetag{1.1}. Then the value of an integral of the
second kind is calculated as follows:
$$
\hskip -2em
I=\pm\int\limits^{\,b}_{a\,}\left(\,\shave{\sum^2_{i=1}}
v_i\,\dot u^i\right)dt.
\mytag{5.3}
$$
This formula reducing the integral of the second kind to the regular
integral over the segment $[a,\,b]$ on the real axis can be taken for
the definition of the integral \mythetag{5.2}. The sign is chosen 
regarding to the direction of the contour on Fig\.~5.1. If $a<b$ and
if when $t$ changes from $a$ to $b$ the corresponding point on the
contour moves along the arrow, we choose plus in \mythetag{5.3}. Otherwise,
we choose minus. Changing the variable $\tilde t=\varphi(t)$ in the 
integral \mythetag{5.3} and choosing the proper sign upon reparametrization
of the contour, one can verify that the value of this integral does not 
depend on the choice of the parametrization on the contour.\par
     Now let's change the curvilinear coordinate system on the surface. 
The derivatives $\dot u^i$ in the integral \mythetag{5.3} under a change
of curvilinear coordinates on the surface are transformed as follows:
$$
\hskip -2em
\dot u^i=\frac{du^i}{dt}=\sum^2_{j=1}\frac{\partial u^i}
{\partial\tilde u^j}\,\frac{d\tilde u^j}{dt}=
\sum^2_{j=1}S^i_j\,\dot{\tilde u}\vphantom{u}^j.
\mytag{5.4}
$$
Substituting \mythetag{5.4} into the formula \mythetag{5.3}, for the integral
$I$ we derive:
$$
\hskip -2em
I=\pm\int\limits^{\,b}_{a\,}\left(\,\shave{\sum^2_{j=1}}
\left(\,\shave{\sum^2_{i=1}}S^i_j\,v_i\right)
\dot{\tilde u}\vphantom{u}^j\right)dt.
\mytag{5.5}
$$
Now let's write the relationship \mythetag{5.3} in coordinates $\tilde u^1,
\,\tilde u^2$. For this purpose we rename $u^i$ to $\tilde u^i$ and $v_i$
to $\tilde v_i$ in the formula \mythetag{5.3}:
$$
\hskip -2em
I=\pm\int\limits^{\,b}_{a\,}\left(\,\shave{\sum^2_{i=1}}
\tilde v_i\,\dot{\tilde u}\vphantom{u}^i\right)dt.
\mytag{5.6}
$$
Comparing the formulas \mythetag{5.5} and \mythetag{5.6}, we see that these
formulas are similar in their structure. For the numeric values of the
integrals \mythetag{5.3} and \mythetag{5.6} to be always equal (irrespective
to the form of the contour $\gamma$ and its parametrization) the quantities
$v_i$ and $\tilde v_i$ should be related as follows:
$$
\xalignat 2
&\tilde v_j=\sum^2_{i=1}S^i_j\,v_i,
&&v_i=\sum^2_{i=1}T^j_i\,\tilde v_j.
\endxalignat
$$
These formulas represent the transformation rule for the components of
a covec\-torial field. Thus, we conclude that any path integral of the
second kind on a surface \mythetag{5.2} is given by some inner covectorial
field on this surface.\par
     Now let's proceed with the integral in the right hand side of the
Green's formula \mythetag{5.1}. Distracting for a while from the particular
integral in this formula, let's consider the following double integral:
$$
\hskip -2em
I=\iint\limits_{\Omega}\,F\,\,du^1du^2.
\mytag{5.7}
$$
A change of curvilinear coordinates can be interpreted as a change of
variables in the integral \mythetag{5.7}. Remember that a change of 
variables in a multiple integral is performed according to the following
formula (see \mycite{2}):
$$
\hskip -2em
\iint\limits_{\vphantom{\tilde\Omega}\Omega}
\,F\,\,du^1du^2=\iint\limits_{\tilde\Omega}\,F\,
|\det J|\,\,d\tilde u^1d\tilde u^2,
\mytag{5.8}
$$
where $J$ is the Jacobi matrix determined by the change of variables: 
$$
\hskip -2em
J=\left\Vert
\vphantom{\vrule height 24pt depth 24pt}
\matrix
\dsize\frac{\partial u^1}{\partial\tilde u^1}&
\dsize\frac{\partial u^1}{\partial\tilde u^2}\\
\vspace{1.2ex}
\dsize\frac{\partial u^2}{\partial\tilde u^1}&
\dsize\frac{\partial u^2}{\partial\tilde u^2}
\endmatrix\right\Vert.
\mytag{5.9}
$$
The Jacobi matrix \mythetag{5.9} coincides with the transition matrix $S$ (see formula \mythetagchapter{2.7}{4} in Chapter 
\uppercase\expandafter{\romannumeral 4}). Therefore, the function $F$ being
integrated in the formula \mythetag{5.7} should obey the transformation rule
$$
\tilde F=|\det S|\,F
\mytag{5.10}
$$
under a change of curvilinear coordinates on the surface. The quantity
$F$ has no indices. However, due to \mythetag{5.10}, this quantity is
neither a scalar nor a pseudoscalar. In order to change this not very
pleasant situation the integral \mythetag{5.7} over a two-dimensional 
domain $\Omega$ on a surface is usually written as
$$
\hskip -2em
I=\iint\limits_{\Omega}\sqrt{\det\bold g}\,\,f\,\,du^1du^2,
\mytag{5.11}
$$
where $\det\bold g$ is the determinant of the first quadratic form. In this
case the quantity $f$ in the formula \mythetag{5.11} is a scalar. This fact
follows from the equality $\det\bold g=(\det T)^2\,\det\tilde\bold g$ that
represent the transformation rule for the determinant of the metric tensor
under a change of coordinate system.\par
     Returning back to the integral in the right hand side of \mythetag{5.1},
we transform it to the form \mythetag{5.11}. For this purpose we use the
above notations $P=v_1$, $Q=v_2$, and remember that $v_1$ and $v_2$ are the
components of the covectorial field. Then
$$
\hskip -2em
\frac{\partial Q}{\partial u^1}-
\frac{\partial P}{\partial u^2}=
\frac{\partial v_2}{\partial u^1}-
\frac{\partial v_1}{\partial u^2}.
\mytag{5.12}
$$
The right hand side of \mythetag{5.12} can be represented in form of the 
contraction with the unit skew-symmetric matrix $d^{ij}$ (see formula
\mythetagchapter{3.6}{4} in Chapter 
\uppercase\expandafter{\romannumeral 4}):
$$
\hskip -2em
\frac{\partial v_2}{\partial u^1}-
\frac{\partial v_1}{\partial u^2}=
\sum^2_{i=1}\sum^2_{j=1} d^{ij}\,\frac{\partial v_j}
{\partial u^i}=
\sum^2_{i=1}\frac{\partial}{\partial u^i}\!
\left(\,\shave{\sum^2_{j=1}}d^{ij}\,v_j\right).
\mytag{5.13}
$$
Note that the quantities $d_{ij}$ with lower indices enter the formula
for the area tensor $\boldsymbol\omega$ (see \mythetagchapter{3.7}{4}
in Chapter \uppercase\expandafter{\romannumeral 4}). Let's raise the
indices of the area tensor by means of the inverse metric tensor:
$$
\omega^{ij}=\sum^2_{p=1}\sum^2_{q=1} g^{ip}\,g^{jq}\,
\omega_{pq}=\sum^2_{p=1}\sum^2_{q=1}\xi_D\,\sqrt{\det\bold g}\,
\,g^{ip}\,g^{jq}\,d_{pq}.
$$
Applying the formula \mythetagchapter{3.7}{4} from Chapter 
\uppercase\expandafter{\romannumeral 4}, we can calculate the components
of the area tensor $\omega^{ij}$ in the explicit form:
$$
\hskip -2em
\omega^{ij}=\xi_D\,\sqrt{\det\bold g^{-1}}\,d^{ij}.
\mytag{5.14}
$$
The formula \mythetag{5.14} expresses $\omega^{ij}$ through $d^{ij}$.
Now we use \mythetag{5.14} in order to express $d^{ij}$ in the formula 
\mythetag{5.13} back through the components of the area tensor:
$$
\frac{\partial v_2}{\partial u^1}-
\frac{\partial v_1}{\partial u^2}=
\sum^2_{i=1}\frac{\partial}{\partial u^i}\!
\left(\,\shave{\sum^2_{j=1}}
\xi_D\,\sqrt{\det\bold g}\,\,\omega^{ij}\,v_j\right).
$$
In order to simplify the further calculations we denote 
$$
\hskip -2em
y^i=\sum^2_{j=1}\omega^{ij}\,v_j.
\mytag{5.15}
$$
Taking into account \mythetag{5.15}, the formula \mythetag{5.13} can be
written as follows:
$$
\hskip -2em
\aligned
\frac{\partial v_2}{\partial u^1}&-
\frac{\partial v_1}{\partial u^2}=
\sum^2_{i=1}\xi_D\,
\frac{\partial\bigl(\sqrt{\det\bold g}\,\,y^i\bigr)}
{\partial u^i}=\\
\vspace{1ex}
&=\xi_D\,\sqrt{\det\bold g}\,\sum^2_{i=1}
\left(\frac{\partial y^i}{\partial u^i}+
\frac{1}{2}\frac{\partial\bigl(\ln\det\bold g\bigr)}
{\partial u^i}\,y^i\right).
\endaligned
\mytag{5.16}
$$
The logarithmic derivative for the determinant of the metric tensor
is calculated by means of the lemma~\mythelemmachapter{7.1}{4} from 
Chapter \uppercase\expandafter{\romannumeral 4}. However, we need 
not repeat these calculations here, since this derivative is already 
calculated (see \mythetagchapter{7.12}{4} and the proof of the 
theorem~\mythetheoremchapter{7.2}{4} in Chapter
\uppercase\expandafter{\romannumeral 4}):
$$
\hskip -2em
\frac{\partial\bigl(\ln\det\bold g\bigr)}{\partial u^i}=
\sum^2_{p=1}\sum^2_{q=1}g^{pq}\,\frac{\partial g_{pq}}
{\partial u^i}=\sum^2_{q=1}2\,\Gamma^q_{iq}.
\mytag{5.17}
$$
With the use of \mythetag{5.17} the formula \mythetag{5.16} is transformed
as follows:
$$
\frac{\partial v_2}{\partial u^1}-
\frac{\partial v_1}{\partial u^2}=
\xi_D\,\sqrt{\det\bold g}\,\sum^2_{i=1}
\left(\frac{\partial y^i}{\partial u^i}+
\shave{\sum^2_{q=1}}\Gamma^q_{qi}\,y^i\right).
$$
In this formula one easily recognizes the contraction of the covariant
differential of the vector field $\bold y$. Indeed, we have 
$$
\hskip -2em
\frac{\partial v_2}{\partial u^1}-
\frac{\partial v_1}{\partial u^2}=
\xi_D\,\sqrt{\det\bold g}\,\sum^2_{i=1}\nabla_iy^i.
\mytag{5.18}
$$\par
     Using the formula \mythetag{5.18}, the notations 
\mythetag{5.15}, and the autoparallelism condition for the area tensor
$\nabla_q\omega^{ij}=0$, we can write the Green's formula as 
$$
\hskip -2em
\oint\limits_{\gamma\,\,}\sum^2_{i=1}v_i\,du^i=
\xi_D\iint\limits_{\Omega}\sum^2_{i=1}\sum^2_{j=1}
\omega^{ij}\nabla_iv_j\,\sqrt{\det\bold g}\,du^1du^2.
\mytag{5.19}
$$
The sign factor $\xi_D$ in \mythetag{5.19} should be especially commented. 
The condition that the domain $\Omega$ should lie to the left of the 
contour $\gamma$ when moving along the arrow is not invariant under an
arbitrary change of coordinates $u^1,\,u^2$ by $\tilde u^1,\,\tilde u^2$.
Indeed, if we set $\tilde u^1=-u^1$ and $\tilde u^2=u^2$, we would have
the mirror image of the domain $\Omega$ and the contour $\gamma$ shown
on Fig\.~5.1. This means that the direction should be assigned to the 
geometric contour $\gamma$ lying on the surface, not to its image in a
chart. Then the sign factor $\xi_D$ in \mythetag{5.19} can be omitted.
\par
     The choice of the direction on a geometric contour outlining a
domain on a surface is closely related to the choice of the normal 
vector on that surface. The normal vector $\bold n$ should be chosen 
so that when observing from the end of the vector $\bold n$ and moving
in the direction of the arrow along the contour $\gamma$ the domain 
$\Omega$ should lie to the left of the contour. The choice of the normal
vector $\bold n$ defines the orientation of the surface thus defining 
the unit pseudoscalar field $\xi_D$.\par
\head
\S~\mysection{6} Gauss-Bonnet theorem.
\endhead
\rightheadtext{\S~6. Gauss-Bonnet theorem.}
     Let's consider again the process of inner parallel translation of
tangent vectors along curves on surfaces. The equation \mythetag{4.6} shows
that from the outer (three-dimensional) point of view this parallel
translation differs substantially from the regular parallel translation:
the vectors being translated do not remain parallel to the fixed direction
in the space --- they change. However, their lengths are preserved, and, if
we translate several vectors along the same curve, the angles between
vectors are preserved (see theorem~\mythetheorem{4.3}).\par
     From the above description, we see that in the process of parallel translation, apart from the motion of the attachment point along the
curve, the rotation of the vectors about the normal vector $\bold n$ 
occurs. Therefore, we have the natural problem --- how to measure the
angle of this rotation\,? We consider this problem just below.\par
     Suppose that we have a surface equipped with the orientation. This
means that the orientation field $\xi_D$ and the area tensor $\boldsymbol\omega$ are defined (see formula \mythetagchapter{3.10}{4}
in Chapter \uppercase\expandafter{\romannumeral 4}). We already know that
$\xi_D$ fixes one of the two possible normal vectors $\bold n$ at each
point of the surface (see formula \mythetagchapter{4.3}{4} in Chapter
\uppercase\expandafter{\romannumeral 4}).
\mytheorem{6.1} The inner tensor field $\Theta$ of the type
$(1,1)$ with the components 
$$
\hskip -2em
\theta^i_j=\sum^2_{k=1}\omega_{jk}\,g^{ki}
\mytag{6.1}
$$
is an operator field describing the counterclockwise rotation in the
tangent plane to the angle $\pi/2=90^\circ$ about the normal vector 
$\bold n$.
\endproclaim
\demo{Proof} Let $\bold a$ be a tangent vector to the surface and let
$\bold n$ be the unit normal vector at the point where $\bold a$ is
attached. Then, in order to construct the vector $\bold b=\Theta(\bold a)$
obtained by rotating $\bold a$ counterclockwise to the angle
$\pi/2=90^\circ$ about the vector $\bold n$ one can use the following vector product:
$$
\hskip -2em
\bold b=\Theta(\bold a)=[\bold n,\,\bold a].
\mytag{6.2}
$$
Let's substitute the expression given by the formula
\mythetagchapter{4.3}{4} from Chapter 
\uppercase\expandafter{\romannumeral 4} for the vector $\bold n$
into \mythetag{6.2}. Then let's expand the vector $\bold a$ in the basis $\bold E_1,\,\bold E_2$:
$$
\hskip -2em
\bold a=a^1\cdot\bold E_1+a^2\cdot\bold E_2.
\mytag{6.3}
$$
As a result for the vector $\bold b$ in the formula \mythetag{6.2} we derive
$$
\hskip -2em
\bold b=\sum^2_{j=1}\xi_D\cdot\frac{[[\bold E_1,\,
\bold E_2],\,\bold E_j]}{\left|[\bold E_1,\,\bold E_2]
\right|}\cdot a^j.
\mytag{6.4}
$$
In order to calculate the denominator in the formula \mythetag{6.4} we use
the well-known formula from the analytical geometry (see \mycite{4}):
$$
\left|[\bold E_1,\,\bold E_2]\right|^2=
\det
\vmatrix
(\bold E_1\,|\,\bold E_1) &(\bold E_1\,|\,\bold E_2) \\
(\bold E_2\,|\,\bold E_1) &(\bold E_2\,|\,\bold E_2)
\endvmatrix
=\det\bold g.
$$
As for the numerator in the formula \mythetag{6.4}, here we use the not 
less known formula for the double vectorial product:
$$
[[\bold E_1,\,\bold E_2],\,\bold E_j]=\bold E_2\cdot
(\bold E_1\,|\,\bold E_j)-\bold E_1\cdot(\bold E_j\,|\,
\bold E_2).
$$
Taking into account these two formulas, we can write \mythetag{6.4} 
as follows:
$$
\hskip -2em
\bold b=\sum^2_{j=1}\xi_D\cdot\frac{g_{1j}\cdot\bold E_2-
g_{2j}\cdot\bold E_1}{\sqrt{\det\bold g}}\cdot a^j.
\mytag{6.5}
$$
Using the components of the area tensor \mythetag{5.14}, no we can rewrite
\mythetag{6.5} in a more compact and substantially more elegant form:
$$
\bold b=\sum^2_{i=1}\left(\,\shave{\sum^2_{j=1}\sum^2_{k=1}}
\omega^{ki}\,g_{kj}\,a^j\right)\cdot\bold E_i.
$$
From this formula it is easy to extract the formula \mythetag{6.1} for 
the components of the linear operator $\Theta$ relating $\bold b$ and 
$\bold a$. The theorem is proved.
\qed\enddemo
     The operator field $\Theta$ is the contraction of the tensor product
of two fields $\boldsymbol\omega$ and $\bold g$. The autoparallelism of 
the latter ones means that $\Theta$ is also an autoparallel field, i\.\,e\.
$\nabla\Theta=0$.\par
     We use the autoparallelism of $\Theta$ in the following way. Let's
choose some parametric curve $\gamma$ on a surface and perform the parallel
translation of some unit vector $\bold a$ along this curve. As a result we
get the vector-valued function $\bold a(t)$ on the curve satisfying the
equation of parallel translation $\nabla_{\!t}\bold a=0$ (see formula
\mythetag{4.8}). Then we define the vector-function $\bold b(t)$ on the 
curve as follows:
$$
\hskip -2em
\bold b(t)=\Theta(\bold a(t)).
\mytag{6.6}
$$
From \mythetag{6.6} we derive $\nabla_{\!t}(\bold b)=\nabla_{\!t}\Theta(\bold
a)+\Theta(\nabla_{\!t}\bold a)=0$. This means that the function \mythetag{6.6}
also satisfies the equation of parallel translation. It follows from the
autoparallelism of $\Theta$ and from the items \therosteritem{2} and
\therosteritem{3} in the theorem~\mythetheorem{4.2}. The vector-functions
$\bold a(t)$ and $\bold b(t)$ determine two mutually perpendicular unit
vectors at each point of the curve. There are the following obvious
relationships for them:
$$
\xalignat 2
&\hskip -2em
\Theta(\bold a)=\bold b,
&&\Theta(\bold b)=-\bold a.
\mytag{6.7}
\endxalignat
$$
Let's remember for the further use that $\bold a(t)$ and $\bold b(t)$ are
obtained by parallel translation of the vectors $\bold a(0)$ and $\bold
b(0)$ along the curve from its initial point.\par
     Now let's consider some inner vector field $\bold x$ on the surface
(it is tangent to the surface in the outer representation). If the field
vectors $\bold x(u^1,u^2)$ are nonzero at each point of the surface, they
can be normalized to the unit length: $\bold x\to\bold x/|\bold x|$.
Therefore, we shall assume $\bold x$ to be a field of unit vectors: 
$|\bold x|=1$. At the points of the curve $\gamma$ this field can be
expanded in the basis of the vectors $\bold a$ and $\bold b$:
$$
\hskip -2em
\bold x=\cos(\varphi)\cdot\bold a+\sin(\varphi)\cdot
\bold b.
\mytag{6.8}
$$
The function $\varphi(t)$ determines the angle between the vector $\bold a$ 
and the field vector $\bold x$ measured from $\bold a$ to $\bold x$ in the
counterclockwise direction. The change of $\varphi$ describes the rotation
of the vectors during their parallel translation along the curve.\par
     Let's apply the covariant differentiation $\nabla_{\!t}$ to the relationship \mythetag{6.8} and take into account that both vectors $\bold a$
and $\bold b$ satisfy the equation of parallel translation:
$$
\hskip -2em
\nabla_{\!t}\bold x=(-\sin(\varphi)\cdot\bold a+
\cos(\varphi)\cdot\bold b)\cdot\dot\varphi.
\mytag{6.9}
$$
Here we used the fact that the covariant derivative $\nabla_{\!t}$ for 
the scalar coincides with the regular derivative with respect to $t$. In
particular, we have $\nabla_{\!t}\varphi=\dot\varphi$. Now we apply the
operator $\Theta$ to both sides of \mythetag{6.8} and take into account
\mythetag{6.7}:
$$
\hskip -2em
\Theta(\bold x)=\cos(\varphi)\cdot\bold b-\sin(\varphi)\cdot
\bold a.
\mytag{6.10}
$$
Now we calculate the scalar product of $\Theta(\bold x)$ from \mythetag{6.10} and  $\nabla_{\!t}\bold x$ from \mythetag{6.9}. Remembering that $\bold a$ and
$\bold b$ are two mutually perpendicular unit vectors, we get
$$
\hskip -2em
(\Theta(\bold x)\,|\,\nabla_{\!t}\bold x)=(\cos^2(\varphi)+
\sin^2(\varphi))\,\dot\varphi=\dot\varphi.
\mytag{6.11}
$$
Let's write the equality \mythetag{6.11} in coordinate form. The
vector-function $\bold x(t)$ on the curve is the restriction of the
vector field $\bold x$, therefore, the covariant derivative $\nabla_{\!t}\bold x$ is the contraction of the covariant differential
$\nabla\bold x$ with the tangent vector of the curve (see formula
\mythetag{4.19}). Hence, we have
$$
\hskip -2em
\dot\varphi=\sum^2_{q=1}\sum^2_{i=1}\sum^2_{j=1}
\bigl(x^i\,\omega_{ij}\,\nabla_qx^j\bigr)\,\dot u^q.
\mytag{6.12}
$$
Here in deriving \mythetag{6.12} from \mythetag{6.11} we used the formula 
\mythetag{6.1} for the components of the operator field $\Theta$.\par
     Let's discuss the role of the field $\bold x$ in the construction
described just above. The vector field $\bold x$ is chosen as a reference
mark relative to which the rotation angle of the vector $\bold a$ is
measured. This way of measuring the angle is relative. Changing the
field $\bold x$, we would change the value of the angle $\varphi$. We have
to admit this inevitable fact since tangent planes to the surface at
different points are not parallel to each other and we have no preferable
direction relative to which we could measure the angles on all of them.
\par
     There is a case where we can exclude the above uncertainty of the
angle. Let's consider a closed parametric contour $\gamma$ on the surface.
Let $[0,1]$ be the range over which the parameter $t$ runs on such contour.
Then $\bold x(0)$ and $\bold x(1)$ do coincide. They represent the same 
field vector at the point with coordinates $u^1(0),u^2(0)$:
$$
\bold x(0)=\bold x(1)=\bold x(u^1(0),u^2(0)).
$$
Unlike $\bold x(t)$, the function $\bold a(t)$ is not the restriction
of a vector field to a curve $\gamma$. Therefore, the vectors $\bold a(0)$
and $\bold a(1)$ can be different. This is an important feature of the
inner parallel translation that differs it from the parallel translation
in the Euclidean space $\Bbb E$.\par
     In the case of a closed contour $\gamma$ the difference $\varphi(1)
-\varphi(0)$ characterizes the angle to which the vector $\bold a$ turns
$\bold a$ as a result of parallel translation along the contour. Note 
that measuring the angle from $\bold x$ to $\bold a$ is opposite to
measuring it from $\bold a$ to $\bold x$ in the formula \mythetag{6.8}. Therefore, taking for positive the angle measured from $\bold x$ in
the counterclockwise direction, we should take for the increment of the
angle gained during the parallel translation along $\gamma$ the following
quantity:
$$
{\ssize\triangle}\varphi=\varphi(0)-\varphi(1)=
-\int\limits^{\,1}_{0\,}\dot\varphi\,dt.
$$
Let's substitute \mythetag{6.12} for $\dot\varphi$ into this formula.
As a result we get
$$
\hskip -2em
{\ssize\triangle}\varphi=-\int\limits^{\,1}_{0\,}\left(
\shave{\sum^2_{q=1}\sum^2_{i=1}\sum^2_{j=1}}
\bigl(x^i\,\omega_{ij}\,\nabla_qx^j\bigr)\,\dot u^q
\right)dt.
\mytag{6.13}
$$
Comparing \mythetag{6.13} with \mythetag{5.3}, we see that \mythetag{6.13}
now can be written in the form of a path integral of the second kind:
$$
\hskip -2em
{\ssize\triangle}\varphi=-\oint\limits_{\gamma\,\,}
\sum^2_{q=1}\sum^2_{i=1}\sum^2_{j=1}
\bigl(x^i\,\omega_{ij}\,\nabla_qx^j\bigr)\,du^q.
\mytag{6.14}
$$\par
     Assume that the contour $\gamma$ outlines some connected and simply
connected fragment $\Omega$ on the surface. Then for this fragment $\Omega$ we can apply to \mythetag{6.14} the Green's formula written in the form of
\mythetag{5.19}:
$$
{\ssize\triangle}\varphi=-\xi_D\iint\limits_{\Omega}
\sum^2_{i=1}\sum^2_{j=1}\sum^2_{p=1}\sum^2_{q=1}
\omega^{ij}\ \nabla_i\bigl(x^p\,\omega_{pq}\,\nabla_jx^q
\bigr)\,\sqrt{\det\bold g}\,du^1du^2.
$$
If the direction of the contour is in agreement with the orientation
of the surface, then the sign factor $\xi_D$ can be omitted:
$$
\hskip -2em
\gathered
{\ssize\triangle}\varphi=-\iint\limits_{\Omega}
\sum^2_{i=1}\sum^2_{j=1}\sum^2_{p=1}\sum^2_{q=1}
\bigr(x^p\,\omega^{ij}\omega_{pq}\,\nabla_i\nabla_jx^q\ +\\
\vspace{1ex}
+\ \nabla_ix^p\,\omega^{ij}\omega_{pq}\,\nabla_jx^q
\bigr)\,\sqrt{\det\bold g}\,du^1du^2.
\endgathered
\mytag{6.15}
$$
Let's show that the term $\nabla_ix^p\,\omega^{ij}\omega_{pq}\,
\nabla_jx^q$ in \mythetag{6.15} yields zero contribution to the
value of the integral. This feature is specific to the two-dimensional
case where we have the following relationship:
$$
\hskip -2em
\omega^{ij}\omega_{pq}=d^{ij}\,d_{pq}=
\delta^i_p\,\delta^j_q-\delta^i_q\,\delta^j_p.
\mytag{6.16}
$$
The proof of the formula \mythetag{6.16} is analogous to the proof of the
formula \mythetagchapter{8.23}{4} in Chapter
\uppercase\expandafter{\romannumeral 4}. It is based on the skew-symmetry
of $d^{ij}$ and $d_{pq}$.\par
     Let's complete the inner vector field $\bold x$ of the surface
by the other inner vector field $\bold y=\Theta(\bold x)$. The vectors
$\bold x$ and $\bold y$ form a pair of mutually perpendicular unit vectors
in the tangents plane. For their components we have 
$$
\xalignat 3
&\hskip -2em
\sum^2_{q=1} x^q\,x_q=1,
&&x_i=\sum^2_{k=1}g_{ik}\,x^k,
&&y_i=\sum^2_{k=1}g_{ik}\,y^k,\qquad
\mytag{6.17}\\
&\hskip -2em
\sum^2_{q=1}\nabla_kx^q\,x_q=0,
&&y_q=\sum^2_{p=1}\omega_{pq}\,x^p,
&&y^i=\sum^2_{j=1}\omega^{ji}\,x_j.\qquad
\mytag{6.18}
\endxalignat
$$
The first relationship \mythetag{6.17} expresses the fact that $|\bold x|=1$,
other two relationships \mythetag{6.17} determine the covariant components
$x_i$ and $y_i$ of $\bold x$ and $\bold y$. The first relationship 
\mythetag{6.18} is obtained by differentiating \mythetag{6.17}, the second 
and the third relationships \mythetag{6.18} express the vectorial relationship
$\bold y=\Theta(\bold x)$.\par
     Let's multiply \mythetag{6.16} by $\nabla_kx^q\,x_j\,x^p$ and 
then sum up over $q$, $p$, and $j$ taking into account the relationships \mythetag{6.17} and \mythetag{6.18}:
$$
\hskip -2em
\nabla_kx^i=\left(\,\shave{\sum^2_{q=1}} y_q\,\nabla_kx^q\right)
y^i=z_k\,y^i.
\mytag{6.19}
$$
Using \mythetag{6.19}, i.\,e\. substituting $\nabla_ix^p=z_i\,y^p$ and
$\nabla_jx^q=z_j\,y^q$ into \mythetag{6.15}, we see that the contribution
of the second term in this formula is zero. Then, applying \mythetag{6.16} 
to \mythetag{6.15}, for the increment ${\ssize\triangle}\varphi$ we derive
$$
{\ssize\triangle}\varphi=-\iint\limits_{\Omega}
\sum^2_{i=1}\sum^2_{j=1}
x^i\bigl(\nabla_i\nabla_jx^j-\nabla_j\nabla_ix^j
\bigr)\,\sqrt{\det\bold g}\,du^1du^2.
$$
Now we apply the relationship \mythetagchapter{8.5}{4} from Chapter
\uppercase\expandafter{\romannumeral 4} to the field $\bold x$. Moreover,
we take into account the formulas \mythetagchapter{8.24}{4} and
\mythetagchapter{9.9}{4} from Chapter 
\uppercase\expandafter{\romannumeral 4}:
$$
{\ssize\triangle}\varphi=\iint\limits_{\Omega}
\sum^2_{i=1}\sum^2_{j=1}\bigl(K\,g_{ij}\,x^i\,x^j
\bigr)\,\sqrt{\det\bold g}\,du^1du^2.
$$
Remember that the vector field $\bold x$ was chosen to be of the unit
length from the very beginning. Therefore, upon summing up over the 
indices $i$ and $j$ we shall have only the Gaussian curvature under the
integration:
$$
\hskip -2em
{\ssize\triangle}\varphi=\iint\limits_{\Omega}
K\,\sqrt{\det\bold g}\,du^1du^2.
\mytag{6.20}
$$\par
\parshape 3 0pt 360pt 0pt 360pt 180pt 180pt 
     Now let's consider some surface on which a connected and simply connected domain $\Omega$ outlined by a piecewise continuously
differentiable contour $\gamma$ is given (see Fig\.~6.1). In other
\vadjust{\vskip 5pt\hbox to 0pt{\kern -5pt
\includegraphics{ris15.eps}\hss}\vskip -5pt}words, 
we have a polygon with curvilinear sides on the surface.
The Green's formula \mythetag{5.1} is applicable to a a piecewise 
continuously differentiable contour, therefore, the for\-mula 
\mythetag{6.20} is valid in this case. The parallel translation of
the vector $\bold a$ along a piecewise continuously 
differentiable
contour is performed step by step. The result of translating 
the vector $\bold a$ along a side of the curvilinear polygon 
$\gamma$ is used as the initial data for the equations of parallel
translation on the succeeding side. Hence, $\varphi(t)$ is a 
continuous func\-tion, though its derivative can be discon\-tinuous at the 
corners of the polygon.\par
\parshape 3 180pt 180pt 180pt 180pt 0pt 360pt 
     Let's introduce the natural paramet\-rization $t=s$ on the sides of
the polygon $\gamma$. Then we have the unit tangent vector $\boldsymbol
\tau$ on them. The vector-function $\boldsymbol\tau(t)$ is a continuous
function on the sides, except for the corners, where $\boldsymbol\tau(t)$
abruptly turns to the angles ${\ssize\triangle}\psi_1,\,
{\ssize\triangle}\psi_2,\,\ldots,\,{\ssize\triangle}\psi_n$
(see Fig\.~6.1). Denote by $\psi(t)$ the angle between the vector
$\boldsymbol\tau(t)$ and the vector $\bold a(t)$ being parallel translated 
along $\gamma$. We measure this angle from $\bold a$ to $\boldsymbol\tau$
taking for positive the counterclockwise direction. The finction $\psi(t)$
is a continuously differentiable function on $\gamma$ except for the
corners. At these points it has jump discontinuities with jumps 
${\ssize\triangle}\psi_1,\,{\ssize\triangle}\psi_2,\,\ldots,\,
{\ssize\triangle}\psi_n$.\par
     Let's calculate the derivative of the function $\psi(t)$ out of 
its discontinuity points. Applying the considerations associated with the
expansions \mythetag{6.8} and \mythetag{6.9} to the vector
$\boldsymbol\tau(t)$, for such derivative we find:
$$
\hskip -2em
\dot\psi=(\Theta(\boldsymbol\tau)\,|\,\nabla_{\!t}\boldsymbol\tau).
\mytag{6.21}
$$
Then let's calculate the components of the vector
$\nabla_{\!t}\boldsymbol\tau$ in the inner representation of the surface
(i\.\,e\. in the basis of the frame vectors $\bold E_1$ and $\bold E_2$):
$$
\hskip -2em
\nabla_{\!t}\tau^k=\ddot u^k+\sum^2_{i=1}\sum^2_{j=1}
\Gamma^k_{ji}\,\dot u^i\,\dot u^j.
\mytag{6.22}
$$
Keeping in mind that $t=s$ is the natural parameter on the sides of the
polygon $\gamma$, we compare \mythetag{6.22} with the formula \mythetag{2.5}
for the geodesic curvature and with the formula \mythetag{2.4}. As a result
we get the equality
$$
\hskip -2em
\nabla_{\!t}\boldsymbol\tau=k\cdot
\bold n_{\,\text{curv}}-k_{\,\text{norm}}\cdot\bold n
=k_{\,\text{geod}}\cdot\bold n_{\,\text{inner}}.
\mytag{6.23}
$$
But $\bold n_{\,\text{inner}}$ is a unit vector in the tangent plane 
perpendicular to the vector $\boldsymbol\tau$. The same is true for
the vector $\Theta(\boldsymbol\tau)$ in the scalar product \mythetag{6.21}.
\pagebreak Hence, the unit vectors $\bold n_{\,\text{inner}}$ and
$\Theta(\boldsymbol\tau)$ are collinear. Let's denote by
$\varepsilon(t)$ the sign factor equal to the scalar product of these
vectors:
$$
\hskip -2em
\varepsilon=(\Theta(\boldsymbol\tau)\,|\,\bold n_{\,\text{inner}})=
\pm 1.
\mytag{6.24}
$$
Now from the formulas \mythetag{6.23} and \mythetag{6.24} we derive:
$$
\hskip -2em
\dot\psi=\varepsilon\,k_{\,\text{geod}}.
\mytag{6.25}
$$
Let's find the increment of the function $\psi(t)$ gained as a result
of round trip along the whole contour. It is composed by two parts: 
the integral of \mythetag{6.25} and the sum jumps at the corners of the polygon $\gamma$:
$$
\hskip -2em
{\ssize\triangle}\psi=\oint\limits_{\gamma\,\,}
\varepsilon\,k_{\,\text{geod}}\,ds+\sum^n_{i=1}
{\ssize\triangle}\psi_i.
\mytag{6.26}
$$\par
     The angle ${\ssize\triangle}\varphi$ is measured from $\bold x$ to
$\bold a$ in the counterclockwise direction, while the angle
${\ssize\triangle}\psi$ is measured from $\bold a$ to $\boldsymbol\tau$ 
in the same direction. Therefore, the sum ${\ssize\triangle}\varphi
+{\ssize\triangle}\psi$ is the total increment of the angle between
$\bold x$ and $\boldsymbol\tau$. It is important to note that the initial
value and the final value of the vector $\boldsymbol\tau$ upon round trip
along the contour do coincide. The same is true for the vector $\bold x$. 
Hence, the sum of increments ${\ssize\triangle}\varphi
+{\ssize\triangle}\psi$ is an integer multiple of the angle
$2\pi=360^{\circ}$:
$$
{\ssize\triangle}\varphi+{\ssize\triangle}\psi=
2\pi\,r.
\mytag{6.27}
$$
Practically, the value of the number $r$ in the formula \mythetag{6.27} 
is equal to unity. Let's prove this fact by means of the following 
considerations: we perform the continuous deformation of the surface
on Fig\.~6.1 flattening it to a plain, then we continuously deform the
contour $\gamma$ to a circle. During such a continuous deformation the
left hand side of the equality \mythetag{6.27} changes continuously, 
while the right hand side can change only in discrete jumps. Therefore,
under the above continuous deformation of the surface and the contour
both sides of \mythetag{6.27} do not change at all. On a circle the total
angle of rotation of the unit tangent vector is calculated explicitly,
it is equal to $2\pi$. Hence, $r=1$. We take into account this
circumstance when substituting \mythetag{6.20} and \mythetag{6.26} into 
the formula \mythetag{6.27}:
$$
\hskip -2em
\iint\limits_{\Omega}K\,\sqrt{\det\bold g}\,du^1du^2+
\oint\limits_{\gamma\,\,}\varepsilon\,k_{\,\text{geod}}\,ds+
\sum^n_{i=1}{\ssize\triangle}\psi_i=2\pi.\quad
\mytag{6.28}
$$
The formula \mythetag{6.28} is the content of the following theorem which 
is known as the {\it Gauss-Bonnet theorem}.
\mytheorem{6.2} The sum of the external angles of a curvilinear
polygon on a surface is equal to $2\pi$ minus two integrals: the area 
integral of the Gaussian curvature over the interior of the polygon
and the integral of the geodesic curvature (taken with the sign factor
$\varepsilon$) over its \pagebreak perimeter.
\endproclaim
     It is interesting to consider the case where the polygon is formed
by geodesic lines on a surface of the constant Gaussian curvature. The second integral in \mythetag{6.28} then is equal to zero, while the first
integral is easily calculated. For the sum of internal angles of a geodesic
triangle in this case we derive 
$$
\alpha_1+\alpha_2+\alpha_3=\pi+K\,S,
$$
where $K\,S$ is the product of the Gaussian curvature of the surface and
the area of the triangle.
\subhead A philosophic remark\endsubhead By measuring the sum of 
angles of some sufficiently big triangle we can decide whether our
world is flat or it is equipped with the curvature. This is not a joke.
The idea of a curved space became generally accepted in the  
modern notions on the structure of the world.
\par
\newpage
\topmatter
\title
References.
\endtitle
\endtopmatter
\document
\setfirstpage
\par\noindent
\myref{1}{Sharipov~R.~A. {\it Course of linear algebra and
multidimensional geometry}, Bashkir State University, Ufa, 1996; 
see on-line \myhref{http://uk.arXiv.org/abs/math/0405323/}{math.HO/0405323/}
in Electronic Archive \myEarXivlink.}
\myref{2}{Kudryavtsev~L.~D. {\it Course of mathematical analysis, 
Vol.~\uppercase\expandafter{\romannumeral 1} and \uppercase\expandafter{\romannumeral 2}}, 
{\tencyr\char '074}Visshaya Shkola{\tencyr\char '076} 
publishers, Moscow, 1985.}
\myref{3}{Kostrikin~A.~I. {\it Introduction to algebra}, 
{\tencyr\char '074}Nauka{\tencyr\char '076} publishers, Moscow, 1977.}
\myref{4}{Beklemishev~D.~V. {\it Course of analytical geometry and linear algebra}, {\tencyr\char '074}Nauka{\tencyr\char '076} publishers, Moscow, 1985.}
\bigskip 
\bigskip 
\centerline{\bf AUXILIARY REFERENCES\,\footnotemark.}
\bigskip 
\bigskip 
\myref{5}{Sharipov~R.~A. {\it Quick introduction to tensor analysis},
free on-line textbook 
\myhref{http://uk.arXiv.org/abs/math.HO/0403252/}{math.HO/0403252}
in Electronic Archive \myEarXivlink.}
\myref{6}{Sharipov~R.~A. {\it Classical electrodynamics and theory 
of relativity}, Bashkir\linebreak State University, Ufa, 1997; see on-line
\myhref{http://uk.arXiv.org/abs/physics/0311011/}{physics/0311011}
in Electronic Archive \myEarXivlink.}
\footnotetext{\ The references \mycite{5} and \mycite{6} are added in 2004.}
\adjustfootnotemark{-1}
\enddocument
\end